%% file: nonrigid-normal-C-2-1.tex
\documentclass[12pt,twoside,leqno,openany]{amsart}

\pdfoutput=1

\usepackage{amsbsy,amscd,amsfonts,amsmath,amssymb,amsthm,color,
fancybox,fancyhdr,footmisc,graphics,graphicx,ifthen,mathrsfs,
multicol,pdfpages,rotating,times,wasysym,pifont,blkarray}
\usepackage[dvipsnames,svgnames,x11names]{xcolor}

\usepackage[all]{xy}
\usepackage[utf8]{inputenc}
\usepackage[T1]{fontenc}
\usepackage{mathtools}
\newtagform{EngelLie}[\scriptstyle]{$}{$}
\makeatletter\newcommand{\leqnomode}{\tagsleft@true}
\newcommand{\reqnomode}{\tagsleft@false}\makeatother

\input macros.tex


\begin{document}

\setcounter{section}{0}

$\:$

\bigskip\bigskip\bigskip\bigskip\bigskip

\begin{center}

{\large\bf On Convergent Poincar\'e-Moser 
Reduction\footnotemark[1]\medskip
\\
for Levi Degenerate Embedded $5$-Dimensional CR Manifolds}
\label{convergent-Poincare-Moser-reduction}

\medskip

\bigskip\bigskip
\bigskip\bigskip

Wei-Guo {\sc Foo}\footnotemark[2],
Joël~{\sc Merker}\footnotemark[3],
The-Anh {\sc Ta}\footnotemark[3]

\footnotetext[1]{\,
---\,\,in the sense of normalization\,\,---}


\footnotetext[2]{\,\,Hua Loo-Keng Center for Mathematical
Sciences, AMSS, Chinese Academy of Sciences, Beijing, China.} 

\footnotetext[3]{\,\,Laboratoire de Mathématiques d'Orsay,
CNRS, Université Paris-Saclay, 91405 Orsay Cedex,
France.}

\end{center}\bigskip

\begin{center}
\begin{minipage}[t]{12.5cm}
\parindent 0.53cm
\footnotesize
\noindent
{\sc Abstract}.
Firstly, 
applying Lie's elementary theory for appropriate prolongations to jet
spaces of orders $1$ and $2$, we show that any
$\mathcal{C}^\omega$ hypersurface $M^5 \subset \C^3$ in the class
$\mathfrak{C}_{2,1}$ carries {\em two} sorts of Cartan-Moser {\sl
chains}, that are of orders $1$ and $2$.

Secondly, integrating and straightening any given order $2$ chain
passing through any point $p \in M$ to be the $v$-axis in coordinates
$(z, \zeta, w = u + i\, v)$ centered at $p$, without setting up the
formal theory in advance, we show that there exists a {\em convergent}
change of complex coordinates $(z, \zeta, w) \longmapsto (z', \zeta',
w')$ fixing the origin in which $\gamma$ is the $v$-axis and in which
$M$ has {\em Poincar\'e-Moser reduced equation} (suppressing primes):
\[
\aligned
u
&
\,=\,
z\overline{z}
+
\tfrac{1}{2}\,{\overline{z}}^2\zeta
+
\tfrac{1}{2}\,{z}^2\overline{\zeta}
+
z\overline{z}\zeta\overline{\zeta}
+
\tfrac{1}{2}\,{\overline{z}}^2\zeta\zeta\overline{\zeta}
+
\tfrac{1}{2}\,{z}^2\overline{\zeta}\zeta\overline{\zeta}
+
z\overline{z}\zeta\overline{\zeta}\zeta\overline{\zeta}
\notag
\\
&
\ \ \ \ \ \ \ \ \ \ \ 
+
2\,\Resmall\,
\Big\{
{z}^3{\overline{\zeta}}^2\,
F_{3,0,0,2}(v)
+
\zeta\overline{\zeta}\,
\big(
3\,{z}^2\overline{z}\overline{\zeta}\,
F_{3,0,0,2}(v)
\big)
\big\}
\notag
\\
&
\ \ \ \ \ \ \ \ \ \ \ \ 
+
2\,\Resmall\,
\Big\{
{z}^5\overline{\zeta}\,
F_{5,0,0,1}(v)
+
{z}^4{\overline{\zeta}}^2\,
F_{4,0,0,2}(v)
+
{z}^3{\overline{z}}^2\overline{\zeta}\,
F_{3,0,2,1}(v)
\\
&
\ \ \ \ \ \ \ \ \ \ \ \ \ \ \ \ \ \ \ \ \ \ \ \ \ \ \ \ \ \ \ \ \ \ \
\ \ \ \ \ \ \ \ \ \ \ \ \ \ \ \ \ \ \ \
+
{z}^3\overline{z}{\overline{\zeta}}^2\,
F_{3,0,1,2}(v)
+
{z}^3{\overline{\zeta}}^3\,
F_{3,0,0,3}(v)
\Big\}
\notag
\\
&
\ \ \ \ \ \ \ \ \ \ \ 
+
{z}^3{\overline{z}}^3\,
{\rm O}_{z,\overline{z}}(1)
+
{\overline{z}}^3\zeta\,
{\rm O}_{z,\zeta,\overline{z}}(3)
+
{z}^3\overline{\zeta}\,
{\rm O}_{z,\overline{z},\overline{\zeta}}(3)
+
\zeta\overline{\zeta}\,
{\rm O}_{z,\zeta,\overline{z},\overline{\zeta}}(5),
\endaligned
\]
where all monomials 
in $\zeta \overline{\zeta} (\cdots)$ gather {\em dependent} 
derivatives on which normalizations act automatically.

Thirdly, starting from an $M$ having preliminary normalized equation:
\[
u
\,=\,
z\overline{z}
+
\tfrac{1}{2}\,\overline{z}^2\zeta
+
\tfrac{1}{2}\,z^2\overline{\zeta}
+
z\overline{z}\zeta\overline{\zeta}
+
{\rm O}_{z,\zeta,\overline{z},\overline{\zeta},v}(5),
\]
assigning weights $[z] := 1$, $[\zeta] := 0$, $[w] := 2$, we show that
a normalizing biholomorphism 
exists and is {\em unique} when it is assumed to be of the
form:
\[
\aligned
z'
\,:=\,
&\,
z
+
f_{\geqslant 2}(z,\zeta,w)
\ \ \ \ \ \ \ \ \ \ \ \ \ \ \ \ \ \ \ \
\zeta'
\,:=\,
\zeta
+
g_{\geqslant 1}(z,\zeta,w),
&
\ \ \ \ \ \ \ \ \ \ \ \ \ \ \ \ \ \ \ \
w'
\,:=\,
&\,
w
+
h_{\geqslant 3}(z,\zeta,w),
\\
0
\,=\,
&\,
f_w(0),
&
\ \ \ \ \ \ \ \ \ \ \ \ \ \ \ \ \ \ \ \
0
\,=\,
&\,
\Imsmall\,h_{ww}(0).
\endaligned
\]

The values at the origin of Pocchiola's two primary 
Cartan-type relative differential invariants are:
\[
\Waux_0
\,=\,
4\,\overline{F_{3,0,0,2}(0)}
\ \ \ \ \ \ \ \ \ \ \ \ \ \ \ \ \ \ \ \
\text{and}
\ \ \ \ \ \ \ \ \ \ \ \ \ \ \ \ \ \ \ \
\Jaux_0
\,=\,
20\,
F_{5,0,0,1}(0).
\]

The proofs are detailed, accessible to non-experts. 
The computer-generated aspects (forthcoming)
have been reduced to a minimum here.
\end{minipage}
\end{center}

\Section{\bf Introduction}
\label{introduction-nonrigid-normal-C-2-1}
\HEAD{{\ref{introduction-nonrigid-normal-C-2-1}}.~{\sf Introduction}
}{
Wei-Guo {\sc Foo}, Joël {\sc Merker}, The-Anh {\sc Ta}}

As explained in the survey introduction
of~{\cite{Chen-Foo-Merker-Ta-2019}}, the appropriate local graphed
model for $2$-nondegenerate constant Levi rank $1$ 
real analytic ($\mathcal{C}^\omega$) hypersurfaces $M^5
\subset \C^3$, generally graphed, in coordinates
$\big(z, \zeta, w = u + i\,v
\big)$ as: 
\[
u
\,=\,
F\big(z,\zeta,\overline{z},\overline{\zeta},v\big),
\]
is the so-called {\sl Gaussier-Merker model:}
\[
u
\,=\,
\frac{z\overline{z}
+\frac{1}{2}\,\overline{z}^2\zeta
+\frac{1}{2}\,z^2\overline{\zeta}}{1-\zeta\overline{\zeta}}
\,\,=:\,\,
\maux\big(z,\zeta,\overline{z},\overline{\zeta}\big).
\]
Fels-Kaup~{\cite{Fels-Kaup-2007}} showed that its (connected)
intersection with $\{ \vert \zeta\vert < 1\}$ is biholomorphic to a
Zariski-open subset of the complex
tube $S_{\sf LC}^2 \times i\,\R^3$ over the
real light cone $(\Re\,z_2)^2 - (\Re\,z_3)^2 = (\Re\,z_1)^2$.  The
light cone $S_{\sf LC}^2 \subset \R^3$ 
is {\em the} maximally symmetric non-flat
parabolic surface, characterized, according
to~{\cite{Chen-Merker-2019}}, by the vanishing of certain 
two differential
invariants.

By applying either Cartan's method of equivalence, or Tanaka's
approach, several recent works
({\cite{Isaev-Zaitsev-2013,
Medori-Spiro-2014, 
Medori-Spiro-2015,
Pocchiola-2013,
Merker-Pocchiola-2018,
Foo-Merker-2019})
have been devoted to construct absolute parallelisms, namely
$10$-dimensional $\{e\}$-structure bundles $P^{10} \longrightarrow
M^5$ for such $M^5 \subset \C^3$, 
invariantly related to biholomorphic equivalences of such
hypersurfaces.

By performing advanced electronic computations,
Merker-Pocchiola~{\cite{Pocchiola-2013, Merker-Pocchiola-2018}} found
that only two primary curvature invariants exist, denoted $\Waux$ and
$\Iaux$.  These intensive computations have been redone manually by
Foo-Merker in~{\cite{Foo-Merker-2019}} all along $\sim 50$ pages.  One
obtains certain `horizontal' (semi-basic) $1$-forms $\big\{ \rho,
\kappa, \zeta, \overline{\kappa}, \overline{\zeta} \big\}$ 
with $\overline{\rho} = \rho$
together
with four `vertical' $1$-forms $\pi^1$, $\pi^2$, $\overline{\pi}^1$,
$\overline{ \pi}^2$ which satisfy `compact' structure equations of the
form:
\[
\aligned
d\rho
&
\,=\,
\big(
\pi^1
+
\overline{\pi}^1
\big)
\wedge\rho
+
i\,\kappa\wedge\overline{\kappa},
\notag
\\
d\kappa
&
\,=\,
\pi^2\wedge\rho
+
\pi^1\wedge\kappa
+
\zeta\wedge\overline{\kappa},
\\
d\zeta
&
\,=\,
\big(\pi^1-\overline{\pi}^1\big)
\wedge\zeta
+
i\,\pi^2\wedge\kappa
\,+
\\
&
\ \ \ \ \
+
\Raux\,
\rho\wedge\zeta
+
i\,
\frac{1}{\overline{\sf c}^3}\,
\overline{\Jaux}_0\,\rho\wedge\overline{\kappa}
+
\frac{1}{{\sf c}}\,
\Waux_0\,
\kappa\wedge\zeta,
\endaligned
\]
conjugate structure equations for 
$d\overline{\kappa}$, $d\overline{\zeta}$
being easily deduced.

In Sections~{\ref{normalizations-origin}}
and~{\ref{end-point-normalization-M5-C3}}, we 
copy the expressions of
the two primary
relative differential invariants 
$\Waux_0 \colon M \longrightarrow \C$ 
and $\overline{\Jaux}_0 \colon M \longrightarrow \C$,
while $\Raux$ is a certain (useless) secondary invariant.

\begin{Theorem}
{\rm {\cite{Pocchiola-2013, Merker-Pocchiola-2018,
Foo-Merker-2019}}}
Only two primary invariants, $\Waux_0$ and $\Jaux_0$, occur for
biholomorphic equivalences of 
$2$-nondegenerate constant Levi rank $1$ real analytic
hypersurfaces $M^5 \subset \C^3$, and:
\[
0
\,\equiv\,
\Waux_0
\,\equiv\,
\Jaux_0
\ \ \ \ \ \ \ 
\Longleftrightarrow
\ \ \ \ \ \ \ 
M\,\,
\text{is equivalent to the
Gaussier-Merker model}.
\]
Furthermore, when either $\Waux_0 \neq 0$ or $\Jaux_0 \neq 0$, the
equivalence problem reduces to a $5$-dimensional $\{e\}$-structure on
$M^5$, and every non-flat $M^5$ has CR automorphisms group of
dimension $\leqslant 5$.\qed
\end{Theorem}

In this article, our motivation is to view again these relative CR
differential invariants by putting the equation of such 
$M^5 \subset \C^3$
into normal form, like Chern-Moser did in~{\cite{Chern-Moser-1974}}.
Generally, the 
Poincaré-Moser normal form~\cite{Chern-Moser-1974} provides a
distinguished choice of local holomorphic coordinates for a
hypersurface, in which its defining equation is approximated 
as far 
as possible by that of the local model,
for instance in $\C^{n+1} \ni (z_1, \dots, z_n, w = u+i\,v)$,
a real hyperquadric:
\[
u
\,=\,
\vert z_1\vert^2
+\cdots+
\vert z_p\vert^2
-
\vert z_{p+1}\vert^2
-\cdots-
\vert z_n\vert^2.
\]
Usually, a biholomorphic transformation bringing a hypersurface to a
normal form at the origin is defined up to
composition with the automorphisms
group of the model.

Two months ago, in~{\cite{Chen-Foo-Merker-Ta-2019}}}, 
joint with Chen, we studied
{\em rigid} $\mathcal{C}^\omega$
hypersurfaces $M^5 \subset \C^3$: 
\[
u
\,=\,
F\big(z,\zeta,\overline{z},\overline{\zeta}\big)
\,=\,
\sum_{a,b,c,d\geqslant 0}\,
z^a\zeta^b\overline{z}^c\overline{\zeta}^d\,
F_{a,b,c,d}
\eqno
{\scriptstyle{(F_{a,b,c,d}\,\in\,\C,\,\,
\overline{F_{c,d,a,b}}\,=\,F_{a,b,c,d})}},
\]
with graphing function $F$ {\em independent of $v$},
which are 
everywhere $2$-nondegenerate and of constant Levi rank $1$,
under the {\sl rigid biholomorphisms group},
a group which consists
of transformations of the form:
\[
(z,\zeta,w)
\,\,\,\longmapsto\,\,\,
\Big(
f(z,\zeta),\,
g(z,\zeta),\,
\rho\,w+h(z,\zeta)
\Big)
\,=:\,
\big(z',\zeta',w'\big),
\]
having nonzero holomorphic Jacobian $f_zg_\zeta -f_\zeta g_z \neq 0$, 
with $\rho \in \R^\ast$.
We established that every such rigid $M^5 \subset \C^3$ 
is {\em rigidly equivalent} to a `perturbation' of
the Gaussier-Merker model:
\[
\aligned
u
&
\,=\,
\frac{z\overline{z}+\frac{1}{2}\,\overline{z}^2\zeta
+\frac{1}{2}\,z^2\overline{\zeta}}{
1-\zeta\overline{\zeta}}
+
2\,\Re\,
\Big\{
F_{4,0,0,1}\,z^4\overline{\zeta}
+
\Re\,F_{3,0,1,1}\,z^3\overline{z}\overline{\zeta}
+
F_{3,0,0,2}\,z^3\overline{\zeta}^2
\Big\}
\\
&
\ \ \ \ \ \ \ \ \ \ \ \ \ \ \
+
z^3\overline{z}^3\,{\rm O}_{z,\overline{z}}(0)
+
2\,\Re\,z^3\overline{\zeta}\,
{\rm O}_{z,\overline{z},\overline{\zeta}}(2)
+
\zeta\overline{\zeta}\,
{\rm O}_{z,\overline{z}}(3)\,
{\rm O}_{z,\zeta,\overline{z},\overline{\zeta}}(1).
\endaligned
\]
Here, by writing $\Re\, F_{3,0,1,1}$, we mean that the
(complex) coefficient $F_{3,0,1,1} \in \C$ has been normalized
to be real.

Furthermore, writing:
\[
\aligned
u
\,=\,
F\big(z,\zeta,\overline{z},\overline{\zeta}\big)
&
\,=\,
\maux\big(z,\zeta,\overline{z},\overline{\zeta}\big)
+
G\big(z,\zeta,\overline{z},\overline{\zeta}\big)
\\
&
\,=\,
\maux\big(z,\zeta,\overline{z},\overline{\zeta}\big)
+
\sum_{a,b,c,d\in\N
\atop
a+c\geqslant 3}\,
G_{a,b,c,d}\,
z^a\zeta^b\overline{z}^c\overline{\zeta}^d,
\endaligned
\]
two such rigid $\mathcal{C}^\omega$ hypersurfaces $M^5
\subset \C^3$ and ${M'}^5 \subset {\C'}^3$, both brought into such a
normal form, are rigidly biholomorphically equivalent if and only if
there exist two constants $\rho \in \R_+^\ast$, $\varphi \in \R$, such
that for all $a$, $b$, $c$, $d$:
\[
G_{a,b,c,d}
\,=\,
G_{a,b,c,d}'\,
\rho^{\frac{a+c-2}{2}}\,
e^{i\varphi(a+2b-c-2d)}.
\]
This means that the normal form is defined only up to the
$2$-dimensional action of the {\em rigid} isotropy group
of the origin:
\[
(z,\zeta,w)
\,\,\,\longmapsto\,\,\,
\big(
\rho^{1/2}\,e^{i\varphi}\,z,\,\,
e^{2i\varphi}\,\zeta,\,\,
\rho\,w
\big)
\eqno
{\scriptstyle{(\rho\,\in\,\R_+^\ast,\,\,
\varphi\,\in\,\R)}},
\]

Before making public this normal form, in~{\cite{Foo-Merker-Ta-2019}}, 
we produced Cartan-type
reduction to an $\{e\}$-structure
for the equivalence problem, under {\em rigid}
(local) biholomorphic transformations, 
of such rigid $M^5$
that are $2$-nondegenerate of constant Levi rank $1$.
We constructed an invariant $7$-dimensional bundle
$P^7 \longrightarrow M^5$ equipped with coordinates:
\[
\big(
z_1,z_2,\overline{z}_1,\overline{z}_2,v,\,
{\sf c},\overline{\sf c}
\big),
\]
with ${\sf c} \in \C$,
together with of seven $1$-forms
generating $T^\ast\!P^7$, denoted:
\[
\big\{
\rho,\,
\kappa,\,
\zeta,\,
\overline{\kappa},\,
\overline{\zeta},\,
\alpha,\,
\overline{\alpha}
\big\}
\eqno
{\scriptstyle{(\overline{\rho}\,=\,\rho)}},
\] 
which satisfy invariant structure equations of the form:
\[
\aligned
d\rho 
&
\,=\,
\big(\alpha+\overline{\alpha}\big)
\wedge
\rho 
+ 
i\,
\kappa
\wedge
\overline{\kappa},
\\
d\kappa 
&
\,=\,
\alpha
\wedge
\kappa 
+ 
\zeta
\wedge
\overline{\kappa},
\\
d\zeta 
&
\,=\,
\big(
\alpha-\overline{\alpha}
\big)
\wedge
\zeta
+
\frac{1}{\sf c}\,\Iaux_{0}\,
\kappa\wedge\zeta
+
\frac{1}{\overline{\sf c}\overline{\sf c}}\,
\Vaux_{0}\,
\kappa\wedge\overline{\kappa},
\\
d\alpha 
&
\,=\, 
\zeta\wedge\overline{\zeta}
-
\frac{1}{\sf c}\,
\Iaux_{0}\,
\zeta\wedge\overline{\kappa}
+
\frac{1}{{\sf c}\overline{\sf c}}\,\Qaux_{0}\,
\kappa
\wedge
\overline{\kappa}
+
\frac{1}{\overline{\sf c}}\,
\overline{\Iaux}_{0}\, 
\overline{\zeta}
\wedge
\kappa.
\endaligned
\]
We refer to~{\cite{Chen-Foo-Merker-Ta-2019}} for
explicit expressions of the two primary invariants $\Iaux_0, \Vaux_0
\colon M \longrightarrow \C$, and of the secondary invariant $\Qaux_0
\colon M \longrightarrow \R$, which is real.  Once $M$ is put into
normal form as above, their values at the origin are:
\[
\Iaux_0
\,=\,
4\,\overline{F_{3,0,0,2}}
\ \ \ \ \ \ \ \ \ \ \ \ \ \ \ \ \ \ \ \
\Vaux_0
\,=\,
-\,8\,\overline{F_{4,0,0,1}}
\ \ \ \ \ \ \ \ \ \ \ \ \ \ \ \ \ \ \ \
\Qaux_0
\,=\,
4\,\Re\,
F_{3,0,1,1}.
\]

The goal of this article is to set up a rigorous
{\em convergent} Poincaré-Moser normal
form for any everywhere $2$-nondegenerate constant Levi rank
$1$ general (nonrigid) $\mathcal{C}^\omega$ 
hypersurface $M^5 \subset \C^3$ under the {\em full}
(not necessarily rigid) biholomorphisms group:
\[
(z,\zeta,w)
\,\,\,\longmapsto\,\,\,
\Big(
f(z,\zeta,w),\,
g(z,\zeta,w),\,
h(z,\zeta,w)
\Big).
\]
Given such an $M^5 \subset \C^3$ 
with $0 \in M$, by examining terms of $F$ up to order $4$,
it is elementary to find a holomorphic system of coordinates
in which it is:
\[
u
\,=\,
F
\,=\,
z\overline{z}
+
\tfrac{1}{2}\,\overline{z}^2\zeta
+
\tfrac{1}{2}\,z^2\overline{\zeta}
+
z\overline{z}\zeta\overline{\zeta}
+
{\rm O}_{z,\zeta,\overline{z},\overline{\zeta},v}(5).
\]
Since the Gaussier-Merker model is invariant under the complex
scalings:
\[
(z,\zeta,w)
\,\,\longmapsto\,\,
\big(
\lambda\,z,\,\,
\tfrac{\lambda}{\overline{\lambda}}\,\zeta,\,\,
\lambda\overline{\lambda}\,w
\big)
\eqno
{\scriptstyle{(\lambda\,\in\,\C^\ast)}},
\]
it is natural to assign the weights:
\[
[z]
\,:=\,
1
\,=:\,
[\overline{z}],
\ \ \ \ \ \ \ \ \ \ \ \ \ \ \ \ \ \ \ \
[\zeta]
\,:=\,
0
\,=:\,
[\overline{\zeta}],
\ \ \ \ \ \ \ \ \ \ \ \ \ \ \ \ \ \ \ \
[w]
\,:=\,
2
\,=:\,
[\overline{w}].
\]
Then by $e_{\geqslant\nu}(z,\zeta,w)$, we will mean a
holomorphic function near the origin all of whose monomials
$z^a \zeta^b w^e$ are of weight $a + 2\,e \geqslant \nu$.

\begin{Theorem}
{\bf [Main]}
There exists a biholomorphism $(z,\zeta,w) \longmapsto (z', \zeta',
w')$ fixing $0$ which maps $(M, 0)$ into $(M',0)$ of normalized
equation (suppressing primes):
\[
\aligned
u
&
\,=\,
\frac{z\overline{z}+\frac{1}{2}\,\overline{z}^2\zeta
+\frac{1}{2}\,z^2\overline{\zeta}}{1-\zeta\overline{\zeta}}
\\
&
\ \ \ \ \ \ \ \ \ \ \ 
+
2\,\Re\,
\Big\{
{z}^3{\overline{\zeta}}^2\,
F_{3,0,0,2}(v)
+
\zeta\overline{\zeta}\,
\big(
3\,{z}^2\overline{z}\overline{\zeta}\,
F_{3,0,0,2}(v)
\big)
\big\}
\notag
\\
&
\ \ \ \ \ \ \ \ \ \ \ \ 
+
2\,\Re\,
\Big\{
{z}^5\overline{\zeta}\,
F_{5,0,0,1}(v)
+
{z}^4{\overline{\zeta}}^2\,
F_{4,0,0,2}(v)
+
{z}^3{\overline{z}}^2\overline{\zeta}\,
F_{3,0,2,1}(v)
\\
&
\ \ \ \ \ \ \ \ \ \ \ \ \ \ \ \ \ \ \ \ \ \ \ \ \ \ \ \ \ \ \ \ \ \ \
\ \ \ \ \ \ \ \ \ \ \ \ \ \ \
+
{z}^3\overline{z}{\overline{\zeta}}^2\,
F_{3,0,1,2}(v)
+
{z}^3{\overline{\zeta}}^3\,
F_{3,0,0,3}(v)
\Big\}
\notag
\\
&
\ \ \ \ \ \ \ \ \ \ \ 
+
{z}^3{\overline{z}}^3\,
{\rm O}_{z,\overline{z}}(1)
+
{\overline{z}}^3\zeta\,
{\rm O}_{z,\zeta,\overline{z}}(3)
+
{z}^3\overline{\zeta}\,
{\rm O}_{z,\overline{z},\overline{\zeta}}(3)
+
\zeta\overline{\zeta}\,
{\rm O}_{z,\overline{z}}(3)\,
{\rm O}_{z,\zeta,\overline{z},\overline{\zeta}}(2).
\endaligned
\]

Furthermore, the map exists and is {\em unique} if it is assumed
to be of the form:
\[
\aligned
z'
\,:=\,
&\,
z
+
f_{\geqslant 2}(z,\zeta,w)
\ \ \ \ \ \ \ \ \ \ \ \ \ \ \ \ \ \ \ \
\zeta'
\,:=\,
\zeta
+
g_{\geqslant 1}(z,\zeta,w),
&
\ \ \ \ \ \ \ \ \ \ \ \ \ \ \ \ \ \ \ \
w'
\,:=\,
&\,
w
+
h_{\geqslant 3}(z,\zeta,w),
\\
0
\,=\,
&\,
f_w(0),
&
\ \ \ \ \ \ \ \ \ \ \ \ \ \ \ \ \ \ \ \
0
\,=\,
&\,
\Im\,h_{ww}(0).
\endaligned
\]
\end{Theorem}

Equivalently, writing:
\[
u
\,=\,
F
\,=\,
\sum_{a,b,c,d\geqslant 0}\,
z^a\zeta^b\overline{z}^c\overline{\zeta}^d\,
F_{a,b,c,d}(v),
\]
the normal form is defined by the general 
{\sl prenormalization conditions:}
\[
\aligned
0
\,\equiv\,
&
F_{a,b,0,0}(v)
\equiv\,
F_{0,0,c,d}(v),
\\
0
\,\equiv\,
&
F_{a,b,1,0}(v)
\,\equiv\,
F_{1,0,c,d}(v),
\\
0
\,\equiv\,
&
F_{a,b,2,0}(v)
\,\equiv\,
F_{2,0,c,d}(v),
\endaligned
\]
with the obvious two exceptions
$F_{1,0,1,0}(v) \equiv 1$ and $F_{0,1,2,0}(v) \equiv
\frac{1}{2} \equiv F_{2,0,0,1}(v)$, together with the
{\sl sporadic normalization conditions}, listed by 
increasing order $4$, $5$, $6$:
\[
\aligned
0
&
\,\equiv\,
F_{3,0,0,1}(v)
\,\equiv\,
F_{0,1,3,0}(v),
&
\ \ \ \ \ \ \ \ \ \ \ \ \ \ \ \ \ \ \ \ \ \ \ \ \ \
&
\\
0
&
\,\equiv\,
F_{4,0,0,1}(v)
\,\equiv\,
F_{0,1,4,0}(v),
&
\ \ \ \ \ \ \ \ \ \ \ \ \ \ \ \ \ \ \ \ \ \ \ \ \ \
0
&
\,\equiv\,
F_{3,0,1,1}(v)
\,\equiv\,
F_{1,1,3,0}(v),
\\
0
&
\,\equiv\,
F_{4,0,1,1}(v)
\,\equiv\,
F_{1,1,4,0}(v),
&
\ \ \ \ \ \ \ \ \ \ \ \ \ \ \ \ \ \ \ \ \ \ \ \ \ \
0
&
\,\equiv\,
F_{3,0,3,0}(v).
\endaligned
\]

Without the above conditions 
$z' = z + f_{\geqslant 2}$, 
$\zeta' = \zeta + g_{\geqslant 1}$,
$w' = w + h_{\geqslant 3}$ guaranteeing uniqueness, 
one can verify that  
a normalizing transformation is unique up to the right action
of the 5-dimensional stability group of the 
Gaussier-Merker model having the finite equations:
\[
\aligned
z'
&
\,:=\,
\lambda\,
\frac{z+i\,\alpha\,z^2+
\big(i\,\alpha\,\zeta-i\,\overline{\alpha}\big)\,w}
{1+2i\,\alpha\,z-\alpha^2z^2
-
\big(\alpha^2\zeta-\alpha\overline{\alpha}+i\,r\big)\,w},
\\
\zeta'
&
\,:=\,
\frac{\lambda}{\overline{\lambda}}\,\,
\frac{\zeta+2i\,\overline{\alpha}\,z
-\big(\alpha\overline{\alpha}+i\,r\big)\,z^2
+\big(\overline{\alpha}^2-i\,r\,\zeta
-\alpha\overline{\alpha}\,\zeta\big)\,w}
{1+2i\,\alpha\,z-\alpha^2z^2
-
\big(\alpha^2\zeta-\alpha\overline{\alpha}+i\,r\big)\,w},
\\
w'
&
\,:=\,
\lambda\overline{\lambda}\,
\frac{w}
{1+2i\,\alpha\,z-\alpha^2z^2
-
\big(\alpha^2\zeta-\alpha\overline{\alpha}+i\,r\big)\,w},
\endaligned
\]
where $\lambda \in \C^\ast$, $\alpha \in \C$, $r \in \R$ are
arbitrary.

Lastly, the values at the origin of Pocchiola's two primary 
Cartan-type relative differential invariants are:
\[
\Waux_0
\,=\,
4\,\overline{F_{3,0,0,2}(0)}
\ \ \ \ \ \ \ \ \ \ \ \ \ \ \ \ \ \ \ \
\text{and}
\ \ \ \ \ \ \ \ \ \ \ \ \ \ \ \ \ \ \ \
\Jaux_0
\,=\,
20\,
F_{5,0,0,1}(0).
\]

However, Poincaré-Moser normal forms or
Cartan-Tanaka reductions to $\{e\}$-structures are  
only a {\em preliminary} towards the understanding
of the biholomorphic equivalence problem 
for embedded $\mathcal{C}^\omega$ CR submanifolds
$M \subset \C^\NN$, quite far from any resolution,
not even to be termed `complete resolution'. 

Indeed, focusing on CR geometry, 
we would like to indicate two `defects' of 
Poincaré-Moser normal forms in comparison to
Cartan-Tanaka principal bundles.

\medskip\noindent$\bullet$\,
Moser-type CR normal forms are in fact {\em incomplete} 
in the sense that their invariants are only {\em relative},
yet defined up to the action of a certain ambiguity (isotropy) group.
 
\medskip\noindent$\bullet$\,
Moser-type CR normal forms hold only at one point, hence
are incapable to fully characterize flatness as 
Cartan's method {\em does}.

\medskip

The main reason why Cartan's method is more powerful is that it
embraces computations {\em at every point} of a given manifold.
Objects manipulated by Cartan's thought are (often very complicated)
rational differential expressions in partial derivatives of
fundamental (graphing) functions. In comparison, objects manipulated
by Moser's method are only plain Taylor coefficients,
hence computations are {\em much more elementary}.

Fortunately, it is known that 
symmetries of a hypersurface can be read off
from subsequently constructed {\em deeper} normal forms,
not touched in the present paper, but forthcoming.

These comments conduct us to at least formulate and raise
a certain number of questions showing
that several mysteries remain. 

\medskip\noindent{\bf Q}\textsuperscript{\ding{192}}
How to get rid of ambiguity in Moser CR-normal forms?
What are the true (absolute) differential invariants?
Can one retrieve Pocchiola's dimension drop $10 \downarrow 5$?
Can one link Moser's punctual invariants with Cartan's
invariants at every point?

\medskip\noindent{\bf Q}\textsuperscript{\ding{193}}
In all possibly existing branches, how to find a minimal
set of generators for the differential algebra of absolute
differential invariants? Using either Moser's or Cartan's method?

\medskip\noindent{\bf Q}\textsuperscript{\ding{194}}
In each branch, what are the differential relations
(syzygies) between differential invariants?

\medskip\noindent{\bf Q}\textsuperscript{\ding{195}}
How to implement the determination of CR-homogeneous
models beyond naive Taylor series manipulations at only one point?
How to employ the theory of Lie? How to view
Cartan's invariants in a Taylor series?

\medskip\noindent{\bf Q}\textsuperscript{\ding{196}}
How to implement, from Moser's side of the bridge,
any sub-branch assumption that requires that an ideal
of differential invariants, or a collection of Taylor coefficients,
vanish (identically)?

\medskip

To close this brief introduction, 
three aspects of the article should be emphasized.

\medskip\noindent{\bf A}\textsuperscript{\!\!\ding{192}}
Analogs of Cartan-Moser chains will be `discovered from scratch'
by applying a method due to Lie, as in~{\cite{Merker-2020}}.

\medskip\noindent{\bf A}\textsuperscript{\!\!\ding{193}}
Detailed proofs for the existence of a {\em convergent} normal form,
missing on arxiv.org, will be offered to the reader.

\medskip\noindent{\bf A}\textsuperscript{\!\!\ding{194}}
The `{\sl formal theory}' will be developped {\em after} the
`{\sl convergent theory}'. 

\medskip\noindent{\bf Acknowledgments.}
Zhangchi Chen provided the Maple figures of 
Sections~{\ref{prolongations-jet-1}}
and~{\ref{prolongations-jet-2}}.

\Section{\bf $\mathfrak{C}_{2,1}$ Hypersurfaces $M^5 \subset \C^3$}
\label{C-2-1-hypersurfaces-M5-C3}
\HEAD{{\ref{C-2-1-hypersurfaces-M5-C3}}.~{\sf $\mathfrak{C}_{2,1}$ 
Hypersurfaces $M^5 \subset \C^3$}
}{
Wei-Guo {\sc Foo}, Joël {\sc Merker}, The-Anh {\sc Ta}}

Our object of study is the collection
of real $\mathcal{C}^\omega$ hypersurfaces $M^5 \subset \C^3$
whose Levi form is of constant rank $1$ at every point
and that are everywhere $2$-nondegenerate
({\em see} below), a {\sl class}
that we will denote as:
\[
\mathfrak{C}_{2,1}.
\]

Pick any point $p \in M$ and adapt affine holomorphic coordinates
$\big(z, \zeta, w = u + i\, v\big) \in \C^3$ 
in which $p$ is the origin, 
so that 
$T_0 M \oplus \R_u = \C^3$.  From any $\mathcal{C}^\omega$ 
real defining
equation for $M$ near $p$, the analytic implicit function theorem
enables to solve for $u$ as:
\[
u
\,=\,
F\big(z,\zeta,\overline{z},\overline{\zeta},v\big),
\]
for some $\mathcal{C}^\omega$ graphing function $F$, the 
{\sl core object} of our study. 
This $F$ is expandable in converging
power series as:
\[
F
\big(
z,\zeta,\overline{z},\overline{\zeta},v
\big)
\,=\,
\sum_{a+b+c+d+e\geqslant 1}\,
F_{a,b,c,d,e}\,
z^a\zeta^b\overline{z}^c\overline{\zeta}^dv^e,
\]
for some infinite collection of complex coefficients
$F_{a,b,c,d,e} \in \C$. Then by conjugating
only complex coefficients, {\sl define:}
\[
\overline{F}
\big(z,\zeta,\overline{z},\overline{\zeta},v\big)
\,:=\,
\sum_{a+b+c+d+e\geqslant 1}\,
\overline{F}_{a,b,c,d,e}\,
z^a\zeta^b\overline{z}^c\overline{\zeta}^dv^e.
\]

The reality $\overline{u} = u$ forces
$\overline{F(z,\zeta,\overline{z},\overline{\zeta},v)} 
= F(z,\zeta,\overline{z},\overline{\zeta},v)$, that is:
\leqnomode\usetagform{default}
\begin{align}
\label{Fbar-F-identity}
\overline{F}
\big(
\overline{z},\overline{\zeta},z,\zeta,v
\big)
\,\equiv\,
F
\big(
z,\zeta,\overline{z},\overline{\zeta},v
\big).
\end{align}
Applying $\frac{1}{a!}\partial_z^a\,
\frac{1}{b!} \partial_\zeta^b\,
\frac{1}{c!} \partial_{\overline{z}}^c\,
\frac{1}{d!} \partial_{\overline{\zeta}}^d\,
\frac{1}{e!} \partial_v^e$ at the origin
$(0, 0, 0, 0, 0)$, we obtain
the (known) condition on the $F_{a,b,c,d,e} \in \C$ 
which guarantees reality of the graphing function:
\[
\overline{F_{c,d,a,b,e}}
\,=\,
F_{a,b,c,d,e}.
\]

Later, we will expand $F$ in powers of $(z, \zeta, \overline{z},
\overline{\zeta})$ only, by introducing:
\[
F\big(z,\zeta,\overline{z},\overline{\zeta},v\big)
\,=\,
\sum_{a,b,c,d}\,
z^a\zeta^b\overline{z}^c\overline{\zeta}^d\,
\sum_e\,
F_{a,b,c,d,e}\,v^e
\,\,=:\,\,
\sum_{a,b,c,d}\,
z^a\zeta^b\overline{z}^c\overline{\zeta}^d\,
F_{a,b,c,d}(v).
\]
The reality of $F$ is then equivalent to:
\leqnomode\usetagform{default}
\begin{align}
\label{reality-F-v}
\overline{F_{c,d,a,b}(v)}
\,=\,
F_{a,b,c,d}(v).
\end{align}

In the literature~{\cite{
Freeman-1977,
Gaussier-Merker-2003,
Merker-2008,
Fels-Kaup-2008,
Merker-Pocchiola-Sabzevari-2013-5-CR-II,
Isaev-Zaitsev-2013,
Medori-Spiro-2014,
Medori-Spiro-2015,
Merker-Nurowski-2020,
Foo-Merker-2019,
Foo-Merker-Ta-2019}}, 
several equivalent definitions of the class $\mathfrak{C}_{2,1}$
exist. We propose a computational formulation
of the two concepts of
constant Levi rank $1$ and of $2$-nondegeneracy,
already shown in~{\cite{Chen-Foo-Merker-Ta-2019}}
when $M$ is {\sl rigid}, namely when $F$ is idenpendent of $v$.

For this, we need the {\sl complex graphed representation}
of any $\mathcal{C}^\omega$ hypersurface $M^5 \subset \C^3$:
\[
w
\,=\,
Q\big(z,\zeta,\overline{z},\overline{\zeta},\overline{w}\big),
\]
with a $\C$-valued 
analytic function $Q$ which is obtained by solving for $w$ in
$\frac{w+\overline{w}}{2} = F \big(z, \zeta, \overline{z},
\overline{\zeta}, \frac{w-\overline{w}}{2i} \big)$, so that:
\[
\tfrac{1}{2}\,
Q\big(z,\zeta,\overline{z},\overline{\zeta},\overline{w}\big)
+
\tfrac{1}{2}\,
\overline{w}
\,\equiv\,
F
\Big(
z,\zeta,\overline{z},\overline{\zeta},\,
\tfrac{1}{2i}\,
Q\big(z,\zeta,\overline{z},\overline{\zeta},\overline{w}\big)
-
\tfrac{1}{2i}\,\overline{w}
\Big).
\]
Such an analytic function $Q$ cannot be arbitrary, it must
satisfy a compatibility condition obtained by replacing
$\overline{w} := \overline{Q}$ in its last argument:
\[
w
\,\equiv\,
Q\Big(
z,\zeta,\overline{z},\overline{\zeta},\,\,
\overline{Q}\big(\overline{z},\overline{\zeta},z,\zeta,w\big)
\Big).
\]

\Section{\bf Two Invariant Determinants}
\label{two-invariant-determinants}
\HEAD{{\ref{two-invariant-determinants}}.~{\sf Two Invariant 
Determinants}
}{
Wei-Guo {\sc Foo}, Joël {\sc Merker}, The-Anh {\sc Ta}}

A local biholomorphism:
\[
(z,\zeta,w)
\,\,\,\longmapsto\,\,\,
\Big(f(z,\zeta,w),\,
g(z,\zeta,w),\,
h(z,\zeta,w)\Big)
\,=:\,
\big(z',\zeta',w'\big),
\]
has nowhere vanishing holomorphic Jacobian determinant:
\[
0
\,\neq\,
\left\vert\!
\begin{array}{ccc}
f_z & g_z & h_z
\\
f_\zeta & g_\zeta & h_\zeta
\\
f_w & g_w & h_w
\end{array}
\!\right\vert.
\]
Suppose that it makes a biholomorphism
between two $\mathcal{C}^\omega$ hypersurfaces
both represented by complex graphing functions:
\[
w
\,=\,
Q\big(z,\zeta,\overline{z},\overline{\zeta},\overline{w}\big)
\ \ \ \ \ \ \ \ \ \ \ \ \ \ \ \ \ \ \ \
\text{and}
\ \ \ \ \ \ \ \ \ \ \ \ \ \ \ \ \ \ \ \
w'
\,=\,
Q'\big(
z',\zeta',\overline{z}',\overline{\zeta}',\overline{w}'
\big).
\]
Plugging the three components of the 
biholomorphism in the target equation,
we get the so-called {\sl fundamental identity:}
\[
h(z,\zeta,w)
\,\,=\,\,
Q'
\Big(
f(z,\zeta,w),\,
g(z,\zeta,w),\,
\overline{f}(\overline{z},\overline{\zeta},\overline{w}),\,
\overline{g}(\overline{z},\overline{\zeta},\overline{w}),\,
\overline{h}(\overline{z},\overline{\zeta},\overline{w})
\Big)
\bigg\vert_{
w\,=\,Q(z,\zeta,\overline{z},\overline{\zeta},\overline{w})},
\]
which holds
identically in the ring of
converging power series $\C\{z, \zeta, \overline{z}, 
\overline{\zeta}, \overline{w} \}$.

By differentiating this identity (exercise!), 
one may express the invariancy of the Levi
form as a relation between the two Levi determinants
defined as:
\[
\left\vert\!
\begin{array}{ccc}
Q_{\overline{z}} & Q_{\overline{\zeta}} 
& 
Q_{\overline{w}}
\\
Q_{z\overline{z}} & Q_{z\overline{\zeta}} 
& 
Q_{z\overline{w}}
\\
Q_{\zeta\overline{z}} & Q_{\zeta\overline{\zeta}} 
& 
Q_{\zeta\overline{w}}
\end{array}
\!\right\vert
\ \ \ \ \ \ \ \ \ \ \ \ \ \ \ \ \ \ \ \
\text{and}
\ \ \ \ \ \ \ \ \ \ \ \ \ \ \ \ \ \ \ \
\left\vert\!
\begin{array}{ccc}
Q_{\overline{z}'}' & Q_{\overline{\zeta}'}' 
& 
Q_{\overline{w}'}'
\\
Q_{z'\overline{z}'}' & Q_{z'\overline{\zeta}'}'
& 
Q_{z'\overline{w}'}'
\\
Q_{\zeta'\overline{z}'}' & Q_{\zeta'\overline{\zeta}'}'
& 
Q_{\zeta'\overline{w}'}'
\end{array}
\!\right\vert.
\]
Indeed, abbreviate:
\[
\mathcal{L}_z
\,:=\,
\frac{\partial}{\partial z}
+
Q_z(z,\zeta,\overline{z},\overline{\zeta},\overline{w})\,
\frac{\partial}{\partial w}
\ \ \ \ \ \ \ \ \ \ \ \ \ \ \ \ \ \ \ \
\text{and}
\ \ \ \ \ \ \ \ \ \ \ \ \ \ \ \ \ \ \ \
\mathcal{L}_\zeta
\,:=\,
\frac{\partial}{\partial\zeta}
+
Q_\zeta(z,\zeta,\overline{z},\overline{\zeta},\overline{w})\,
\frac{\partial}{\partial w}.
\]

\begin{Proposition}
\label{Prp-Levi-determinant}
Through any biholomorphism between
real hypersurfaces 
$\{w = Q\} \subset \C^3$ and $\{w' = Q'\} \subset {\C'}^3$, one has:
\[
\left\vert\!
\begin{array}{ccc}
Q_{\overline{z}'}' & Q_{\overline{\zeta}'}' 
& 
Q_{\overline{w}'}'
\\
Q_{z'\overline{z}'}' & Q_{z'\overline{\zeta}'}'
& 
Q_{z'\overline{w}'}'
\\
Q_{\zeta'\overline{z}'}' & Q_{\zeta'\overline{\zeta}'}'
& 
Q_{\zeta'\overline{w}'}'
\end{array}
\!\right\vert
\,\,=\,\,
\frac{
\left\vert\!
\begin{array}{ccc}
f_z & f_\zeta & f_w
\\
g_z & g_\zeta & g_w
\\
h_z & h_\zeta & h_w
\end{array}
\!\right\vert^3
}{
\left\vert\!
\begin{array}{ccc}
\overline{f}_{\overline{z}} & \overline{f}_{\overline{\zeta}} &
\overline{f}_{\overline{w}}
\\
\overline{g}_{\overline{z}} & \overline{g}_{\overline{\zeta}} &
\overline{g}_{\overline{w}}
\\
\overline{h}_{\overline{z}} & \overline{h}_{\overline{\zeta}} &
\overline{h}_{\overline{w}}
\end{array}
\!\right\vert^1}\,\,
\frac{1}{
\left\vert\!
\begin{array}{cc}
\mathcal{L}_z(f) & \mathcal{L}_{\zeta}(f)
\\
\mathcal{L}_z(g) & \mathcal{L}_{\zeta}(g)
\end{array}
\!\right\vert^4
}\,\,
\left\vert\!
\begin{array}{ccc}
Q_{\overline{z}} & Q_{\overline{\zeta}} 
& 
Q_{\overline{w}}
\\
Q_{z\overline{z}} & Q_{z\overline{\zeta}} 
& 
Q_{z\overline{w}}
\\
Q_{\zeta\overline{z}} & Q_{\zeta\overline{\zeta}} 
& 
Q_{\zeta\overline{w}}
\end{array}
\!\right\vert.
\eqno\qed
\]
\end{Proposition}

Consequently, the property that the Levi form is of constant
rank $1$ is biholomorphically invariant.
The $2$-nondegeneracy property~{\cite{Merker-Nurowski-2020}} 
then expresses as the nonvanishing of: 
\[
\left\vert\!
\begin{array}{ccc}
Q_{\overline{z}} & Q_{\overline{\zeta}} 
& 
Q_{\overline{w}}
\\
Q_{z\overline{z}} & Q_{z\overline{\zeta}} 
& 
Q_{z\overline{w}}
\\
Q_{zz\overline{z}} & Q_{zz\overline{\zeta}} 
& 
Q_{zz\overline{w}}
\end{array}
\!\right\vert
\ \ \ \ \ \ \ \ \ \ \ \ \ \ \ \ \ \ \ \
\text{and}
\ \ \ \ \ \ \ \ \ \ \ \ \ \ \ \ \ \ \ \
\left\vert\!
\begin{array}{ccc}
Q_{\overline{z}'}' & Q_{\overline{\zeta}'}' 
& 
Q_{\overline{w}'}'
\\
Q_{z'\overline{z}'}' & Q_{z\overline{\zeta}'}'
& 
Q_{z'\overline{w}'}'
\\
Q_{z'z'\overline{z}'}' & Q_{z'z'\overline{\zeta}'}'
& 
Q_{z'z'\overline{w}'}'
\end{array}
\!\right\vert.
\]

\begin{Proposition}
\label{Prp-2ndg-determinant}
When the Levi form is of constant rank $1$, 
through any biholomorphism between real hypersurfaces
$\{w=Q\} \subset \C^3$ and 
$\{w'=Q'\} \subset {\C'}^3$, one has:
\[
\frac{\left\vert\!
\begin{array}{ccc}
Q_{\overline{z}'}' & Q_{\overline{\zeta}'}' 
& 
Q_{\overline{w}'}'
\\
Q_{z'\overline{z}'}' & Q_{z'\overline{\zeta}'}'
& 
Q_{z'\overline{w}'}'
\\
Q_{z'z'\overline{z}'}' & Q_{z'z'\overline{\zeta}'}'
& 
Q_{z'z'\overline{w}'}'
\end{array}
\!\right\vert
}{
\left\vert\!
\begin{array}{ccc}
Q_{\overline{z}} & Q_{\overline{\zeta}} 
& 
Q_{\overline{w}}
\\
Q_{z\overline{z}} & Q_{z\overline{\zeta}} 
& 
Q_{z\overline{w}}
\\
Q_{zz\overline{z}} & Q_{zz\overline{\zeta}} 
& 
Q_{zz\overline{w}}
\end{array}
\!\right\vert
}
\,\,=\,\,
\frac{
\left\vert\!
\begin{array}{ccc}
f_z & f_\zeta & f_w
\\
g_z & g_\zeta & g_w
\\
h_z & h_\zeta & h_w
\end{array}
\!\right\vert^3
}{
\left\vert\!
\begin{array}{ccc}
\overline{f}_{\overline{z}} & \overline{f}_{\overline{\zeta}} &
\overline{f}_{\overline{w}}
\\
\overline{g}_{\overline{z}} & \overline{g}_{\overline{\zeta}} &
\overline{g}_{\overline{w}}
\\
\overline{h}_{\overline{z}} & \overline{h}_{\overline{\zeta}} &
\overline{h}_{\overline{w}}
\end{array}
\!\right\vert^1}\,\,
\frac{
\left(
\mathcal{L}_\zeta(g)
\left\vert\!
\begin{array}{cc}
Q_{\overline{z}} & Q_{\overline{w}}
\\
Q_{z\overline{z}} & Q_{z\overline{w}}
\end{array}
\!\right\vert
-
\mathcal{L}_z(g)
\left\vert\!
\begin{array}{cc}
Q_{\overline{z}} & Q_{\overline{w}}
\\
Q_{\zeta\overline{z}} & Q_{\zeta\overline{w}}
\end{array}
\!\right\vert
\right)^3
}{
\left\vert\!
\begin{array}{cc}
\mathcal{L}_z(f) & \mathcal{L}_\zeta(f)
\\
\mathcal{L}_z(g) & \mathcal{L}_\zeta(g)
\end{array}
\!\right\vert^6
\,\,\,\,\,
\left\vert\!
\begin{array}{cc}
Q_{\overline{z}} & Q_{\overline{w}}
\\
Q_{z\overline{z}} & Q_{z\overline{w}}
\end{array}
\!\right\vert^3}.
\eqno\qed
\]
\end{Proposition}

Recall that we denote the class of 
(local) hypersurfaces $M^5 \subset \C^3$
passing through the origin $0 \in M$ 
that are $2$-nondegenerate and whose Levi form
has constant rank $1$ as:
\[
\mathfrak{C}_{2,1}.
\]

Repeatedly, we shall use the real expression of the 
{\sl Levi determinant:}
\leqnomode\usetagform{default}
\begin{align}
\label{Levi-determinant-F}
\Levi(F)
\,:=\,
\left\vert\!
\def\arraystretch{1.5}
\begin{array}{cccc}
0 & F_z & F_\zeta & -\frac{1}{2}+\frac{1}{2i}F_v
\\
F_{\overline{z}} & F_{z\overline{z}} & F_{\zeta\overline{z}} &
\frac{1}{2i}F_{\overline{z}v}
\\
F_{\overline{\zeta}} & F_{z\overline{\zeta}} & 
F_{\zeta\overline{\zeta}} &
\frac{1}{2i}F_{\overline{\zeta}v}
\\
-\frac{1}{2}-\frac{1}{2i}F_v & -\frac{1}{2i}F_{zv} &
-\frac{1}{2i}F_{\zeta v} & \frac{1}{4}F_{vv}
\end{array}
\!\right\vert.
\end{align}

The next (known) statement applies to $\rho := - u + F$.

\begin{Lemma}
{\rm {\cite{Foo-Merker-Ta-2018}}}
If $M^5 \subset \C^3$ is implicitly defined by $\rho \big(z, \zeta, w,
\overline{z}, \overline{\zeta}, \overline{w} \big) = 0$ with a
$\mathcal{C}^\omega$ real function $\rho = \overline{\rho}$ satisfying
$\rho_w \neq 0$, and if $w = Q \big( z, \zeta, \overline{z},
\overline{\zeta}, \overline{w} \big)$ is its associated complex
graphing function, then:
\[
\left\vert\!
\begin{array}{cccc}
0 & \rho_z & \rho_\zeta & \rho_w
\\
\rho_{\overline{z}} & \rho_{z\overline{z}} & 
\rho_{\zeta\overline{z}} &
\rho_{w\overline{z}}
\\
\rho_{\overline{\zeta}} & \rho_{z\overline{\zeta}} & 
\rho_{\zeta\overline{\zeta}} &
\rho_{w\overline{\zeta}}
\\
\rho_{\overline{w}} & \rho_{z\overline{w}} &
\rho_{\zeta\overline{w}} & \rho_{w\overline{w}}
\end{array}
\!\right\vert
\,\,=\,\,
\rho_w^4\,
\left\vert\!
\begin{array}{ccc}
Q_{\overline{z}} & Q_{\overline{\zeta}} 
& 
Q_{\overline{w}}
\\
Q_{z\overline{z}} & Q_{z\overline{\zeta}} 
& 
Q_{z\overline{w}}
\\
Q_{\zeta\overline{z}} & Q_{\zeta\overline{\zeta}} 
& 
Q_{\zeta\overline{w}}
\end{array}
\!\right\vert.
\eqno
\qed
\]
\end{Lemma}

We leave as an exercise to find some invariant determinant
expressed in terms of $F$ which corresponds to the
$2$-nondegeneracy determinant of 
Proposition~{\ref{Prp-2ndg-determinant}} in terms of $Q$.

\Section{\bf Infinitesimal CR Automorphisms}
\label{infinitesimal-CR-automorphisms}
\HEAD{{\ref{infinitesimal-CR-automorphisms}}.~{\sf 
Infinitesimal CR Automorphisms}
}{
Wei-Guo {\sc Foo}, Joël {\sc Merker}, The-Anh {\sc Ta}}

In the class $\mathfrak{C}_{2,1}$, 
the appropriate homogeneous model, named $M_{\sf LC}$,
was set up by Gaussier-Merker
in~{\cite{Gaussier-Merker-2003}} and
Fels-Kaup in~{\cite{Fels-Kaup-2007}}, {\em see} 
also~{\cite{Chen-Foo-Merker-Ta-2019}}:
\[
M_{\sf LC}
\colon
\ \ \ \ \
u
\,=\,
\frac{z\overline{z}+\frac{1}{2}\,z^2\overline{\zeta}
+\frac{1}{2}\overline{z}^2\zeta}{1-\zeta\overline{\zeta}}
\,\,=:\,\,
\maux\big(z,\zeta,\overline{z},\overline{\zeta}\big).
\]
The letter $\maux$ here stands for $\maux$odel.

The $10$-dimensional simple Lie algebra of its infinitesimal 
CR automorphisms:
\[
\mathfrak{g}
\,:=\,
\mathfrak{aut}_{CR}\big(M_{\sf LC}\big)
\,\cong\,
\mathfrak{so}_{2,3}(\R),
\]
has $10$ natural generators $X_1, \dots, X_{10}$, which are $(1,0)$
vector fields in $\C^3$ having holomorphic coefficients with $X_\sigma
+ \overline{X}_\sigma$ tangent to $M_{\sf LC}$.

It is natural to assign the following weights to variables and to
vector fields:
\leqnomode\usetagform{default}
\begin{align}
\label{weighting-z-zeta-w}
[z]
\,:=\,
1
\ \ \ \ \ \ \ \ 
[\zeta]
\,:=\,
0,
\ \ \ \ \ \ \ \ 
[w]
\,:=\,
2
\ \ \ \ \ \ \ \ 
\big[\partial_z\big]
\,:=\,
-\,1
\ \ \ \ \ \ \ \ 
\big[\partial_\zeta\big]
\,:=\,
0
\ \ \ \ \ \ \ \ 
\big[\partial_w\big]
\,:=\,
-\,2.
\end{align}
The Lie algebra $\mathfrak{g} = \mathfrak{aut}_{CR}(M_{\sf LC})$ can
be graded as:
\[
\mathfrak{g}
\,=\,
\mathfrak{g}_{-2}
\oplus
\mathfrak{g}_{-1}
\oplus
\mathfrak{g}_0
\oplus
\mathfrak{g}_1
\oplus
\mathfrak{g}_2,
\]
where, as shown in~{\cite{Gaussier-Merker-2003,
Foo-Merker-Ta-2019}}:
\[
\aligned
\mathfrak{g}_{-2}
&
\,:=\,
\Span\,
\big\{
i\,\partial_w
\big\},
\\
\mathfrak{g}_{-1}
&
\,:=\,
\Span\,
\big\{
(\zeta-1)\,\partial_z
-
2z\,\partial_w,
\ \ \
(i+i\zeta)\,\partial_z
-
2iz\,\partial_w
\big\},
\endaligned
\]
where $\mathfrak{g}_0 = \mathfrak{g}_0^{\sf trans} 
\oplus \mathfrak{g}_0^{\sf iso}$:
\[
\aligned
\mathfrak{g}_0^{\sf trans}
&
\,:=\,
\Span\,
\Big\{
z\zeta\,\partial_z
+
(\zeta^2-1)\,\partial_\zeta
-
z^2\,\partial_w,
\ \ \
iz\zeta\,\partial_z
+
(i+i\zeta^2)\,\partial_\zeta
-
iz^2\,\partial_w
\Big\},
\\
\mathfrak{g}_0^{\sf iso}
&
\,:=\,
\Span\,
\big\{
z\,\partial_z
+
2w\,\partial_w,
\ \ \ 
iz\,\partial_z
+
2i\zeta\,\partial_\zeta
\big\},
\endaligned
\]
while:
\[
\aligned
\mathfrak{g}_1
&
\,:=\,
\Span\,
\big\{
\big(z^2-\zeta w-w)\,\partial_z
+
\big(2z\zeta+2z\big)\,\partial_\zeta
+
2zw\,\partial_w,
\\
&
\ \ \ \ \ \ \ \ \ \ \ \ \ \ \ \ \ \ \ \
\ \ \ \ \ 
\big(-iz^2+i\zeta w-iw\big)\,\partial_z
+
\big(-2iz\zeta+2iz\big)\,\partial_\zeta
-
2izw\,\partial_w
\big\},
\\
\mathfrak{g}_2
&
\,:=\,
\Span\,
\big\{
izw\,\partial_z
-
iz^2\,\partial_\zeta
+
iw^2\,\partial_w
\big\}.
\endaligned
\]

Calling these $X_1, \dots, X_{10}$ in order of appearance,
the five $X_\sigma + \overline{X}_\sigma$ for $\sigma = 1, 2, 3, 4, 5$
span $TM^5$ while those for $\sigma = 6, 7, 8, 9, 10$
generate the isotropy subgroup of the origin.

\smallskip

In fact, we will use the alternative names for the
$5$ generators of the isotropy subroup:
\[
\aligned
{\sf D}
&
\,:=\,
z\,\partial_z
+
2w\,\partial_w,
\\
{\sf R}
&
\,:=\,
iz\,\partial_z
+
2i\,\zeta\,\partial_\zeta,
\\
{\sf I}_1
&
\,:=\,
\big(z^2-\zeta w-w\big)\,
\partial_z
+
\big(2\,z\zeta+2\,z\big)\,\partial_\zeta
+
2\,zw\,\partial_w,
\\
{\sf I}_2
&
\,:=\,
\big(-i\,z^2+i\,\zeta w-i\,w\big)\,\partial_z
+
\big(-2i\,z\zeta+2i\,z\big)\,\partial_\zeta
-
2i\,zw\,\partial_w,
\\
{\sf J}
&
\,:=\,
i\,zw\,\partial_z
-
i\,z^2\,\partial_\zeta
+
i\,w^2\,\partial_w,
\endaligned
\]
having commutator table:
\begin{center}
\begin{tabular} [t] { l | l l l l l }
& ${\sf D}$ & ${\sf R}$ & ${\sf I}_1$ &
${\sf I}_2$ & ${\sf
J}$
\\
\hline
${\sf D}$ & $0$ & $0$ & ${\sf I}_1$ & ${\sf I}_2$ & $2\,{\sf J}$
\\
${\sf R}$ & $*$ & $0$ & $-{\sf I}_2$ & ${\sf I}_1$ & $0$
\\
${\sf I}_1$ & $*$ & $*$ & $0$ & $4\,{\sf J}$ & $0$
\\
${\sf I}_2$ & $*$ & $*$ & $*$ & $0$ & $0$
\\
${\sf J}$ & $*$ & $*$ & $*$ & $*$ & $0$ %
\end{tabular}
\end{center}

\Section{\bf Fractional Representation of the Isotropy Group}
\label{fractional-representation-isotropy-group}
\HEAD{{\ref{fractional-representation-isotropy-group}}.~{\sf 
Fractional Representation of the Isotropy Group}
}{
Wei-Guo {\sc Foo}, Joël {\sc Merker}, The-Anh {\sc Ta}}

By integrating iterated flows of 
${\sf D}$, ${\sf R}$, ${\sf I}_1$, ${\sf I}_2$, ${\sf J}$, 
it can be shown (exercise) that the isotropy subgroup
of the origin $0 \in M_{\sf LC}$ in the Gaussier-Merker model 
has the finite equations:
\[
\aligned
z'
&
\,:=\,
\lambda\,
\frac{z+i\,\alpha\,z^2+
\big(i\,\alpha\,\zeta-i\,\overline{\alpha}\big)\,w}
{1+2i\,\alpha\,z-\alpha^2z^2
-
\big(\alpha^2\zeta-\alpha\overline{\alpha}+i\,r\big)\,w},
\\
\zeta'
&
\,:=\,
\frac{\lambda}{\overline{\lambda}}\,\,
\frac{\zeta+2i\,\overline{\alpha}\,z
-\big(\alpha\overline{\alpha}+i\,r\big)\,z^2
+\big(\overline{\alpha}^2-i\,r\,\zeta
-\alpha\overline{\alpha}\,\zeta\big)\,w}
{1+2i\,\alpha\,z-\alpha^2z^2
-
\big(\alpha^2\zeta-\alpha\overline{\alpha}+i\,r\big)\,w},
\\
w'
&
\,:=\,
\lambda\overline{\lambda}\,
\frac{w}
{1+2i\,\alpha\,z-\alpha^2z^2
-
\big(\alpha^2\zeta-\alpha\overline{\alpha}+i\,r\big)\,w},
\endaligned
\]
where $\lambda \in \C^\ast$, $\alpha \in \C$, $r \in \R$ are
arbitrary.

The Taylor expansions up to respective weighted orders
$5$, $4$, $6$, will soon be useful:
\[
\!\!\!\!\!\!\!\!\!\!\!\!\!\!\!\!\!\!\!\!
\footnotesize
\aligned
z'
&
\,=\,
\lambda\,z
\\
&
\ \ \ \ \ 
-\,
i\,\lambda\alpha\,z^2
-
i\,\lambda\,\overline{\alpha}\,w
\\
&
\ \ \ \ \
-\,\lambda\alpha^2\,z^3
+
\big(-3\,\lambda\alpha\overline{\alpha}
+
i\,\lambda r\big)\,
zw
+
i\,\lambda\alpha\,\zeta w
\\
&
\ \ \ \ \
+
i\,\lambda\alpha^3\,z^4
+
\big(
6i\,\lambda\alpha^2\overline{\alpha}
+
3\,\lambda\alpha r
\big)\,z^2w
+
\big(
\lambda\overline{\alpha}r
+
i\,\lambda\overline{\alpha}^2\alpha
\big)\,
w^2
+
3\,\lambda\alpha^2\,z\zeta w
\\
&
\ \ \ \ \
+
\lambda\alpha^4\,z^5
+
\big(
-6i\,\alpha^2\lambda r
+
10\,\lambda\alpha^3\overline{\alpha}
\big)\,
z^3w
-
6i\,\lambda\alpha^3\,
z^2\zeta w
+
\big(
5\,\lambda\alpha^2\overline{\alpha}^2
-
6i\,\lambda\alpha\overline{\alpha}r
-
\lambda r^2
\big)\,
zw^2
+
\big(
-2i\,\lambda\alpha^2\overline{\alpha}
-
\lambda\alpha r
\big)\,
\zeta w^2,
\endaligned
\]
\[
\footnotesize
\aligned
\zeta'
&
\,=\,
2i\,
\frac{\lambda}{\overline{\lambda}}\,
\overline{\alpha}\,z
+
\frac{\lambda}{\overline{\lambda}}\,
\zeta
\\
&
\ \ \ \ \
+
\Big(
-i\,
\frac{\lambda}{\overline{\lambda}}\,
r
+
3\,
\frac{\lambda}{\overline{\lambda}}\,
\alpha\overline{\alpha}
\Big)\,
z^2
-
2i\,
\frac{\lambda}{\overline{\lambda}}\,
\alpha\,z\zeta
+
\frac{\lambda}{\overline{\lambda}}\,
\overline{\alpha}^2\,w
\\
&
\ \ \ \ \
+
\Big(
-4i\,
\frac{\lambda}{\overline{\lambda}}\,
\alpha^2\overline{\alpha}
-
2\,
\frac{\lambda}{\overline{\lambda}}\,
\alpha r
\Big)\,
z^3
-
3\,
\frac{\lambda}{\overline{\lambda}}\,
\alpha^2\,
z^2\zeta
+
\Big(
-4i\,
\frac{\lambda}{\overline{\lambda}}\,
\alpha\overline{\alpha}^2
-
2\,
\frac{\lambda}{\overline{\lambda}}\,
\overline{\alpha} r
\Big)\,
zw
-
2\,
\frac{\lambda}{\overline{\lambda}}\,
\alpha\overline{\alpha}\,
\zeta w
\\
&
\ \ \ \ \
+
\Big(
-5\,
\frac{\lambda}{\overline{\lambda}}\,
\alpha^3\overline{\alpha}
+
3i\,
\frac{\lambda}{\overline{\lambda}}\,
\alpha^2 r
\Big)\,
z^4
+
4i\,
\frac{\lambda}{\overline{\lambda}}\,
\alpha^3\,z^3\zeta
+
\Big(
-10\,
\frac{\lambda}{\overline{\lambda}}\,
\alpha^2\overline{\alpha}^2
+
8i\,
\frac{\lambda}{\overline{\lambda}}\,
\alpha\overline{\alpha}r
+
\frac{\lambda}{\overline{\lambda}}\,
r^2
\Big)\,
z^2w
\\
&
\ \ \ \ \ \ \ \ \ \ \ \ \ \ \ \ \ \ \ \ \ \ \ \ \ \ \ \ \ \ \ \ \ \ \ 
\ \ \ \ \ \ \ \ \ \ \ \ \ \ \ \ \ \ \ \ \ \ \ \ \ \ \ \ \ \ \ \ \ \ \ 
+
\Big(
8i\,
\frac{\lambda}{\overline{\lambda}}\,
\alpha^2\overline{\alpha}
+
2\,
\frac{\lambda}{\overline{\lambda}}\,
\alpha r
\Big)\,
z\zeta w
+
\frac{\lambda}{\overline{\lambda}}\,
\alpha^2\,\zeta^2w
+
\Big(
i\,
\frac{\lambda}{\overline{\lambda}}\,
\overline{\alpha}^2r
-
\frac{\lambda}{\overline{\lambda}}\,
\alpha\overline{\alpha}^3
\Big)\,
w^2,
\endaligned
\]
\[
\footnotesize
\aligned
w'
&
\,=\,
0
\\
&
\ \ \ \ \
+
\lambda\overline{\lambda}\,
w
\\
&
\ \ \ \ \
-\,2i\,\lambda\overline{\lambda}\,\alpha\,
zw
\\
&
\ \ \ \ \
-\,
3\,\lambda\overline{\lambda}\alpha^2\,
z^2w
+
\big(
i\,\lambda\overline{\lambda}r
-
\lambda\overline{\lambda}\alpha\overline{\alpha}
\big)\,
w^2
\\
&
\ \ \ \ \
+
4i\,\lambda\overline{\lambda}\alpha^3\,
z^3
+
\big(
4i\,\lambda\overline{\lambda}\alpha^2\overline{\alpha}
+
4\,\lambda\overline{\lambda}\,\alpha r
\big)\,
zw^2
+
\lambda\overline{\lambda}\alpha^2\,
\zeta w^2
\\
&
\ \ \ \ \
+
5\,\lambda\overline{\lambda}\alpha^4\,
z^4w
+
\big(
10\,\lambda\overline{\lambda}\alpha^3\overline{\alpha}
-
10i\,\lambda\overline{\lambda}\,\alpha^2r
\big)\,
z^2w^2
-
4i\,\lambda\overline{\lambda}\alpha^3\,
z\zeta w^2
+
\big(
-\lambda\overline{\lambda}r^2
-
2i\,\lambda\overline{\lambda}\alpha\overline{\alpha}r
+
\lambda\overline{\lambda}\alpha^2\overline{\alpha}^2
\big)\,
w^3.
\endaligned
\]

\Section{\bf Lie Jet Theory}
\label{Lie-jet-theory}
\HEAD{{\ref{Lie-jet-theory}}.~{\sf Lie Jet Theory}
}{
Wei-Guo {\sc Foo}, Joël {\sc Merker}, The-Anh {\sc Ta}}

To apply Lie's theory similarly as in~{\cite{Merker-2020}}, we must
work with the five {\em intrinsic}, {\em real}, coordinates
$(x,y,s,t,v)$ on $M^5$, where:
\[
z
\,=\,
x+i\,y,
\ \ \ \ \ \ \ \ \ \ \ \ \ \ \ \ \ \ \ \
\zeta
\,=\,
s+i\,t,
\ \ \ \ \ \ \ \ \ \ \ \ \ \ \ \ \ \ \ \
w
\,=\,
u+i\,v.
\]

As in~{\cite{Merker-2020}}, we consider 
parametrized local real $\mathcal{C}^\omega$ 
curves passing by the origin
\[
\tau
\,\,\longmapsto\,\,
\big(x(\tau),y(\tau),s(\tau),t(\tau),\,\tau\big).
\]
with $v(\tau) \equiv \tau$ guaranteeing that the curve is  
{\em not} CR-tangential. We then use the parameter-letter $v$ instead
of $\tau$. 

The eight independent coordinates corresponding
to $\dot{x}(v)$, $\dot{y}(v)$, $\dot{s}(v)$, $\dot{t}(v)$,
$\ddot{x}(v)$, $\ddot{y}(v)$, $\ddot{s}(v)$, $\ddot{t}(v)$
will be denoted:
\[
\big(
v,\,x,y,s,t,\,
x_1,y_1,s_1,t_1,\,
x_2,y_2,s_2,t_2
\big).
\]
The first jet space is $J_{1,4}^1 \equiv \R^{1+4+4}$, and
the second jet space is $J_{1,4}^2 \equiv \R^{1+4+4+4}$.

Any diffeomorphism $(v,x,y,s,t) \longmapsto (v',x',y',s',t')$ 
lifts to jet
spaces of any order. 
Because the formulas rapidly become complicated
{\cite{Olver-1995, Merker-2008, Chen-Merker-2019}},
Lie linearized the action of diffeomorphisms.

As in~{\cite{Merker-2020}}, we will 
{\em apply} Lie's formulas. Start from a general vector field:
\[
\vec{\bf v}
\,:=\,
\xi(v,x,y,s,t)\,
\frac{\partial}{\partial v}
+
\varphi(v,x,y,s,t)\,
\frac{\partial}{\partial x}
+
\psi(v,x,y,s,t)\,
\frac{\partial}{\partial y}
+
\lambda(v,x,y,s,t)\,
\frac{\partial}{\partial s}
+
\mu(v,x,y,s,t)\,
\frac{\partial}{\partial t}.
\]
Introduce the {\sl total differentiation
operator:}
\[
\!\!\!\!\!\!\!\!\!\!\!\!\!\!\!
{\sf D}_v
\,:=\,
\frac{\partial}{\partial v}
+
x_1\,
\frac{\partial}{\partial x}
+
y_1\,
\frac{\partial}{\partial y}
+
s_1\,
\frac{\partial}{\partial s}
+
t_1\,
\frac{\partial}{\partial t}
+
x_2\,
\frac{\partial}{\partial x_1}
+
y_2\,
\frac{\partial}{\partial y_1}
+
s_2\,
\frac{\partial}{\partial s_1}
+
t_2\,
\frac{\partial}{\partial t_1}
+
x_3\,
\frac{\partial}{\partial x_2}
+
y_3\,
\frac{\partial}{\partial y_2}
+
s_3\,
\frac{\partial}{\partial s_2}
+
t_3\,
\frac{\partial}{\partial t_2}.
\]
Then the second prolongation of $\vec{\bf v}$:
\[
\aligned
\vec{\bf v}^{(2)}
&
\,=\,
\vec{\bf v}
+
\varphi_1\,
\frac{\partial}{\partial x_1}
+
\psi_1\,
\frac{\partial}{\partial y_1}
+
\lambda_1\,
\frac{\partial}{\partial s_1}
+
\mu_1\,
\frac{\partial}{\partial t_1}
\\
&
\ \ \ \ \ \ \ \ \ 
+
\varphi_2\,
\frac{\partial}{\partial x_2}
+
\psi_2\,
\frac{\partial}{\partial y_2}
+
\lambda_2\,
\frac{\partial}{\partial s_2}
+
\mu_2\,
\frac{\partial}{\partial t_2},
\endaligned
\]
has coefficients
({\cite{Lie-Merker-2015,
Olver-1995, Merker-2008, Chen-Merker-2019}}):
\[
\!\!\!\!\!\!\!\!\!\!\!\!\!\!\!\!\!\!\!\!\!\!\!\!\!\!\!\!\!\!
\footnotesize
\aligned
\varphi_1
&
\,:=\,
{\sf D}_v
\big(
\varphi
-
\xi\,x_1
\big)
+
\xi\,
x_2,
\ \ \ \ \
\psi_1
\,:=\,
{\sf D}_v
\big(
\psi
-
\xi\,y_1
\big)
+
\xi\,
y_2,
\ \ \ \ \
\lambda_1
\,:=\,
{\sf D}_v
\big(
\lambda
-
\xi\,s_1
\big)
+
\xi\,
s_2,
\ \ \ \ \
\mu_1
\,:=\,
{\sf D}_v
\big(
\mu
-
\xi\,t_1
\big)
+
\xi\,
t_2,
\\
\varphi_2
&
\,:=\,
{\sf D}_v{\sf D}_v
\big(
\varphi
-
\xi\,x_1
\big)
+
\xi\,
x_3,
\ \ \ \ \
\psi_2
\,:=\,
{\sf D}_v{\sf D}_v
\big(
\psi
-
\xi\,y_1
\big)
+
\xi\,
y_3,
\ \ \ \ \
\lambda_2
\,:=\,
{\sf D}_v{\sf D}_v
\big(
\lambda
-
\xi\,s_1
\big)
+
\xi\,
s_3,
\ \ \ \ \
\mu_2
\,:=\,
{\sf D}_v{\sf D}_v
\big(
\mu
-
\xi\,t_1
\big)
+
\xi\,
t_3.
\endaligned
\]
\Section{\bf Intrinsic Isotropy Automorphisms of the Gaussier-Merker
Model}
\label{intrinsic-automorphisms-GM-model}
\HEAD{{\ref{intrinsic-automorphisms-GM-model}}.~{\sf Intrinsic 
Isotropy Automorphisms of the Gaussier-Merker Model}
}{
Wei-Guo {\sc Foo}, Joël {\sc Merker}, The-Anh {\sc Ta}}

We want to apply Lie's prolongation formulas within the {\em first}
jet space to our $5$ vector fields ${\sf X} = {\sf D}$, ${\sf R}$,
${\sf I}_1$, ${\sf I}_2$, ${\sf J}$.  But these {\em holomorphic}
$(1,0)$ fields were {\em extrinsic}, defined in $\C^3$.  We
must therefore write up the five fields ${\sf X} + \overline{\sf X}$
in the {\em intrinsic} coordinates $(x,y,s,t,v) \in M_{\sf LC}^5$.  By
slight abuse, we keep the notation ${\sf X}$ instead of ${\sf X} +
\overline{\sf X}$:
\[
\aligned
{\sf D}
&
\,=\,
x\,\partial_x
+
y\,\partial_y
+
2v\,\partial_v,
\\
{\sf R}
&
\,=\,
-\,y\,\partial_x
+
x\,\partial_y
-
2t\,\partial_s
+
2s\,\partial_t,
\endaligned
\]
\[
\footnotesize
\aligned
{\sf I}_1
&
\,=\,
\Big[
\frac{
2\,x^2s^2-2\,y^2s^2+2\,y^2+2\,xyt+x^2t^2-y^2t^2
+2\,xyst-tv+s^2tv+t^3v+2\,x^2s
}{
-1+s^2+t^2}
\Big]\,
\frac{\partial}{\partial x}
\\
&
\ \ \ \ \
+
\Big[
\frac{-\,y^2t-x^2st-v+s^2v+t^2v-sv+s^3v+st^2v-x^2t-4\,xyt^2
+y^2st+2\,xy-2\,xys^2}
{1-s^2+t^2}
\Big]\,
\frac{\partial}{\partial y}
\\
&
\ \ \ \ \
+
\big[
2\,x-2\,yt+2\,xs
\big]\,
\frac{\partial}{\partial s}
+
\big[
2\,y+2\,ys+2\,xt
\big]\,
\frac{\partial}{\partial t}
\\
&
\ \ \ \ \
+
\Big[
\frac{
-\,4\,xy^2t-2\,x^2ys-2\,x^2y+2\,y^3s-2\,xv+2\,xs^2v
+2\,xt^2v-2\,y^3
}{-1+s^2+t^2}
\Big]\,
\frac{\partial}{\partial v},
\endaligned
\]
\[
\footnotesize
\aligned
\\
{\sf I}_2
&
\,=\,
\Big[
\frac{
-\,y^2t-x^2st-4\,xyt^2+y^2st-sv+s^3v+st^2v+2\,xy
-2\,xys^2+v-s^2v-t^2v-x^2t
}{
-1+s^2+t^2}
\Big]\,
\frac{\partial}{\partial x}
\\
&
\ \ \ \ \
+
\Big[
\frac{
-\,2\,x^2+2\,x^2s^2+x^2t^2-2\,xyt-2\,y^2s^2-y^2t^2
+2\,xyst-tv+s^2tv+t^3v+2\,y^2s
}{1-s^2+t^2}
\Big]\,
\frac{\partial}{\partial y}
\\
&
\ \ \ \ \
+
\big[
2\,xt-2\,y+2\,ys
\big]\,
\frac{\partial}{\partial s}
+
\big[
-\,2\,xs+2\,x+2\,yt
\big]\,
\frac{\partial}{\partial t}
\\
&
\ \ \ \ \
+
\Big[
\frac{
-\,2\,xy^2s+2\,xy^2+4\,x^2yt+2\,x^3s+2\,x^3
-2\,yv+2\,ys^2v+2\,yt^2v
}{-1+s^2+t^2}
\Big]\,
\frac{\partial}{\partial v},
\endaligned
\]
\[
\footnotesize
\aligned
{\sf J}
&
\,=\,
\Big[
\frac{
-\,2\,xy^2t-x^2ys-x^2y+y^3s-xv+xs^2v+xt^2v-y^3
}{
-1+s^2+t^2}
\Big]\,
\frac{\partial}{\partial x}
\\
&
\ \ \ \ \
+
\Big[
\frac{
-\,xy^2s+xy^2+2\,x^2yt+x^3s+x^3-yv+ys^2v+yt^2v
}{1-s^2+t^2}
\Big]\,
\frac{\partial}{\partial y}
\\
&
\ \ \ \ \
+
\big[
2\,xy
\big]\,
\frac{\partial}{\partial s}
+
\big[
-\,x^2+y^2
\big]\,
\frac{\partial}{\partial t}
\\
&
\ \ \ \ \
+
\Big[
\frac{
\big(
v-s^2v-t^2v-x^2s-x^2-2\,xyt+y^2s-y^2
\big)
\big(
-\,v+s^2v+t^2v-x^2-x^2s-2\,xyt+y^2s-y^2
\big)
}{\big(1-s^2-t^2\big)^2}
\Big]\,
\frac{\partial}{\partial v}.
\endaligned
\] 

\Section{\bf Prolongation to the Jet Space of Order $1$}
\label{prolongations-jet-1}
\HEAD{{\ref{prolongations-jet-1}}.~{\sf Prolongation
to the Jet Space of Order $1$}
}{
Wei-Guo {\sc Foo}, Joël {\sc Merker}, The-Anh {\sc Ta}}

As said, we work above the origin $0 \in M_{\sf LC}$. 

\begin{center}
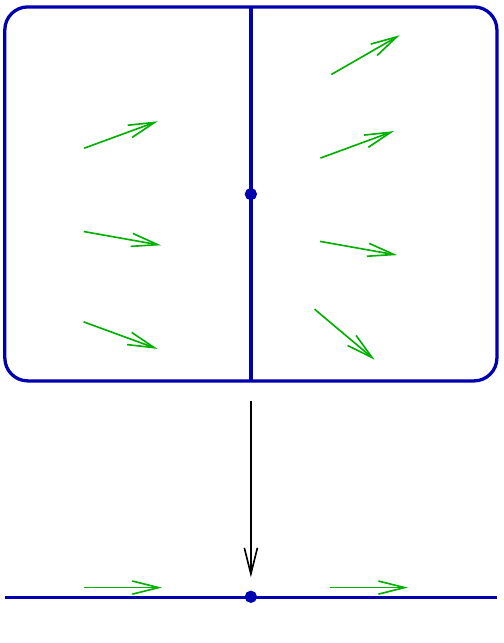
\end{center}

By Lie's theory, any vector field $\vec{\bf v}$ on the base $M$ lifts
as a vector field $\vec{\bf v}^{(1)}$ on the first jet space
$J_{1,4}^1 = \R^{1+4+4}$.

Because our five intrinsic vector fields ${\sf D}$, ${\sf R}$, ${\sf
I}_1$, ${\sf I}_2$, ${\sf J}$ vanish at $v = x = y = s = t= 0$,
their prolongations will automatically be tangent to the fiber $\big\{
(0, 0, 0, 0, 0, x_1, y_1, s_1, t_1) \big\}$ above 
$(0,0,0,0,0)$ in the first jet space.

Lie's formulas yield the very simple values of these first
prolongations above the origin
$v=x=y=s=t=0$:
\[
\def\arraystretch{1.25}
\begin{array}{ccccc}
& \partial_{x_1} & \partial_{y_1} & \partial_{s_1} & \partial_{t_1}
\\
{\sf D}^{(1)} & -x_1 & -y_1 & -2s_1 & -2t_1
\\
{\sf R}^{(1)} & -y_1 & x_1 & -2t_1 & 2s_1
\\
{\sf I}_1^{(1)} & 0 & -1 & 2x_1 & 2y_1
\\
{\sf I}_2^{(1)} & 1 & 0 & -2y_1 & 2y_1
\\
{\sf J}^{(1)} & 0 & 0 & 0 & 0
\end{array}
\]

\begin{center}
\includegraphics[scale=0.50]{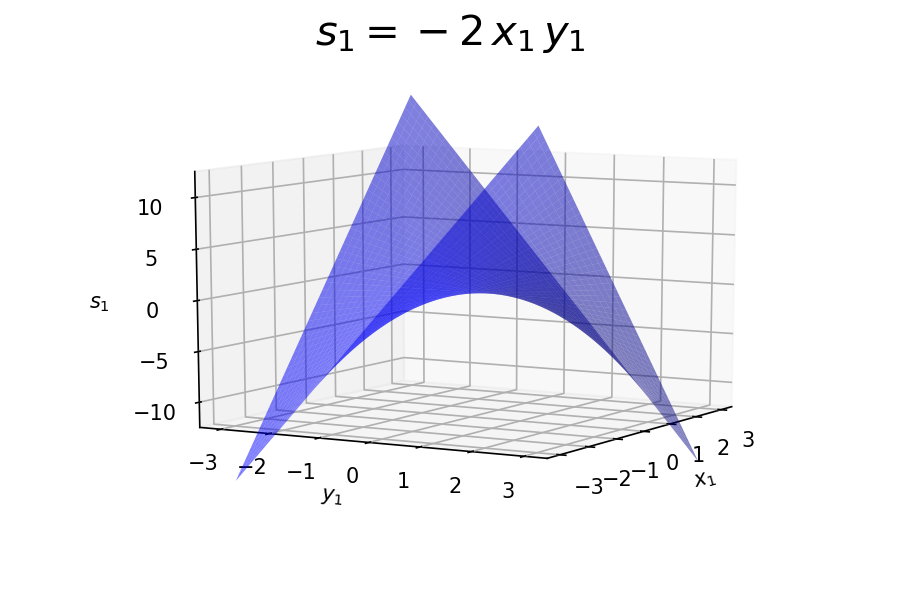}
\!\!\!\!\!\!\!\!\!\!
\includegraphics[scale=0.50]{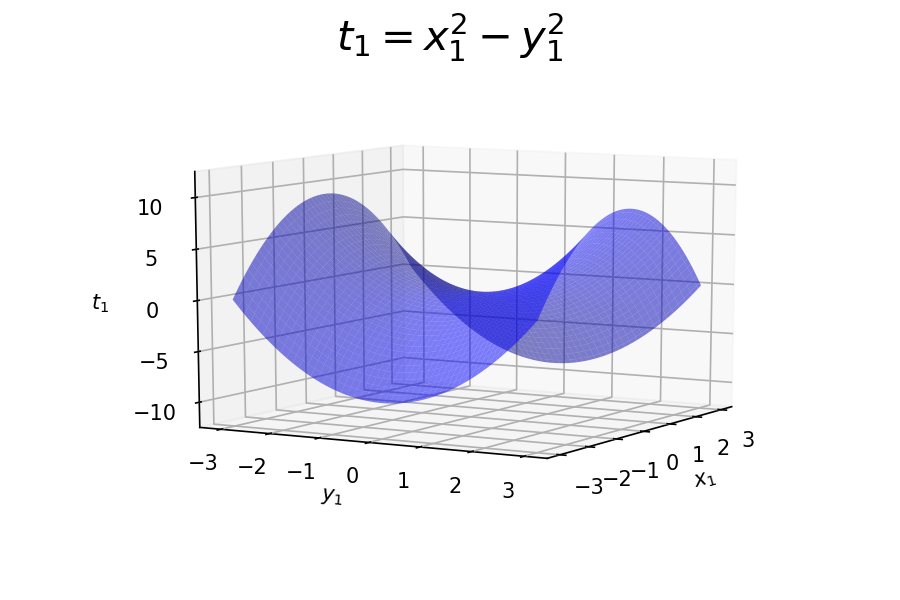}
\end{center}

\begin{Observation}
\label{Obs-Sigma-0-jets-1}
On $\R^4 = \R_{x_1, y_1, s_1, t_1}^4$, 
there exists a unique 
$\big\{ {\sf D}^{(1)}, {\sf R}^{(1)},
{\sf I}_1^{(1)}, {\sf I}_2^{(1)}, {\sf J}^{(1)} 
\big\}$-invariant $2$-dimensional submanifold
$\Sigma_0^1 \subset \R^4$, algebraic, graphed as:
\[
\left[
\aligned
s_1
&
\,=\,
-\,2\,x_1y_1,
\\
t_1
&
\,=\,
x_1^2-y_1^2,
\endaligned\right.
\]
Moreover, the complement $\R^4 \backslash \Sigma_0^1$ is a unique
(transitive)
orbit under ${\sf D}^{(1)}$, ${\sf R}^{(1)}$,
${\sf I}_1^{(1)}$, ${\sf I}_2^{(1)}$, ${\sf J}^{(1)}$.
\end{Observation}

\proof
We can drop the fifth line of ${\sf J}^{(1)}$
containing only zeros. With
$a_1$ and $b_1$ being parameters, any point of $\R^4$
can be written as $(x_1, y_1, s_1, t_1)$ with:
\[
s_1
\,:=\,
-\,2\,x_1y_1
+
a_1,
\ \ \ \ \ \ \ \ \ \ \ \ \ \ \ \ \ \ \ \
t_1
\,:=\,
x_1^2-y_1^2
+
b_1.
\]
Then replacing $s_1$ and $t_1$:
\[
\left(\!
\begin{array}{cccc}
-x_1 & -y_1 & -2s_1 & -2t_1
\\
-y_1 & x_1 & -2t_1 & 2s_1
\\
0 & -1 & 2x_1 & 2y_1
\\
1 & 0 & -2y_1 & 2x_1
\end{array}
\!\right)
\ \ \ \ \
\xrightarrow[{\rule[0pt]{50pt}{0pt}}]{\text{\sf Gauss-pivot}}
\ \ \ \ \
\left(\!
\begin{array}{cccc}
0 & 0 & -2a_1 & -2b_1
\\
0 & 0 & -2b_1 & 2a_1
\\
0 & -1 & 2x_1 & 2y_1
\\
1 & 0 & -2y_1 & 2x_1
\end{array}
\!\right).
\]
This matrix has determinant $-4a_1^2 - 4b_1^2$, hence is of
rank $4$ when $(a_1, b_1) \neq (0, 0)$. In the
corresponding locus, namely in $\R^4 \big\backslash \Sigma_0^1$,
the five prolonged vector fields
${\sf D}^{(1)}$, ${\sf R}^{(1)}$, ${\sf I}_1^{(1)}$,
${\sf I}_2^{(1)}$, ${\sf J}^{(1)}$ have everywhere rank $4$, 
hence generate locally open orbits, so that 
$\R^4 \big\backslash \Sigma_0^1$ is a single 
orbit under their action.
 
When $a_1 = b_1 = 0$, the above matrix has rank $2$. 
In this $2$-dimensional graphed locus, the rank of
${\sf D}^{(1)}$, ${\sf R}^{(1)}$, ${\sf I}_1^{(1)}$,
${\sf I}_2^{(1)}$, $J^{(1)}$ is everywhere equal to $2$,
whence $\Sigma_0^1$ is a single orbit under their action.
\endproof

Thus, the model $M_{\sf LC}$ has an invariant cone:
\[
s_1
+
i\,t_1
\,=\,
i\,
\big(
x_1
+
i\,y_1
\big)^2,
\]
namely a cone invariant under the action of
${\sf D}^{(1)}$, ${\sf R}^{(1)}$, ${\sf I}_1^{(1)}$,
${\sf I}_2^{(1)}$, ${\sf J}^{(1)}$.
Soon, we will see that {\em every} $M^5 \subset \C^3$ in the class
$\mathfrak{C}_{2,1}$ also possesses an invariant cone
at {\em any} of its points $p \in M^5$.

\Section{\bf Prolongation to the Jet Space of Order $2$}
\label{prolongations-jet-2}
\HEAD{{\ref{prolongations-jet-2}}.~{\sf Prolongation
to the Jet Space of Order $2$}
}{
Wei-Guo {\sc Foo}, Joël {\sc Merker}, The-Anh {\sc Ta}}

Next, we increment the jet order by one unit.
The second order Lie prolongations
${\sf D}^{(2)}$, ${\sf R}^{(2)}$, ${\sf I}_1^{(2)}$, 
${\sf I}_2^{(2)}$, ${\sf J}^{(2)}$
have the following coefficients above the origin, 
$v = x = y = s = t = 0$:
\[
\def\arraystretch{1.25}
\begin{array}{ccccccccc}
& \partial_{x_1} & \partial_{y_1} & \partial_{s_1} & \partial_{t_1}
& \partial_{x_2} & \partial_{y_2} & \partial_{s_2} & \partial_{t_2}
\\
{\sf D}^{(2)} & -x_1 & -y_1 & -2s_1 & -2t_1 &
-3x_2 & -3y_2 & -4s_2 & -4t_2
\\
{\sf R}^{(2)} & -y_1 & x_1 & -2t_1 & 2s_1 &
-y_2 & x_2 & -2t_2 & 2s_2
\\
{\sf I}_1^{(2)} & 0 & -1 & 2x_1 & 2y_1 &
2t_1-4x_1^2-4y_1^2 & -2s_1 & 2x_2-4y_1t_1 & 2y_2+4x_1s_1
\\
{\sf I}_2^{(2)} & 1 & 0 & -2y_1 & 2y_1 &
-2s_1 & -2t_1-4x_1^2-4y_1^2 & -2y_2+4x_1t_1 & 2x_2-4x_1s_1
\\
{\sf J}^{(2)} & 0 & 0 & 0 & 0 &
0 & 0 & 2s_1+4x_1y_1 & 2t_1-2x_1^2+2y_1^2
\end{array}
\]

Of course, we pull this matrix back to $\Sigma_0^1$,
hence the last line becomes null. Keeping only the
first $4$ lines, and performing a Gauss pivot, we get:
\[
\footnotesize
\aligned
\left(\!
\def\arraystretch{1.25}
\begin{array}{ccccccccc}
0 & 0 & 0 & 0 &
6x_1^2y_1+6y_1^3-3x_2 & -6x_1y_1^2-6x_1^3-3y_2 &
\underset{-2x_1y_2+4x_1^4-4s_2}{{\scriptstyle{-2x_2y_1-4y_1^4}}} & 
\underset{+2x_1x_2+8x_1^3y_1-4t_2}{{\scriptstyle{-2y_1y_2+8x_1y_1^3}}}
\\
0 & 0 & 0 & 0 &
-2x_1^3-2x_1y_1^2-y_2 & -2x_1^2y_1-2y_1^3+x_2 &
2x_1x_2-2y_1y_2-2t_2 & 2x_1y_2+2y_1x_2+2s_2
\\
0 & -1 & 2x_1 & 2y_1 & 
-2x_1^2-6y_1^2 & 4x_1y_1 & 2x_2-4x_1^2y_1+4y_1^3 &
2y_2-8x_1y_1^2
\\
1 & 0 & -2y_1 & 2x_1 & 4x_1y_1 & -6x_1^2-2y_1^2 &
-2y_2+4x_1^3-4x_1y_1^2 & 2x_2+8x_1^2y_1
\end{array}
\!\right).
\endaligned
\]
The upper $2 \times 4$ block, having $8$ entries,
then shows that $x_2$, $y_2$, $s_2$, $t_2$ can be 
uniquely and consistently defined in terms of $x_1$, $y_1$,
so that they define an invariant surface under the action of
${\sf D}^{(2)}$, ${\sf R}^{(2)}$, ${\sf I}_1^{(2)}$,
${\sf I}_2^{(2)}$, ${\sf J}^{(2)}$.

\begin{center}
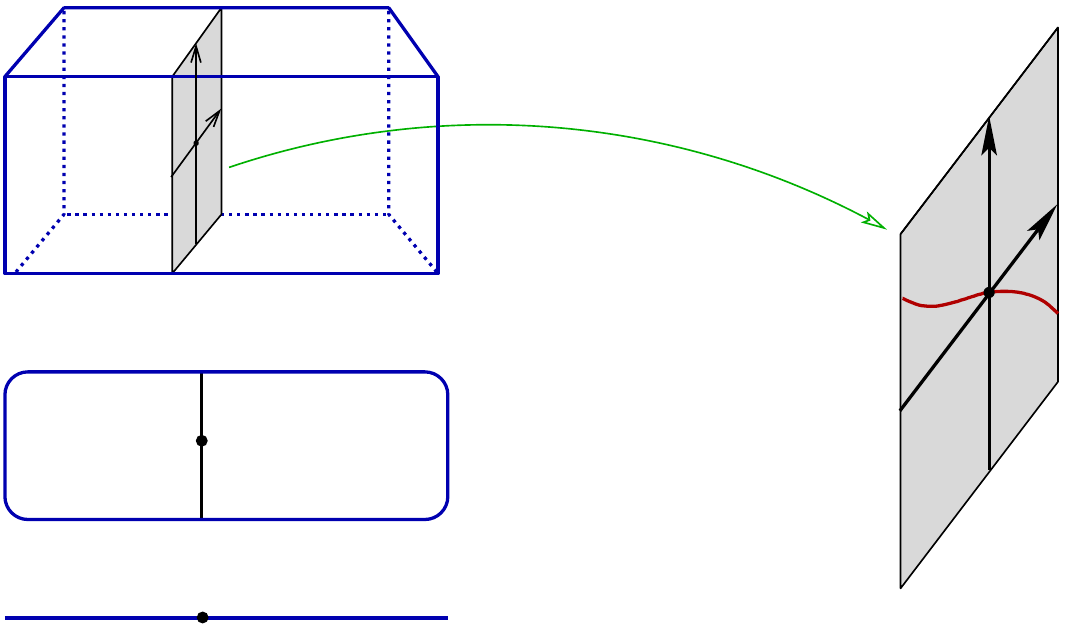
\end{center}

\begin{Observation}
\label{Obs-Sigma-0-jets-2}
On $\R^8 = \R_{x_1, y_1, s_1, t_1}^4 \times 
\R_{x_2, y_2, s_2, t_2}^4$, there
exists a unique $\big\{ {\sf D}^{(2)}, {\sf R}^{(2)},
{\sf I}_1^{(2)}, {\sf I}_2^{(2)}, {\sf J}^{(2)} 
\big\}$-invariant
$2$-dimensional submanifold
$\Sigma_0^2 \subset \R^8$, algebraic, graphed as:
\[
\left[
\aligned
s_1
&
\,=\,
-\,2\,x_1y_1,
\\
t_1
&
\,=\,
x_1^2-y_1^2,
\endaligned\right.
\ \ \ \ \ \ \ \ \ \ \ \ \ \ \ \ \ \ \ \ \ \ \ \ \ \
\left[
\aligned
x_2
&
\,=\,
2\,x_1^2y_1
+
2\,y_1^3,
\\
y_2
&
\,=\,
-\,2\,x_1^3
-
2\,x_1y_1^2,
\\
s_2
&
\,=\,
-\,2\,y_1^4+2\,x_1^4,
\\
t_2
&
\,=\,
4\,x_1^3y_1
+
4\,x_1y_1^3.
\endaligned\right.
\]
Moreover, the complement $\R^8 \backslash \Sigma_0^2$ is a unique
orbit under the transitive action of ${\sf D}^{(2)}$, ${\sf R}^{(2)}$,
${\sf I}_1^{(2)}$, ${\sf I}_2^{(2)}$, ${\sf J}^{(2)}$.
\end{Observation}

\[
\xymatrix{
J_{1,4}^2
\ar[d]
&
\supset
&
\Sigma_0^2
\ar[d]
\\
J_{1,4}^1
\ar[d]
&
\supset
&
\Sigma_0^1
\ar[d]
\\
M
&
\ni
&
0.
}
\]

\proof
As said, we pull everything back to $\Sigma_0^1$ having equations
$s_1 = - 2 x_1y_1$, $t_1 = x_1^2 - y_1^2$. 
With $a_2$, $b_2$, $c_2$, $d_2$ being parameters, any point
of $\R_{x_2, y_2, s_2, t_2}^4$ can be written as:
\[
\left[
\aligned
x_2
&
\,=\,
2\,x_1^2y_1
+
2\,y_1^3
+
a_2,
\\
y_2
&
\,=\,
-\,2\,x_1^3
-
2\,x_1y_1^2
+
b_2,
\endaligned\right.
\ \ \ \ \ \ \ \ \ \ \ \ \ \ \ \ \ \ \ \ \ \ \ \ \ \
\left[
\aligned
s_2
&
\,=\,
-\,2\,y_1^4+2\,x_1^4
+
c_2,
\\
t_2
&
\,=\,
4\,x_1^3y_1
+
4\,x_1y_1^3
+
d_2.
\endaligned\right.
\]

Replacing $x_2$, $y_2$ without replacing $s_2$, $t_2$,
the upper right $2 \times 4$ block becomes:
\[
\left(\!
\begin{array}{cccc}
-3a_2 & -3b_2 & -8y_1^4+8x_1^4-4s_2-2y_1a_2-2x_1b_2 & 
16x_1^3y_1+16x_1y_1^3-4t_2-2y_1b_2+2x_1a_2
\\
-b_2 & a_2 & 8x_1^3y_1+8x_1y_1^3-2t_2+2x_1a_2-2y_1b_2 &
-4x_1^4+4y_1^4+2s_2+2x_1b_2+2y_1a_2
\end{array}
\!\right).
\]
Visibly, it is of rank $2$ whenever $(a_2, b_2) \neq (0, 0)$.

Thus, put in it $a_2 := 0$ and $b_2 := 0$:
\[
\left(\!
\begin{array}{cccc}
0 & 0 & -8y_1^4+8x_1^4-4s_2 & 
16x_1^3y_1+16x_1y_1^3-4t_2
\\
0 & 0 & 8x_1^3y_1+8x_1y_1^3-2t_2 &
-4x_1^4+4y_1^4+2s_2
\end{array}
\!\right),
\]
and now replace $s_2$, $t_2$, to get:
\[
\left(\!
\begin{array}{cccc}
0 & 0 & -4c_2 & -4d_2
\\
0 & 0 & -2c_2 & 2d_2
\end{array}
\!\right),
\]
a submatrix which has maximal rank $2$ 
if and only if $(c_2, d_2) \neq (0, 0)$.
This concludes.
\endproof

We have therefore shown that, to every (fixed) $1$-jet
at the origin $0 \in M_{\sf LC}$
of the form:
\[
j_0^1 
\,=\, 
\big(
x_1,y_1,\,
-2x_1y_1,\,
x_1^2-y_1^2
\big)
\] 
is associated a unique
second order jet at the origin:
\[
j_0^2
\,=\,
\Big(
x_1,y_1,\,
-2x_1y_1,\,
x_1^2-y_1^2,\,\,
2x_1^2y_1+2y_1^3,\,
-2x_1^3-2x_1y_1^2,\,
-2y_1^4+2x_1^4,\,
4x_1^3y_1+4x_1y_1^3
\Big),
\]
and since $\Sigma_0^2$ is invariant under the
action of the stability group
of the Gaussier-Merker model, 
this association is invariant.

\begin{center}
\includegraphics[scale=0.50]{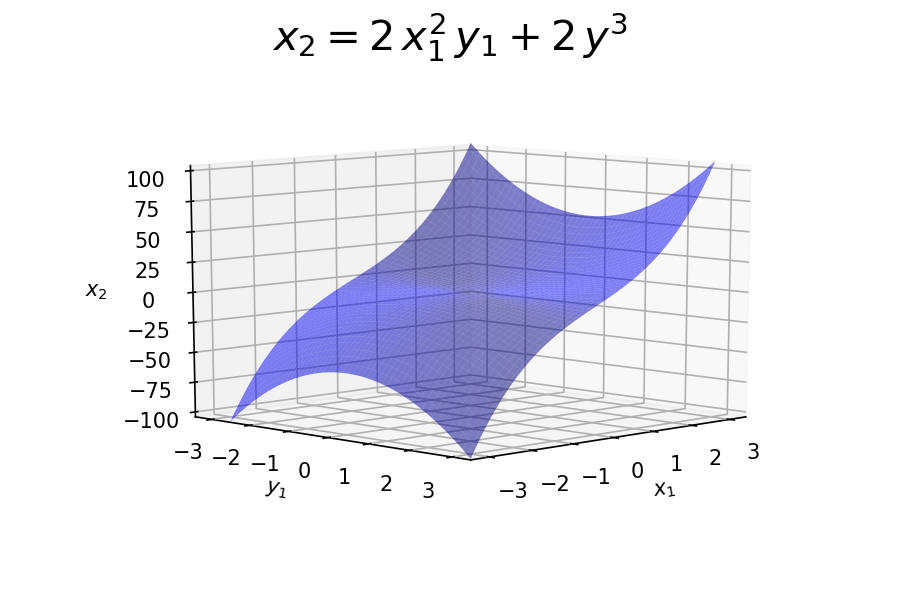}
\!\!\!\!\!\!\!\!\!\!
\includegraphics[scale=0.50]{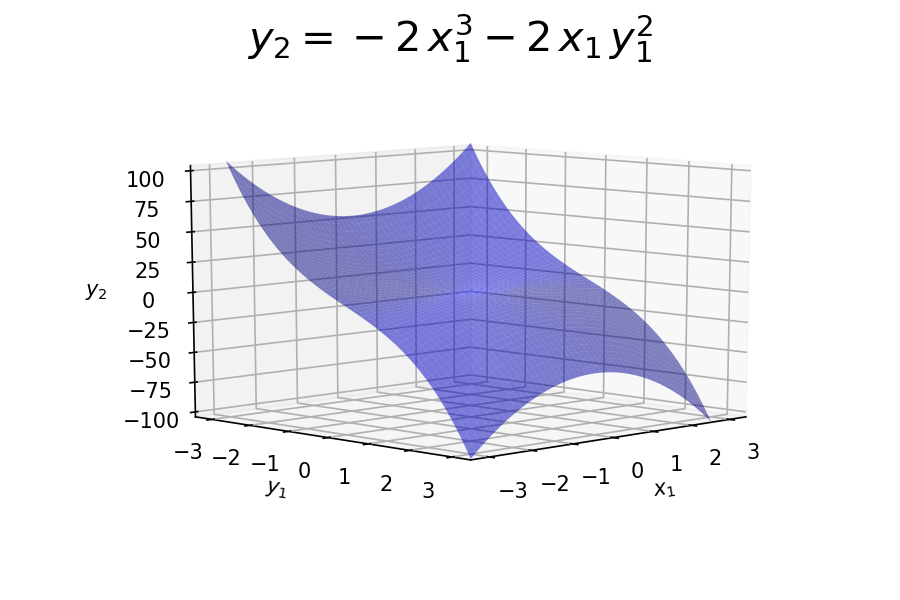}
\end{center}

\begin{center}
\includegraphics[scale=0.50]{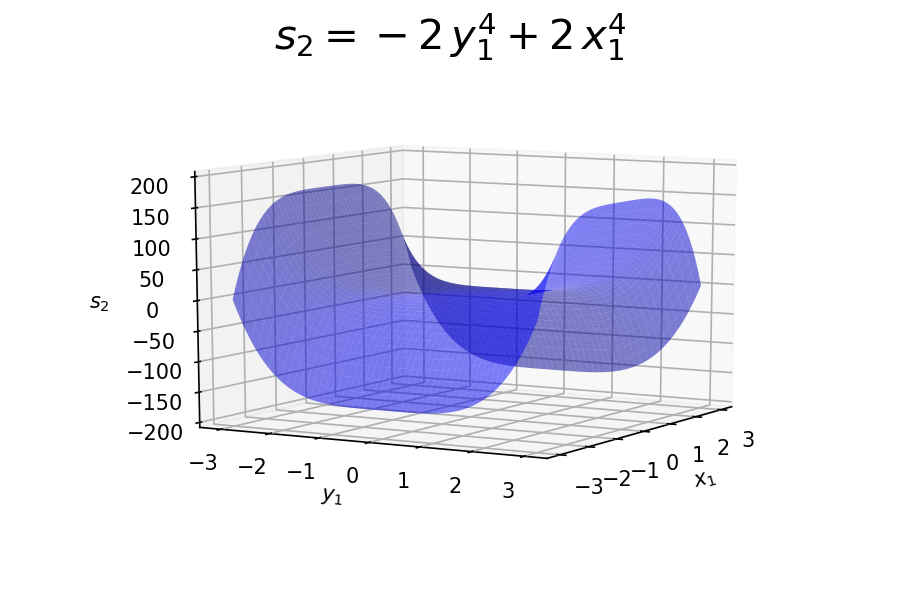}
\!\!\!\!\!\!\!\!\!\!
\includegraphics[scale=0.50]{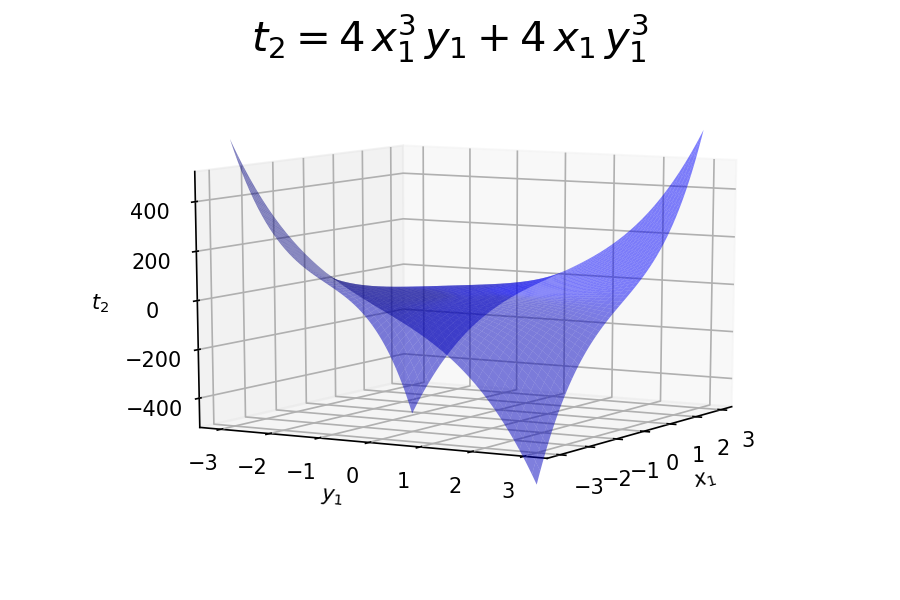}
\end{center}

Our next goal will be to transfer this invariancy property
to {\em any} $M^5 \in \mathfrak{C}_{2,1}$. But subtleties 
will spice up our job.

\Section{\bf Road Map to Convergent Normal Form}
\label{road-map-convergent-normal-form}
\HEAD{{\ref{road-map-convergent-normal-form}}.~{\sf Road Map to 
Convergent Normal Form}
}{
Wei-Guo {\sc Foo}, Joël {\sc Merker}, The-Anh {\sc Ta}}

A certain {\sl Lie-theoretic}
construction of Cartan-Moser chains for
Levi nondegenerate hypersurfaces $M^3 \subset \C^2$ was set up
in~{\cite{Merker-2020}} in order to be imitated 
when studying hypersurfaces
$M^5 \subset \C^3$ in the class $\mathfrak{C}_{2,1}$, 
in the present memoir.
However, we will encounter not only analogies,
but also differences.

Recall that any Levi nondegenerate $M^3 \subset \C^2$, taken
at any point $p \in M$, can be brought, in local coordinates
$(z, w = u + iv)$ vanishing at $p$, to the preliminary normal
form~{\cite[Prp.~2.2]{Merker-2020}}:
\[
v
\,=\,
z\overline{z}
+
{\rm O}(6),
\]
where 
the remainder is {\sl weighted} according to
$[z] := 1$, $[w] := 2$.
Furthermore, the {\sl ambiguity} of such a
punctual preliminary normalization, namely {\em any}
map:
\[
z'
\,=\,
f_1+f_2+f_3+f_4
+
{\rm O}(5),
\ \ \ \ \ \ \ \ \ \ \ \ \ \ \ \ \ \ \ \ \ \ \ \ \ \
w'
\,=\,
g_1+g_2+g_3+g_4+g_5
+
{\rm O}(6),
\]
which preserves this normalization, {\em i.e.} which sends
$v = z\overline{z} + {\rm O}(6)$ to 
$v' = z' \overline{z}' + {\rm O}(6)$, can be shown
to be necessarily of the form~{\cite[Prp.~2.4]{Merker-2020}}:
\[
\aligned
z'
&
\,:=\,
\lambda\,z
+
2i\lambda\overline{\alpha}\,
z^2
+
\big(
-4\lambda\overline{\alpha}^2
\big)\,
z^3
+
\big(
-8i\lambda\overline{\alpha}^3
\big)\,
z^4
\\
&
\ \ \ \ \ \ \ \ \ \ \ \
+
\lambda\alpha\,
w
+
\big(
3i\lambda\alpha\overline{\alpha}
+
\lambda r
\big)\,
zw
+
\big(
-8\lambda\alpha\overline{\alpha}^2
+
4i\overline{\alpha}\lambda r
\big)\,
z^2w
\notag
\\
&
\ \ \ \ \ \ \ \ \ \ \ \ \ \ \ \ \ \ \ \ \ \ \ \ \ \ \ \ \ \ \ \ \ \ 
\ \ \ \ \ \ \ \ \ \ \ \ \ \ \ \ \ \ \ \ \ \ \ \ \ \
+
\big(
\lambda\alpha r
+
i\lambda\alpha^2\overline{\alpha}
\big)\,
w^2
+
{\rm O}(5),
\\
w'
&
\,=\,
\lambda\overline{\lambda}\,
w
+
2i\lambda\overline{\lambda}\overline{\alpha}\,
zw
+
\big(
-4\lambda\overline{\lambda}\overline{\alpha}^2
\big)\,
z^2w
+
\big(
-8i\lambda\overline{\lambda}\overline{\alpha}^3
\big)\,
z^3w
\notag
\\
&
\ \ \ \ \ \ \ \ \ \ \ \ \ \ \ \ \ \ \ \
+
\big(
i\lambda\overline{\lambda}\alpha\overline{\alpha}
+
\lambda\overline{\lambda}r
\big)\,
w^2
+
\big(
4i\lambda\overline{\lambda}\overline{\alpha}r
-
4\lambda\overline{\lambda}\overline{\alpha}^2\alpha
\big)\,
zw^2
+
{\rm O}(6),
\endaligned
\]
and this form 
coincides exactly with the Taylor expansion, up to 
weighted orders $4$, $5$, of the general stability group
of the {\em model} $\{v = z \overline{z}\} \longrightarrow 
\{ v' = z'\overline{z}' \}$, which is well know to be:
\[
z'
\,=\,
\frac{\lambda\,(z+\alpha\,w)}{
1-2i\overline{\alpha}\,z
-
(r+i\alpha\overline{\alpha})\,w},
\ \ \ \ \ \ \ \ \ \ \ \ \ \ \ \ 
w'
\,=\,
\frac{\lambda\overline{\lambda}\,w}{
1-2i\overline{\alpha}\,z
-
(r+i\alpha\overline{\alpha})\,w},
\]
with arbitrary $\lambda \in \C^\ast$, $\alpha \in \C$, $r \in \R$.

One could then figure out that precisely similar statements hold for 
$M^5 \in \mathfrak{C}_{2,1}$.
However, some `{\sl discrepancies}', which we will overcome, will 
occur. Indeed, let us briefly describe some differences,
as a preliminary view on the technical road we will 
drive into the forest. 

Taking the weights $[z] := 1$, $[\zeta] := 1$, 
$[w] := 2$, starting with $u = F(z, \zeta, \overline{z}, 
\overline{\zeta}, v)$ passing through the origin,
by progressively normalizing the power series expansion 
of $F$, it is not difficult to show that any
$M^5 \in \mathfrak{C}_{2,1}$ can be brought to the form:
\[
u
\,=\,
z\overline{z}
+
\tfrac{1}{2}\,
\overline{z}^2\zeta
+
\tfrac{1}{2}\,
z^2\overline{\zeta}
+
{\rm O}_{z,\zeta,\overline{z},\overline{\zeta},v}(4).
\]
As we know from 
Section~{\ref{fractional-representation-isotropy-group}},
the isotropy group of the Gaussier-Merker model
is {\sl also} parametrized by $5$ real constants
$\lambda \in \C^\ast$, $\alpha \in \C$, $r \in \R$,
and an expansion of the concerned fractional formulas 
was provided there. 

However, one can verify (exercise) that the stability 
group of the above punctual normalization up to order
$3$ happens to be:
\[
\aligned
z'
&
\,:=\,
\lambda\,
z
+
\Big(
\frac{\delta}{\overline{\lambda}}
-
\tfrac{1}{2}\,
\frac{\lambda^2}{\overline{\lambda}}\,
\overline{\beta}
\Big)\,
z^2
-
\tfrac{1}{2}\,
\frac{\overline{\delta}}{\overline{\lambda}}
w,
\\
\zeta'
&
\,:=\,
\frac{\lambda}{\overline{\lambda}}\,
\zeta
+
\beta\,
z,
\\
w'
&
\,:=\,
\lambda\overline{\lambda}\,
w
+
\delta\,
zw,
\endaligned
\]
with arbitrary $\lambda \in \C^\ast$, $\beta \in \C$, 
$\delta \in \C$. This looks different from 
the stability group of 
the model, shown in
Section~{\ref{fractional-representation-isotropy-group}} 
and truncated to orders $2$, $1$, $3$.

Next, it can be shown (and we will do it) that that any $M^5 \in 
\mathfrak{C}_{2,1}$ can be brought to the form:
\[
u
\,=\,
z\overline{z}
+
\tfrac{1}{2}\,
\overline{z}^2\zeta
+
\tfrac{1}{2}\,
z^2\overline{\zeta}
+
z\overline{z}\zeta\overline{\zeta}
+
{\rm O}_{z,\zeta,\overline{z},\overline{\zeta},v}(5).
\] 
Lemma~{\ref{Lm-stability-order-4}} will
show that the stability of this equation reads as:
\[
\aligned
z'
&
\,:=\,
\lambda\,z
-
i\,\lambda\alpha\,z^2
-
i\,\lambda\overline{\alpha}\,w
-
\frac{\lambda^2}{\overline{\lambda}}\,
\overline{\beta}\,
z^3
+
\Big(
i\,\lambda r
-
\tfrac{3}{2}\,
\lambda\alpha\overline{\alpha}
-
\tfrac{1}{4}\,
\frac{\lambda^2}{\overline{\lambda}}\,
\overline{\varepsilon}
-
\tfrac{1}{4}\,
\overline{\lambda}\varepsilon
\Big)\,
zw
+
i\,\lambda\alpha\,
\zeta w,
\\
\zeta'
&
\,:=\,
\frac{\lambda}{\overline{\lambda}}\,\zeta
+
2i\,\frac{\lambda}{\overline{\lambda}}\,\overline{\alpha}\,
z
+
\varepsilon\,z^2
-
2i\,\frac{\lambda}{\overline{\lambda}}\,\alpha\,
z\zeta
+
\beta\,w,
\\
w'
&
\,:=\,
\lambda\overline{\lambda}\,w
-
2i\,\lambda\overline{\lambda}\alpha\,
zw
-
\big(
2\,\lambda\overline{\lambda}\alpha^2
+
\lambda^2\overline{\beta}
\big)\,
z^2w
+
\big(
-\lambda\overline{\lambda}\alpha\overline{\alpha}
+
i\,\lambda\overline{\lambda}\,r
\big)\,w^2,
\endaligned
\]
where $\lambda \in \C^\ast$, $\alpha \in \C$, $r \in \R$,
$\beta \in \C$, $\varepsilon \in \C$ are arbitrary
parameters. Thus, in comparison with the isotropy of the GM-model, 
shown in
Section~{\ref{fractional-representation-isotropy-group}} and
truncated to orders $3$, $2$, $4$,
there are two `extra' complex parameters,
namely $\beta$, $\varepsilon$.

Also, in Proposition~{\ref{Prp-F-origin-order-5}}
we will normalize, still at the origin only:
\[
\aligned
u
&
\,=\,
z\overline{z}
+
\tfrac{1}{2}\,\overline{z}^2\zeta
+
\tfrac{1}{2}\,z^2\overline{\zeta}
+
z\overline{z}\zeta\overline{\zeta}
+
\tfrac{1}{2}\,\overline{z}^2\zeta\zeta\overline{\zeta}
+
\tfrac{1}{2}\,z^2\overline{\zeta}\zeta\overline{\zeta}
\notag
\\
&
\ \ \ \ \ \ \ \ \ \ \
\ \ \ \ \ \ \ \ \ \ \ \ \ \ \ \ \ \ \ \ \ \ \ \ \ \
+
z^3\overline{\zeta}^2\,F_{3,0,0,2,0}
+
\overline{z}^3\zeta^2\,\overline{F_{3,0,0,2,0}}
+
{\rm O}_{z,\zeta,\overline{z},\overline{\zeta},v}(6),
\endaligned
\]
and in Lemma~{\ref{Lm-isotropy-order-5}}, we will see that the
stability group of this normal form is:

\[
\footnotesize
\aligned
z'
&
\,:=\,
\lambda\,z
-
i\,\lambda\alpha\,z^2
-
i\,\lambda\overline{\alpha}\,w
-
\lambda\alpha^2\,
z^3
+
\Big(
i\,\lambda r
-
3\,\lambda\alpha\overline{\alpha}
+
2i\,\lambda\alpha\,
F_{3,0,0,2,0}
-
2i\,\lambda\overline{\alpha}\,\overline{F_{3,0,0,2,0}}
\Big)\,
zw
+
i\,\lambda\alpha\,
\zeta w
\\
&
\ \ \ \ \ \ \ \ \ \ \ \ \ \ \
+ 
i\,\lambda\alpha^3\,
z^4
+
\Big(
8i\,\lambda\alpha^2\overline{\alpha}
+
\tfrac{1}{2}\,\frac{\lambda^2}{\overline{\lambda}}\,\overline{\gamma}
+
4\,\frac{\lambda}{\overline{\lambda}}\,\overline{\tau}
+
4\,\lambda\alpha^2\,F_{3,0,0,2,0}
-
8\,\lambda\alpha\overline{\alpha}\,\overline{F_{3,0,0,2,0}}
\Big)\,
z^2w
+
3\,\lambda\alpha^2\,
z\zeta w
+
\tau\,w^2,
\\
\zeta'
&
\,:=\,
\frac{\lambda}{\overline{\lambda}}\,\zeta
+
2i\,\frac{\lambda}{\overline{\lambda}}\,\overline{\alpha}\,
z
+
\Big(
3\,\frac{\lambda}{\overline{\lambda}}\,\alpha\overline{\alpha}\,
-
i\,\frac{\lambda}{\overline{\lambda}}\,r
-
2i\,\frac{\lambda}{\overline{\lambda}}\,\alpha\,
F_{3,0,0,2,0}
+
6i\,\frac{\lambda}{\overline{\lambda}}\,\overline{\alpha}\,
\overline{F_{3,0,0,2,0}}
\Big)\,
z^2
-
2i\,\frac{\lambda}{\overline{\lambda}}\,\alpha\,
z\zeta
+
\frac{\lambda}{\overline{\lambda}}\,\overline{\alpha}^2\,
w
\\
&
\ \ \ \ \ \ \ \ \ \ \ \ \ \ \
+
\Big(
2\,\frac{\lambda}{\overline{\lambda}}\,\alpha r
-
4i\,\frac{\lambda}{\overline{\lambda}}\,\alpha^2\overline{\alpha}
-
2\,\frac{\lambda^2}{\overline{\lambda}^2}\,\overline{\gamma}
-
8\,\frac{\lambda}{\overline{\lambda}^2}\,\overline{\tau}
+
12\,\frac{\lambda}{\overline{\lambda}}\,\alpha^2\,
F_{3,0,0,2,0}
+
4\,\frac{\lambda}{\overline{\lambda}}\,\alpha\overline{\alpha}\,
\overline{F_{3,0,0,2,0}}
\Big)\,z^3
-
3\,\frac{\lambda}{\overline{\lambda}}\,\alpha^2\,
z^2\zeta
+
\gamma\,zw
\\
&
\ \ \ \ \ \ \ \ \ \ \ \ \ \ \
+
\Big(
-2\,\frac{\lambda}{\overline{\lambda}}\,
\alpha\overline{\alpha}
+
4i\,\frac{\lambda}{\overline{\lambda}}\,\alpha\,
F_{3,0,0,2,0}
-
4i\,\frac{\lambda}{\overline{\lambda}}\,\overline{\alpha}\,
\overline{F_{3,0,0,2,0}}
\Big)\,
\zeta w,
\\
w'
&
\,:=\,
\lambda\overline{\lambda}\,w
-
2i\,\lambda\overline{\lambda}\alpha\,
zw
-
3\,\lambda\overline{\lambda}\alpha^2\,
z^2w
+
\big(
-\lambda\overline{\lambda}\alpha\overline{\alpha}
+
i\,\lambda\overline{\lambda}\,r
\big)\,w^2
+
4i\,\lambda\overline{\lambda}\alpha^3\,
z^3w
\\
&
\ \ \ \ \ \ \ \ \ \ \ \ \ \ \
+
\Big(
6i\,\lambda\overline{\lambda}\alpha^2\overline{\alpha}
+
2\,\lambda\overline{\lambda}\alpha r
+
2\,\lambda\overline{\tau}
+
4\,\lambda\overline{\lambda}\alpha^2\,
F_{3,0,0,2,0}
-
4\,\lambda\overline{\lambda}\alpha\overline{\alpha}\,
\overline{F_{3,0,0,2,0}}
\Big)\,
zw^2
+
\lambda\overline{\lambda}\alpha^2\,\zeta w^2.
\endaligned
\]
where $\lambda \in \C^\ast$, $\alpha \in \C$, $r \in \R$,
$\gamma \in \C$, $\tau \in \C$ are arbitrary.
Thus, there are again two `extra' complex parameters,
namely $\gamma$, $\tau$.

To realize a Moser-like normal form for hypersurfaces 
$M^5 \in \mathfrak{C}_{2,1}$ and to define analogs
of Cartan-Moser chains, we will therefore have to adapt a bit
our ideas. Let us give a quick summary.

To start with, we will pick any curve
$0 \in \gamma \subset M$ 
which is {\sl CR-transversal} in the sense
that $\dot{\gamma} \not\in T^cM$.
It is well known that one can always straighten it
to be $\gamma = \{ (0, 0, iv) \} \subset M$, the $v$-axis. 
It is also well known that, after an appropriate 
biholomorphism, one can make the 
graphing function
$F(z, \zeta, \overline{z}, \overline{\zeta}, v)$
to have {\em no} pluriharmonic terms, in the sense 
that $F(z, \zeta, 0, 0, v) \equiv 0$.

In Section~{\ref{chain-straightening-harmonic-killing}} 
to~{\ref{normalization-F-3-0-1-1}}, 
we will continue to {\em prenormalize} and even
start to {\em normalize} $F$ 
further, {\em without touching $\gamma$}, namely
by always stabilizing $\{ (0, 0, iv) \} \subset M$.

However, at some moment of the normalization process,
exactly as what occurs~{\cite{Chern-Moser-1974,
Jacobowitz-1990}}
for Levi nondegenerate
$M^3 \subset \C^2$, one is `{\sl forced}' to 
perform additional normalizations which
{\em bend} the $v$-axis, hence destroy what was preserved
up to this point. This fact confirms
that it was inappropriate to choose
at the beginning any CR-transversal curve $0 \in \gamma \subset M$,
`{\sl at random}'. 

It is at this crucial moment that the Cartan-Moser chains
start to appear to eyes. By appropriately interpreting
the algebraic or geometric normalization conditions
that force to change the $v$-axis, one realizes that
certain CR-transversal curves are {\em invariant} 
under biholomorphisms of $\C^2$. Our goal is to 
view something similar and new about $M^5 \in \mathfrak{C}_{2,1}$.
We will do it.

The Lie-theoretical path taken in~{\cite{Merker-2020}}
consisted in normalizing the equation
of $M$ at only one point, only up to order $5$,
which is quite elementary, can be done by hand or on a computer,
and does not employ (at all) the implicit function theorem.
In this memoir, we will conduct essentially the same method
as in~{\cite{Merker-2020}} but with two differences.
Firstly, we will prenormalize the equation of $M$
not only at $0$ but all along the $v$-axis $\gamma \subset M$
(chosen at random) and reach
Proposition~{\ref{Prp-normalization-order-5-v-axis}}, 
until we come to the point where chains start to appear to eyes.
Then we will work only at $0$, with power series expansions
of orders $5$, $6$, $7$, and `discover' that the chains
are the same as stated by
Observations~{\ref{Obs-Sigma-0-jets-1}}
and~{\ref{Obs-Sigma-0-jets-2}} for the Gaussier-Merker model,
notwithstanding the presence of extra complex parameters.

Once chains are known, we will go back to the starting point,
and choose the CR-transversal $\gamma \subset M$ to {\em be a chain},
then we will plainly apply all what was done for a random
$\gamma$, and
we will deduce that
two normalizations of certain coefficients
$F_{a,b,c,d}(v)$ 
realize themselves gratuitously thanks to chains,
and lastly, we will obtain a complete 
Moser-like normal form.

To terminate our mathematical work and get some
uniqueness property, we will work out 
the formal theory of the normal form only at the end of the paper.

\Section{\bf Chain Straightening and Harmonic Killing}
\label{chain-straightening-harmonic-killing}
\HEAD{{\ref{chain-straightening-harmonic-killing}}.~{\sf Chain 
Straightening and Harmonic Killing}
}{
Wei-Guo {\sc Foo}, Joël {\sc Merker}, The-Anh {\sc Ta}}

Start with any $\mathfrak{C}_{2,1}$ hypersurface $M \subset {\C}^3$,
passing by the origin $0\in M$.  Since $T_{0}^c M \cong \C^2$, we
can assume after a $\C$-linear transformation that $T_{0}^c M =
\C_{z} \times \C_{\zeta} \times \{0\}$, in coordinates $(z, \zeta,
w) \in {\C}^3$.

The `{\sl game}' is to transform $M$ progressively into more and more
{\sl normalized} hypersurfaces. Each (partial) normalization step can
represented by means of a biholomorphism fixing the origin
as:
\[
\aligned
{\C}^3
\,\supset\,
({M}^5,0)
\ \ \ \ \ \ \ \ \ \ 
&
\xrightarrow[{\rule[0pt]{50pt}{0pt}}]{{\sf normalize}}
\ \ \ \ \ \ \ \ \ \ 
({M'}^5,0)
\,\subset\,
{\C'}^3,
\\
\big(z,\zeta,w\big)
\ \ \ \ \ \ \ \ \ \ 
&
\xrightarrow[{\rule[0pt]{50pt}{0pt}}]{}
\ \ \ \ \ \ \ \ \ \
\big(
f(z,\zeta,w),\,
g(z,\zeta,w),\,
h(z,\zeta,w)
\big)
\\
&
\ \ \ \ \ \ \ \ \ \ \ \ \ \ \ \ \ \ \ \ \ \ \ 
\,=:\,
\big(z',\zeta',w'\big).
\endaligned
\]
Without loss of generality, both hypersurfaces will
be assumed, with $w = u + iv$ and $w' = u' +
iv'$, to be $\mathcal{C}^\omega$-graphed as:
\[
u
\,=\,
F\big(z,\zeta,\overline{z},\overline{\zeta},v\big)
\ \ \ \ \ \ \ \ \ \ \ \ \ \ \ \ \ \ \ \
\text{and}
\ \ \ \ \ \ \ \ \ \ \ \ \ \ \ \ \ \ \ \
u'
\,=\,
F'\big(z',\zeta',\overline{z}',\overline{\zeta}',v'\big).
\]
We may assume that $T_{0}^c M = \{ w = 0\}$
is left untouched, so that $T_{0'}^c M' = \{ w' = 0\}$ too.

In fact, step by step, all previously achieved normalizations
will be conserved while performing any further normalization.
Once $M$ has been partly normalized to some new $M'$,
we will erase primes to the obtained $M' =: M$,
normalize once more, and so on.

Now, the hypothesis that the biholomorphism establishes a
CR-diffeomorphism $M \overset{\sim}{\longrightarrow} M'$, expresses as
saying that $u' = F'$ when $u = F$, namely:
\[
0
\,=\,
-\,
\Re\,
h(z,\zeta,w)
+
F'
\Big(
f(z,\zeta,w),\,
g(z,\zeta,w),\,
\overline{f}(\overline{z},\overline{\zeta},\overline{w}),\,
\overline{g}(\overline{z},\overline{\zeta},\overline{w}),\,\,
\Im\,
h(z,\zeta,w)
\Big)
\bigg\vert_{w\,=\,F(z,\zeta,\overline{z},\overline{\zeta},v)+iv}.
\]
Performing the indicated replacement
$w = F + i\,v$ yields 

\begin{Lemma}
\label{Lm-fundamental-identity}
{\bf [Fundamental identity]}
The map $(z', \zeta', w') = (f,g,h)$ sends $M = \{ u = F \}$ 
to $M' = \{u' = F'\}$ if and only if:
\[
\footnotesize
\aligned
0
&
\,\equiv\,
-\,
{\textstyle{\frac{1}{2}}}\,
h
\big(
z,\zeta,\,
F(z,\zeta,\overline{z},\overline{\zeta},v)+iv
\big)
-
{\textstyle{\frac{1}{2}}}\,
\overline{h}
\big(
\overline{z},\overline{\zeta},\,
F(z,\zeta,\overline{z},\overline{\zeta},v)-iv
\big)
\,+
\notag
\\
&
\ \ \ \ \
+
F'
\Big(
f
\big(
z,\zeta,\,
F(z,\zeta,\overline{z},\overline{\zeta},v)+iv
\big),\,\,
g
\big(
z,\zeta,\,
F(z,\zeta,\overline{z},\overline{\zeta},v)+iv
\big),\,\,
\overline{f}
\big(
\overline{z},\overline{\zeta},\,
F(z,\zeta,\overline{z},\overline{\zeta},v)-iv
\big),\,\,
\notag
\\
&
\ \ \ \ \ \ \ \ \ \ \ \ \ \ \ \ \
\overline{g}
\big(
\overline{z},\overline{\zeta},\,
F(z,\zeta,\overline{z},\overline{\zeta},v)-iv
\big),\,\,
{\textstyle{\frac{1}{2i}}}\,
h
\big(
z,\zeta,\,
F(z,\zeta,\overline{z},\overline{\zeta},v)+iv
\big)
-
{\textstyle{\frac{1}{2i}}}\,
\overline{h}
\big(
\overline{z},\overline{\zeta},\,
F(z,\zeta,\overline{z},\overline{\zeta},v)-iv
\big)
\Big),
\endaligned
\]
holds identically in $\C \{ z, \zeta, \overline{z}, 
\overline{\zeta}, v \}$.\qed
\end{Lemma}

Although this equation
looks complicated, it must be dealt with.  Progressive normalizations
will make it more tractable.

One of the first tasks is to annihilate all pluriharmonic monomials
$F_{a,b,0,0,e}\, {z}^a {\zeta}^b {v}^e$ in $(z, \zeta)$, 
and their conjugates as
well. For completeness, we explain in details how to do this known
normalization.  We proceed in two
steps. 

As already explained in 
Section~{\ref{road-map-convergent-normal-form}},
a CR-transversal curve with $0 \in \gamma \subset M$ 
is now at first chosen `at random',
while a better choice will be made later, when
the normalization process will reach a certain deeper point.

\begin{Lemma}
Let $\gamma \colon \R \longrightarrow M$ be any local
$\mathcal{C}^\omega$ curve with $\gamma(0) = 0 \in M$ and
$\dot{\gamma}(0) \not\in T_{0}^cM = \{w = 0\}$.  Then there exists
a biholomorphism $(z, \zeta, w) \longmapsto (z', w', \zeta')$
sending (stabilizing) $T_0^c M = \{w=0\}$ to
$T_{0'}^c M' = \{ w' = 0\}$ 
which sends $\gamma$ to the curve
$\gamma(t) = (0, 0, it)$ straightened along the $v$-axis.
\end{Lemma}

Notice that the CR-transversal direction 
$\dot{\gamma}'(0) \in T_{0'} M' \big\backslash
T_{0'}^c M'$ together with $T_{0'}^c M' = \{ w' = 0\}$ 
implies $T_0 M' = \{ u' = 0\}$.

\proof
Write the curve as:
\[
\gamma(t)
\,=\,
\big(
\varphi(t),\,
\psi(t),\,
\chi(t)
\big),
\]
with some complex-valued analytic functions
$\varphi$, $\psi$, $\chi$. 
By assumption, $\dot{\chi}(0) \neq 0$. 
This guarantees invertibility of the {\em inverse} holomorphic 
change of coordinates:
\[
z
\,:=\,
z'
+
\varphi\big(-iw'\big),
\ \ \ \ \ \ \ \ \ \ \ \ \ \ \ \ \ \ \ \
\zeta
\,:=\,
\zeta'
+
\psi\big(-iw'\big),
\ \ \ \ \ \ \ \ \ \ \ \ \ \ \ \ \ \ \ \
w
\,:=\,
\chi\big(-iw'\big).
\]

Similarly, the target (transformed)
curve can be written $\gamma'(t) = 
\big( \varphi'(t), \psi'(t), i \chi'(t) \big)$\,\,---\,\,note
the $i$ factor\,\,---, and the 
pointwise correspondence
between curves writes as:
\[
\varphi(t)
\,\equiv\,
\varphi'(t)
+
\varphi\big(
-i(i\chi'(t))\big),
\ \ \ \ \ \ \ \ \ \ \ \ \ \ \ \ \ \ \ \
\psi(t)
\,\equiv\,
\psi'(t)
+
\psi\big(
-i(i\chi'(t))\big),
\ \ \ \ \ \ \ \ \ \ \ \ \ \ \ \ \ \ \ \
\chi(t)
\,\equiv\,
\chi\big(
-i(i\chi'(t))\big).
\]
This last identity yields $t \equiv \chi'(t)$ thanks to 
$0 \neq \dot{\chi}(0)$. Replacing then $\chi'(t) := t$ inside
the first two identities 
concludes that $0 \equiv \varphi'(t) \equiv \psi'(t)$.
\endproof

Consequently, the graphing function of the transformed
hypersurface writes, after erasing primes:
\[
M
\colon\ \ \ \ \
u
\,=\,
F\big(z,\zeta,\overline{z},\overline{\zeta},v\big)
\]
with $F = {\rm O}(2)$ and also $F(0,0,0,0,v) \equiv 0$.
This last condition is
technically needed for the next second elementary normalization.

\begin{Lemma}
Starting from $F = {\rm O}(2)$ with $F(0,0,0,0,v) \equiv 0$,
there exists a biholomorphism of the form:
\[
z'
\,:=\,
z,
\ \ \ \ \ \ \ \ \ \ \ \ \ \ \ \ \ \ \ \
\zeta'
\,:=\,
\zeta,
\ \ \ \ \ \ \ \ \ \ \ \ \ \ \ \ \ \ \ \
w'
\,:=\,
w+h(z,\zeta,w),
\]
with $h = {\rm O}(2)$ and $h(0,0,w) \equiv 0$ which transforms 
$\{ u = F\}$ to $\{ u' = F' \}$ satisfying:
\[
0
\,\equiv\,
F'\big(z',\zeta',0,0,v'\big)
\,\equiv\,
F'(0,0,\overline{z}',\overline{\zeta}',v'\big).
\]
\end{Lemma}

The second vanishing identity
is a consequence of the first
by conjugation, thanks 
to~({\ref{Fbar-F-identity}}).
Equivalently, $F_{a,b,0,0,e}' = 0 = F_{0,0,c,d,e}'$ for all
integer indices. Notice that
$F'(0,0,0,0,v') \equiv 0$ still holds. 

\proof
If such a biholomorphism exists, the fundamental identity 
of Lemma~{\ref{Lm-fundamental-identity}}
shows that:
\leqnomode\usetagform{default}
\begin{footnotesize}
\begin{align}
\label{identity-before-harmonic-killing}
0
&
\,\equiv\,
-\,
F(z,\zeta,\overline{z},\overline{\zeta},v)
-
{\textstyle{\frac{1}{2}}}\,
h
\big(
z,\zeta,\,
F(z,\zeta,\overline{z},\overline{\zeta},v)
+
iv
\big)
-
{\textstyle{\frac{1}{2}}}\,
\overline{h}
\big(
\overline{z},\overline{\zeta},\,
F(z,\zeta,\overline{z},\overline{\zeta},v)
-
iv
\big)
\,+
\notag
\\
&
\ \ \ \ \
+
F'
\Big(
z,\zeta,\overline{z},\overline{\zeta},\,
v
+
{\textstyle{\frac{1}{2i}}}\,
h
\big(
z,\zeta,\,
F(z,\zeta,\overline{z},\overline{\zeta},v)
+
iv
\big)
-
{\textstyle{\frac{1}{2i}}}\,
\overline{h}
\big(
\overline{z},\overline{\zeta},\,
F(z,\zeta,\overline{z},\overline{\zeta},v)
-
iv
\big)
\Big).
\end{align}
\end{footnotesize}
Our goal is to make $F'(z', \zeta', 0, 0, v) \equiv 0$.

If this vanishing identity would hold, putting  
$\overline{z} := 0 =: \overline{\zeta}$ 
in~({\ref{identity-before-harmonic-killing}}) 
we would deduce:
\leqnomode\usetagform{default}
\begin{align}
\label{would-deduce-z-zeta-0}
0
\,\equiv\,
-\,F(z,\zeta,0,0,v)
-
{\textstyle{\frac{1}{2}}}\,
h\big(
z,\zeta,\,
F(z,\zeta,0,0,v)
+
i\,v
\big)
-
{\textstyle{\frac{1}{2}}}\,
\overline{h}
\big(
0,0,\,
F(z,\zeta,0,0,v)
-
i\,v
\big)
+
0.
\end{align}
We claim that such an identity can be employed in order
to define $h(z, \zeta, w)$ uniquely, with the supplementary
condition that the last term $-\frac{1}{2}\, \overline{h}$
of~({\ref{would-deduce-z-zeta-0}})
is zero.

Indeed, thanks to $F = {\rm O}(2)$, we may apply the 
implicit function theorem to invert:
\[
F(z,\zeta,0,0,v)
+
i\,v
\,=:\,
\omega
\ \ \ \ \ \ \ \ \ \ \ \ \ \ \ \ \ \ \ \
\Longleftrightarrow
\ \ \ \ \ \ \ \ \ \ \ \ \ \ \ \ \ \ \ \
v
\,=\,
\TT(z,\zeta,\omega)
\,=\,
-\,i\,\omega
+
{\rm O}(2).
\]
Define therefore $h(z, \zeta, w)$ accordingly:
\[
0
\,\equiv\,
-\,
F\big(z,\zeta,\,\TT(z,\zeta,\omega)\big)
-
{\textstyle{\frac{1}{2}}}\,
h(z,\zeta,\omega)
-
{\textstyle{\frac{1}{2}}}
\cdot
0.
\]
Now, because $F(0, 0, 0, 0, v) \equiv 0$ by hypothesis,
it comes $0 \equiv h(0, 0, \omega)$, just by putting
$z := 0 =: \zeta$ in~({\ref{would-deduce-z-zeta-0}}).

Consequently, the identity~({\ref{would-deduce-z-zeta-0}})
is indeed realized with $- \frac{1}{2}\, \overline{h} = 0$.
Finally, coming back 
to~({\ref{identity-before-harmonic-killing}})$\big\vert_{
\overline{z}=\overline{\zeta}=0}$, we get in conclusion
what we want:
\[
0
\,\equiv\,
0
+
F'
\Big(
z,\zeta,0,0,\,
v
+
{\textstyle{\frac{1}{2i}}}\,
h\big(
z,\zeta,\,
F(z,\zeta,0,0,v)
+
i\,v
\big)
-
0
\Big).
\qedhere
\]
\endproof

Thus, erasing primes, we have obtained the preliminary normalization:
\[
u
\,=\,
F
\,=\,
\sum_{a+b\geqslant 1
\atop
c+d\geqslant 1}\,
z^a\zeta^b\overline{z}^c\overline{\zeta}^d\,
F_{a,b,c,d}(v)
\ \ \ \ \ \ \ \ \ 
\text{with}
\ \ \ \ \ \ \ \ \ \ \
F_{a,b,c,d}(v)
\,:=\,
\sum_{e\geqslant 1}\,
F_{a,b,c,d,e}\,
v^e.
\]
In the sequel, we shall perform normalizing biholomorphisms
which stabilize this form. 

\Section{\bf Prenormalization: Step~I}
\label{prenormalization-step-I}
\HEAD{{\ref{prenormalization-step-I}}.~{\sf 
Prenormalization: Step~I}
}{
Wei-Guo {\sc Foo}, Joël {\sc Merker}, The-Anh {\sc Ta}}

To start with, let us expand:
\[
u
\,=\,
z\overline{z}\,
F_{1,0,1,0}(v)
+
z\overline{\zeta}\,
F_{1,0,0,1}(v)
+
\overline{z}\zeta\,
F_{0,1,1,0}(v)
+
\zeta\overline{\zeta}\,
F_{0,1,0,1}(v)
+
{\rm O}_{z,\zeta,\overline{z},\overline{\zeta}}(3).
\]
By assumption, the Levi
matrix of $F$ has rank $1$ everywhere,
hence in particular at the origin.
We compute this matrix:
\[
\Levi(F)
\,=\,
\left(\!
\begin{array}{cccc}
0 & {\rm O}(1) & {\rm O}(1) & -\frac{1}{2}+{\rm O}(2)
\\
{\rm O}(1) & F_{1,0,1,0}(0)+{\rm O}(1) & 
F_{0,1,1,0}(0)+{\rm O}(1) & {\rm O}(1)
\\
{\rm O}(1) & F_{1,0,0,1}(0)+{\rm O}(1) &
F_{0,1,0,1}(0)+{\rm O}(1) & {\rm O}(1)
\\
-\frac{1}{2}+{\rm O}(2) & {\rm O}(1) & {\rm O}(1) & {\rm O}(1)
\end{array}
\!\right),
\]
where ${\rm O}(\NN) = {\rm O}_{z,\zeta, \overline{z}, 
\overline{\zeta}, v}(\NN)$ for any integer $\NN \in \N$.
Hence at the origin $(z, \zeta, \overline{z}, \overline{\zeta}, v)
= (0, 0, 0, 0, 0)$:
\[
1
\,=\,
\rank\,
\left(\!
\begin{array}{cc}
F_{1,0,1,0}(0) & F_{0,1,1,0}(0)
\\
F_{1,0,0,1}(0) & F_{0,1,0,1}(0)
\end{array}
\!\right).
\]

After a $\C$-linear invertible transformation in the
$(z, \zeta)$-space, we can assume:
\leqnomode\usetagform{default}
\begin{align}
\label{1-rank-origin-Levi-form}
1
\,=\,
F_{1,0,1,0}(0)
\ \ \ \ \ \ \
\text{and}
\ \ \ \ \ \ \
0
\,=\,
F_{1,0,0,1}(0)
\,=\,
F_{0,1,1,0}(0)
\,=\,
F_{0,1,0,1}(0),
\end{align}
so that:
\[
u
\,=\,
z\overline{z}
+
{\rm O}_{z,\zeta,\overline{z},\overline{\zeta},v}(3).
\]

\begin{Lemma}
There exists a biholomorphism of the form:
\[
z'
\,:=\,
z\,\varphi(w),
\ \ \ \ \ \ \ \ \ \ \ \ \ \ \ \ \ \ \ \ \ \ \ \ \ \
\zeta'
\,:=\,
\zeta,
\ \ \ \ \ \ \ \ \ \ \ \ \ \ \ \ \ \ \ \ \ \ \ \ \ \
w'
\,:=\,
w,
\]
which transforms $M = \{u = F\}$ into $M'$ of equation:
\[
u'
\,=\,
z'\overline{z}'
+
\sum_{(a,b,c,d)\neq(1,0,1,0)
\atop
a+b\geqslant1,\,
c+d\geqslant 1}\,\,
{z'}^a{\zeta'}^b{\overline{z}'}^c{\overline{w}'}^d\,
F_{a,b,c,d}'(v').
\]
\end{Lemma}

\proof
We write the source hypersurface as:
\[
u
\,=\,
F
\,=\,
z\overline{z}\,F_{1,0,1,0}(v)
+
\zeta\,\big(\cdots\big)
+
\overline{\zeta}\,\big(\cdots\big)
+
{\rm O}_{z,\zeta,\overline{z},\overline{\zeta}}(3),
\]
and similarly for the target:
\[
u'
\,=\,
F'
\,=\,
{z'}{\overline{z}'}\,
F_{1,0,1,0}'(v')
+
\zeta'\,\big(\cdots\big)
+
\overline{\zeta}'\,\big(\cdots\big)
+
{\rm O}_{z',\zeta',\overline{z}',\overline{\zeta}'}(3).
\]
Through any map of the form being considered, since $z' = z\,
(\cdots)$ and $\zeta' = \zeta$, 
it is clear that the remainders correspond to one another:
\[
\zeta'\,\big(\cdots\big)
\,=\,
\zeta\,\big(\cdots\big),
\ \ \ \ \ \ \ \ \ \ \ \ \ \ \ \ \ \ \ \
{\rm O}_{z',\zeta',\overline{z}',\overline{\zeta}'}(3)
\,=\,
{\rm O}_{z,\zeta,\overline{z},\overline{\zeta}}(3).
\]

Since $u = u'$, the fundamental 
identity~({\ref{Lm-fundamental-identity}}) writes:
\[
\aligned
0
&
\,\equiv\,
-\,F\big(z,\zeta,\overline{z},\overline{\zeta}\big)
\\
&
\ \ \ \ \
+
F'
\Big(
z\,\varphi\big(F+iv\big),\,\,
\zeta,\,\,
\overline{z}\,\overline{\varphi}\big(F-iv\big),\,\,
\overline{\zeta},\,\,v
\Big),
\endaligned
\]
which implies, after taking account of the fact that remainders
are the same and that $v = v'$:
\[
\aligned
0
&
\,\equiv\,
-\,z\overline{z}\,F_{1,0,1,0}(v)
\\
&
\ \ \ \ \
+
z\overline{z}\,
\varphi\big(F+iv\big)\,
\overline{\varphi}\big(F-iv\big)\,
F_{1,0,1,0}'(v)
+
\zeta\,\big(\cdots\big)
+
\overline{\zeta}\,\big(\cdots\big)
+
{\rm O}_{z,\zeta,\overline{z},\overline{\zeta}}(3).
\endaligned
\]

Next, by Taylor expanding at $i\,v$, we get:
\[
\varphi\big(iv+F\big)
\,=\,
\varphi(iv)
+
F\,\big(\cdots\big)
\,=\,
\varphi(iv)
+
{\rm O}_{z,\zeta,\overline{z},\overline{\zeta}}(2),
\]
and by inserting this above, we obtain:
\[
\aligned
0
&
\,\equiv\,
-\,z\overline{z}\,F_{1,0,1,0}(v)
\\
&
\ \ \ \ \
+
z\overline{z}\,
\varphi(iv)\,\overline{\varphi}(-iv)\,
F_{1,0,1,0}'(v)
+
{\rm O}_{z,\zeta,\overline{z},\overline{\zeta}}(4)
+
\zeta\,\big(\cdots\big)
+
\overline{\zeta}\,\big(\cdots\big)
+
{\rm O}_{z,\zeta,\overline{z},\overline{\zeta}}(3).
\endaligned
\]

Identifying the coefficients of $z\overline{z}$ yields:
\[
0
\,\equiv\,
-\,F_{1,0,1,0}(v)
+
\varphi(iv)\,\overline{\varphi}(iv)\,
F_{1,0,1,0}'(v).
\]
We can normalize $F_{1,0,1,0}'(v) \equiv 1$ provided
$\varphi$ satisfies:
\[
\varphi(iv)\,\overline{\varphi}(-iv)
\,\equiv\,
F_{1,0,1,0}(v).
\]
Observing that $\overline{F_{1,0,1,0}(v)} = F_{1,0,1,0}(v)$ 
by the reality condition~({\ref{reality-F-v}}),
it suffices to set:
\[
\varphi(w)
\,:=\,
\sqrt{
F_{1,0,1,0}\big(-i\,w\big)},
\]
a function which is holomorphic thanks to $F_{1,0,1,0}(0) = 1$.
\endproof

So, erasing primes, we have obtained:
\leqnomode\usetagform{default}
\begin{align}
\label{normalization-end-step-I}
u
\,=\,
z\overline{z}
+
\sum_{(a,b,c,d)\neq(1,0,1,0)
\atop
a+b\geqslant 1,\,c+d\geqslant 1}\,
z^a\zeta^b\overline{z}^c\overline{\zeta}^d\,
F_{a,b,c,d}(v).
\end{align}

\Section{\bf Dependent and Independent Jets}
\label{dependent-jets-independent-jets}
\HEAD{{\ref{dependent-jets-independent-jets}}.~{\sf 
Dependent and Independent Jets}
}{
Wei-Guo {\sc Foo}, Joël {\sc Merker}, The-Anh {\sc Ta}}

Now, the assumption of Levi degeneracy states as the 
vanishing identity:
\[
0
\,\equiv\,
\Levi(F)
\,:=\,
\left\vert\!
\def\arraystretch{1.25}
\begin{array}{cccc}
0 & F_z & F_\zeta & -\frac{1}{2}+\frac{1}{2i}F_v
\\
F_{\overline{z}} & F_{z\overline{z}} & F_{\zeta\overline{z}} &
\frac{1}{2i}F_{\overline{z}v}
\\
F_{\overline{\zeta}} & F_{z\overline{\zeta}} & 
F_{\zeta\overline{\zeta}} &
\frac{1}{2i}F_{\overline{\zeta}v}
\\
-\frac{1}{2}-\frac{1}{2i}F_v & -\frac{1}{2i}F_{zv} &
-\frac{1}{2i}F_{\zeta v} & \frac{1}{4}F_{vv}
\end{array}
\!\right\vert.
\]
But the Levi form is {\em not} assumed to be identically zero,
it is assumed to be constantly of rank $1$. With $F = z \overline{z}
+ {\rm O}(3)$ in~({\ref{normalization-end-step-I}}), 
this assumption
expresses as the nonvanishing of the minor:
\[
0
\,\neq\,
\Levi_1(F)
\,:=\,
\left\vert\!
\def\arraystretch{1.25}
\begin{array}{ccc}
0 & F_z & -\frac{1}{2}+\frac{1}{2i}F_v
\\
F_{\overline{z}} & F_{z\overline{z}} & \frac{1}{2i}F_{\overline{z}v}
\\
-\frac{1}{2}-\frac{1}{2i}F_v & -\frac{1}{2i}F_{zv} & \frac{1}{4}F_{vv}
\end{array}
\!\right\vert.
\]

Expanding $\Levi_2(F)$ along its third column gives:
\[
\aligned
F_{\zeta\overline{\zeta}}
\cdot
\Levi_1(F)
&
\,\,\equiv\,\,
-\,F_\zeta\,
\left\vert\!
\def\arraystretch{1.25}
\begin{array}{cccc}
F_{\overline{z}} & F_{z\overline{z}} &
\frac{1}{2i}F_{\overline{z}v}
\\
F_{\overline{\zeta}} & F_{z\overline{\zeta}} & 
\frac{1}{2i}F_{\overline{\zeta}v}
\\
-\frac{1}{2}-\frac{1}{2i}F_v & -\frac{1}{2i}F_{zv} &
\frac{1}{4}F_{vv}
\end{array}
\!\right\vert
\\
&
\ \ \ \ \
+
F_{\zeta\overline{z}}\,
\left\vert\!
\def\arraystretch{1.25}
\begin{array}{cccc}
0 & F_z & -\frac{1}{2}+\frac{1}{2i}F_v
\\
F_{\overline{\zeta}} & F_{z\overline{\zeta}} & 
\frac{1}{2i}F_{\overline{\zeta}v}
\\
-\frac{1}{2}-\frac{1}{2i}F_v & -\frac{1}{2i}F_{zv} &
\frac{1}{4}F_{vv}
\end{array}
\!\right\vert
-
\tfrac{1}{2i}\,
F_{\zeta v}\,
\left\vert\!
\def\arraystretch{1.25}
\begin{array}{cccc}
0 & F_z & -\frac{1}{2}+\frac{1}{2i}F_v
\\
F_{\overline{z}} & F_{z\overline{z}} &
\frac{1}{2i}F_{\overline{z}v}
\\
F_{\overline{\zeta}} & F_{z\overline{\zeta}} & 
\frac{1}{2i}F_{\overline{\zeta}v}
\end{array}
\!\right\vert.
\endaligned
\]

Expanding $\Levi_1(F)$ and
dividing,
we get a rational expression:
\[
F_{\zeta\overline{\zeta}}
\,\equiv\,
\frac{\mathcal{P}\big(
F_z,\,
F_\zeta,\,
F_{\overline{z}},\,
F_{\overline{\zeta}},\,
F_v,\,
F_{z\overline{z}},\,
F_{z\overline{\zeta}},\,
F_{\zeta\overline{z}},\,
F_{zv},\,
F_{\zeta v},\,
F_{\overline{z}v},\,
F_{\overline{\zeta}v},\,
F_{vv}
\big)}{
F_{z\overline{z}}
+
F_vF_vF_{z\overline{z}}
+
i\,F_{\overline{z}}F_{zv}
-
i\,F_zF_{\overline{z}v}
+
F_zF_{\overline{z}}F_{vv}
-
F_vF_{\overline{z}}F_{zv}
-
F_zF_vF_{\overline{z}v}
},
\]
whose numerator $\mathcal{P}$ is a certain
universal polynomial, not depending on $F$.
By assumption, the denominator is nonvanishing (locally).

Differentiating this identity and successively performing 
appropriate replacements (exercise), we obtain

\begin{Proposition}
\label{Prp-dependent-derivatives}
For all integers $a,b,c,d,e \in \N$ with $b \geqslant 1$ and
$d \geqslant 1$, there exist a polynomial $\mathcal{P}_{a,b,c,d,e}$
and an exponent $\NN_{a,b,c,d,e} \in \N_{\geqslant 1}$ such that:
\[
\scriptsize
\!\!\!\!\!\!\!\!\!\!\!\!\!\!\!\!\!\!\!\!
\aligned
F_{z^a\zeta^b\overline{z}^c\overline{\zeta}^dv^e}
\,\equiv\,
\frac{\mathcal{P}_{a,b,c,d,e}\Big(
\big\{
F_{z^{a'}\overline{z}^{c'}v^{e'}}
\big\}_{a'+c'+e'\leqslant a+b+c+d+e},\,\,
\big\{
F_{z^{a'}\zeta^{b'}\overline{z}^{c'}v^{e'}}
\big\}_{a'+b'+c'+e'\leqslant a+b+c+d+e}^{b'\geqslant 1},\,\,
\big\{
F_{z^{a'}\overline{z}^{c'}\overline{\zeta}^{d'}v^{e'}}
\big\}_{a'+c'+d'+e'\leqslant a+b+c+d+e}^{d'\geqslant 1}
\Big)
}{
\big(
F_{z\overline{z}}
+
F_vF_vF_{z\overline{z}}
+
i\,F_{\overline{z}}F_{zv}
-
i\,F_zF_{\overline{z}v}
+
F_zF_{\overline{z}}F_{vv}
-
F_vF_{\overline{z}}F_{zv}
-
F_zF_vF_{\overline{z}v}
\big)^{\NN_{a,b,c,d,e}}
}
\endaligned
\]
\end{Proposition}

Accordingly, as in~{\cite{Chen-Merker-2019}}, we will term:
\[
\aligned
\text{\sl Dependent derivatives}
&
\,:=\,
\Big\{
F_{z^a\zeta^b\overline{z}^c\overline{\zeta}^dv^e}
\Big\}_{a,\,b,\,c,\,d,\,e\,\geqslant\,0}^{b\,\geqslant\,1,\,\,
d\,\geqslant\,1},
\\
\text{\sl Independent derivatives}
&
\,:=\,
\big\{
F_{z^a\overline{z}^cv^e}
\big\}_{a,\,c,\,e\,\geqslant\,0}
\medcup
\big\{
F_{z^a\zeta^b\overline{z}^cv^e}
\big\}_{a,\,c,\,e\,\geqslant\,0}^{b\,\geqslant\,1}
\medcup
\big\{
F_{z^a\overline{z}^c\overline{\zeta}^dv^e}
\big\}_{a,\,c,\,e\,\geqslant\,0}^{d\,\geqslant\,1}.
\endaligned
\]

At the origin when we will progressively normalize 
the power series $F$, any modification of the
values of the {\em independent} derivatives
of $F$ at $0$ will automatically transfer to the
{\em dependent} derivatives of $F$ at $0$
{\em via} the formulas of 
Proposition~{\ref{Prp-dependent-derivatives}}.
Thus, freedom of normalization concerns only 
{\em independent} derivatives:
\[
\tfrac{1}{a!}\,\tfrac{1}{b!}\,\tfrac{1}{c!}\,
\tfrac{1}{d!}\,\tfrac{1}{e!}\,
\partial_z^a\,\partial_\zeta^b\,
\partial_{\overline{z}}^c\,
\partial_{\overline{\zeta}}^d\,
\partial_v^e\,
F\big(0,0,0,0,0\big)
\,\,=\,\,
F_{a,b,c,d,e}
\eqno
{\scriptstyle{(b+d\,\leqslant\,1)}}.
\]

For this reason, we will often write:
\[
\aligned
u
\,=\,
F
&
\,=\,
z\overline{z}
+
\sum_{a+c\geqslant 3
\atop
a\geqslant 1,\,c\geqslant 1}\,
z^a\overline{z}^c\,
F_{a,0,c,0}(v)
+
\sum_{b\geqslant 1}\,
z^a\zeta^b\overline{z}^c\,
F_{a,b,c,0}(v)
+
\sum_{d\geqslant 1}\,
z^a\overline{z}^c\overline{\zeta}^d\,
F_{a,0,c,d}(v)
\\
&
\ \ \ \ \ \ \ \ \ \ \ \ \ \ 
+
\zeta\overline{\zeta}\,
\big(\cdots\big),
\endaligned
\]
pointing out that all terms behind $\zeta \overline{\zeta}\,
(\cdots)$ are sorts of `{\sl remainder terms}'.
However, some information will be needed 
about these remainders anyway while normalizing the
main independent derivatives. Indeed, 
regularly, we will come back to the
Levi determinant~({\ref{Levi-determinant-F}}).

\Section{\bf Prenormalization: Step~II}
\label{prenormalization-step-II}
\HEAD{{\ref{prenormalization-step-II}}.~{\sf 
Prenormalization: Step~II}
}{
Wei-Guo {\sc Foo}, Joël {\sc Merker}, The-Anh {\sc Ta}}

Now, we come back to~({\ref{normalization-end-step-I}}),
which we rewrite by selecting monomials 
having $\overline{z}^1$
as single antiholomorphic component:
\[
\aligned
u
&
\,=\,
z\overline{z}
+
\sum_{a+b\geqslant 1
\atop
(a,b)\neq(1,0)}\,
z^a\zeta^b\overline{z}^1\,
F_{a,b,1,0}(v)
+
\sum_{a+b\geqslant 1
\atop
c\geqslant 2}\,
z^a\zeta^b\overline{z}^c\,
F_{a,b,c,0}(v)
+
\sum_{a+b\geqslant 1
\atop
d\geqslant 1}\,
z^a\zeta^b\overline{z}^c\overline{\zeta}^d\,
F_{a,b,c,d}(v)
\\
&
\,=\,
\overline{z}\,
\Big(
z
+
\sum_{a+b\geqslant 1
\atop
(a,b)\neq(1,0)}\,
z^a\zeta^b\,
F_{a,b,1,0}(v)
\Big)
+
\overline{z}^2\,
\big(\cdots\big)
+
\overline{\zeta}\,
\big(\cdots\big).
\endaligned
\]

\begin{Lemma}
There exists a biholomorphism of the form:
\[
z'
\,:=\,
z
+
\Lambda(z,\zeta,w)
\,=\,
z
+
{\rm O}_{z,\zeta,w}(2),
\ \ \ \ \ \ \ \ \ \ \ \ \ \ \ \ \ \ \ \
\zeta'
\,:=\,
\zeta,
\ \ \ \ \ \ \ \ \ \ \ \ \ \ \ \ \ \ \ \
w'
\,:=\,
w,
\]
which transforms $M = \{ u = F\}$ into $M'$ of equation:
\[
u'
\,=\,
z'\overline{z}'
+
{\overline{z}'}^2\,
\big(\cdots\big)
+
\overline{\zeta}'\,
\big(\cdots\big).
\]
\end{Lemma}

\proof
Set:
\[
\Lambda(z,\zeta,w)
\,:=\,
\sum_{a+b\geqslant 1
\atop
(a,b)\neq(1,0)}\,
z^a\zeta^b\,
F_{a,b,1,0}\big(-i\,w\big)
\,\,=\,\,
z^2\,\big(\cdots\big)
+
\zeta\,\big(\cdots\big).
\]
Since $F_{0,1,1,0}(0) = 0$ 
by~({\ref{1-rank-origin-Levi-form}}),
we indeed have $\Lambda = 
{\rm O}_{z,\zeta,w}(2)$. Thus the equation of $M$ writes:
\[
u
\,=\,
\overline{z}\,
\big(
z
+
\Lambda(z,\zeta,v)
\big)
+
\overline{z}^2\,\big(\cdots\big)
+
\overline{\zeta}\,\big(\cdots\big).
\]

Restricting $z' = z + \Lambda (z, \zeta, -iw)$ 
to $M$, Taylor expanding at
$(z, \zeta, v)$, and using 
$0 \equiv F(z,\zeta, 0, 0, v)$
we obtain:
\[
\aligned
z'
\,=\,
z
+
\Lambda\big(z,\,\zeta,\,v-iF\big)
&
\,=\,
z
+
\Lambda(z,\zeta,v)
+
F\,\big(\cdots\big)
\\
&
\,=\,
z
+
\Lambda(z,\zeta,v)
+
\overline{z}\,
\big(\cdots\big)
+
\overline{\zeta}\,
\big(\cdots\big),
\endaligned
\]
hence replacing $z + \Lambda(z,\zeta,v) = 
z' - \overline{z}(\cdots) - \overline{\zeta}(\cdots)$ 
and replacing $\zeta := \zeta'$:
\[
\aligned
u'
\,=\,
u
&
\,=\,
\overline{z}
\Big(
z'
-
\overline{z}\,\big(\cdots\big)
-
\overline{\zeta}\,\big(\cdots\big)
\Big)
+
\overline{z}^2\,\big(\cdots\big)
+
\overline{\zeta}\,\big(\cdots\big)
\\
&
\,=\,
\overline{z}\,z'
+
\overline{z}^2\,\big(\cdots\big)
+
\overline{\zeta}'\,\big(\cdots\big).
\endaligned
\]
Now, an inversion gives:
\[
\aligned
z
+
\Lambda
\,=\,
z
+
z^2\,\big(\cdots\big)
+
\zeta\,\big(\cdots\big)
\,=\,
z'
\ \ \ \ \ \ \ \ \ \ \ \ \ \ \ \ \ \ \ \
&
\Longleftrightarrow
\ \ \ \ \ \ \ \ \ \ \ \ \ \ \ \ \ \ \ \
z
\,=\,
z'
+
{z'}^2\,\big(\cdots\big)
+
\zeta'\,\big(\cdots\big)
\\
&\,\,
\Longrightarrow
\ \ \ \ \ \ \ \ \ \ \ \ \ \ \ \ \ \ 
\overline{z}^2
\,=\,
\overline{z}'}^2\,\big(\cdots\big)
+
{\overline{\zeta}'\,\big(\cdots\big),
\endaligned
\]
which concludes:
\[
u'
\,=\,
z'\overline{z}'
+
{\overline{z}'}^2\,\big(\cdots\big)
+
\overline{\zeta}'\,\big(\cdots\big).
\qedhere
\]
\endproof

Erasing primes, and using the fact that the graphing
function is real, we obtain

\begin{Corollary}
\label{Cor-Lfm-0-intermediate}
Any $\mathcal{C}^\omega$ hypersurface $0 \in M^5 \subset
\C^3$ whose Levi form is
of rank $1$ at the origin can be brought to the form:
\[
u
\,=\,
z\overline{z}
+
z^2\overline{z}^2\,
\big(\cdots\big)
+
\overline{z}^2\zeta\,
\big(\cdots\big)
+
z^2\overline{\zeta}\,
\big(\cdots\big)
+
\zeta\overline{\zeta}\,
\big(\cdots\big).
\eqno\qed
\]
\end{Corollary}

Next, as said, we need more information about the appearing
dependent derivatives
in the remainder $\zeta \overline{\zeta} (\cdots)$.
We start to really use the assumption that the Levi form of 
$M \in \mathfrak{C}_{2,1}$ has constant rank $1$.

\begin{Lemma}
\label{Lm-dependent-jets-order-2}
Any $\mathcal{C}^\omega$ hypersurface $0 \in M^5 \subset
\C^3$ whose Levi form is
of {\em constant rank $1$}
around the origin can be brought to the form:
\[
u
\,=\,
z\overline{z}
+
z^2\overline{z}^2\,
{\rm O}_{z,\overline{z}}(0)
+
\overline{z}^2\zeta\,
{\rm O}_{z,\zeta,\overline{z}}(0)
+
z^2\overline{\zeta}\,
{\rm O}_{z,\overline{z},\overline{\zeta}}(0)
+
\zeta\overline{\zeta}\,
{\rm O}_{z,\zeta,\overline{z},\overline{\zeta}}(2).
\]
\end{Lemma}

\proof
Indeed, from the equation
of Corollary~{\ref{Cor-Lfm-0-intermediate}}, 
rewritten by emphasizing the remainder $R$, which is
{\em real}, as:
\[
u
\,=\,
z\overline{z}
+
z^2\overline{z}^2\,
\big(\cdots\big)
+
\overline{z}^2\zeta\,
\big(\cdots\big)
+
z^2\overline{\zeta}\,
\big(\cdots\big)
+
\zeta\overline{\zeta}\,R,
\]
the Levi
determinant~({\ref{Levi-determinant-F}}) writes:
\[
0
\,\equiv\,
\left\vert\!
\begin{array}{cccc}
0 & \overline{z}+{\rm O}(2) & {\rm O}(1) &
-\frac{1}{2}+{\rm O}(2)
\\
z+{\rm O}(2) & 1 + {\rm O}(2) & {\rm O}(1) & {\rm O}(2)
\\
{\rm O}(1) & {\rm O}(1) & 
\big[\zeta\overline{\zeta}R\big]_{\zeta\overline{\zeta}} & {\rm O}(1)
\\
-\frac{1}{2}+{\rm O}(2) & {\rm O}(2) & {\rm O}(1) & {\rm O}(2)
\end{array}
\!\right\vert,
\]
where, for abbreviation, we denote shortly ${\rm O} (\NN)$
in the places of ${\rm O}_{z,\zeta,\overline{z},\overline{\zeta}}
(\NN)$, with $\NN \in \N$. Expanding the determinant
along its first column and computing
modulo ${\rm O}(2)$, we get:
\[
\aligned
0
&
\,\equiv\,
-\,\big(z+{\rm O}(2)\big)\,
\left\vert\!
\begin{array}{ccc}
\overline{z}+{\rm O}(2) & {\rm O}(1) &
-\frac{1}{2}+{\rm O}(2)
\\
{\rm O}(1) & 
\big[\zeta\overline{\zeta}R\big]_{\zeta\overline{\zeta}} & {\rm O}(1)
\\
{\rm O}(2) & {\rm O}(1) & {\rm O}(2)
\end{array}
\!\right\vert
+
{\rm O}(1)\,
\left\vert\!
\begin{array}{ccc}
\overline{z}+{\rm O}(2) & {\rm O}(1) &
-\frac{1}{2}+{\rm O}(2)
\\
1 + {\rm O}(2) & {\rm O}(1) & {\rm O}(2)
\\
{\rm O}(2) & {\rm O}(1) & {\rm O}(2)
\end{array}
\!\right\vert
\\
&
\ \ \ \ \ \ \ \ \ \ \ \ \ \ \ \ \ \ \ \ \ \ \ \ \ \ \ \ \ \ \ \ \ \ \ 
\ \ \ \ \ \ \ \ \ \ \ \ \ \ \ \ \ \ \ \ \ \ \ \ \ \ \
-\,
\big(-\tfrac{1}{2}+{\rm O}(2)\big)\,\left\vert\!
\begin{array}{ccc}
\overline{z}+{\rm O}(2) & {\rm O}(1) &
-\frac{1}{2}+{\rm O}(2)
\\
1 + {\rm O}(2) & {\rm O}(1) & {\rm O}(2)
\\
{\rm O}(1) & 
\big[\zeta\overline{\zeta}R\big]_{\zeta\overline{\zeta}} & {\rm O}(1)
\end{array}
\!\right\vert
\\
&
\,=\,
{\rm O}(2)
+
{\rm O}(2)
-
\tfrac{1}{4}\,
\big[
\zeta\overline{\zeta}\,R
\big]_{\zeta\overline{\zeta}}
+
{\rm O}(2),
\endaligned
\]
whence:
\[
R
+
\zeta\,R_{\zeta}
+
\overline{\zeta}\,R_{\overline{\zeta}}
\,=\,
{\rm O}(2).
\]
Then certainly $R = {\rm O}_{z, \zeta, 
\overline{z}, \overline{\zeta}} (1)$. Since $\overline{R} = R$ 
is real:
\[
R
\,=\,
z\,A(v)
+
\zeta\,B(v)
+
\overline{z}\,\overline{A}(v)
+
\overline{\zeta}\,\overline{B}(v)
+
{\rm O}_{z,\zeta,\overline{z},\overline{\zeta}}(2),
\]
and replacing $R$, $R_\zeta$, $R_z$ above
yields $0 \equiv A(v) \equiv 2\,B(v)$, so
$R = {\rm O}_{z, \zeta, \overline{z}, \overline{\zeta}}(2)$.
\endproof

\Section{\bf Expression of the Assumption of $2$-Nondegeneracy
at the Origin}
\label{expression-2-ndg-origin}
\HEAD{{\ref{expression-2-ndg-origin}}.~{\sf Expression of the 
Assumption of $2$-Nondegeneracy at the Origin}
}{
Wei-Guo {\sc Foo}, Joël {\sc Merker}, The-Anh {\sc Ta}}

Consequently, abbreviating $\alpha := F_{2,0,0,1,0} \in \C$, we may
show cubic terms:
\[
u
\,=\,
z\overline{z}
+
\alpha\,
z^2\overline{\zeta}
+
\overline{\alpha}\,
\overline{z}^2\zeta
+
{\rm O}_{z,\zeta,\overline{z},\overline{\zeta},v}(4).
\]
Writing $u = \frac{1}{2} w + \frac{1}{2} \overline{w}$,
and solving for $w$, we get:
\[
w
\,=\,
Q\big(z,\zeta,\overline{z},\overline{\zeta},\overline{w}\big)
\,=\,
-\,\overline{w}
+
2\,
z\overline{z}
+
2\,\alpha\,z^2\overline{\zeta}
+
2\,\overline{\alpha}\,
\overline{z}^2\zeta
+
{\rm O}_{z,\zeta,\overline{z},\overline{\zeta},\overline{w}}(4).
\]
Inserting this in the $3 \times 3$ invariant 
determinant of Proposition~{\ref{Prp-2ndg-determinant}}, we get, 
with ${\rm O}(\NN)$ abbreviating ${\rm O}_{z, \zeta, 
\overline{z}, \overline{\zeta}, \overline{w}} (\NN)$:
\[
0
\,\neq\,
\left\vert\!
\begin{array}{ccc}
Q_{\overline{z}} & Q_{\overline{\zeta}} 
& 
Q_{\overline{w}}
\\
Q_{z\overline{z}} & Q_{z\overline{\zeta}} 
& 
Q_{z\overline{w}}
\\
Q_{zz\overline{z}} & Q_{zz\overline{\zeta}} 
& 
Q_{zz\overline{w}}
\end{array}
\!\right\vert
\,\,=\,\,
\left\vert\!
\begin{array}{ccc}
2z+{\rm O}(2) & 2\alpha z^2+{\rm O}(3) & -1+{\rm O}(3)
\\
2+{\rm O}(2) & 4\alpha z+{\rm O}(2) & {\rm O}(2)
\\
{\rm O}(1) & 4\alpha+{\rm O}(1) & {\rm O}(1)
\end{array}
\!\right\vert. 
\]
Expanding along the last column and computing modulo ${\rm O}(1)$:
\[
0
\,\neq\,
-\,8\,\alpha
+
{\rm O}(1).
\]
So the assumption of $2$-nondegeneracy at the origin means that
$\alpha \neq 0$. After the dilation $\zeta \longmapsto 
\frac{1}{2\alpha}\, \zeta$, we obtain:
\[
u
\,=\,
z\overline{z}
+
\tfrac{1}{2}\,
z^2\overline{\zeta}
+
\tfrac{1}{2}\,
\overline{z}^2\zeta
+
{\rm O}_{z,\zeta,\overline{z},\overline{\zeta},v}(4).
\]

\Section{\bf Prenormalization: Step~III}
\label{prenormalization-step-III}
\HEAD{{\ref{prenormalization-step-III}}.~{\sf 
Prenormalization: Step~III}
}{
Wei-Guo {\sc Foo}, Joël {\sc Merker}, The-Anh {\sc Ta}}

Thus, we have obtained the partial normalization:
\[
\aligned
u
\,=\,
z\overline{z}
+
\overline{z}^2\zeta\,
F_{0,1,2,0}(v)
+
z^2\overline{\zeta}\,
F_{2,0,0,1}(v)
+
z^2\overline{z}^2\,
{\rm O}_{z,\overline{z}}(0)
+
\overline{z}^2\zeta\,
{\rm O}_{z,\zeta,\overline{z}}(1)
+
z^2\overline{\zeta}\,
{\rm O}_{z,\overline{z},\overline{\zeta}}(1)
+
\zeta\overline{\zeta}\,
{\rm O}_{z,\zeta,\overline{z},\overline{\zeta}}(2),
\endaligned
\]
with $F_{0,1,2,0}(0) = \frac{1}{2} = F_{2,0,0,1}(0)$.

\begin{Lemma}
There exists a biholomorphism of the form:
\[
z'
\,:=\,
z,
\ \ \ \ \ \ \ \ \ \ \ \ \ \ \ \ \ \ \ \
\zeta'
\,:=\,
\zeta\,\psi(w),
\ \ \ \ \ \ \ \ \ \ \ \ \ \ \ \ \ \ \ \
w'
\,:=\,
w,
\]
with $\psi(0) \neq 0$, which 
normalizes $F_{0,1,2,0}'(v') \equiv \frac{1}{2} \equiv
F_{2,0,0,1}'(v')$:
\[
u'
\,=\,
z'\overline{z}'
+
\tfrac{1}{2}\,
{\overline{z}'}^2\zeta'
+
\tfrac{1}{2}\,
{z'}^2\overline{\zeta}'
+
{z'}^2{\overline{z}'}^2\,
{\rm O}_{z',\overline{z}'}(0)
+
{\overline{z}'}^2\zeta'\,
{\rm O}_{z',\zeta',\overline{z}'}(1)
+
{z'}^2\overline{\zeta}'\,
{\rm O}_{z',\overline{z}',\overline{\zeta}'}(1)
+
\zeta'\overline{\zeta}'\,
{\rm O}_{z',\zeta',\overline{z}',\overline{\zeta}'}(2).
\]
\end{Lemma}

\proof
It is obvious that ${\rm O}_{z', \zeta', \overline{z}', 
\overline{\zeta}'} (\NN) = {\rm O}_{z, \zeta, \overline{z},
\overline{\zeta}}(\NN)$. 

From the source equation:
\[
u
\,=\,
z\overline{z}
+
\overline{z}^2\zeta\,
F_{0,1,2,0}(v)
+
z^2\overline{\zeta}\,
F_{2,0,0,1}(v)
+
{\rm O}_{z,\zeta,\overline{z},\overline{\zeta}}(4),
\]
with $F_{0,1,2,0}(0) = \frac{1}{2} = F_{2,0,0,1}(0)$,
the target equation will be of a similar form:
\[
u'
\,=\,
z'\overline{z}'
+
{\overline{z}'}^2\zeta'\,
F_{0,1,2,0}'(v')
+
{z'}^2\overline{\zeta}'\,
F_{2,0,0,1}'(v')
+
{\rm O}_{z',\zeta',\overline{z}',\overline{\zeta}'}(4).
\]

Since $u = F$ and $u = u' = F'$, 
the fundamental equation writes:
\[
\aligned
0
&
\,\equiv\,
-\,F\big(z,\zeta,\overline{z},\overline{\zeta},v\big)
\\
&
\ \ \ \ \
+
F'
\Big(
z,\,\zeta\,\psi\big(F+iv\big),\,\,
\overline{z},\,\,
\zeta\,\overline{\psi}\big(F-iv\big),\,\,
v
\Big),
\endaligned
\]
that is:
\[
\aligned
0
&
\,\equiv\,
-\,z\overline{z}
-
\overline{z}^2\zeta\,
F_{0,1,2,0}(v)
-
z^2\overline{\zeta}\,
F_{2,0,0,1}(v)
-
{\rm O}_{z,\zeta,\overline{z},\overline{\zeta}}(4)
\\
&
\ \ \ \ \
+
z\overline{z}
+
\overline{z}^2\zeta\,
\psi\big(F+iv\big)\,
F_{0,1,2,0}'(v)
+
z^2\overline{\zeta}\,
\overline{\psi}\big(F-iv\big)\,
F_{2,0,0,1}'(v)
+
{\rm O}_{z,\zeta,\overline{z},\overline{\zeta}}(4).
\endaligned
\]

Next, by Taylor expanding at $iv$:
\[
\psi\big(F+iv\big)
\,=\,
\psi(iv)
+
F\,\big(\cdots\big)
\,=\,
\psi(iv)
+
{\rm O}_{z,\zeta,\overline{z},\overline{\zeta}}(2),
\]
we get:
\[
0
\,\equiv\,
-\,
\overline{z}^2\zeta\,
\Big(
F_{0,1,2,0}(v)
-
\psi(iv)\,
F_{0,1,2,0}'(v)
\Big)
-
z^2\overline{\zeta}\,
\Big(
F_{2,0,0,1}(v)
-
\overline{\psi}\big(-iv\big)\,
F_{2,0,0,1}'(v)
\Big)
+
{\rm O}_{z,\zeta,\overline{z},\overline{\zeta}}(4).
\]

Thus, to normalize $F_{0,1,2,0}'(v) \equiv \frac{1}{2} \equiv
F_{2,0,0,1}'(v)$, it suffices to set:
\[
\psi(w)
\,:=\,
2\,F_{0,1,2,0}\big(-i\,w\big).
\qedhere
\]
\endproof

So erasing primes, we have normalized:
\leqnomode\usetagform{default}
\begin{align}
\label{zbar-2-zeta-1-over-2}
u
\,=\,
z\overline{z}
+
\tfrac{1}{2}\,
\overline{z}^2\zeta
+
\tfrac{1}{2}\,
z^2\overline{\zeta}
&
+
z^2\overline{z}^2\,
F_{2,0,2,0}(v)
+
z^2\overline{z}^2\,
{\rm O}_{z,\overline{z}}(1)
\notag
\\
&\,
+
\overline{z}^2\zeta\,
{\rm O}_{z,\zeta,\overline{\zeta}}(1)
+
z^2\overline{\zeta}\,
{\rm O}_{z,\overline{z},\overline{\zeta}}(1)
+
\zeta\overline{\zeta}\,
{\rm O}_{z,\zeta,\overline{z},\overline{\zeta}}(2).
\end{align}
Our next goal is to eliminate $F_{2,0,2,0}(v)$.

\begin{Lemma}
There exists a biholomorphism of the form:
\[
z'
\,:=\,
z,
\ \ \ \ \ \ \ \ \ \ \ \ \ \ \ \ \ \ \ \
\zeta'
\,:=\,
\zeta
+
z^2\,\psi(w),
\ \ \ \ \ \ \ \ \ \ \ \ \ \ \ \ \ \ \ \
w'
\,:=\,
w,
\]
which normalizes $F_{2,0,2,0}'(v') \equiv 0$:
\[
\aligned
u'
\,=\,
z'\overline{z}'
+
\tfrac{1}{2}\,
{\overline{z}'}^2\zeta'
+
\tfrac{1}{2}\,
{z'}^2\overline{\zeta}'
&
+
{z'}^2{\overline{z}'}^2\,
{\rm O}_{z',\overline{z}'}(1)
\\
&\,
+
{\overline{z}'}^2\zeta'\,
{\rm O}_{z',\zeta',\overline{z}'}(1)
+
{z'}^2{\overline{\zeta}'}\,
{\rm O}_{z',\overline{z}',\overline{\zeta}'}(1)
+
\zeta'\overline{\zeta}'\,
{\rm O}_{z',\zeta',\overline{z}',\overline{\zeta}'}(2).
\endaligned
\]
\end{Lemma}

\proof
In~({\ref{zbar-2-zeta-1-over-2}}), extract the real term
$F_{2,0,2,0}(v)$ and split it:
\leqnomode\usetagform{default}
\begin{align}
\label{split-F-2-0-2-0}
u
\,=\,
z\overline{z}
+
\tfrac{1}{2}\,
\overline{z}^2\,
\big(
\zeta
+
z^2\,F_{2,0,2,0}(v)
\big)
+
\tfrac{1}{2}\,
z^2\,
\big(
\overline{\zeta}
+
\overline{z}^2\,F_{2,0,2,0}(v)
\big)
&
+
z^2\overline{z}^2\,
{\rm O}_{z,\overline{z}}(1)
\notag
\\
&\,
+
\overline{z}^2\zeta\,
{\rm O}_{z,\zeta,\overline{z}}(1)
+
z^2\overline{\zeta}\,
{\rm O}_{z,\overline{z},\overline{\zeta}}(1)
+
\zeta\overline{\zeta}\,
{\rm O}_{z,\zeta,\overline{z},\overline{\zeta}}(2).
\end{align}

We claim that the biholomorphism which works is:
\[
z'
\,:=\,
z,
\ \ \ \ \ \ \ \ \ \ \ \ \ \ \ \ \ \ \ \
\zeta'
\,:=\,
\zeta
+
z^2\,
F_{2,0,2,0}\big(-i\,w\big),
\ \ \ \ \ \ \ \ \ \ \ \ \ \ \ \ \ \ \ \
w'
\,:=\,
w.
\]
The inverse is:
\[
\zeta
\,=\,
\zeta'
-
{z'}^2\,
F_{2,0,2,0}\big(-i\,w'\big)
\,\,=\,\,
\zeta'
+
{z'}^2\,\big(\cdots\big).
\]
We verify first that all remainders correspond to one another:
\[
\footnotesize
\aligned
z^2\overline{z}^2\,
{\rm O}_{z,\overline{z}}(1)
&
\,=\,
{z'}^2{\overline{z}'}^2\,
{\rm O}_{z',\overline{z}'}(1),
\\
\overline{z}^2\zeta\,
{\rm O}_{z,\zeta,\overline{z}}(1)
&
\,=\,
{\overline{z}'}^2\,
\big(
\zeta'
+
{z'}^2\,(\cdots)
\big)\,
{\rm O}_{z',\zeta',\overline{z}'}(1)
\\
&
\,=\,
{\overline{z}'}^2\zeta'\,
{\rm O}_{z',\zeta',\overline{z}'}(1)
+
{z'}^2{\overline{z}'}^2\,
\big[
{\rm O}_{z',\overline{z}'}(1)
+
\zeta'\,
{\rm O}_{z',\zeta',\overline{z}'}(0)
\big]
\\
&
\,=\,
{\overline{z}'}^2\zeta'\,
{\rm O}_{z',\zeta,\overline{z}'}(1)
+
{z'}^2{\overline{z}'}^2\,
{\rm O}_{z',\overline{z}'}(1),
\\
\zeta\overline{\zeta}\,
{\rm O}_{z,\zeta,\overline{z},\overline{\zeta}}(2)
&
\,=\,
\big(
\zeta'
+
{z'}^2\,(\cdots)
\big)\,
\big(
\overline{\zeta}'
+
{\overline{z}'}^2\,(\cdots)
\big)\,
{\rm O}_{z',\zeta',\overline{z}',\overline{\zeta}'}(2)
\\
&
\,=\,
\zeta'\overline{\zeta}'\,
{\rm O}_{z',\zeta',\overline{z}',\overline{\zeta}'}(2)
+
\zeta'{\overline{z}'}^2\,
\big[
{\rm O}_{z',\zeta',\overline{z}'}(2)
+
\overline{\zeta}'\,
{\rm O}_{z',\zeta',\overline{z}',\overline{\zeta}'}(1)
\big]
\\
&
\ \ \ \ \ \ \ \ \ \ \ \ \ \ \ \ \ \ \ \ \ \ \ \ \ \ \ \ \ \ \ \ \ 
+
\overline{\zeta}'{z'}^2\,
\big[
{\rm O}_{z',\overline{z}',\overline{\zeta}'}(2)
+
\zeta'\,{\rm O}_{z',\zeta',\overline{z}',\overline{\zeta}'}(1)
\big]
\\
&
\ \ \ \ \ \ \ \ \ \ \ \ \ \ \ \ \ \ \ \ \ \ \ \ \ \ \ \ \ \ \ \ \ 
+
{z'}^2{\overline{z}'}^2\,
\big[
{\rm O}_{z',\overline{z}'}(2)
+
\zeta'\,{\rm O}_{z',\zeta',\overline{z}'}(1)
+
\overline{\zeta}'\,
{\rm O}_{z',\overline{z}',\overline{\zeta}'}(1)
+
\zeta'\overline{\zeta}'\,
{\rm O}_{z',\zeta',\overline{z}',\overline{\zeta}'}(0)
\big]
\\
&
\,=\,
{z'}^2{\overline{z}'}^2\,
{\rm O}_{z',\overline{z}'}(1)
+
{\overline{z}'}^2\zeta'\,
{\rm O}_{z',\zeta',\overline{z}'}(1)
+
{z'}^2\overline{\zeta}'\,
{\rm O}_{z',\overline{z}',\overline{\zeta}'}(1)
+
\zeta'{\overline{\zeta}'}\,
{\rm O}_{z',\zeta',\overline{z}',\overline{\zeta}'}(2).
\endaligned
\]

Next, using $0 \equiv F\big(0, 0, \overline{z}, \overline{\zeta}, 
v\big)$, and Taylor expanding at $v'$, we can write:
\[
\aligned
\zeta
&
\,=\,
\zeta'
-
{z'}^2\,F_{2,0,2,0}\big(v'-iF\big)
\\
&
\,=\,
\zeta'
-
{z'}^2\,
F_{2,0,2,0}(v')
-
{z'}^2\,F\,\big(\cdots\big)
\\
&
\,=\,
\zeta'
-
{z'}^2\,
F_{2,0,2,0}(v')
-
{z'}^2\,
\big[
z\,(\cdots)
+
\zeta\,(\cdots)
\big]
\\
&
\,=\,
\zeta'
-
{z'}^2\,
F_{2,0,2,0}(v')
-
{z'}^2\,
\big[
z'\,(\cdots)
+
\zeta'(\cdots)
\big].
\endaligned
\]

Lastly, replacing $z$, $\zeta$,  $\overline{z}$,
$\overline{\zeta}$, $u$, $v$ in terms of 
$z'$, $\zeta'$, $\overline{z}'$, $\overline{\zeta}'$,
$u'$, $v'$ in~({\ref{split-F-2-0-2-0}}), 
we obtain what was asserted:
\begin{align*}
u'
&
\,=\,
z'\overline{z}'
+
\tfrac{1}{2}\,
{\overline{z}'}^2\,
\Big(
\zeta'
-
\zero{
{z'}^2\,F_{2,0,2,0}(v')}
-
{z'}^3\,\big(\cdots\big)
-
{z'}^2\zeta'\,\big(\cdots\big)
+
\zero{
{z'}^2\,F_{2,0,2,0}(v')}
\Big)
\\
&
\ \ \ \ \ \ \ \ \ \ \ \ \,
+
\tfrac{1}{2}\,
{z'}^2\,
\Big(
\overline{\zeta}'
-
\zerozero{
{\overline{z}'}^2\,F_{2,0,2,0}(v')}
-
{\overline{z}'}^3\,\big(\cdots\big)
-
{\overline{z}'}^2\overline{\zeta}'\,\big(\cdots\big)
+
\zerozero{
{\overline{z}'}^2\,F_{2,0,2,0}(v')}
\Big)
\\
&
\ \ \ \ \ \ \ \ \ \ \ \ \,
+
{z'}^2{\overline{z}'}^2\,
{\rm O}_{z',\overline{z}'}(1)
+
{\overline{z}'}^2\zeta'\,
{\rm O}_{z',\zeta',\overline{z}'}(1)
+
{z'}^2\overline{\zeta}'\,
{\rm O}_{z',\overline{z}',\overline{\zeta}'}(1)
+
\zeta'\overline{\zeta}'\,
{\rm O}_{z',\zeta',\overline{z}',\overline{\zeta}'}(2)
\\
&
\,=\,
z'\overline{z}'
+
\tfrac{1}{2}\,
{\overline{z}'}^2\zeta'
+
\tfrac{1}{2}\,
{z'}^2\overline{\zeta}'
\\
&
\ \ \ \ \ \ \ \ \ \ \ \ \,
+
{z'}^2{\overline{z}'}^2\,
{\rm O}_{z',\overline{z}'}(1)
+
{\overline{z}'}^2\zeta'\,
{\rm O}_{z',\zeta',\overline{z}'}(1)
+
{z'}^2\overline{\zeta}'\,
{\rm O}_{z',\overline{z}',\overline{\zeta}'}(1)
+
\zeta'\overline{\zeta}'\,
{\rm O}_{z',\zeta',\overline{z}',\overline{\zeta}'}(2).
\qedhere
\end{align*}
\endproof

Thus, dropping primes, we have reached the following
normalization, where we 
show all monomials in $F$ which have
$\overline{z}^2$ as only antiholomorphic part:
\[
\aligned
u
&
\,=\,
z\overline{z}
+
\tfrac{1}{2}\,
\overline{z}^2\zeta
+
\tfrac{1}{2}\,
z^2\overline{\zeta}
+
\sum_{a+c\geqslant 5
\atop
a\geqslant 2,\,c\geqslant 2}\,
z^a\overline{z}^c\,
F_{a,0,c,0}(v)
+
\sum_{a+b+c\geqslant 4
\atop
b\geqslant 1,\,c\geqslant 2}\,
z^a\zeta^b\overline{z}^c\,
F_{a,b,c,0}(v)
\\
&
+
\sum_{a+c+d\geqslant 4
\atop
a\geqslant 2,\,d\geqslant 1}\,
z^a\overline{z}^c\overline{\zeta}^d\,
F_{a,0,c,d}(v)
+
\sum_{a+b+c+d\geqslant 4
\atop
b\geqslant 1,\,d\geqslant 1}\,
z^a\zeta^b\overline{z}^c\overline{\zeta}^d\,
F_{a,b,c,d}(v).
\endaligned
\]

Now, we will work modulo $\overline{z}^3\, (\cdots) + 
\overline{\zeta}\, (\cdots)$, so the last two sums above
disappear and many terms in the first
two sums as well, so that it remains:
\leqnomode\usetagform{default}
\begin{align}
\label{u-a-3-a-b-2}
u
&
\,=\,
z\overline{z}
+
\notag
\tfrac{1}{2}\,\overline{z}^2\,
\bigg[
\zeta
+
2\,\sum_{a\geqslant 3}\,
z^a\,F_{a,0,2,0}(v)
+
2\,\sum_{a+b\geqslant 2
\atop \ \ \ \ 
b\geqslant 1}\,
z^a\zeta^b\,
F_{a,b,2,0}(v)
\bigg]
\\
&
\ \ \ \ \ \ \ \ \ \ \ \ \ \ \ \ \ \ \ \ \ \ \ \ \ \ \ \ \ \ \ \ \ \ \
\ \ \ \ \ \ \ \ \ \ \ \ \ \ \ \ \ \ \ \ \ \ \ \ \ \ \ \ \ \ \ \ \ \ \
+
\overline{z}^3\,\big(\cdots\big)
+
\overline{\zeta}\,\big(\cdots\big).
\end{align}

\begin{Lemma}
The biholomorphism:
\[
\aligned
z'
\,:=\,
z,
\ \ \ \ \ \ \ \ \ \ \ \ \ \ \ \ \ \ \ \
\zeta'
&
\,:=\,
\zeta
+
2\,\sum_{a\geqslant 3}\,
z^a\,F_{a,0,2,0}\big(-i\,w\big)
+
2\,\sum_{a+b\geqslant 2
\atop\ \ \ \ 
b\geqslant 1}\,
z^a\zeta^b\,
F_{a,b,2,0}\big(-i\,w\big),
\\
w'
&
\,:=\,
w,
\endaligned
\]
transforms $M$ into $M'$ of equation:
\[
u
'\,=\,
z'\overline{z}'
+
\tfrac{1}{2}\,{\overline{z}'}^2\zeta'
+
\tfrac{1}{2}\,{z'}^2\overline{\zeta}'
+
{\overline{z}'}^3\,\big(\cdots\big)
+
\overline{\zeta}'\,\big(\cdots\big).
\]
\end{Lemma}

\proof
As in~{\cite{Foo-Merker-Ta-2019}}, we write:
\[
\zeta'
\,:=\,
\zeta
+
\tau(z,w)
+
\zeta\,\omega(z,\zeta,w),
\]
where:
\[
\tau
\,=\,
z^3\,\big(\cdots\big)
\ \ \ \ \ \ \ \ \ \ \ \ \ \ \ \ \ \ \ \
\text{and}
\ \ \ \ \ \ \ \ \ \ \ \ \ \ \ \ \ \ \ \
\omega
\,=\,
{\rm O}_{z,\zeta,w}(1).
\]
The inverse is certainly of the form
$\zeta = \zeta' + {\rm O}_{z', \zeta', w'}(2)$, hence:
\[
\zeta
\,=\,
\zeta'
+
\tau'(z',w')
+
\zeta'\,
\omega'(z',\zeta',w'),
\]
with $\tau' = {\rm O}_{z',w'}(2)$ and
$\omega' = {\rm O}_{z',\zeta',w'}(1)$. We claim that
$\tau' = {z'}^3\, (\cdots)$.

Indeed, replacing $\zeta' = \tau(z,w) + \zeta\, 
[1 + \omega(z, \zeta, \omega)]$ into 
$\zeta = \tau'(z', w') + \zeta'\, [1 + \omega'(z', 
\zeta', w')]$, the following identity must hold
in $\C\{ z, \zeta, w\}$:
\[
\zeta
\,\equiv\,
\tau'(z,w)
+
\big(
\tau(z,w)
+
\zeta\,[1+\omega(z,\zeta,w)]
\big)\,
\Big[
1
+
\omega'\big(
z,\,
\tau(z,w)
+
\zeta\,[1+\omega(z,\zeta,w)],\,
w
\big)
\Big].
\]
Putting $\zeta := 0$, it comes:
\[
0
\,\equiv\,
\tau'(z,w)
+
\tau(z,w)\,
\big[
1
+
{\rm O}_{z,w}(1)
\big]
\,\,\equiv\,\,
\tau'(z,w)
+
z^3\,(\cdots)\,
\big[
1
+
{\rm O}_{z,w}(1)
\big].
\]
Thus $\zeta = \zeta'\, (\cdots) + 
{z'}^3\, (\cdots)$, which enables us to verify that remainders 
correspond as follows:
\[
\aligned
\overline{\zeta}\,\big(\cdots\big)
&
\,=\,
\overline{\zeta}'\,\big(\cdots\big)
+
{\overline{z}'}^3\,
\big(\cdots\big),
\\
\overline{z}^3\,\big(\cdots\big)
&
\,=\,
{\overline{z}'}^3\,\big(\cdots\big).
\endaligned
\]

Next, using $0 \equiv F(z, \zeta, 0, 0, 0)$, so that
$F = \overline{z} (\cdots) + \overline{\zeta} (\cdots) = 
\overline{z}' (\cdots) + \overline{\zeta}' (\cdots)$,
we have:
\[
\aligned
\zeta'
&
\,=\,
\zeta
+
2\,\sum_{a\geqslant 3}\,
z^a\,F_{a,0,2,0}\big(v-iF\big)
+
2\,\sum_{a+b\geqslant 2
\atop\ \ \ \
b\geqslant 1}\,
z^a\zeta^b\,
F_{a,b,2,0}\big(v-iF\big)
\\
&
\,=\,
\zeta
+
2\,\sum_{a\geqslant 3}\,
z^a\,F_{a,0,2,0}(v)
+
F\,\big(\cdots\big)
+
2\,\sum_{a+b\geqslant 2
\atop\ \ \ \
b\geqslant 1}\,
z^a\zeta^b\,F_{a,b,2,0}(v)
+
F\,\big(\cdots\big).
\endaligned
\]

Lastly, coming back to~({\ref{u-a-3-a-b-2}}), we conclude:
\begin{align*}
u'
\,=\,
u
&
\,=\,
z'\overline{z}'
+
\tfrac{1}{2}\,
{\overline{z}'}^2\,
\big[
\zeta'
-
\overline{\zeta}\,(\cdots)
-
\overline{z}\,(\cdots)
\big]
+
{\overline{z}'}^3\,\big(\cdots\big)
+
{\overline{\zeta}'}\,\big(\cdots\big)
\\
&
\,=\,
z'\overline{z}'
+
\tfrac{1}{2}\,
{\overline{z}'}^2\zeta'
+
{\overline{z}'}^3\,\big(\cdots\big)
+
\overline{\zeta}'\,\big(\cdots\big).
\qedhere
\end{align*}
\endproof

Erasing primes, and using the fact that the graphing function
is real, we obtain:
\[
u
\,=\,
z\overline{z}
+
\tfrac{1}{2}\,
\overline{z}^2\zeta
+
\tfrac{1}{2}\,
z^2\overline{\zeta}
+
z^3\overline{z}^3\,\big(\cdots\big)
+
\overline{z}^3\zeta\,\big(\cdots\big)
+
z^3\overline{\zeta}\,\big(\cdots\big)
+
\zeta\overline{\zeta}\,\big(\cdots\big).
\]
It remains only to analyze the dependent-derivatives remainder
$\zeta\overline{\zeta} \big( \cdots \big)$. For this, 
we must extract the single $4^{\rm th}$ order monomial
$z\overline{z} \zeta\overline{\zeta}$
in the GM-model $\maux(z, \zeta, \overline{z}, \overline{\zeta})$.
Then we realize that behind $\zeta \overline{\zeta} (\cdots)$,
there must be order $3$ terms only.

\begin{Proposition}
\label{Prp-prenormalization}
{\bf [Prenormalization]}
Any hypersurface $M^5 \in \mathfrak{C}_{2,1}$
can be brought to the prenormal form:
\[
\aligned
u
&
\,=\,
z\overline{z}
+
\tfrac{1}{2}\,\overline{z}^2\zeta
+
\tfrac{1}{2}\,z^2\overline{\zeta}
+
z\overline{z}\zeta\overline{\zeta}
\\
&
\ \ \ \ \ 
+
z^3\overline{z}^3\,{\rm O}_{z,\overline{z}}(0)
+
\overline{z}^3\zeta\,{\rm O}_{z,\zeta,\overline{z}}(0)
+
z^3\overline{\zeta}\,{\rm O}_{z,\overline{z},\overline{\zeta}}(0)
+
\zeta\overline{\zeta}\,
{\rm O}_{z,\zeta,\overline{z},\overline{\zeta}}(3).
\endaligned
\]
\end{Proposition}

\proof
We write:
\[
u
\,=\,
\overline{z}z
+
\tfrac{1}{2}\,
\overline{z}^2\zeta
+
\tfrac{1}{2}\,
z^2\overline{\zeta}
+
z\overline{z}\zeta\overline{\zeta}
+
z^3\overline{z}^3\,\big(\cdots\big)
+
\overline{z}^3\zeta\,\big(\cdots\big)
+
z^3\overline{\zeta}\,\big(\cdots\big)
+
\zeta\overline{\zeta}\,R.
\]
From Lemma~{\ref{Lm-dependent-jets-order-2}}, we already
know that $R = {\rm O}_{z, \zeta, \overline{z}, 
\overline{\zeta}} (2)$. 

To get more, we look at the Levi determinant:
\[
0
\,\equiv\,
\left\vert\!
\begin{array}{cccc}
0 & \overline{z}+z\overline{\zeta}+{\rm O}(3) &
\frac{1}{2}\overline{z}^2+{\rm O}(3) &
-\frac{1}{2}+{\rm O}(4)
\\
z+\overline{z}\zeta+{\rm O}(3) & 1+{\rm O}(2) &
\overline{z}+{\rm O}(2) & {\rm O}(3)
\\
\frac{1}{2}z^2+{\rm O}(3) & z+{\rm O}(2) &
z\overline{z}+\big[\zeta\overline{\zeta}R\big]_{\zeta\overline{\zeta}}
& {\rm O}(3)
\\
-\frac{1}{2}+{\rm O}(4) & {\rm O}(3) & {\rm O}(3) & {\rm O}(4)
\end{array}
\!\right\vert.
\]
Computing modulo ${\rm O}(3)$, so that the entries
$(2,4)$, $(3,4)$, $(4,2)$, $(4,3)$, $(4,4)$ are
`zero', we get:
\[
0
\,\equiv\,
-\,\big(-\tfrac{1}{2}\big)\,\big(-\tfrac{1}{2}\big)\,
\left\vert\!
\begin{array}{cc}
1+{\rm O}(2) & \overline{z}+{\rm O}(2)
\\
z+{\rm O}(2) & z\overline{z}+
\big[\zeta\overline{\zeta}R\big]_{\zeta\overline{\zeta}}
\end{array}
\!\right\vert
+
{\rm O}(3).
\]
that is:
\[
\big[\zeta\overline{\zeta}R\big]_{\zeta\overline{\zeta}}
\,\equiv\,
{\rm O}(3).
\]

Thanks to the already known $R = {\rm O}(2)$:
\[
R
\,=\,
A\,zz+B\,z\zeta+C\,z\overline{z}+D\,z\overline{\zeta}
+
E\,\zeta\zeta+\overline{D}\,\zeta\overline{z}+G\,\zeta\overline{\zeta}
+
\overline{A}\,\overline{z}z+\overline{B}\,\overline{z}\overline{\zeta}
+
\overline{E}\overline{\zeta}\overline{\zeta}
+
{\rm O}_{z,\zeta,\overline{z},\overline{\zeta}}(3),
\]
with both $\overline{C} = C$ and $\overline{G} = G$ real, hence:
\[
\aligned
{\rm O}(3)
&
\,\equiv\,
R
+
\zeta\,R_\zeta
+
\overline{\zeta}\,
R_{\overline{\zeta}}
+
\zeta\overline{\zeta}\,
R_{\zeta\overline{\zeta}}
\\
&
\,\equiv\,
A\,zz
+
2B\,z\zeta
+
\big(\overline{A}+C\big)\,
z\overline{z}
+
2D\,z\overline{\zeta}
+
3E\,\zeta\zeta
+
2\overline{D}\,\zeta\overline{z}
+
4G\,\zeta\overline{\zeta}
+
2\overline{B}\,\overline{z}\overline{\zeta}
+
3\overline{E}\,\overline{\zeta}\overline{\zeta},
\endaligned
\]
and this forces $A = B = C = D = E = G = 0$, whence
$R = {\rm O}_{z, \zeta, \overline{z}, \overline{\zeta}}(3)$.
\endproof

\Section{\bf Normalization $F_{3,0,0,1}(v) = 0$}
\label{normalization-F-3-0-0-1}
\HEAD{{\ref{normalization-F-3-0-0-1}}.~{\sf 
Normalization $F_{3,0,0,1}(v) = 0$}
}{
Wei-Guo {\sc Foo}, Joël {\sc Merker}, The-Anh {\sc Ta}}

Now, we specify the unique term of order $4$ in $(z, \zeta, 
\overline{z}, \overline{\zeta})$:
\[
u
\,=\,
z\overline{z}
+
\tfrac{1}{2}\,
\overline{z}^2\zeta
+
\tfrac{1}{2}\,
z^2\overline{\zeta}
+
z\overline{z}\zeta\overline{\zeta}
+
z^3\overline{\zeta}\,
F_{3,0,0,1}(v)
+
\overline{z}^3\zeta\,
\overline{F_{3,0,0,1}(v)}
+
{\rm O}_{z,\zeta,\overline{z},\overline{\zeta}}(5).
\]
Abbreviate:
\[
\varphi(v)
\,:=\,
F_{3,0,0,1}(v).
\]

\begin{Lemma}
The biholomorphism:
\[
\aligned
z'
&
\,:=\,
z
+
z^2\,\varphi(-iw)
+
2\,z^3\,\varphi(-iw)\,\varphi(-iw),
\\
\zeta'
&
\,:=\,
\zeta
-
2\,z\,\overline{\varphi}(-iw)
+
4\,z\zeta\,
\varphi(-iw)
-
5\,z^2\,\varphi(-iw)\,
\overline{\varphi}(-iw),
\\
w'
&
\,:=\,
w,
\endaligned
\]
transforms $M$ into $M'$ of equation:
\[
u'
\,=\,
z'\overline{z}'
+
\tfrac{1}{2}\,{\overline{z}'}^2\zeta'
+
\tfrac{1}{2}\,{z'}^2\overline{\zeta}'
+
z'\overline{z}'\zeta'\overline{\zeta}'
+
{\rm O}_{z',\zeta',\overline{z}',\overline{\zeta}'}(5).
\]

\end{Lemma}

\proof
On restriction to $M$ where $-iw = v - iF$:
\[
\aligned
z'
&
\,:=\,
z
+
z^2\,\varphi\big(v-iF\big)
+
2\,z^3\,\varphi\big(v-iF\big)\,\varphi\big(v-iF\big),
\\
\zeta'
&
\,:=\,
\zeta
-
2\,z\,\overline{\varphi}\big(v-iF\big)
+
4\,z\zeta\,
\varphi\big(v-iF\big)
-
5\,z^2\,\varphi\big(v-iF\big)\,
\overline{\varphi}\big(v-iF\big),
\endaligned
\]
hence Taylor expanding at $v$ and using $F = {\rm O}(2)$:
\[
\aligned
z'
&
\,=\,
z
+
z^2\,\varphi(v)
+
2\,z^3\,\varphi(v)\,\varphi(v)
+
{\rm O}_{z,\zeta,\overline{z},\overline{\zeta}}(4),
\\
\zeta'
&
\,=\,
\zeta
-
2\,z\,\overline{\varphi}(v)
+
4\,z\zeta\,\varphi(v)
-
5\,z^2\,\varphi(v)\,\overline{\varphi}(v)
+
{\rm O}_{z,\zeta,\overline{z},\overline{\zeta}}(3).
\endaligned
\]
An expansion concludes:
\begin{align*}
{}
&
z'\overline{z}'
+
\tfrac{1}{2}\,{\overline{z}'}^2\zeta'
+
\tfrac{1}{2}\,{z'}^2\overline{\zeta}'
+
z'\overline{z}'\zeta'\overline{\zeta}'
+
{\rm O}_{z',\zeta',\overline{z}',\overline{\zeta}'}(5)
\,=\,
\notag
\\
&
\,=\,
\Big(
z
+
z^2\,\varphi(v)
+
2\,z^3\,\varphi(v)\,\varphi(v)
\Big)\,
\Big(
\overline{z}
+
\overline{z}^2\,\overline{\varphi}(v)
+
2\,\overline{z}^3\,\overline{\varphi}(v)\,\overline{\varphi}(v)
\Big)
+
{\rm O}_{z,\zeta,\overline{z},\overline{\zeta}}(5)
\notag
\\
&
\ \ \ \ \
+
\tfrac{1}{2}\,
\Big(
\overline{z}
+
\overline{z}^2\,\overline{\varphi}(v)
\Big)^2\,
\Big(
\zeta
-
2\,z\,\overline{\varphi}(v)
+
4\,z\zeta\,\varphi(v)
-
5\,z^2\,\varphi(v)\,\overline{\varphi}(v)
\Big)
+
{\rm O}_{z,\zeta,\overline{z},\overline{\zeta}}(5)
\notag
\\
&
\ \ \ \ \
+
\tfrac{1}{2}\,
\Big(
z
+
z^2\,\varphi(v)
\Big)^2\,
\Big(
\overline{\zeta}
-
2\,\overline{z}\,\varphi(v)
+
4\,\overline{z}\overline{\zeta}\,\overline{\varphi}(v)
-
5\,\overline{z}^2\,\overline{\varphi}(v)\,\varphi(v)
\Big)
+
{\rm O}_{z,\zeta,\overline{z},\overline{\zeta}}(5)
\notag
\\
&
\ \ \ \ \
+
z\overline{z}\,
\big(
\zeta
-
2\,z\,\overline{\varphi}(v)
\big)\,
\big(
\overline{\zeta}
-
2\,\overline{z}\,\varphi(v)
\big)
+
{\rm O}_{z,\zeta,\overline{z},\overline{\zeta}}(5)
\\
&
\,=\,
z\overline{z}
+
\tfrac{1}{2}\,\overline{z}^2\zeta
+
\tfrac{1}{2}\,z^2\overline{\zeta}
+
z\overline{z}\zeta\overline{\zeta}
+
z^3\overline{\zeta}\,\varphi(v)
+
\overline{z}^3\zeta\,\overline{\varphi}(v)
+
{\rm O}_{z,\zeta,\overline{z},\overline{\zeta}}(5).
\qedhere
\end{align*}
\endproof

After this, 
although $F_{a,b,0,0}(v) \equiv 0$ for all $(a,b)$,
it is not necessarily still true that 
prenormalization holds:
\reqnomode\usetagform{EngelLie}
\begin{align}
0
\overset{\text{\bf ?}}{\,\,\equiv\,\,}
F_{a,b,1,0}(v)
\tag{(\forall\,(a,b)\,\neq\,(1,0)),}
\\
0
\overset{\text{\bf ?}}{\,\,\equiv\,\,}
F_{a,b,2,0}(v)
\tag{(\forall\,(a,b)\,\neq\,(0,1)).}
\end{align}

\Section{\bf Repetition of Prenormalization}
\label{repetition-prenormalization}
\HEAD{{\ref{repetition-prenormalization}}.~{\sf 
Repetition of Prenormalization}
}{
Wei-Guo {\sc Foo}, Joël {\sc Merker}, The-Anh {\sc Ta}}

Fortunately, we can repeat the prenormalization. Indeed, let us write:
\[
u
\,=\,
z\overline{z}
+
\tfrac{1}{2}\,\overline{z}^2\zeta
+
\tfrac{1}{2}\,z^2\overline{\zeta}
+
z\overline{z}\zeta\overline{\zeta}
+
\sum_{a+b+c+d\geqslant 5
\atop
a+b\geqslant 1,\,c+d\geqslant 1}\,
z^a\zeta^b\overline{z}^c\overline{\zeta}^d\,
F_{a,b,c,d}(v).
\]
We will perform two biholomorphisms of the form:
\[
z'
\,:=\,
z
+
{\rm O}_{z,\zeta}(4),
\ \ \ \ \ \ \ \ \ \ \ \ \ \ \ \ \ \ \ \
\zeta'
\,:=\,
\zeta
+
{\rm O}_{z,\zeta}(3),
\ \ \ \ \ \ \ \ \ \ \ \ \ \ \ \ \ \ \ \
w'
\,=\,
w,
\]
so that normalizations of terms up to order $4$ included
will be stabilized and preserved.

Starting from:
\[
u
\,=\,
\overline{z}\,
\Big(
z
+
\sum_{a+b\geqslant 4}\,
z^a\zeta^b\,F_{a,b,1,0}(v)
\Big)
+
\overline{z}^2\,\big(\cdots\big)
+
\overline{\zeta}\,\big(\cdots\big),
\]
we perform the following first biholomorphism, 
with $z' := z + {\rm O}_{z, \zeta}(4)$,
$\zeta' := \zeta$, $w' := w$,
which we restrict to
$M$, using $F = \overline{z} (\cdots) + \overline{\zeta} (\cdots)$:
\[
\aligned
z'
\,:=\,
&\,
z
+
\sum_{a+b\geqslant 4}\,
z^a\zeta^b\,
F_{a,b,1,0}(-iw)
\\
\,=\,
&\,
z
+
\sum_{a+b\geqslant 4}\,
z^a\zeta^b\,
\Big[
F_{a,b,1,0}(v)
+
F\,\big(\cdots\big)
\Big]
\\
\,=\,
&\,
z
+
\sum_{a+b\geqslant 4}\,
z^a\zeta^b\,
F_{a,b,1,0}(v)
+
\overline{z}\,\big(\cdots\big)
+
\overline{\zeta}\,\big(\cdots\big),
\endaligned
\]
hence:
\[
z'
-
\overline{z}'\,\big(\cdots\big)
-
\overline{\zeta}'\,\big(\cdots\big)
\,\,=\,\,
z
+
\sum_{a+b\geqslant 4}\,
z^a\zeta^b\,
F_{a,b,1,0}(v),
\]
so we can replace, using $z' = z + z^4 (\cdots) + \zeta (\cdots)$
which gives by inversion $z = z' + {z'}^4 (\cdots) + \zeta' (\cdots)$:
\[
\aligned
u'
\,=\,
u
&
\,=\,
\big(
\overline{z}'
+
{\overline{z}'}^4\,(\cdots)
+
\overline{\zeta}'\,(\cdots)
\big)\,
\big(
z'
-
\overline{z}'\,(\cdots)
-
\overline{\zeta}\,(\cdots)
\big)
+
{\overline{z}'}^2\,\big(\cdots\big)
+
\overline{\zeta}'\,\big(\cdots\big)
\\
&
\,=\,
\overline{z}'z'
+
{\overline{z}'}^2\,\big(\cdots\big)
+
\overline{\zeta}'\,\big(\cdots\big).
\endaligned
\]

Next, erase primes, specify terms having $\overline{z}^2$ as only
antiholomorphic part:
\[
u
\,=\,
z\overline{z}
+
\tfrac{1}{2}\,
\Big(
\zeta
+
2\,\sum_{a+b\geqslant 3}\,
z^a\zeta^b\,F_{a,b,2,0}(v)
\Big)\,
\overline{z}^2
+
\overline{z}^3\,\big(\cdots\big)
+
\overline{\zeta}\,\big(\cdots\big),
\]
and perform the second biholomorphism:
\[
z'
\,:=\,
z,
\ \ \ \ \ \ \ \ \ \ \ \ \ \ \ \ \ \ \ \
\zeta'
\,:=\,
\zeta
+
2\,\sum_{a+b\geqslant 3}\,
z^a\zeta^b\,F_{a,b,2,0}(-iw),
\ \ \ \ \ \ \ \ \ \ \ \ \ \ \ \ \ \ \ \
w'
\,:=\,
w.
\]
Since $-iw = v - iF$ on $M$, using $F = \overline{z} (\cdots) +
\overline{\zeta} (\cdots)$, we have:
\[
\aligned
\zeta'
&
\,=\,
\zeta
+
2\,\sum_{a+b\geqslant 3}\,
z^a\zeta^b\,
F_{a,b,2,0}\big(v-iF\big)
\\
&
\,=\,
\zeta
+
2\,\sum_{a+b\geqslant 3}\,
z^a\zeta^b\,
F_{a,b,2,0}(v)
+
\overline{z}\,\big(\cdots\big)
+
\overline{\zeta}\,\big(\cdots\big),
\endaligned
\]
hence after an inversion:
\[
\zeta'
-
\overline{z}'\,\big(\cdots\big)
-
\overline{\zeta}'\,\big(\cdots\big)
\,\,=\,\,
\zeta
+
2\,\sum_{a+b\geqslant 3}\,
z^a\zeta^b\,F_{a,b,2,0}(v).
\]
So using $\zeta' = \zeta + z^3 (\cdots) + \zeta\,
{\rm O}_{z,\zeta}(2)$
which gives after inversion $\zeta = \zeta' + {z'}^3 (\cdots)
+ \zeta' {\rm O}_{z', \zeta'}(2)$, 
and observing that remainders correspond
to one another, we can replace:
\[
\aligned
u'
\,=\,
u
&
\,=\,
z'\overline{z}'
+
\tfrac{1}{2}\,
{\overline{z}'}^2\,
\big(
\zeta'
-
\overline{z}'\,(\cdots)
-
\overline{\zeta}'\,(\cdots)
\big)
+
{\overline{z}'}^3\,
\big(\cdots\big)
+
\overline{\zeta}'\,\big(\cdots\big)
\\
&
\,=\,
z'\overline{z}'
+
\tfrac{1}{2}\,
{\overline{z}'}^2\zeta'
+
{\overline{z}'}^3\,\big(\cdots\big)
+
\overline{\zeta}'\,\big(\cdots\big).
\endaligned
\]

Since terms are unchanged up to order $5$, 
and since the right-hand side is real, 
we have reached:
\[
u'
\,=\,
z'\overline{z}'
+
\tfrac{1}{2}\,{\overline{z}'}^2\zeta'
+
\tfrac{1}{2}\,{z'}^2\overline{\zeta}'
+
z'\overline{z}'\zeta'\overline{\zeta}'
+
{z'}^3{\overline{z}'}^3\,{\rm O}_{z',\overline{z}'}(0)
+
{\overline{z}'}^3\zeta'\,{\rm O}_{z',\zeta',\overline{z}'}(1)
+
{z'}^3\overline{\zeta}'\,
{\rm O}_{z',\overline{z}',\overline{\zeta}'}(1)
+
\zeta'\overline{\zeta}'\,
{\rm O}_{z',\zeta',\overline{z}',\overline{\zeta}'}(3).
\] 

\begin{Lemma}
Starting from:
\[
u
\,=\,
z\overline{z}
+
\tfrac{1}{2}\,
\overline{z}^2\zeta
+
\tfrac{1}{2}\,
z^2\overline{\zeta}
+
z\overline{z}\zeta\overline{\zeta}
+
\sum_{a+b+c+d\geqslant 5
\atop
a+b\geqslant 1,\,c+d\geqslant 1}\,
z^a\zeta^b\overline{z}^c\overline{\zeta}^d\,
F_{a,b,c,d}(v),
\]
there exists a biholomorphism of the form:
\[
z'
\,=\,
z
+
{\rm O}_{z,\zeta}(4),
\ \ \ \ \ \ \ \ \ \ \ \ \ \ \ \ \ \ \ \
\zeta'
\,=\,
\zeta
+
{\rm O}_{z,\zeta}(3),
\ \ \ \ \ \ \ \ \ \ \ \ \ \ \ \ \ \ \ \
w'
\,:=\,
w,
\]
which transforms $M$ into $M'$ of equation:
\[
u'
\,=\,
z'\overline{z}'
+
\tfrac{1}{2}\,{\overline{z}'}^2\zeta'
+
\tfrac{1}{2}\,{z'}^2\overline{\zeta}'
+
z'\overline{z}'\zeta'\overline{\zeta}'
+
{z'}^3{\overline{z}'}^3\,{\rm O}_{z',\overline{z}'}(0)
+
{\overline{z}'}^3\zeta'\,{\rm O}_{z',\zeta',\overline{z}'}(1)
+
{z'}^3\overline{\zeta}'\,
{\rm O}_{z',\overline{z}',\overline{\zeta}'}(1)
+
\zeta'\overline{\zeta}'\,
{\rm O}_{z',\zeta',\overline{z}',\overline{\zeta}'}(4).
\]
\end{Lemma}

\proof
The only modification is the information about the dependent jets
remainder being an ${\rm O}(4)$ after
$\zeta' \overline{\zeta}'$, which improves 
the previous ${\rm O}(3)$. The proof consists in 
examining the Levi determinant, and proceeds similarly
as at the end of the proof of
Proposition~{\ref{Prp-prenormalization}}.
\endproof

\Section{\bf Normalization $F_{3,0,1,1}(v) = 0$}
\label{normalization-F-3-0-1-1}
\HEAD{{\ref{normalization-F-3-0-1-1}}.~{\sf 
Normalization $F_{3,0,1,1}(v) = 0$}
}{
Wei-Guo {\sc Foo}, Joël {\sc Merker}, The-Anh {\sc Ta}}

Including order $5$ terms from $z^3\overline{\zeta} 
{\rm O}_{z, \overline{z}, \overline{\zeta}}(1)$, three
new terms appear:
\leqnomode\usetagform{default}
\begin{align}
\label{gather-O-6}
u
&
\,=\,
z\overline{z}
+
\tfrac{1}{2}\,\overline{z}^2\zeta
+
\tfrac{1}{2}\,z^2\overline{\zeta}
+
z\overline{z}\zeta\overline{\zeta}
+
\tfrac{1}{2}\,\overline{z}^2\zeta\zeta\overline{\zeta}
+
\tfrac{1}{2}\,z^2\overline{\zeta}\zeta\overline{\zeta}
\notag
\\
&
\ \ \ \ \
+
2\,\Re\,
\Big\{
z^3\overline{z}\overline{\zeta}\,F_{3,0,1,1}(v)
+
z^4\overline{\zeta}\,F_{4,0,0,1}(v)
+
z^3\overline{\zeta}^2\,F_{3,0,0,2}(v)
\Big\}
+
{\rm O}_{z,\zeta,\overline{z},\overline{\zeta}}(6),
\end{align}
and we gather all remainder terms as an ${\rm O}(6)$.

\begin{Lemma}
There exists a biholomorphism of the form:
\[
z'
\,:=\,
z,
\ \ \ \ \ \ \ \ \ \ \ \ \ \ \ \ \ \ \ \
\zeta'
\,:=\,
\zeta
+
i\,\varphi(-iw)\,z^2,
\ \ \ \ \ \ \ \ \ \ \ \ \ \ \ \ \ \ \ \
w'
\,:=\,
w,
\]
with $\varphi(v) \in \R$ for $v \in \R$, which normalizes:
\[
\Im\,
F_{3,0,1,1}'(v')
\,\equiv\,
0.
\]
\end{Lemma}

\proof
On restriction to $M$, the inverse writes:
\[
\aligned
\zeta
&
\,=\,
\zeta'
-
i\,\varphi(-iw)\,{z'}^2
\\
&
\,=\,
\zeta'
-
i\,\varphi\big(v-iF\big)\,{z'}^2
\\
&
\,=\,
\zeta'
-
i\,\varphi(v)\,{z'}^2
+
{z'}^2\,F\,\big(\cdots\big)
\\
&
\,=\,
\zeta'
-
i\,\varphi(v')\,{z'}^2
+
{\rm O}_{z',\zeta',\overline{z}',\overline{\zeta}'}(4).
\endaligned
\]
So we insert in~({\ref{gather-O-6}}) and we conclude:
\begin{footnotesize}
\begin{align*}
u'
\,=\,
u
&
\,=\,
z'\overline{z}'
+
\tfrac{1}{2}\,
{\overline{z}'}^2\,
\big(
\zeta'
-
i\,\varphi(v')\,{z'}^2
+
{\rm O}(4)
\big)
+
\tfrac{1}{2}\,
{z'}^2\,
\big(
\overline{\zeta}'
+
i\,\overline{\varphi}(v')\,{\overline{z}'}^2
+
{\rm O}(4)
\big)
\\
&
\ \ \ \ \
+
z'\overline{z}'\,
\big(
\zeta'
-
i\,\varphi(v')\,{z'}^2
\big)\,
\big(
\overline{\zeta}'
+
i\,\overline{\varphi}(v)\,
{\overline{z}'}^2
\big)
+
\tfrac{1}{2}\,
{\overline{z}'}^2\zeta'\zeta'\overline{\zeta}'
+
\tfrac{1}{2}\,{z'}^2\overline{\zeta}'\zeta'\overline{\zeta}'
\\
&
\ \ \ \ \
+
2\,\Re\,
\Big\{
{z'}^3\overline{z}'\overline{\zeta}'\,
F_{3,0,1,1}(v')
+
{z'}^4\overline{\zeta}'\,
F_{4,0,0,1}(v')
+
{z'}^3{\overline{\zeta}'}^2\,F_{3,0,0,2}(v')
\Big\}
+
{\rm O}_{z',\zeta',\overline{z}',\overline{\zeta}'}(6)
\\
&
\,=\,
z'\overline{z}'
+
\tfrac{1}{2}\,{\overline{z}'}^2\zeta'
+
\tfrac{1}{2}\,{z'}^2\overline{\zeta}'
+
z'\overline{z}'\zeta'\overline{\zeta}'
+
\tfrac{1}{2}\,{\overline{z}'}^2\zeta'\zeta'\overline{\zeta}'
+
\tfrac{1}{2}\,{z'}^2\overline{\zeta}'\zeta'\overline{\zeta}'
\\
&
\ \ \ \ \
+
{z'}^2{\overline{z}'}^2\,
\big[
-
\tfrac{i}{2}\,\varphi(v')
+
\tfrac{i}{2}\,\overline{\varphi}(v')
\big]
\\
&
\ \ \ \ \
+
2\,\Re\,
\Big\{
{z'}^3\overline{z}'\overline{\zeta}'\,
\big[
F_{3,0,1,1}(v')
-
i\,\varphi(v')
\big]
+
{z'}^4\overline{\zeta}'\,
F_{4,0,0,1}(v')
+
{z'}^3{\overline{\zeta}'}^2\,F_{3,0,0,2}(v')
\Big\}
\\
&
\ \ \ \ \ 
+
{\rm O}_{z',\zeta',\overline{z}',\overline{\zeta}'}(6).
\qedhere
\end{align*}
\end{footnotesize}
\endproof

Breaking routine, we do not erase primes.

\begin{Lemma}
There exists a biholomorphism whose inverse is of the form:
\[
z'
\,:=\,
z\,e^{i\varphi(-iw)},
\ \ \ \ \ \ \ \ \ \ \ \ \ \ \ \ \ \ \ \
\zeta'
\,:=\,
\zeta\,e^{2i\varphi(-iw)}
+
\psi(-iw)\,z^2,
\ \ \ \ \ \ \ \ \ \ \ \ \ \ \ \ \ \ \ \
w'
\,:=\,
w,
\]
with $\varphi(v) \in \R$ for $v \in \R$, which normalizes
$u' = F'$ above to $u = F$ of the same shape, but with:
\[
\Re\,F_{3,0,1,1}(v)
\,\equiv\,
0.
\]
\end{Lemma}

\proof
Start with:
\[
\aligned
z'\overline{z}'
&
\,=\,
z\overline{z}\,
e^{i[\varphi(v-iF)-\overline{\varphi}(v+iF)]}
\\
&
\,=\,
z\overline{z}\,
e^{i[\varphi(v)+\varphi_v(v)(-iF)+F^2(\cdots)
-\overline{\varphi}(v)-\overline{\varphi}_v(v)(iF)-F^2(\cdots)]}
\\
&
\,=\,
z\overline{z}\,
e^{2\varphi_v(v)\,F+F^2(\cdots)}
\\
&
\,=\,
z\overline{z}\,
\big(
1
+
2\,\varphi_v(v)\,F
+
{\rm O}(4)
\big)
\\
&
\,=\,
z\overline{z}
+
2\,\varphi_v(v)\,z^2\overline{z}^2
+
\varphi_v(v)\,z\zeta\overline{z}^3
+
\varphi_v(v)\,z^3\overline{z}\overline{\zeta}
+
{\rm O}(6).
\endaligned
\]

Next:
\[
\aligned
\Re\,\big({\overline{z}'}^2\zeta'\big)
&
\,=\,
\Re\,
\Big(
\overline{z}^2\,e^{-2i\overline{\varphi}(i\overline{w})}\,
\big[
\zeta\,e^{2i\varphi(-iw)}
+
\psi(-iw)\,z^2
\big]
\Big)
\\
&
\,=\,
\Re\,
\Big(
\overline{z}^2\zeta\,
e^{2i[-\overline{\varphi}(v+iF)+\varphi(v-iF)]}
+
z^2\overline{z}^2\,e^{-2i\overline{\varphi}(v+iF)}\,
\psi(v-iF)
\Big)
\\
&
\,=\,
\Re\,
\Big(
\overline{z}^2\zeta\,
e^{2i[-\overline{\varphi}(v)-\overline{\varphi}_v(v)(iF)-F^2(\cdots)
+\varphi(v)+\varphi_v(v)(-iF)+F^2(\cdots)]}
\Big)
+
z^2\overline{z}^2\,\psi(v)
+
{\rm O}(6)
\\
&
\,=\,
\Re\,
\Big(
\overline{z}^2\zeta\,
e^{2[\overline{\varphi}_v(v)+\varphi_v(v)]\,F+F^2(\cdots)}
\Big)
+
z^2\overline{z}^2\,\psi(v)
+
{\rm O}(6)
\\
&
\,=\,
\Re\,
\Big(
\overline{z}^2\zeta\,
\big[
1
+
4\,\varphi_v(v)\,
\big(
z\overline{z}+{\rm O}(3)
\big)
+
{\rm O}(4)
\big]
\Big)
+
z^2\overline{z}^2\,\psi(v)
+
{\rm O}(6)
\\
&
\,=\,
\tfrac{1}{2}\,
\overline{z}^2\zeta
+
\tfrac{1}{2}\,
z^2\overline{\zeta}
+
z^2\overline{z}^2\,\psi(v)
+
2\,\varphi_v(v)\,
\overline{z}^3z\zeta
+
2\,\varphi_v(v)\,z^3\overline{z}\overline{\zeta}
+
{\rm O}(6).
\endaligned
\]

Lastly:
\[
\aligned
z'\overline{z}'\zeta'\overline{\zeta}'
&
\,=\,
\big(
z\overline{z}
+
{\rm O}(4)
\big)\,
\big(
\zeta\,e^{2i\varphi(v)+F(\cdots)}
+
(\psi(v)+F(\cdots))\,z^2
\big)\,
\big(
\overline{\zeta}\,
e^{-2i\overline{\varphi}(v)+F(\cdots)}
+
(\overline{\psi}(v)+F(\cdots))\,\overline{z}^2
\big)
\\
&
\,=\,
z\overline{z}\,
\big(
\zeta
+
\psi(v)\,z^2
\big)\,
\big(
\overline{\zeta}
+
\overline{\psi}(v)\,
\overline{z}^2
\big)
+
{\rm O}(6)
\\
&
\,=\,
z\overline{z}\zeta\overline{\zeta}
+
z\overline{z}\zeta\overline{z}^2\,\overline{\psi}(v)
+
z\overline{z}\,\psi(v)\,z^2\overline{\zeta}
+
{\rm O}(6).
\endaligned
\]

Summing, we conclude by taking $\psi(v) := -2\, \varphi_v(v)$
and $\varphi_v(v) := - F_{3,0,1,1}'(v)$:
\begin{align*}
u'
&
\,=\,
z'\overline{z}'
+
\tfrac{1}{2}\,{\overline{z}'}^2\zeta'
+
\tfrac{1}{2}\,{z'}^2\overline{\zeta}'
+
z'\overline{z}'\zeta'\overline{\zeta}'
+
\tfrac{1}{2}\,{\overline{z}'}^2\zeta'\zeta'\overline{\zeta}'
+
\tfrac{1}{2}\,{z'}^2\overline{\zeta}'\zeta'\overline{\zeta}'
\\
&
\ \ \ \ \
+
2\,\Re\,
\Big\{
F_{4,0,0,1}'(v')\,
{z'}^3\overline{z}'\overline{\zeta}'
+
F_{3,0,0,1}'{z'}^3\overline{z}'\overline{\zeta}'
+
F_{3,0,0,2}'(v')\,{z'}^3{\overline{\zeta}'}^2
\Big\}
+
{\rm O}_{z',\zeta',\overline{z}',\overline{\zeta}'}(6)
\\
&
\,=\,
z\overline{z}
+
\tfrac{1}{2}\,\overline{z}^2\zeta
+
\tfrac{1}{2}\,z^2\overline{\zeta}
+
z\overline{z}\zeta\overline{\zeta}
+
\tfrac{1}{2}\,\overline{z}^2\zeta\zeta\overline{\zeta}
+
\tfrac{1}{2}\,z^2\overline{\zeta}\zeta\overline{\zeta}
\\
&
\ \ \ \ \
+
z^2\overline{z}^2\,
\big[
2\,\varphi_v(v)
+
\psi(v)
\big]
\\
&
\ \ \ \ \
+
2\,\Re\,
\Big\{
2\,\varphi_v(v)
+
\psi(v)
+
\varphi_v(v)
+
F_{3,0,1,1}'(v)
\Big\}
\\
&
\ \ \ \ \
+
2\,\Re\,
\Big\{
F_{4,0,0,1}'(v)\,
z^4\overline{\zeta}
+
F_{3,0,0,2}'(v)\,
z^3\overline{\zeta}^2
\Big\}
+
{\rm O}_{z,\zeta,\overline{z},\overline{\zeta}}(6).
\qedhere
\end{align*}
\endproof

\begin{Proposition}
\label{Prp-normalization-order-5-v-axis}
For every hypersurface $M^5 \in \mathfrak{C}_{2,1}$, 
at any point $p \in M$, 
given any CR-transversal curve $p \in \gamma \subset M$,
there exist holomorphic coordinates $(z, \zeta, w) \in \C^3$
vanishing at $p$ in which $\gamma$ is the $v$-axis and in which
$M$ has equation:
\[
\aligned
u
&
\,=\,
z\overline{z}
+
\tfrac{1}{2}\,\overline{z}^2\zeta
+
\tfrac{1}{2}\,z^2\overline{\zeta}
+
z\overline{z}\zeta\overline{\zeta}
+
\tfrac{1}{2}\,\overline{z}^2\zeta\zeta\overline{\zeta}
+
\tfrac{1}{2}\,z^2\overline{\zeta}\zeta\overline{\zeta}
\\
&
\ \ \ \ \ \ \ \ \ \ \ 
+
z^3\overline{z}^3\,
{\rm O}_{z,\overline{z}}(0)
\\
&
\ \ \ \ \ \ \ \ \ \ \ 
+
2\,\Re\,
\Big\{
0
+
z^4\overline{\zeta}\,
F_{4,0,0,1}(v)
+
z^3\overline{\zeta}^2\,
F_{3,0,0,2}(v)
\Big\}
\\
&
\ \ \ \ \ \ \ \ \ \ \ 
+
\overline{z}^3\zeta\,
{\rm O}_{z,\zeta,\overline{\zeta}}(2)
+
z^3\overline{\zeta}\,
{\rm O}_{z,\overline{z},\overline{\zeta}}(2)
+
\zeta\overline{\zeta}\,
{\rm O}_{z,\zeta,\overline{z},\overline{\zeta}}(4).
\endaligned
\]
\end{Proposition}

\proof
The annihilation of $F_{3,0,1,1}(v) \equiv 0$ has been 
performed
above. After that, it is necessary to
repeat prenormalization, as was done in
Section~{\ref{repetition-prenormalization}}, and this
does not perturb the normalizations done
up to order $5$ in $(z, \zeta, \overline{z}, \overline{\zeta})$.

Lastly, 
it remains to justify the vanishing order $4$
of the dependent-derivatives remainder
$\zeta \overline{\zeta} \big( \cdots \big)$.
This can be done by examining the Levi 
determinant~({\ref{Levi-determinant-F}}),
similarly as was done in {\em e.g.} the proof
of Proposition~{\ref{Prp-prenormalization}}.
\endproof

\Section{\bf Normalizations at the Origin}
\label{normalizations-origin}
\HEAD{{\ref{normalizations-origin}}.~{\sf 
Normalizations at the Origin}
}{
Wei-Guo {\sc Foo}, Joël {\sc Merker}, The-Anh {\sc Ta}}

Now, we work at the origin.
Expanding now in terms of all five variables
$(z, \zeta, \overline{z}, \overline{\zeta}, v)$, 
and working modulo {\em weighted} order $6$ terms,
for the weights $[z] = 1$, $[\zeta] = 1$, $[w] = 2$,
we have obtained:
\[
\aligned
u
&
\,=\,
z\overline{z}
+
\tfrac{1}{2}\,\overline{z}^2\zeta
+
\tfrac{1}{2}\,z^2\overline{\zeta}
+
z\overline{z}\zeta\overline{\zeta}
+
\tfrac{1}{2}\,\overline{z}^2\zeta\zeta\overline{\zeta}
+
\tfrac{1}{2}\,z^2\overline{\zeta}\zeta\overline{\zeta}
\\
&
\ \ \ \ \ \ \ \ \ \ \
+
z^4\overline{\zeta}\,F_{4,0,0,1,0}
+
\overline{z}^4\zeta\,\overline{F_{4,0,0,1,0}}
+
z^3\overline{\zeta}^2\,F_{3,0,0,2,0}
+
\overline{z}^3\zeta^2\,\overline{F_{3,0,0,2,0}}
+
{\rm O}_{z,\zeta,\overline{z},\overline{\zeta},v}(6).
\endaligned
\]
To normalize further, we can assume that the target hypersurface
has already been normalized in the same way:
\[
\aligned
u'
&
\,=\,
z'\overline{z}'
+
\tfrac{1}{2}\,{\overline{z}'}^2\zeta'
+
\tfrac{1}{2}\,{z'}^2\overline{\zeta}'
+
z'\overline{z}'\zeta'\overline{\zeta}'
+
\tfrac{1}{2}\,{\overline{z}'}^2\zeta'\zeta'\overline{\zeta}'
+
\tfrac{1}{2}\,{z'}^2\overline{\zeta}'\zeta'\overline{\zeta}'
\\
&
\ \ \ \ \ \ \ \ \ \ \
+
{z'}^4\overline{\zeta}'\,F_{4,0,0,1,0}'
+
{\overline{z}'}^4\zeta'\,\overline{F_{4,0,0,1,0}'}
+
{z'}^3{\overline{\zeta}'}^2\,F_{3,0,0,2,0}'
+
{\overline{z}'}^3{\zeta'}^2\,\overline{F_{3,0,0,2,0}'}
+
{\rm O}_{z',\zeta',\overline{z}',\overline{\zeta}',v'}(6).
\endaligned
\]

But then, it is necessary to stabilize the normalization obtained
up to order $4$. With the help of a computer,
one can prove the following:

\begin{Lemma}
\label{Lm-stability-order-4}
Any biholomorphic map of the form:
\[
z'
\,:=\,
f_1+f_2+f_3,
\ \ \ \ \ \ \ \ \ \ \ \ \ \ \ \ \ \ \ \
\zeta'
\,:=\,
g_1+g_2,
\ \ \ \ \ \ \ \ \ \ \ \ \ \ \ \ \ \ \ \
w'
\,:=\,
h_1+h_2+h_3+h_4,
\]
where $f_1$, $f_2$, $f_3$, $g_1$, $g_2$, $h_1$, $h_2$, 
$h_3$, $h_4$ are weighted homogeneous polynomials
in $(z, \zeta, w)$ of degrees equal to their
indices,
which stabilizes the normalization up to order $4$:
\[
z\overline{z}
+
\tfrac{1}{2}\,\overline{z}^2\zeta
+
\tfrac{1}{2}\,z^2\overline{\zeta}
+
z\overline{z}\zeta\overline{\zeta}
+
{\rm O}_{z,\zeta,\overline{z},\overline{\zeta},v}(5)
\ \ \ \ \ \ \ \ \ \ \ \ \ \ \ \ \ \ \ \
\longrightarrow
\ \ \ \ \ \ \ \ \ \ \ \ \ \ \ \ \ \ \ \
z'\overline{z}'
+
\tfrac{1}{2}\,{\overline{z}'}^2\zeta'
+
\tfrac{1}{2}\,{z'}^2\overline{\zeta}'
+
z'\overline{z}'\zeta'\overline{\zeta}'
+
{\rm O}_{z',\zeta',\overline{z}',\overline{\zeta}',v'}(5)
\]
is of the form:
\[
\aligned
z'
&
\,:=\,
\lambda\,z
-
i\,\lambda\alpha\,z^2
-
i\,\lambda\overline{\alpha}\,w
-
\frac{\lambda^2}{\overline{\lambda}}\,
\overline{\beta}\,
z^3
+
\Big(
i\,\lambda r
-
\tfrac{3}{2}\,
\lambda\alpha\overline{\alpha}
-
\tfrac{1}{4}\,
\frac{\lambda^2}{\overline{\lambda}}\,
\overline{\varepsilon}
-
\tfrac{1}{4}\,
\overline{\lambda}\varepsilon
\Big)\,
zw
+
i\,\lambda\alpha\,
\zeta w,
\\
\zeta'
&
\,:=\,
\frac{\lambda}{\overline{\lambda}}\,\zeta
+
2i\,\frac{\lambda}{\overline{\lambda}}\,\overline{\alpha}\,
z
+
\varepsilon\,z^2
-
2i\,\frac{\lambda}{\overline{\lambda}}\,\alpha\,
z\zeta
+
\beta\,w,
\\
w'
&
\,:=\,
\lambda\overline{\lambda}\,w
-
2i\,\lambda\overline{\lambda}\alpha\,
zw
-
\big(
2\,\lambda\overline{\lambda}\alpha^2
+
\lambda^2\overline{\beta}
\big)\,
z^2w
+
\big(
-\lambda\overline{\lambda}\alpha\overline{\alpha}
+
i\,\lambda\overline{\lambda}\,r
\big)\,w^2,
\endaligned
\]
where $\lambda \in \C^\ast$, $\alpha \in \C$, $r \in \R$,
$\beta \in \C$, $\varepsilon \in \C$ are arbitrary
parameters.\qed
\end{Lemma}

Compared to the expansions to orders
$3$, $2$, $4$ of the components of the isotropy
group of the Gaussier-Merker model shown in
Section~{\ref{fractional-representation-isotropy-group}},
{\em two new parameters appear}, namely $\beta$ and $\varepsilon$.
This causes little trouble to define {\sl chains}
for $M^5 \in \mathfrak{C}_{2,1}$, 
analogous to the Cartan-Moser chains for 
Levi nondegenerate $M^3 \subset \C^2$ redefined
in~{\cite{Merker-2020}}, 
because the linearization of the above collection of maps
(in fact a group) is:
\[
\aligned
z'
&
\,:=\,
\lambda\,z
-
i\,\lambda\overline{\alpha}\,w,
\\
\zeta'
&
\,:=\,
\frac{\lambda}{\overline{\lambda}}\,
\zeta
+
2i\,\frac{\lambda}{\overline{\lambda}}\,\overline{\alpha}\,
z
+
\beta\,w,
\\
w'
&
\,:=\,
\lambda\overline{\lambda}\,
w,
\endaligned
\]
and this action, parametrized by $6$ variables
$\lambda$, $\overline{\lambda}$, 
$\alpha$, $\overline{\alpha}$, 
$\beta$, $\overline{\beta}$, is {\em transitive}
on $1$-jets at the origin (exercise),
contrary to the linearization of the action of the
isotropy group of the Gaussier-Merker model:
\[
\aligned
z'
&
\,:=\,
\lambda\,z
-
i\,\lambda\overline{\alpha}\,w,
\\
\zeta'
&
\,:=\,
\frac{\lambda}{\overline{\lambda}}\,
\zeta
+
2i\,\frac{\lambda}{\overline{\lambda}}\,\overline{\alpha}\,
z
+
\frac{\lambda}{\overline{\lambda}}\,\overline{\alpha}^2\,
w,
\\
w'
&
\,:=\,
\lambda\overline{\lambda}\,
w,
\endaligned
\]
in which $\beta = \frac{\lambda}{\overline{\lambda}}\,
\overline{\alpha}^2$ is a {\em dependent} parameter. This is why
we obtained an invariant
submanifold $\Sigma_0^1$ in 
Observation~{\ref{Obs-Sigma-0-jets-1}}.

To resolve this little discrepancy, we must normalize to higher
order at the origin.

\smallskip

So to normalize further, we will employ maps of the form:
\[
\aligned
z'
&
\,:=\,
\lambda\,z
-
i\,\lambda\alpha\,z^2
-
i\,\lambda\overline{\alpha}\,w
-
\frac{\lambda^2}{\overline{\lambda}}\,
\overline{\beta}\,
z^3
+
\Big(
i\,\lambda r
-
\tfrac{3}{2}\,
\lambda\alpha\overline{\alpha}
-
\tfrac{1}{4}\,
\frac{\lambda^2}{\overline{\lambda}}\,
\overline{\varepsilon}
-
\tfrac{1}{4}\,
\overline{\lambda}\varepsilon
\Big)\,
zw
+
i\,\lambda\alpha\,
\zeta w
\\
&
\ \ \ \ \ 
+
\sum_{a+b+2e=4}\,
f_{a,b,e}\,
z^a\zeta^bw^e,
\\
\zeta'
&
\,:=\,
\frac{\lambda}{\overline{\lambda}}\,\zeta
+
2i\,\frac{\lambda}{\overline{\lambda}}\,\overline{\alpha}\,
z
+
\varepsilon\,z^2
-
2i\,\frac{\lambda}{\overline{\lambda}}\,\alpha\,
z\zeta
+
\beta\,w
\\
&
\ \ \ \ \ 
+
\sum_{a+b+2e=3}\,
g_{a,b,e}\,
z^a\zeta^bw^e,
\\
w'
&
\,:=\,
\lambda\overline{\lambda}\,w
-
2i\,\lambda\overline{\lambda}\alpha\,
zw
-
\big(
2\,\lambda\overline{\lambda}\alpha^2
+
\lambda^2\overline{\beta}
\big)\,
z^2w
+
\big(
-\lambda\overline{\lambda}\alpha\overline{\alpha}
+
i\,\lambda\overline{\lambda}\,r
\big)\,w^2
\\
&
\ \ \ \ \ 
+
\sum_{a+b+2e=5}\,
h_{a,b,e}\,
z^a\zeta^bw^e.
\endaligned
\]
Still on a computer, we verify

\begin{Assertion}
Whatever map is chosen, one has:
\[
F_{3,0,0,2,0}'
\,=\,
\frac{1}{\overline{\lambda}}\,
F_{3,0,0,2,0}.
\eqno\qed 
\]
\end{Assertion}

Furthermore, the map:
\[
\aligned
z'
&
\,:=\,
z
+
2\,F_{4,0,0,1,0}\,z^3
-
2\,F_{4,0,0,1,0}\,z\zeta w,
\\
\zeta'
&
\,:=\,
\zeta
-
2\,\overline{F_{4,0,0,1,0}}\,
w
+
10\,z^2\zeta\,
F_{4,0,0,1,0},
\\
w'
&
\,:=\,
w
+
2\,z^2w\,
F_{4,0,0,1,0},
\endaligned
\]
normalizes $F_{4,0,0,1,0}' := 0$ (exercise). What we have proved
so far deserved to be stated as a 

\begin{Proposition}
\label{Prp-F-origin-order-5}
At every point $p \in M^5$ of a
hypersurface $M^5 \subset \C^3$ in the class $\mathfrak{C}_{2,1}$, 
there exist 
holomorphic coordinates $(z, \zeta, w) \in \C^3$ 
centered at $p = (z_p, \zeta_p, w_p) 
= (0, 0, 0)$ in which $M$ has equation:
\reqnomode\usetagform{EngelLie}
\begin{align}
u
&
\,=\,
z\overline{z}
+
\tfrac{1}{2}\,\overline{z}^2\zeta
+
\tfrac{1}{2}\,z^2\overline{\zeta}
+
z\overline{z}\zeta\overline{\zeta}
+
\tfrac{1}{2}\,\overline{z}^2\zeta\zeta\overline{\zeta}
+
\tfrac{1}{2}\,z^2\overline{\zeta}\zeta\overline{\zeta}
\notag
\\
&
\ \ \ \ \ \ \ \ \ \ \
\ \ \ \ \ \ \ \ \ \ \ \ \ \ \ \ \ \ \ \ \ \ \ \ \ \
+
z^3\overline{\zeta}^2\,F_{3,0,0,2,0}
+
\overline{z}^3\zeta^2\,\overline{F_{3,0,0,2,0}}
+
{\rm O}_{z,\zeta,\overline{z},\overline{\zeta},v}(6).
\tag{\qed}
\end{align}
\end{Proposition}

By applying the technique of
Chen-Foo-Merker-Ta~{\cite[Sections~9, 10]{Chen-Foo-Merker-Ta-2019}},
one can realize, after rather hard computations,
that there corresponds to this Taylor coefficient
$F_{3,0,0,2,0}$, the relative invariant
$\Waux_0$ of Pocchiola, presented in~{\cite{Pocchiola-2013, 
Merker-Pocchiola-2018, Foo-Merker-2019}}:
\[
\aligned
\Waux_0
&
\,:=\,
-\,\frac{1}{3}\,
\frac{\mathcal{K}\big(\overline{\mathcal{L}}_1\big(
\overline{\mathcal{L}}_1(\kaux)\big)\big)}{
\overline{\mathcal{L}}_1(\kaux)^2}
+
\frac{1}{3}\,
\frac{\mathcal{K}\big(\overline{\mathcal{L}}_1(\kaux)\big)\,\,
\overline{\mathcal{L}}_1\big(\overline{\mathcal{L}}_1(\kaux)\big)}{
\overline{\mathcal{L}}_1(\kaux)^3}
\,+
\\
&
\ \ \ \ \
+
\frac{2}{3}\,
\frac{\mathcal{L}_1\big(\mathcal{L}_1(\overline{\kaux})\big)}{
\mathcal{L}_1(\overline{\kaux})}
+
\frac{2}{3}\,
\frac{\mathcal{L}_1\big(\overline{\mathcal{L}}_1(\kaux)\big)}{
\overline{\mathcal{L}}_1(\kaux)}
+
\frac{i}{3}\,
\frac{\mathcal{T}(\kaux)}{\overline{\mathcal{L}}_1(\kaux)},
\endaligned
\]
Much more simply, 
by plugging this normalized $F$ into 
this formula, we obtain its value
{\em only at one point}, namely at the origin:
\[
\Waux_0
\,=\,
4\,
\overline{F_{3,0,0,2,0}}.
\]

Next, we determine the isotropy of this normalization.

\begin{Lemma}
\label{Lm-isotropy-order-5}
Any biholomorphic map of the form:
\[
z'
\,:=\,
f_1+f_2+f_3+f_4,
\ \ \ \ \ \ \ \ \ \ \ \ \ \ \ \ \ \ \ \
\zeta'
\,:=\,
g_1+g_2+g_3,
\ \ \ \ \ \ \ \ \ \ \ \ \ \ \ \ \ \ \ \
w'
\,:=\,
h_1+h_2+h_3+h_4+h_5,
\]
where $f_1$, $f_2$, $f_3$, $f_4$, $g_1$, $g_2$, $g_3$,
$h_1$, $h_2$, $h_3$, $h_4$, $h_5$, are weighted homogeneous
polynomials in $(z, \zeta, w)$ of degrees equal to their
indices,
which stabilizes the normalization up to order $5$ included:
\[
\footnotesize
\aligned
{}
&
z\overline{z}
+
\tfrac{1}{2}\,\overline{z}^2\zeta
+
\tfrac{1}{2}\,z^2\overline{\zeta}
+
z\overline{z}\zeta\overline{\zeta}
+
\tfrac{1}{2}\,\overline{z}^2\zeta\zeta\overline{\zeta}
+
\tfrac{1}{2}\,z^2\overline{\zeta}\zeta\overline{\zeta}
+
F_{3,0,0,2,0}\,z^3\overline{\zeta}^2
+
\overline{F_{3,0,0,2,0}}\,\overline{z}^3\zeta^2
+
{\rm O}_{z,\zeta,\overline{z},\overline{\zeta},v}(6)
\\
\longrightarrow
\ \ \ \ \
&
z'\overline{z}'
+
\tfrac{1}{2}\,{\overline{z}'}^2\zeta'
+
\tfrac{1}{2}\,{z'}^2\overline{\zeta}'
+
z'\overline{z}'\zeta'\overline{\zeta}'
+
\tfrac{1}{2}\,{\overline{z}'}^2\zeta'\zeta'\overline{\zeta}'
+
\tfrac{1}{2}\,{z'}^2\overline{\zeta}'\zeta'\overline{\zeta}'
+
F_{3,0,0,2,0}'\,{z'}^3{\overline{\zeta}'}^2
+
\overline{F_{3,0,0,2,0}'}\,{\overline{z}'}^3{\zeta'}^2
+
{\rm O}_{z',\zeta',\overline{z}',\overline{\zeta}',v'}(6),
\endaligned
\]
is of the form:
\[
\footnotesize
\aligned
z'
&
\,:=\,
\lambda\,z
-
i\,\lambda\alpha\,z^2
-
i\,\lambda\overline{\alpha}\,w
-
\lambda\alpha^2\,
z^3
+
\Big(
i\,\lambda r
-
3\,\lambda\alpha\overline{\alpha}
+
2i\,\lambda\alpha\,
F_{3,0,0,2,0}
-
2i\,\lambda\overline{\alpha}\,\overline{F_{3,0,0,2,0}}
\Big)\,
zw
+
i\,\lambda\alpha\,
\zeta w
\\
&
\ \ \ \ \ \ \ \ \ \ \ \ \ \ \
+ 
i\,\lambda\alpha^3\,
z^4
+
\Big(
8i\,\lambda\alpha^2\overline{\alpha}
+
\tfrac{1}{2}\,\frac{\lambda^2}{\overline{\lambda}}\,\overline{\gamma}
+
4\,\frac{\lambda}{\overline{\lambda}}\,\overline{\tau}
+
4\,\lambda\alpha^2\,F_{3,0,0,2,0}
-
8\,\lambda\alpha\overline{\alpha}\,\overline{F_{3,0,0,2,0}}
\Big)\,
z^2w
+
3\,\lambda\alpha^2\,
z\zeta w
+
\tau\,w^2,
\\
\zeta'
&
\,:=\,
\frac{\lambda}{\overline{\lambda}}\,\zeta
+
2i\,\frac{\lambda}{\overline{\lambda}}\,\overline{\alpha}\,
z
+
\Big(
3\,\frac{\lambda}{\overline{\lambda}}\,\alpha\overline{\alpha}\,
-
i\,\frac{\lambda}{\overline{\lambda}}\,r
-
2i\,\frac{\lambda}{\overline{\lambda}}\,\alpha\,
F_{3,0,0,2,0}
+
6i\,\frac{\lambda}{\overline{\lambda}}\,\overline{\alpha}\,
\overline{F_{3,0,0,2,0}}
\Big)\,
z^2
-
2i\,\frac{\lambda}{\overline{\lambda}}\,\alpha\,
z\zeta
+
\frac{\lambda}{\overline{\lambda}}\,\overline{\alpha}^2\,
w
\\
&
\ \ \ \ \ \ \ \ \ \ \ \ \ \ \
+
\Big(
2\,\frac{\lambda}{\overline{\lambda}}\,\alpha r
-
4i\,\frac{\lambda}{\overline{\lambda}}\,\alpha^2\overline{\alpha}
-
2\,\frac{\lambda^2}{\overline{\lambda}^2}\,\overline{\gamma}
-
8\,\frac{\lambda}{\overline{\lambda}^2}\,\overline{\tau}
+
12\,\frac{\lambda}{\overline{\lambda}}\,\alpha^2\,
F_{3,0,0,2,0}
+
4\,\frac{\lambda}{\overline{\lambda}}\,\alpha\overline{\alpha}\,
\overline{F_{3,0,0,2,0}}
\Big)\,z^3
-
3\,\frac{\lambda}{\overline{\lambda}}\,\alpha^2\,
z^2\zeta
+
\gamma\,zw
\\
&
\ \ \ \ \ \ \ \ \ \ \ \ \ \ \
+
\Big(
-2\,\frac{\lambda}{\overline{\lambda}}\,
\alpha\overline{\alpha}
+
4i\,\frac{\lambda}{\overline{\lambda}}\,\alpha\,
F_{3,0,0,2,0}
-
4i\,\frac{\lambda}{\overline{\lambda}}\,\overline{\alpha}\,
\overline{F_{3,0,0,2,0}}
\Big)\,
\zeta w,
\\
w'
&
\,:=\,
\lambda\overline{\lambda}\,w
-
2i\,\lambda\overline{\lambda}\alpha\,
zw
-
3\,\lambda\overline{\lambda}\alpha^2\,
z^2w
+
\big(
-\lambda\overline{\lambda}\alpha\overline{\alpha}
+
i\,\lambda\overline{\lambda}\,r
\big)\,w^2
+
4i\,\lambda\overline{\lambda}\alpha^3\,
z^3w
\\
&
\ \ \ \ \ \ \ \ \ \ \ \ \ \ \
+
\Big(
6i\,\lambda\overline{\lambda}\alpha^2\overline{\alpha}
+
2\,\lambda\overline{\lambda}\alpha r
+
2\,\lambda\overline{\tau}
+
4\,\lambda\overline{\lambda}\alpha^2\,
F_{3,0,0,2,0}
-
4\,\lambda\overline{\lambda}\alpha\overline{\alpha}\,
\overline{F_{3,0,0,2,0}}
\Big)\,
zw^2
+
\lambda\overline{\lambda}\alpha^2\,\zeta w^2.
\endaligned
\]
where $\lambda \in \C^\ast$, $\alpha \in \C$, $r \in \R$,
$\gamma \in \C$, $\tau \in \C$ are arbitrary
parameters.\qed
\end{Lemma}

In comparison to the normalization up to order $4$, observe
that the previous
two supplementary parameters {\em have now been normalized}:
\[
\aligned
\beta
&
\,:=\,
\frac{\lambda}{\overline{\lambda}}\,
\overline{\alpha}^2,
\\
\varepsilon
&
\,:=\,
-\,2i\,\frac{\lambda}{\overline{\lambda}}\,\alpha\,
F_{3,0,0,2,0}
+
6i\,\frac{\lambda}{\overline{\lambda}}\,\overline{\alpha}\,
\overline{F_{3,0,0,2,0}}
+
3\,\frac{\lambda}{\overline{\lambda}}\,\alpha\overline{\alpha}
-
i\,\frac{\lambda}{\overline{\lambda}}\,r.
\endaligned
\]
With this, the {\em linearized} isotropy has become the 
{\em same} as the one of the GM-model written above:
\leqnomode\usetagform{default}
\begin{align}
\label{linearized-GM-order-1}
z'
&
\,:=\,
\lambda\,z
-
i\,\lambda\overline{\alpha}\,w,
\notag
\\
\zeta'
&
\,:=\,
\frac{\lambda}{\overline{\lambda}}\,
\zeta
+
2i\,\frac{\lambda}{\overline{\lambda}}\,\overline{\alpha}\,
z
+
\frac{\lambda}{\overline{\lambda}}\,\overline{\alpha}^2\,
w,
\\
w'
&
\,:=\,
\lambda\overline{\lambda}\,
w.
\notag
\end{align}
This key fact will enable us to define,
at {\em every} point of {\em any}
$\mathfrak{C}_{2,1}$ hypersurface $M^5 \subset \C^3$,
a CR-invariant $1$-jet
locus $\Sigma_p^1 \subset J_{M,p}^1$ in the bundle
of CR-transversal $1$-jets of $\mathcal{C}^\omega$ 
curves $\gamma \subset M$.

We will follow the guide~{\cite{Merker-2020}}, 
which was prepared in advance on this purpose.

\Section{\bf Point Translations of $\mathcal{C}^\omega$
Hypersurfaces $M^5 \subset \C^3$}
\label{point-translations-M5-C3}
\HEAD{{\ref{point-translations-M5-C3}}.~{\sf 
Point Translations of $\mathcal{C}^\omega$
Hypersurfaces $M^5 \subset \C^3$}
}{
Wei-Guo {\sc Foo}, Joël {\sc Merker}, The-Anh {\sc Ta}}

Consider as before a local $\mathcal{C}^\omega$
hypersurface $M^5 \subset \C^3$ which is
$2$-nondegenerate and of constant Levi rank $1$,
namely belongs to the class $\mathfrak{C}_{2,1}$.

In coordinates $(z, \zeta, w) = (x+iy,\, s+it,\, u+iv)$,
assume that $M$ is locally graphed as $u = F(z, \zeta,
\overline{z}, \overline{\zeta}, v)$. At all points
$p = (z_p, \zeta_p, w_p) \in M$ with $u_p = 
F\big( z_p, \zeta_p, \overline{z}_p, \overline{\zeta}_p,
v_p \big)$, let us expand up to weighted order $5$:
\[
\!\!\!\!\!\!\!\!\!\!\!\!\!\!\!
\aligned
u
\,=\,
F\big(
z,\zeta,\overline{z},\overline{\zeta},v
\big)
\,=\,
\sum_{a+b+c+d+2e\leqslant 5}\,
\tfrac{(z-z_p)^a}{a!}\,
\tfrac{(\zeta-\zeta_p)^b}{b!}\,
\tfrac{(\overline{z}-\overline{z}_p)^c}{c!}\,
\tfrac{(\overline{\zeta}-\overline{\zeta}_p)^d}{d!}\,
\tfrac{(v-v_p)^e}{e!}\,\,
F_{z^a\zeta^b\overline{z}^c\overline{\zeta}^dv^e}
\big(
z_p,\zeta_p,\overline{z}_p,\overline{\zeta}_p,v_p
\big)
+
{\rm O}(6),
\endaligned
\]
subtract $u - u_p$, translate coordinates $z := z-z_p$,
$\zeta := \zeta - \zeta_p$, $w := w - w_p$,
and get a family of hypersurfaces $M^p \subset \C^3$,
parametrized by $p \in M$ and passing through the origin:
\[
u
\,=\,
F^p\big(
z,\zeta,\overline{z},\overline{\zeta},v
\big)
\,\,=\,\,
\sum_{1\leqslant a+b+c+d+2e\leqslant 5}\,
z^a\zeta^b\overline{z}^c\overline{\zeta}^dv^e\,
F_{a,b,c,d,e}^p
+
{\rm O}_{z,\zeta,\overline{z},\overline{\zeta},v}(6),
\]
namely with $F^p(0,0,0,0,0) = 0$, whose graphing
function has coefficients:
\[
F_{a,b,c,d,e}^p
\,:=\,
\tfrac{1}{a!}\,
\tfrac{1}{b!}\,
\tfrac{1}{c!}\,
\tfrac{1}{d!}\,
\tfrac{1}{e!}\,
F_{z^a\zeta^b\overline{z}^c\overline{\zeta}^dv^e}
\big(
z_p,\zeta_p,\overline{z}_p,\overline{\zeta}_p,v_p
\big),
\]
analytically parametrized by $p \in M$. 
Thanks to this, working at {\em only one} point,
namely at the origin, we will treat {\em all} points
$p \in M$.

\begin{Question}
{\sl Are there analogs, on hypersurfaces $M^5 \in \mathfrak{C}_{2,1}$,
of Cartan-Moser chains~{\cite{Cartan-1932-I, Cartan-1932-I, 
Jacobowitz-1990, Merker-2020}}
for Levi nondegenerate 
hypersurfaces $M^3 \subset \C^2$\text{\bf ?}}
\end{Question}

Thanks to Lemma~{\ref{Lm-isotropy-order-5}}, we will construct,
at each point $p \in M$, 
an invariant surface in the bundle of $1$-jets of
CR-transversal curves
in $M$.
So there will be an important difference with
Cartan-Moser chains for
Levi nondegenerate $M^3 \subset \C^2$: the phenomenon
that there exists a CR-transversal
invariant object which is of {\em order $1$}.

To view this object, similarly as in~{\cite{Merker-2020}},
we need to introduce bundles $J_M^1$ and $J_M^2$ of
$1$-jets and $2$-jets of CR-transversal curves
$\gamma \colon \R \longrightarrow M$ with
$\dot{\gamma} \not\in T_\gamma^cM$ nowhere complex-tangential.

\Section{\bf CR-Invariant $1$-Jets $2$-codimensional
Submanifold $\Sigma^1 \subset
J_M^1 \cong M^5 \times \R^4$}
\label{CR-invariant-1-jet-surface}
\HEAD{{\ref{CR-invariant-1-jet-surface}}.~{\sf CR-Invariant $1$-Jets 
$2$-Codimensional Submanifold 
$\Sigma^1 \subset J_M^1 \cong M^5 \times \R^4$}
}{
Wei-Guo {\sc Foo}, Joël {\sc Merker}, The-Anh {\sc Ta}}

In local coordinates for which $M$ is locally graphed as 
$u = F (z, \zeta, \overline{z}, \overline{\zeta}, v)$,
at any point $p \in M$, the CR-transversal curves can be
parametrized as:
\[
v
\,\,\longmapsto\,\,
\big(
x(v),\,
y(v),\,
s(v),\,
t(v),\,
v
\big)
\,\in\,
\R_{x,y,z,t,v}^5
\]
with $\gamma(0) = p = (x_p, y_p, s_p, t_p, v_p)$.

The $4 + 4 = 8$ independent coordinates 
corresponding to 
the first derivatives $\big( \dot{x}(v), \dot{y}(v), \dot{s}(v),
\dot{t}(v) \big)$ and to
the second derivatives $\big( \ddot{x}(v),
\ddot{y}(v), \ddot{s}(v), \ddot{t}(v) \big)$
will be denoted as follows:
\[
\aligned
J_M^1
&
\,:=\,
\big\{
\big(
x_p,y_p,s_p,t_p,v_p,\,
x_p^1,y_p^1,s_p^1,t_p^1,v_p^1
\big)
\big\}
\,\,=\,\,
\R^{5+4},
\\
J_M^2
&
\,:=\,
\big\{
\big(
x_p,y_p,s_p,t_p,v_p,\,
x_p^1,y_p^1,s_p^1,t_p^1,v_p^1,\,
x_p^2,y_p^2,s_p^2,t_p^2,v_p^2
\big)
\big\}
\,\,=\,\,
\R^{5+4+4}.
\endaligned
\]

Now, denote the translation map as:
\[
\tau_p
\colon\ \ \ \ \
(z,\zeta,w)
\,\,
\xrightarrow[{\rule[0pt]{50pt}{0pt}}]{}
\,\,
\big(
z-z_p,\,\,
\zeta-\zeta_p,\,\,
w-w_p
\big)
\,\,=:\,\,
(z,\zeta,w),
\]
so that:
\[
\tau_p
\big(M,p\big)
\,=:\,
\big(M^p,0\big).
\]

\begin{center}
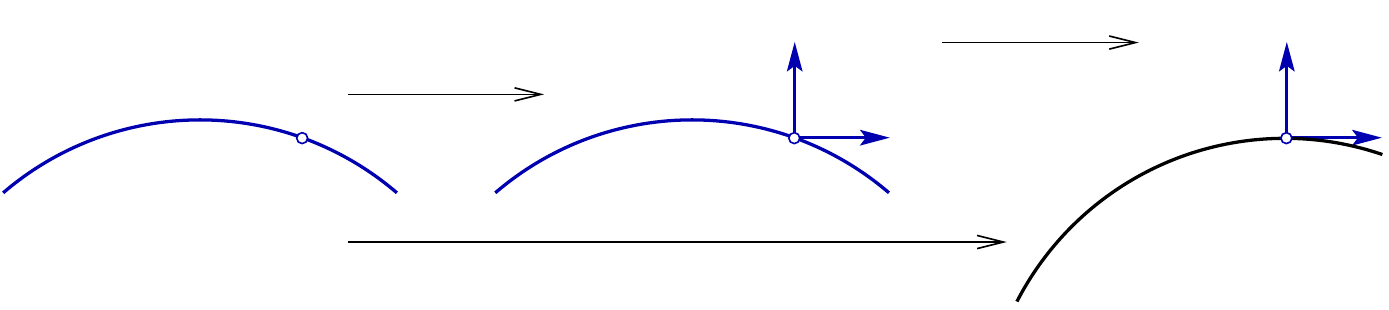
\end{center}

Also, let the punctual (at the origin) normalization map
constructed up to now, by Proposition~{\ref{Prp-F-origin-order-5}},
be denoted by:
\[
\aligned
\Phi_p
\colon\ \ \ \ \
&
(M^p,0)
\,=\,
\Big\{
u
\,=\,
\sum_{1\leqslant a+b+c+d+2e\leqslant 5}\,
F_{a,b,c,d,e}^p\,
z^a\zeta^b\overline{z}^c\overline{\zeta}^dv^e
+
{\rm O}(6)
\Big\}
\\
&
\xrightarrow[{\rule[0pt]{50pt}{0pt}}]{}
\ \ \ \ \
(N^p,0)
\,=\,
\Big\{
u
\,=\,
z\overline{z}
+
\tfrac{1}{2}\,\overline{z}^2\zeta
+
\tfrac{1}{2}\,z^2\overline{\zeta}
+
z\overline{z}\zeta\overline{\zeta}
+
\tfrac{1}{2}\,\overline{z}^2\zeta\zeta\overline{\zeta}
+
\tfrac{1}{2}\,z^2\overline{\zeta}\zeta\overline{\zeta}
\\
&
\ \ \ \ \ \ \ \ \ \ \ \ \ \ \ \ \ \ \ \ \ \ \ \ \ \ \ \ \ \ \ \ \ \ \
\ \ \ \ \ \ \ \ \ \ \ \ \ \ \ \ \ \ \ \ \ \ \ \ \ \ \ \ \ \ \ \ \ \ \
+
z^3\overline{\zeta}^2\,
F_{3,0,0,2,0}^p
+
\overline{z}^3\zeta^2\,
\overline{F_{3,0,0,2,0}^p}
+
{\rm O}(6)
\Big\}.
\endaligned
\]
According to the 
constructions done in
Sections~{\ref{chain-straightening-harmonic-killing}}
to~{\ref{normalization-F-3-0-1-1}}
and according to Proposition~{\ref{Prp-F-origin-order-5}},
we know that $\Phi_p$ depends analytically on $p$.

Abbreviate:
\[
\varphi
\,:=\,
\Phi_p\circ\tau_p,
\]
and consider the diagram:
\[
\xymatrix{
J_{M,p}^1
\ar[rr]^{\varphi^{(1)}}
\ar[d]
&
&
J_{N^p,0}^1
\ar[d]
\\
(M,p)
\ar[rr]_\varphi
&
&
(N^p,0).
}
\]
As in Observation~{\ref{Obs-Sigma-0-jets-1}},
in the $1$-jet fiber above $0 \in N^p$,
introduce the surface:
\[
\Sigma_0^1
\,:=\,
\big\{
(x_1,y_1,s_1,t_1)
\in
J_{N^p,0}^1
\colon\,\,
s_1
=
-2\,x_1y_1,\,\,\,
t_1
=
x_1^2-y_1^2
\big\}.
\]
Using the first prolongation $\varphi^{(1)}$, define the
$2$-dimensional submanifold of $J_{M,p}^1$:
\[
\Sigma_p^1
\,:=\,
{\varphi^{(1)}}^{-1}
\big(
\Sigma_0^1
\big).
\]
Since $\varphi^{(1)}$ is a diffeomorphism $J_{M,p}^1 
\overset{\sim}{\longrightarrow} J_{N^p,0}^1$,
this $\Sigma_p^1$ is also graphed, say of the form:
\[
s_p^1
\,=\,
A\big(x_p^1,y_p^1\big),
\ \ \ \ \ \ \ \ \ \ \ \ \ \ \ \ \ \ \ \ \ \ \ \ \ \
t_p^1
\,=\,
B\big(x_p^1,y_p^1\big),
\]
with two $\mathcal{C}^\omega$ functions 
$A$, $B$ which depend on $p$, and depend {\em also}
a priori on the normalizing map $\varphi$.
\[
\xymatrix{
\Sigma^1
\ar@{^{(}->}[r]
&
J_{M,p}^1
\ar[rr]^{\varphi^{(1)}}
\ar[d]
& &
J_{N^p,0}^1
\ar[d]
&
\ar@{^{(}->}[l]_f
\Sigma^1
\ar@/_3pc/[llll]_{{\varphi^{(1)}}^{-1}}
\\
&
(M,p)
\ar[rr]^{\varphi}
& &
(N^p,0)
&
}
\]

The union:
\[
\underset{p\in M}{\medcup}\,
\Sigma_p^1
\,=:\,
\Sigma^1
\,\,\subset\,\,
J_M^1
\]
is a $\mathcal{C}^\omega$ submanifold of dimension $5 + 2$ within
$J_M^1$ which has dimension $5 + 4$.

\begin{Assertion}
This graphed surface $\Sigma_p^1 \subset J_{M,p}^1 
\cong \R^4$ is {\em independent} of the map
$\varphi = \Phi_p \circ \tau_p$ normalizing the initial
hypersurface $M$
of equation $u = F(z, \zeta, \overline{z}, \overline{\zeta}, v)$
near any of its points $p \in M$, to:
\[
\aligned
u
&
\,=\,
z\overline{z}
+
\tfrac{1}{2}\,\overline{z}^2\zeta
+
\tfrac{1}{2}\,z^2\overline{\zeta}
+
z\overline{z}\zeta\overline{\zeta}
+
\tfrac{1}{2}\,\overline{z}^2\zeta\zeta\overline{\zeta}
+
\tfrac{1}{2}\,z^2\overline{\zeta}\zeta\overline{\zeta}
\\
&
\ \ \ \ \ \ \ \ \ \ \
\ \ \ \ \ \ \ \ \ \ \ \ \ \ \ \ \ \ \ \ \ \ \ \ \ \
+
z^3\overline{\zeta}^2\,F_{3,0,0,2,0}^p
+
\overline{z}^3\zeta^2\,\overline{F_{3,0,0,2,0}^p}
+
{\rm O}_{z,\zeta,\overline{z},\overline{\zeta},v}(6).
\endaligned
\]
\end{Assertion}

\proof
Suppose another such normalizing map is given:
\[
\xymatrix{
& &
(N^p,0)
\ar[ddr]^{\psi\,:=\,\varphi_\prime\circ\varphi^{-1}}
\\
(M,p)
\ar[rru]^{\varphi}
\ar[rrrd]^{\varphi_\prime}
&
&
&
\\
& & &
(N_\prime^p,0).
}
\]
By Lemma~{\ref{Lm-isotropy-order-5}}
which holds for maps stabilizing the origin,
$\psi$ has linear terms exactly equal to the linear
terms of the isotropy group of the GM-model, for
which we already know, thanks to 
Observation~{\ref{Obs-Sigma-0-jets-1}}, that:
\[
\psi^{(1)}
\big(
\Sigma_0^1
\big)
\,=\,
\Sigma_0^{\prime1}.
\]
Hence in conclusion:
\begin{align*}
\Sigma_p^{\prime1}
\,=\,
{\varphi_\prime^{(1)}}^{-1}
\big(
\Sigma_0^{\prime1}
\big)
\,=\,
{\varphi_\prime^{(1)}}^{-1}
\Big(
\psi^{(1)}
\big(
\Sigma_0^1
\big)
\Big)
&
\,=\,
{\varphi_\prime^{(1)}}^{-1}
\Big(
\big(
\varphi_\prime
\circ
\varphi^{-1}
\big)^{(1)}
\big(
\Sigma_0^1
\big)
\Big)
\\
&
\,=\,
{\varphi^{(1)}}^{-1}
\big(
\Sigma_0^1
\big)
\,\,=\,\,
\Sigma_p^1.
\qedhere
\end{align*}
\endproof

So at each point $p \in M$, there exists a CR-invariant,
or biholomorphically invariant,
surface $\Sigma_p^1 \subset J_{M,p}^1$.
Therefore, it is natural to select only
CR-transversal curves $\gamma \colon \R \longrightarrow
M$, $\gamma(0) = p$, such that $\dot{\gamma}(\tau) 
\in \Sigma_{\gamma(\tau)}^1$ for every $\tau \in \R$.

But the `discovery' of this CR-invariant submanifold
$\Sigma_M^1 \subset J_M^1$ 
does not suffice, because the linear action:
\[
\aligned
z'
&
\,:=\,
\lambda\,z
-
i\,\lambda\overline{\alpha}\,w,
\\
\zeta'
&
\,:=\,
\frac{\lambda}{\overline{\lambda}}\,
\zeta
+
2i\,\frac{\lambda}{\overline{\lambda}}\,\overline{\alpha}\,
z
+
\frac{\lambda}{\overline{\lambda}}\,\overline{\alpha}^2\,
w,
\\
w'
&
\,:=\,
\lambda\overline{\lambda}\,
w,
\endaligned
\]
happens to be {\em transitive} on the invariant surface
$\Sigma_0^1 \subset \R^4$ of $1$-jets, according to the fact
that the prolonged
symmetry vector fields 
${\sf D}^{(1)}$, ${\sf R}^{(1)}$, ${\sf I}_1^{(1)}$,
${\sf I}_2^{(1)}$, ${\sf J}^{(1)}$,
shown in Section~{\ref{prolongations-jet-1}},
are of rank $2 = \dim\, \Sigma_0^1$ everywhere.

Remind from~{\cite{Cartan-1932-I, Cartan-1932-II, 
Jacobowitz-1990, Merker-2020}}
that Cartan-Moser chains were 
strictly of {\em second order}.
Hence, we need to explore deeper, 
and to normalize further, still at $0 \in M^p$.
We will realize that to each $1$-jet
$j_p^1 \in \Sigma_p^1$,
there is associated a unique invariant $2$-jet
$j_p^2 = j_p^2 ( j_p^1)$, as we already saw when
studying the GM-model in 
Section~{\ref{prolongations-jet-2}}.

\Section{\bf Order $1$ Chains in $\mathfrak{C}_{2,1}$ 
Hypersurfaces $M^5 \subset \C^3$}
\label{order-1-chains}
\HEAD{{\ref{order-1-chains}}.~{\sf Order $1$ Chains 
in $\mathfrak{C}_{2,1}$ Hypersurfaces $M^5 \subset \C^3$
}
}{
Wei-Guo {\sc Foo}, Joël {\sc Merker}, The-Anh {\sc Ta}}

So far, at the origin,  
we have constructed a normalizing map $\Phi_p$,
composed with a translation map $\tau_p$:
\[
\varphi
\colon\ \ \ \ \
(M,p)
\,\,\,\xrightarrow[{\rule[0pt]{50pt}{0pt}}]{\tau_p}\,\,\,
(M^p,0)
\,\,\,\xrightarrow[{\rule[0pt]{50pt}{0pt}}]{\Phi_p}\,\,\,
(N^p,0),
\]
which brings $(M, p)$ to $(N^p, 0)$ at the origin of equation
fully normalized up to order $5$ included:
\[
\aligned
u
&
\,=\,
z\overline{z}
+
\tfrac{1}{2}\,\overline{z}^2\zeta
+
\tfrac{1}{2}\,z^2\overline{\zeta}
+
z\overline{z}\zeta\overline{\zeta}
+
\tfrac{1}{2}\,\overline{z}^2\zeta\zeta\overline{\zeta}
+
\tfrac{1}{2}\,z^2\overline{\zeta}\zeta\overline{\zeta}
\\
&
\ \ \ \ \ \ \ \ \ \ \ \ \ \ \ \ \ \ \ \ \ \ \ \ \ \ \ \ \ \ \ \ \ \ \
+
2\,\Re\,
\Big\{
0
+
0
+
z^3\overline{\zeta}^2\,
F_{3,0,0,2,0}^p
\big\}
+
{\rm O}_{z,\zeta,\overline{z},\overline{\zeta},v}(6),
\endaligned
\]
namely with $0 = F_{3,0,1,1,0}^p = F_{4,0,0,1,0}^p$,
knowing that $F_{3,0,0,2,0}^p$ is a relative invariant. 

The differential $\varphi_\ast$ establishes isomorphisms:
\[
\aligned
T_pM
&
\,\,\,\overset{\sim}{\longrightarrow}\,\,\,
T_0N^p,
\\
T_p^cM
&
\,\,\,\overset{\sim}{\longrightarrow}\,\,\,
T_0^cN^p,
\\
K_p^cM
&
\,\,\,\overset{\sim}{\longrightarrow}\,\,\,
K_0^cN^p,
\endaligned
\]
where $K^cM \subset T^cM$ is the Levi-kernel 
subbundle~{\cite{Merker-Pocchiola-Sabzevari-2013-5-CR-II}}.
It follows that $\varphi_\ast$ establishes an isomorphism
between the $3$-dimensional real quotient bundles:
\[
T_pM
\Big/
\big(
T_p^cM\big/K_p^cM
\big)
\,\,\,\overset{\sim}{\longrightarrow}\,\,\,
T_0N^p
\Big/
\big(
T_0^cN^p\big/K_0^cN^p
\big).
\]
By definition, on these bundles $T^c / K^c$, the Levi form 
of $M$ is nondegenerate,
of maximal possible rank $1$.

In a neighborhood of some reference point $p_0 \in M$,
we can take coordinates $(z, w, \zeta)$ with $z = x + iy$,
$\zeta = s + it$, $w = u + iv$, so that
$M$ is locally graphed as $u = F(z, \zeta, \overline{z}, 
\overline{\zeta}, v)$, with $(v, x, y, s, t) \in M^5$
being intrinsic coordinates, so that the
Levi form of $M$ is {\em nonzero near $p_0$} along the  
intrinsic $(1,0)$ vector field:
\[
\mathcal{L}
\,:=\,
\frac{\partial}{\partial z}
-
i\,\frac{F_z}{1+i\,F_v}\,
\frac{\partial}{\partial v}.
\]
We will let $p \sim p_0$ vary in a neighborhood of $p_0$.

Taking jet coordinates $(x_1, y_1, s_1, t_1)$ near $p_0$ so that:
\[
J_M^1
\,=\,
\big\{
(v,x,y,s,t,\,x_1,y_1,s_1,t_1)
\big\},
\]
it follows from the above isomorphisms and from the definition
of $\Sigma_0^1 \subset J_{M^p,0}^1$ 
that $\Sigma^1 \subset J_M^1$ is locally defined near $p_0$
as a graph:
\[
s_1
\,=\,
A\big(v,x,y,s,t,\,x_1,y_1\big),
\ \ \ \ \ \ \ \ \ \ \ \ \ \ \ \ \ \ \ \ \ \ \ \ \ \
t_1
\,=\,
B\big(v,x,y,s,t,\,x_1,y_1\big),
\]
in terms of certain two $\mathcal{C}^\omega$ functions $A$, $B$,
which vanish for $x_1 = y_1 = 0$.
In this respect, the first two coordinates $(x_p^1, y_p^1)$ of a
$1$-jet $j_p^1$ at some point $p = (v_p, x_p, y_p, s_p, t_p)
\in M$ near $p_0$ should be thought of
as being {\sl horizontal}, and
the last two coordinates $(s_p^1, t_p^1)$ as being {\sl vertical}.

An alternative presentation of CR-invariant
CR-transversal $1$-jets on hypersurfaces $M^5 \subset \C^3$
will be useful in a moment.

\begin{Definition}
\label{Def-order-1-chains}
A $1$-jet $j_p^1 \in J_{M,p}^1$ is said to be the jet of an {\sl order
$1$ chain} at a point $p \in M$, or to {\em belong} to
the invariant surface $\Sigma_p^1 \subset J_{M,p}^1$,
if, given any punctual normalizing map 
from $(M, p)$ to $(N^p, 0)$ up to
order $5$ as in 
Proposition~{\ref{Prp-F-origin-order-5}}:
\[
\aligned
u
&
\,=\,
z\overline{z}
+
\tfrac{1}{2}\,\overline{z}^2\zeta
+
\tfrac{1}{2}\,z^2\overline{\zeta}
+
z\overline{z}\zeta\overline{\zeta}
+
\tfrac{1}{2}\,\overline{z}^2\zeta\zeta\overline{\zeta}
+
\tfrac{1}{2}\,z^2\overline{\zeta}\zeta\overline{\zeta}
\\
&
\ \ \ \ \ \ \ \ \ \ \ \ \ \ \ \ \ \ \ \ \ \ \ \ \ \ \ \ \ \ \ \ \ \ \
+
2\,\Re\,
\Big\{
0
+
0
+
z^3\overline{\zeta}^2\,
F_{3,0,0,2,0}^p
\big\}
+
{\rm O}_{z,\zeta,\overline{z},\overline{\zeta},v}(6),
\endaligned
\]
which sends $j_p^1$ to a $1$-jet at $0 \in N^p$ having
{\em vanishing} horizontal part:
\[
\varphi^{(1)}\big(j_p^1\big)
\,=\,
\big(0,0,\,s_0^1,t_0^1\big),
\]
then in fact $j_p^1$ is the inverse image of the {\em flat
$1$-jet} at the origin:
\[
j_p^1
\,=\,
{\varphi^{(1)}}^{-1}
\big(0,0,\,0,0\big),
\]
or equivalently $s_0^1 = t_0^1 = 0$.
\end{Definition}

This definition does not depend on the normalizing map $\Phi_p$
in $\varphi = \tau_p \circ \Phi_p$, because if another
$\Phi_p^\prime$ is chosen, which leads to the diagram:
\[
\xymatrix{
& &
(N^p,0)
\ar[ddr]^{\psi\,:=\,\varphi_\prime\circ\varphi}
\\
(M,p)
\ar[rru]^{\varphi}
\ar[rrrd]^{\varphi_\prime}
&
&
&
\\
& & &
(N_\prime^p,0),
}
\]
with $(N_\prime^p, 0)$ having an equation similar to the one
of $(N^p, 0)$ above, then the ambiguity map $\psi := 
\varphi_\prime \circ \varphi$ should stabilize the
flat $1$-jet,
and for this to hold, we already know 
from the formulas~({\ref{linearized-GM-order-1}})
that this forces $\alpha = 0$. 

\smallskip

We will now employ this definition in two ways. 
It is clear that the graphed equations of $\Sigma^1 \subset 
J_M^1$ lead to a system of two first-order ordinary
differential equations:
\[
\dot{s}
\,=\,
A
\big(
v,x,y,s,t,\dot{x},\dot{y}
\big),
\ \ \ \ \ \ \ \ \ \ \ \ \ \ \ \ \ \ \ \ \ \ \ \ \ \
\dot{t}
\,=\,
B
\big(
v,x,y,s,t,\dot{x},\dot{y}
\big),
\]
the time parameter being $v$. 
For any choice of any two functions $(x(v), y(v))$ 
with $(x(0), y(0)) = (x_p, y_p)$, 
with $(\dot{x}(0), \dot{y}(0)) \neq (0,0)$, 
and with $(s(0), t(0)) = (s_p, t_p)$,
there exists a unique local $\mathcal{C}^\omega$
solution to this system passing through $p$ at
`time' $v = 0$,
which is a CR-transversal curve having tangents
in $\Sigma_M^1$. 

\begin{Terminology}
Such a curve will be called an {\sl order $1$ chain}.
\end{Terminology}

Later, when passing to {\sl order $2$ chains},
we will see that the large freedom 
in the choice of arbitrary functions $(x(v), y(v))$
will drop.

\smallskip

Once order $1$ chains are known, it is natural
to restart the whole process of prenormalization
and of partial normalization which begun in 
Section~{\ref{chain-straightening-harmonic-killing}},
by assuming that the CR-transversal
curve $p \in \gamma \subset M$
(not anymore chosen at random)
{\em is} an order $1$ chain.

\smallskip

Then, coming back to
Proposition~{\ref{Prp-normalization-order-5-v-axis}},
but viewed at the origin up to order $6$
in {\em all} variables $(z, \zeta, \overline{z}, \overline{\zeta},
v)$,  
we remember that we have constructed a normalizing map $\Phi_p$,
composed with a translation map $\tau_p$:
\[
\varphi
\colon\ \ \ \ \
(M,p)
\,\,\,\xrightarrow[{\rule[0pt]{50pt}{0pt}}]{\tau_p}\,\,\,
(M^p,0)
\,\,\,\xrightarrow[{\rule[0pt]{50pt}{0pt}}]{\Phi_p}\,\,\,
(N^p,0),
\]
which brings $(M, p)$ to $(N^p, 0)$ at the origin of equation:
\[
\aligned
u
&
\,=\,
z\overline{z}
+
\tfrac{1}{2}\,\overline{z}^2\zeta
+
\tfrac{1}{2}\,z^2\overline{\zeta}
+
z\overline{z}\zeta\overline{\zeta}
+
\tfrac{1}{2}\,\overline{z}^2\zeta\zeta\overline{\zeta}
+
\tfrac{1}{2}\,z^2\overline{\zeta}\zeta\overline{\zeta}
\\
&
\ \ \ \ \ \ \ \ \ \ \ \ \ \ \ \ \ \ \ \ \ \ \ \ \ \ \ \ \ \ \ \ \ \ \
+
2\,\Re\,
\Big\{
0
+
z^4\overline{\zeta}\,
F_{4,0,0,1,0}^p
+
z^3\overline{\zeta}^2\,
F_{3,0,0,2,0}^p
\big\}
+
{\rm O}_{z,\zeta,\overline{z},\overline{\zeta},v}(6),
\endaligned
\]
{\em without changing} the CR-transversal curve 
$0 \in \gamma \subset M$ being the $v$-axis,
hence having {\em flat} $1$-jet at the origin.

\begin{Assertion}
\label{Assertion-F-40010-zero}
Then $F_{4,0,0,1,0}^p = 0$ holds automatically, without
having the needs to perform any further biholomorphism.
\end{Assertion}

\proof
Indeed, we already know that one can continue to normalize
and make $F_{4,0,0,1,0}^p = 0$ by means of the map:
\[
\aligned
z'
&
\,:=\,
z
+
2\,F_{4,0,0,1,0}^p\,z^3
-
2\,F_{4,0,0,1,0}^p\,z\zeta w,
\\
\zeta'
&
\,:=\,
\zeta
-
2\,\overline{F_{4,0,0,1,0}^p}\,
w
+
10\,z^2\zeta\,
F_{4,0,0,1,0}^p,
\\
w'
&
\,:=\,
w
+
2\,z^2w\,
F_{4,0,0,1,0}^p,
\endaligned
\]
which we may call $\Psi \colon (N^p, 0) \longrightarrow
(N_\prime^p, 0)$.
We then reason as in~{\cite[9.5]{Merker-2020}}.

If $F_{4,0,0,1,0}^p \neq 0$ would be nonzero, due to the
presence in $\zeta'$ of the linear term 
$2\, \overline{F_{4,0,0,1,0}^p}\, w$, 
this map $\Psi$ would {\sl not} stabilize the flat order $1$ jet
$j_0^1 = (0, 0, 0, 0)$, and so,
this would contradict
Definition~{\ref{Def-order-1-chains}} applied
to $(M,p) := (N^p, 0)$, to $\varphi := \Psi$,
and to $(N^p, 0) := (N_\prime^p, 0)$.
\endproof

Lastly, coming again back to
Proposition~{\ref{Prp-normalization-order-5-v-axis}}, 
we remember that we have constructed a normalizing map
which brings $M$ near $0 \in M$ to the equation:
\[
\aligned
u
&
\,=\,
z\overline{z}
+
\tfrac{1}{2}\,\overline{z}^2\zeta
+
\tfrac{1}{2}\,z^2\overline{\zeta}
+
z\overline{z}\zeta\overline{\zeta}
+
\tfrac{1}{2}\,\overline{z}^2\zeta\zeta\overline{\zeta}
+
\tfrac{1}{2}\,z^2\overline{\zeta}\zeta\overline{\zeta}
\\
&
\ \ \ \ \ \ \ \ \ \ \ 
+
z^3\overline{z}^3\,
{\rm O}_{z,\overline{z}}(0)
\\
&
\ \ \ \ \ \ \ \ \ \ \ 
+
2\,\Re\,
\Big\{
0
+
z^4\overline{\zeta}\,
F_{4,0,0,1}(v)
+
z^3\overline{\zeta}^2\,
F_{3,0,0,2}(v)
\Big\}
\\
&
\ \ \ \ \ \ \ \ \ \ \ 
+
\overline{z}^3\zeta\,
{\rm O}_{z,\zeta,\overline{z}}(2)
+
z^3\overline{\zeta}\,
{\rm O}_{z,\overline{z},\overline{\zeta}}(2)
+
\zeta\overline{\zeta}\,
{\rm O}_{z,\zeta,\overline{z},\overline{\zeta}}(4).
\endaligned
\]
{\em without changing} any starting CR-transversal curve 
$0 \in \gamma \subset M$. We now realize that $F_{4,0,0,1}(v)
\equiv 0$
vanishes for free.

\begin{Proposition}
\label{Prp-full-normalization-order-5-v-axis}
For every hypersurface $M^5 \in \mathfrak{C}_{2,1}$, 
at any point $p \in M$, 
given any CR-transversal curve $p \in \gamma \subset M$
which is an order $1$ chain,
there exist holomorphic coordinates $(z, \zeta, w) \in \C^3$
vanishing at $p$ in which $\gamma$ is the $v$-axis and in which
$M$ has equation:
\[
\aligned
u
&
\,=\,
z\overline{z}
+
\tfrac{1}{2}\,\overline{z}^2\zeta
+
\tfrac{1}{2}\,z^2\overline{\zeta}
+
z\overline{z}\zeta\overline{\zeta}
+
\tfrac{1}{2}\,\overline{z}^2\zeta\zeta\overline{\zeta}
+
\tfrac{1}{2}\,z^2\overline{\zeta}\zeta\overline{\zeta}
\\
&
\ \ \ \ \ \ \ \ \ \ \ 
+
z^3\overline{z}^3\,
{\rm O}_{z,\overline{z}}(0)
\\
&
\ \ \ \ \ \ \ \ \ \ \ 
+
2\,\Re\,
\Big\{
0
+
0
+
z^3\overline{\zeta}^2\,
F_{3,0,0,2}(v)
\Big\}
\\
&
\ \ \ \ \ \ \ \ \ \ \ 
+
\overline{z}^3\zeta\,
{\rm O}_{z,\zeta,\overline{z}}(2)
+
z^3\overline{\zeta}\,
{\rm O}_{z,\overline{z},\overline{\zeta}}(2)
+
\zeta\overline{\zeta}\,
{\rm O}_{z,\zeta,\overline{z},\overline{\zeta}}(4).
\endaligned
\]
\end{Proposition}

\proof
What was done an instant ago
by Assertion~{\ref{Assertion-F-40010-zero}} 
at the origin $(z, \zeta, w) = (0, 0, 0)$
applies in fact at {\em every} point $(0, 0, iv)$ along the
$v$-axis, thanks to the fact that the
(pre)normalizations
of Propositions~{\ref{Prp-prenormalization}}
and~{\ref{Prp-normalization-order-5-v-axis}}
were achieved {\em all along} the $v$-axis.
\endproof

Because we know the existence of a CR-invariant
surface $\Sigma_p^1 \subset J_{M,p}^1$ on which
the isotropy is transitive, we will assume
that, starting with any fixed $1$-jet $j_p^1 \in \Sigma_p^1$,
the partial normalization map performed up to now
sends $j_p^1$ to the {\em flat} $1$-jet at
$0 \in M^p$, namely to $j_0^1 = (0, 0, 0, 0)$. 
We will assume that subsequent normalizations
{\em stabilize} this invariant flat $1$-jet.
For this, at the very beginning, we have to assume that
the CR-transversal curve used in 
Section~{\ref{chain-straightening-harmonic-killing}},
whose choice was left free, has $1$-jet at the origin $0$
equal to the flat $1$-jet.
By surveying all normalizations done up to now,
one realizes that the $v$-axis was always stabilized,
contained in $M$, hence the flat $1$-jet was
always preserved (implicitly). 

Preserving the flat $1$-jet at $0$
corresponds to making $\alpha := 0$ 
in the formulas
of Section~{\ref{prolongations-jet-2}}
and of Lemma~{\ref{Lm-isotropy-order-5}}.
We state this explicitly as a

\begin{Corollary}
\label{Cor-order-5-alpha-zero}
The biholomorphic maps of Lemma~{\ref{Lm-isotropy-order-5}}
which stabilize punctual normalizations of $(M^p, 0)$
at the origin up
to order $5$ {\em and} which stabilize {\em also}
the flat $1$-jet $j_0^1 = (0, 0, 0, 0) \in \Sigma_0^1$ read,
with $\alpha := 0$ and $\theta := \gamma$, as:
\reqnomode\usetagform{EngelLie}
\begin{align}
z'
&
\,:=\,
\lambda\,z
+
i\,\lambda\,r\,
zw
+
\Big(
\tfrac{1}{2}\,
\frac{\lambda^2}{\overline{\lambda}}\,
\overline{\theta}
+
4\,\frac{\lambda}{\overline{\lambda}}\,
\overline{\tau}
\Big)\,
z^2w
+
\tau\,
w^2,
\notag
\\
\zeta'
&
\,:=\,
\frac{\lambda}{\overline{\lambda}}\,
\zeta
-
i\,
\frac{\lambda}{\overline{\lambda}}\,r\,
z^2
+
\Big(
-2\,\frac{\lambda^2}{\overline{\lambda}^2}\,
\overline{\theta}
-
8\,\frac{\lambda}{\overline{\lambda}^2}\,\overline{\tau}
\Big)\,
z^3
+
\theta\,
zw,
\notag
\\
w'
&
\,:=\,
\lambda\overline{\lambda}\,w
+
i\,\lambda\overline{\lambda}\,r\,
w^2
+
2\,\lambda\overline{\tau}\,
zw^2.
\tag{\qed}
\end{align}
\end{Corollary}

\Section{\bf End of Point Normalization of $\mathcal{C}^\omega$
Hypersurfaces $M^5 \subset \C^3$}
\label{end-point-normalization-M5-C3}
\HEAD{{\ref{end-point-normalization-M5-C3}}.~{\sf 
End of Point Normalization of $\mathcal{C}^\omega$
Hypersurfaces $M^5 \subset \C^3$}
}{
Wei-Guo {\sc Foo}, Joël {\sc Merker}, The-Anh {\sc Ta}}

Thus, we have to look at 6\textsuperscript{th} order terms
in the currently normalized equation of $(M^p, 0)$,
which, taking account of the vanishing of the Levi determinant,
are of the form (exercise):
\[
\aligned
u
&
\,=\,
z\overline{z}
+
\tfrac{1}{2}\,
\overline{z}^2\zeta
+
\tfrac{1}{2}\,
z^2\overline{\zeta}
+
z\overline{z}\zeta\overline{\zeta}
+
\tfrac{1}{2}\,\overline{z}^2\zeta\zeta\overline{\zeta}
+
\tfrac{1}{2}\,z^2\overline{\zeta}\zeta\overline{\zeta}
+
z\overline{z}\zeta\overline{\zeta}\zeta\overline{\zeta}
\\
&
\ \ \ \ \
+
F_{3,0,0,2,0}^p\,z^3\overline{\zeta}^2
+
\overline{F_{3,0,0,2,0}^p}\,
\overline{z}^3\zeta^2
+
\zeta\overline{\zeta}\,
\Big(
3\,F_{3,0,0,2,0}^p\,
z^2\overline{z}\overline{\zeta}
+
3\,\overline{F_{3,0,0,2,0}^p}\,z\zeta\overline{z}^2
\Big)
\\
&
\ \ \ \ \
+
z^3\overline{z}^3\,
F_{3,0,3,0,0}^p
\\
&
\ \ \ \ \
+
2\,\Re\,\Big\{
z^5\overline{\zeta}\,
F_{5,0,0,1,0}^p
+
z^4\overline{z}\overline{\zeta}\,
F_{4,0,1,1,0}^p
+
z^4\overline{\zeta}^2\,
F_{4,0,0,2,0}^p
\\
&
\ \ \ \ \ \ \ \ \ \ \ \ \ \ \ \ \ \ \ \ \ \ \ \ \ \ \ \ \ \ \ \ \ \ \
\ \ \ \ 
+
z^3\overline{z}^2\overline{\zeta}\,
F_{3,0,2,1,0}^p
+
z^3\overline{z}\overline{\zeta}^2\,
F_{3,0,1,2,0}^p
\\
&
\ \ \ \ \ \ \ \ \ \ \ \ \ \ \ \ \ \ \ \ \ \ \ \ \ \ \ \ \ \ \ \ \ \ \
\ \ \ \ \ \ \ \ \ \ \ \ \ \ \ \ \ \ \ \ \ \ \ \ \ \ \ 
\ \ \ \ \
+
z^3\overline{\zeta}^3\,
F_{3,0,0,3,0}^p
\Big\}
+
{\rm O}_{z,\zeta,\overline{z},\overline{\zeta},v}(7).
\endaligned
\]

To normalize further order $6$ terms, it is natural 
to assume that the normalizations up to order $5$ included
are stabilized, {\em and also} that the flat $1$-jet
at the origin is stabilized
as well. Thus we will employ maps of the form:
\[
\aligned
z'
&
\,:=\,
\lambda\,z
+
i\,\lambda\,r\,
zw
+
\Big(
\tfrac{1}{2}\,
\frac{\lambda^2}{\overline{\lambda}}\,
\overline{\theta}
+
4\,\frac{\lambda}{\overline{\lambda}}\,
\overline{\tau}
\Big)\,
z^2w
+
\tau\,
w^2
+
\sum_{a+b+2e=5}\,
f_{a,b,e}\,
z^a\zeta^bw^e,
\\
\zeta'
&
\,:=\,
\frac{\lambda}{\overline{\lambda}}\,
\zeta
-
i\,
\frac{\lambda}{\overline{\lambda}}\,r\,
z^2
+
\Big(
-2\,\frac{\lambda^2}{\overline{\lambda}^2}\,
\overline{\theta}
-
8\,\frac{\lambda}{\overline{\lambda}^2}\,\overline{\tau}
\Big)\,
z^3
+
\theta\,
zw
+
\sum_{a+b+2e=4}\,
g_{a,b,e}\,
z^a\zeta^bw^e,
\\
w'
&
\,:=\,
\lambda\overline{\lambda}\,w
+
i\,\lambda\overline{\lambda}\,r\,
w^2
+
2\,\lambda\overline{\tau}\,
zw^2
+
\sum_{a+b+2e=6}\,
h_{a,b,e}\,
z^a\zeta^bw^e.
\endaligned
\]

\begin{Lemma}
\label{Lm-kill-F-30300-F40110}
One can annihilate:
\[
F_{3,0,3,0,0}^p
\,=\,
0
\ \ \ \ \ \ \ \ \ \ \ \ \ \ \ \ \ \ \ \
\text{and}
\ \ \ \ \ \ \ \ \ \ \ \ \ \ \ \ \ \ \ \
\Big(
\text{either}\ \ \ \ \
F_{4,0,1,1,0}^p
\,=\,
0
\ \ \ \ \ 
\text{or}
\ \ \ \ \ 
F_{3,0,2,1,0}^p
\,=\,
0
\Big).
\]
\end{Lemma}

\proof
By hand or on a computer, one verifies that the map:
\[
\aligned
z'
&
\,:=\,
z
+
\tfrac{3}{4}\,
F_{3,0,3,0,0}^p\,
zw^2,
\\
\zeta'
&
\,:=\,
\zeta,
\\
w'
&
\,:=\,
w
+
\big(
\tfrac{1}{4}\,
F_{3,0,3,0,0}^p
+
\overline{F_{3,0,3,0,0}^p}
\big)\,
w^3,
\endaligned
\]
makes $F_{3,0,3,0,0}^{p\,\prime} = 0$. It is visible (eyes exercise)
that this map stabilizes the flat $1$-jet $j_0^1 = (0, 0, 0, 0)$.

Next, assuming that $F_{3,0,3,0,0}^p = 0 = F_{3,0,3,0,0}^{p\,\prime}$,
the map parametrized by $\tau \in \C$:
\[
\aligned
z'
&
\,:=\,
z
+
2\,\overline{\tau}\,z^2w
+
\tau\,w^2
-
\overline{\tau}\,\zeta w^2,
\\
\zeta'
&
\,:=\,
\zeta
-
4\,\tau\,zw
+
4\,\overline{\tau}\,z\zeta w,
\\
w'
&
\,:=\,
w
+
2\,\overline{\tau}\,zw^2,
\endaligned
\]
also stabilizes the flat $1$-jet $j_0^1 = (0, 0, 0, 0)$, and it
transforms as follows the six remaining coefficients:
\[
\aligned
F_{5,0,0,1,0}^{p\,\prime}
\,=\,
F_{5,0,0,1,0}^p
\ \ \ \ \ \ \ \ \ \ \ \ \ \ \ \ \ \ \ \
F_{4,0,1,1,0}^{p\,\prime}
\,=\,
F_{4,0,1,1,0}^p
-
2\,\overline{\tau},
\ \ \ \ \ \ \ \ \ \ \ \ \ \ \ \ \ \ \ \
F_{4,0,0,2,0}^{p\,\prime}
\,=\,
F_{4,0,0,2,0}^p,
\\
F_{3,0,2,1,0}^{p\,\prime}
\,=\,
F_{3,0,2,1,0}^p
+
2\,\tau,
\ \ \ \ \ \ \ \ \ \ \ \ \ \ \ \ \ \ \ \
F_{3,0,1,2,0}^{p\,\prime}
\,=\,
F_{3,0,1,2,0}^p,
\\
F_{3,0,0,3,0}^{p\,\prime}
\,=\,
F_{3,0,0,3,0}^p.
\endaligned
\]
So one of the two mentioned coefficients can be normalized.
\endproof

A choice must be made. We then determine the stability group
for both choices of normalizations, again
with the constraint of stabilizing the flat $1$-jet $j_0^1$.
Both choices lead to the same 
stability group (exercise on a computer).

\begin{Lemma}
\label{Lm-isotropy-order-6}
Any biholomorphic map of the form:
\[
z'
\,:=\,
f_1+f_2+f_3+f_4+f_5,
\ \ \ \ \ \ \ \ \ \ \ \ \ \ \ \ \ \ \ \
\zeta'
\,:=\,
g_1+g_2+g_3+g_4,
\ \ \ \ \ \ \ \ \ \ \ \ \ \ \ \ \ \ \ \
w'
\,:=\,
h_1+h_2+h_3+h_4+h_5+h_6,
\]
where $f_1$, $f_2$, $f_3$, $f_4$, $f_5$, $g_1$, $g_2$, $g_3$, $g_4$,
$h_1$, $h_2$, $h_3$, $h_4$, $h_5$, $h_6$, are weighted homogeneous,
which stabilizes the normalization up to order $6$ included:
\[
\aligned
u
&
\,=\,
z\overline{z}
+
\tfrac{1}{2}\,\overline{z}^2\zeta
+
\tfrac{1}{2}\,z^2\overline{\zeta}
+
z\overline{z}\zeta\overline{\zeta}
+
\tfrac{1}{2}\,\overline{z}^2\zeta\zeta\overline{\zeta}
+
\tfrac{1}{2}\,z^2\overline{\zeta}\zeta\overline{\zeta}
+
z\overline{z}\zeta\overline{\zeta}\zeta\overline{\zeta}
\\
&
\ \ \ \ \
+
F_{3,0,0,2,0}^p\,z^3\overline{\zeta}^2
+
\overline{F_{3,0,0,2,0}^p}\,\overline{z}^3\zeta^2
+
\zeta\overline{\zeta}\,
\Big(
3\,F_{3,0,0,2,0}^p\,
z^2\overline{z}\overline{\zeta}
+
3\,\overline{F_{3,0,0,2,0}^p}\,
z\zeta\overline{z}^2
\Big)
\\
&
\ \ \ \ \
+
0
+
2\,\Re\,\Big\{
z^5\overline{\zeta}\,
F_{5,0,0,1,0}^p
+
\ \ \ \ \ \ \ \ \ \
0
\ \ \ \ \ \ \ \ \ \ \,
+
z^4\overline{\zeta}^2\,
F_{4,0,0,2,0}^p
\\
&
\ \ \ \ \ \ \ \ \ \ \ \ \ \ \ \ \ \ \ \ \ \ \ \ \ \ \ \ \ \ \ \ \ \ \
\ \ \ \ \ \ \ \ \ \ \ 
+
z^3\overline{z}^2\overline{\zeta}\,
F_{3,0,2,1,0}^p
+
z^3\overline{z}\overline{\zeta}^2\,
F_{3,0,1,2,0}^p
\\
&
\ \ \ \ \ \ \ \ \ \ \ \ \ \ \ \ \ \ \ \ \ \ \ \ \ \ \ \ \ \ \ \ \ \ \
\ \ \ \ \ \ \ \ \ \ \ \ \ \ \ \ \ \ \ \ \ \ \ \ \ \ \ 
\ \ \ \ \ \ \ \ \ \ \ \,
+
z^3\overline{\zeta}^3\,
F_{3,0,0,3,0}^p
\Big\}
+
{\rm O}_{z,\zeta,\overline{z},\overline{\zeta},v}(7),
\endaligned
\]
and which stabilizes the flat $1$-jet at the origin, 
is of the form:
\[
\aligned
z'
&
\,:=\,
\lambda\,z
+
i\,\lambda\,r\,
zw
+
2\,\frac{\lambda^2}{\overline{\lambda}}\,\overline{\chi}\,
z^3w
+
\psi\,
zw^2,
\\
\zeta'
&
\,:=\,
\frac{\lambda}{\overline{\lambda}}\,
\zeta
-
i\,
\frac{\lambda}{\overline{\lambda}}\,r\,
z^2
-
4\,\frac{\lambda^2}{\overline{\lambda}^2}\,\overline{\chi}\,
z^4
+
\Big(
-\tfrac{8}{3}\,\frac{\psi}{\overline{\lambda}}
+
\tfrac{4}{3}\,\frac{\lambda}{\overline{\lambda}^2}\,
\overline{\psi}
-
\tfrac{1}{3}\,\frac{\lambda}{\overline{\lambda}}\,r^2
\Big)\,
z^2w
+
\chi\,w^2,
\\
w'
&
\,:=\,
\lambda\overline{\lambda}\,w
+
i\,\lambda\overline{\lambda}\,r\,
w^2
+
\lambda^2\overline{\chi}\,
z^2w^2
+
\Big(
-\tfrac{1}{3}\,\lambda\overline{\lambda}\,r^2
+
\tfrac{1}{3}\,\overline{\lambda}\psi
+
\tfrac{1}{3}\,\lambda\overline{\psi}
\Big)\,
w^3.
\endaligned
\]
where $\lambda \in \C^\ast$, $r \in \R$,
$\psi \in \C$, $\chi \in \C$ are arbitrary
parameters.\qed
\end{Lemma}

Furthermore, with this map, 
if one stabilizes the normalization $F_{4,0,1,1,0} = 0 = 
F_{4,0,1,1,0}^{p\,\prime}$, 
the other coefficients
transform as:
\[
\footnotesize
\aligned
F_{5,0,0,1,0}^{p\,\prime}
&
\,=\,
\tfrac{1}{\lambda^3}\,
F_{5,0,0,1,0}^p
\ \ \ \ \ \ \ \ \ \ \ \ \ \ \ \ \ \ \ \ \ \ \ \ \ \ \ \ \ \ \ \ \ \ 
\ \   
0
\,=\,
0,
\ \ \ \ \ \ \ \ \ 
F_{4,0,0,2,0}^{p\,\prime}
\,=\,
\tfrac{1}{\lambda\overline{\lambda}}\,
F_{4,0,0,2,0}^p,
\\
F_{3,0,2,1,0}^{p\,\prime}
&
\,=\,
\tfrac{1}{\lambda\overline{\lambda}^2}
F_{3,0,2,1,0}^p
-
2i\,\tfrac{1}{\lambda\overline{\lambda}^2}
F_{3,0,0,2,0}^p,
\ \ \ \ \ \ \ \ \ \ \ \ \ \ \ \ \ \ \ \ \ \ \ \ \ 
F_{3,0,1,2,0}^{p\,\prime}
\,=\,
\tfrac{1}{\overline{\lambda}^2}\,
F_{3,0,1,2,0}^p,
\\
&
\ \ \ \ \ \ \ \ \ \ \ \ \ \ \ \ \ \ \ \ \ \ \ \ \ \ \ \ \ \ \ \ \ \
\ \ \ \ \ \ \ \ \ \ \ \ \ \ \ \ \ \ \ \ \ \ \ \ \ \ \ \ \ \ \ \ \ \
\ \ \ \ \ \ \ \ \ \ \ \ \ \ \ \ 
F_{3,0,0,3,0}^{p\,\prime}
\,=\,
\tfrac{\lambda}{\overline{\lambda}^2}\,
F_{3,0,0,3,0}^p,
\endaligned
\]
while if one stabilizes the normalization $F_{3,0,2,1,0} = 0 = 
F_{3,0,2,1,0}^{p\,\prime}$, the other coefficients
transform as:
\[
\aligned
F_{5,0,0,1,0}^{p\,\prime}
\,=\,
\tfrac{1}{\lambda^3}\,
F_{5,0,0,1,0}^p
\ \ \ \ \ \ \ \ \ 
F_{4,0,1,1,0}^{p\,\prime}
\,=\,
\tfrac{1}{\lambda^2\overline{\lambda}^2}\,
F_{4,0,1,1,0}^p
-
2\,\overline{\tau},
\ \ \ \ \ \ \ \ \
F_{4,0,0,2,0}^{p\,\prime}
\,=\,
\tfrac{1}{\lambda\overline{\lambda}}\,
F_{4,0,0,2,0}^p,
\\
0
\,=\,
2i\,\lambda\overline{\lambda}\,r
F_{3,0,0,2,0}^p,
\ \ \ \ \ \ \ \ \ \ \ \ \ \ \ \, 
F_{3,0,1,2,0}^{p\,\prime}
\,=\,
\tfrac{1}{\overline{\lambda}^2}\,
F_{3,0,1,2,0}^p,
\\
F_{3,0,0,3,0}^{p\,\prime}
\,=\,
\tfrac{\lambda}{\overline{\lambda}^2}\,
F_{3,0,0,3,0}^p.
\endaligned
\]
This second choice happens to be less natural than the first one,
because it forces to discuss the dichotomy branching:
\[
\xymatrix{
&&
F_{3,0,0,2,0}^p\,=\,0,
\\
F_{3,0,0,2,0}^p
\ar[urr]
\ar[drr]
&&
\\
&&
F_{3,0,0,2,0}^p\,\neq\,0,
}
\]
and when $F_{3,0,0,2,0}^p\,\neq\,0$,
it leads to normalize the parameter $r$, which
belongs to the isotropy of the GM-model, 
and such a normalization is too early to be done.

Therefore, we choose the normalization 
$F_{4,0,1,1,0}^{p\,\prime} = 0$.

\smallskip

By applying the technique of
Chen-Foo-Merker-Ta~{\cite[Sections~9, 10]{Chen-Foo-Merker-Ta-2019}},
one can realize, after rather hard computations,
that there corresponds to the Taylor coefficient
$F_{5,0,0,1,0}$, the relative invariant
$\Jaux_0$ of Pocchiola, presented in~{\cite{Pocchiola-2013, 
Merker-Pocchiola-2018, Foo-Merker-2019}}:
\[
\aligned
\overline{\Jaux}_0
&
\,:=\,
\frac{1}{6}\,
\frac{\overline{\mathcal{L}}_1\big(
\overline{\mathcal{L}}_1\big(
\overline{\mathcal{L}}_1\big(
\overline{\mathcal{L}}_1(\kaux)\big)\big)\big)}{
\overline{\mathcal{L}}_1(\kaux)}
-
\frac{5}{6}\,
\frac{\overline{\mathcal{L}}_1\big(
\overline{\mathcal{L}}_1\big(
\overline{\mathcal{L}}_1(\kaux)\big)\big)\,\,
\overline{\mathcal{L}}_1\big(
\overline{\mathcal{L}}_1(\kaux)\big)
}{
\overline{\mathcal{L}}_1(\kaux)^2}
-
\frac{1}{6}\,
\frac{\overline{\mathcal{L}}_1\big(
\overline{\mathcal{L}}_1\big(
\overline{\mathcal{L}}_1(\kaux)\big)\big)
}{
\overline{\mathcal{L}}_1(\kaux)}\,
\overline{\Paux}
\,+
\\
&
\ \ \ \ \
+
\frac{20}{27}\,
\frac{\overline{\mathcal{L}}_1\big(\overline{\mathcal{L}}_1
(\kaux)\big)^3}{
\overline{\mathcal{L}}_1(\kaux)^3}
+
\frac{5}{18}\,
\frac{\overline{\mathcal{L}}_1\big(
\overline{\mathcal{L}}_1(\kaux)\big)^2}{
\overline{\mathcal{L}}_1(\kaux)^2}\,
\overline{\Paux}
+
\frac{1}{6}\,
\frac{\overline{\mathcal{L}}_1\big(
\overline{\mathcal{L}}_1(\kaux)\big)\,\,
\overline{\mathcal{L}}_1\big(\overline{\Paux}\big)}{
\overline{\mathcal{L}}_1(\kaux)}
-
\frac{1}{9}\,
\frac{\overline{\mathcal{L}}_1\big(
\overline{\mathcal{L}}_1(\kaux)\big)}{
\overline{\mathcal{L}}_1(\kaux)}\,\,
\overline{\Paux}\,\overline{\Paux}
\,-
\\
&
\ \ \ \ \
-
\frac{1}{6}\,
\overline{\mathcal{L}}_1\big(
\overline{\mathcal{L}}_1\big(
\overline{\Paux}\big)\big)
+
\frac{1}{3}\,
\overline{\mathcal{L}}_1\big(\overline{\Paux}\big)\,
\overline{\Paux}
-
\frac{2}{27}\,
\overline{\Paux}\,
\overline{\Paux}\,
\overline{\Paux}.
\endaligned
\]
Much more simply, 
by plugging this normalized $F$ into 
this formula, we obtain its value
{\em only at one point}, namely at the origin:
\[
\overline{\Jaux}_0
\,=\,
20\,
\overline{F_{5,0,0,1,0}}.
\]

\Section{\bf Order $2$ Chains in $\mathfrak{C}_{2,1}$ Hypersurfaces 
$M^5 \subset \C^3$}
\label{order-2-chains-C-2-1-hypersurfaces-M5-C3}
\HEAD{{\ref{order-2-chains-C-2-1-hypersurfaces-M5-C3}}.~{\sf Order 
$2$ Chains in $\mathfrak{C}_{2,1}$ Hypersurfaces $M^5 \subset \C^3$}
}{
Wei-Guo {\sc Foo}, Joël {\sc Merker}, The-Anh {\sc Ta}}

In Lemma~{\ref{Lm-isotropy-order-6}}, the presence of the 
free parameter
$\chi \in \C$ in the last term $\chi\, w^2$, of order $4$, of
$\zeta' = \frac{\lambda}{\overline{\lambda}}\,
\zeta + \cdots + \chi\, w^2$, shows that
the flat second jet $j_0^2 = (0, 0, 0, 0,\,\, 0, 0, 0, 0)$
is {\em not} invariant by transformations which 
stabilize the normalizations achieved up to now at order $6$.

To define chains as in
Definition~8.4 of~{\cite{Merker-2020}}, 
we need then to explore a bit further the normalizations.

As we already know thanks to 
Proposition~{\ref{Prp-prenormalization}}, it is possible,
by some punctual normalization, to also make, at order $7$:
\reqnomode\usetagform{EngelLie}
\begin{align}
0
\,=\,
F_{a,b,0,0,e}^p
\tag{(a+b+2e=7),}
\\
0
\,=\,
F_{a,b,1,0,e}^p
\tag{(a+b+2e=6),}
\\
0
\,=\,
F_{a,b,2,0,e}^p
\tag{(a+b+2e=5).}
\end{align}
Once these normalizations are done, the condition that they
are preserved forces $\chi = 0$ (exercise).

We therefore come to maps which express the `ambiguity' of punctual
normalizations being of the form:
\[
\aligned
z'
&
\,:=\,
\lambda\,z
+
i\,\lambda\,r\,
zw
+
\psi\,
zw^2,
\\
\zeta'
&
\,:=\,
\frac{\lambda}{\overline{\lambda}}\,
\zeta
-
i\,
\frac{\lambda}{\overline{\lambda}}\,r\,
z^2
+
\Big(
-\tfrac{8}{3}\,\frac{\psi}{\overline{\lambda}}
+
\tfrac{4}{3}\,\frac{\lambda}{\overline{\lambda}^2}\,
\overline{\psi}
-
\tfrac{1}{3}\,\frac{\lambda}{\overline{\lambda}}\,r^2
\Big)\,
z^2w,
\\
w'
&
\,:=\,
\lambda\overline{\lambda}\,w
+
i\,\lambda\overline{\lambda}\,r\,
w^2
+
\Big(
-\tfrac{1}{3}\,\lambda\overline{\lambda}\,r^2
+
\tfrac{1}{3}\,\overline{\lambda}\psi
+
\tfrac{1}{3}\,\lambda\overline{\psi}
\Big)\,
w^3.
\endaligned
\]
Then such maps have the property that 
they send curves
$\R_v^1 \longrightarrow \R_{x,y,s,t}^4$ of the form:
\[
x
\,=\,
{\rm O}_v(2),
\ \ \ \ \ \ \ \ \ \ \ \
y
\,=\,
{\rm O}_v(2),
\ \ \ \ \ \ \ \ \ \ \ \
s
\,=\,
{\rm O}_v(2),
\ \ \ \ \ \ \ \ \ \ \ \
t
\,=\,
{\rm O}_v(2),
\]
to curves of the similar form:
\[
x'
\,=\,
{\rm O}_{v'}(2),
\ \ \ \ \ \ \ \ \ \ \ \
y
\,=\,
{\rm O}_{v'}(2),
\ \ \ \ \ \ \ \ \ \ \ \
s
\,=\,
{\rm O}_{v'}(2),
\ \ \ \ \ \ \ \ \ \ \ \
t
\,=\,
{\rm O}_{v'}(2),
\]
hence they {\em stabilize the flat $2$-jet $j_0^2 = 
(0, 0, 0, 0,\, 0, 0, 0, 0)$}.

In conclusion, we have reached
a point at which we can state an analog of 
Definition~8.4 in~{\cite{Merker-2020}}.

\begin{Definition}
\label{Def-2-jet-chain-pullback-flat-M5-C3}
Given a hypersurface $M^5 \subset \C^3$
in the class $\mathfrak{C}_{2,1}$, a point $p \in M$,
a $1$-jet $j_p^1 \in \Sigma_p^1$ at $p$, 
given 
the translation map $\tau_p \colon (M,p) 
\longrightarrow (M^p,0)$, 
and using {\em any} normalizing map $\Phi_p \colon M^p 
\longrightarrow N^p$ which sends $(M^p, 0)$ to
a hypersurface
$(N^p, 0)$ of equation:
\[
\aligned
{}
&
z\overline{z}
+
\tfrac{1}{2}\,\overline{z}^2\zeta
+
\tfrac{1}{2}\,z^2\overline{\zeta}
+
z\overline{z}\zeta\overline{\zeta}
+
\tfrac{1}{2}\,\overline{z}^2\zeta\zeta\overline{\zeta}
+
\tfrac{1}{2}\,z^2\overline{\zeta}\zeta\overline{\zeta}
\\
&
\ \ \ \ \
+
2\,\Re\,
\Big\{
0
+
0
+
F_{3,0,0,2,0}^p\,z^3\overline{\zeta}^2
+
\zeta\overline{\zeta}\,
\big(
3\,z^2\overline{z}\overline{\zeta}\,
F_{3,0,0,2,0}^p
\big)
\Big\}
\\
&
\ \ \ \ \
+
0
+
2\,\Re\,\Big\{
z^5\overline{\zeta}\,
F_{5,0,0,1,0}^p
+
\ \ \ \ \ \ \ \ \ \
0
\ \ \ \ \ \ \ \ \ \ \,
+
z^4\overline{\zeta}^2\,
F_{4,0,0,2,0}^p
\\
&
\ \ \ \ \ \ \ \ \ \ \ \ \ \ \ \ \ \ \ \ \ \ \ \ \ \ \ \ \ \ \ \ \ \ \
\ \ \ \ \ \ \ \ \ \ \ 
+
z^3\overline{z}^2\overline{\zeta}\,
F_{3,0,2,1,0}^p
+
z^3\overline{z}\overline{\zeta}^2\,
F_{3,0,1,2,0}^p
\\
&
\ \ \ \ \ \ \ \ \ \ \ \ \ \ \ \ \ \ \ \ \ \ \ \ \ \ \ \ \ \ \ \ \ \ \
\ \ \ \ \ \ \ \ \ \ \ \ \ \ \ \ \ \ \ \ \ \ \ \ \ \ \ 
\ \ \ \ \ \ \ \ \ \ \ \,
+
z^3\overline{\zeta}^3\,
F_{3,0,0,3,0}^p
\Big\}
+
{\rm O}_{z,\zeta,\overline{z},\overline{\zeta},v}(7),
\endaligned
\]
with in addition:
\reqnomode\usetagform{EngelLie}
\begin{align}
0
\,=\,
F_{a,b,0,0,e}^p
\tag{(a+b+2e=7),}
\\
0
\,=\,
F_{a,b,1,0,e}^p
\tag{(a+b+2e=6),}
\\
0
\,=\,
F_{a,b,2,0,e}^p
\tag{(a+b+2e=5),}
\end{align}
and which {\em also sends} $j_p^1$ to the flat $1$-jet
$j_0^1 = (0, 0, 0, 0)$ at $0 \in N^p$, 
assign the $2$-jet $j_p^2$
of the {\sl chain} at $p \in M$ associated with $j_p^1$ 
to be
the inverse image of the {\em flat} $2$-jet at $0 \in N^p$:
\[
j_p^2
\,:=\,
{\big(\Phi_p\circ\tau_p\big)^{(2)}}^{-1}
\big(0,0,0,0,\,\,0,0,0,0\big).
\]
\end{Definition}

Thanks to the preceding reasonings, the result $j_p^2$ is
independent of the normalizing map 
$\Phi_p \circ \tau_p$ 
satisfying ${(\Phi_p \circ \tau_p)}^{(1)} (j_p^1) = (0,0,0,0)$,
the flat $1$-jet at $0 \in N^p$. 

Furthermore, there are $\mathcal{C}^\omega$ functions
$A$, $B$, $C$, $D$, $E$, $F$, which can
be made explicit in terms of $\big\{ F_{a,b,c,d,e} \big\}_{
1 \leqslant a+b+c+d+2e \leqslant 6}$,
such that equations of chains are, with time parameter $v$:
\[
\aligned
\dot{s}
&
\,=\,
A\big(v,x,y,s,t,\dot{x},\dot{y}\big),
\\
\dot{t}
&
\,=\,
B\big(v,x,y,s,t,\dot{x},\dot{y}\big),
\endaligned
\ \ \ \ \ \ \ \ \ \ \ \ \ \ \ \ \ \ \ \ \ \ \ \ \ \
\aligned
\ddot{x}
&
\,=\,
C\big(v,x,y,s,t,\dot{x},\dot{y}\big),
\\
\ddot{y}
&
\,=\,
D\big(v,x,y,s,t,\dot{x},\dot{y}\big),
\\
\ddot{s}
&
\,=\,
E\big(v,x,y,s,t,\dot{x},\dot{y}\big),
\\
\ddot{t}
&
\,=\,
F\big(v,x,y,s,t,\dot{x},\dot{y}\big).
\endaligned
\]
Integrability follows from the fact that $\Sigma_0^2$ {\em is}
a surface.

After that order $2$ chains are known, it is natural
to restart once more the whole process of prenormalization
and of partial normalization which begun in 
Section~{\ref{chain-straightening-harmonic-killing}},
by assuming that the CR-transversal
curve $p \in \gamma \subset M$
(not anymore chosen at random)
{\em is} an order $2$ chain.
In fact, to have a second order chain at a point
$p \in M$, it suffices
to prescribe two real constants, the initial values
$\dot{x}(0)$, $\dot{y}(0)$.

\smallskip

Then, coming back to
Proposition~{\ref{Prp-full-normalization-order-5-v-axis}},
but viewed at the origin up to order $6$
in {\em all} variables $(z, \zeta, \overline{z}, \overline{\zeta},
v)$,  
we remember that we have constructed a normalizing map $\Phi_p$,
composed with a translation map $\tau_p$:
\[
\varphi
\colon\ \ \ \ \
(M,p)
\,\,\,\xrightarrow[{\rule[0pt]{50pt}{0pt}}]{\tau_p}\,\,\,
(M^p,0)
\,\,\,\xrightarrow[{\rule[0pt]{50pt}{0pt}}]{\Phi_p}\,\,\,
(N^p,0),
\]
which brings $(M, p)$ to $(N^p, 0)$ at the origin of equation:
\[
\aligned
u
&
\,=\,
z\overline{z}
+
\tfrac{1}{2}\,\overline{z}^2\zeta
+
\tfrac{1}{2}\,z^2\overline{\zeta}
+
z\overline{z}\zeta\overline{\zeta}
+
\tfrac{1}{2}\,\overline{z}^2\zeta\zeta\overline{\zeta}
+
\tfrac{1}{2}\,z^2\overline{\zeta}\zeta\overline{\zeta}
+
z\overline{z}\zeta\overline{\zeta}\zeta\overline{\zeta}
\\
&
\ \ \ \ \ \ \ \ \ \ \ \ 
+
z^3\overline{z}^3\,F_{3,0,3,0,0}^p
+
z^3\overline{z}^3\,
{\rm O}_{z,\overline{z},v}(1)
\\
&
\ \ \ \ \ \ \ \ \ \ \ \ 
+
2\,\Re\,
\Big\{
0
+
0
+
z^3\overline{\zeta}^2\,
F_{3,0,0,2,0}^p
+
\zeta\overline{\zeta}\,
\big(
3\,z^2\overline{z}\overline{\zeta}\,
F_{3,0,0,2,0}^p
\big)
\big\}
\\
&
\ \ \ \ \ \ \ \ \ \ \ \ 
+
2\,\Re\,
\Big\{
z^5\overline{\zeta}\,
F_{5,0,0,1,0}^p
+
z^4\overline{z}\overline{\zeta}\,
F_{4,0,1,1,0}^p
+
z^4\overline{\zeta}^2\,
F_{4,0,0,2,0}^p
\\
&
\ \ \ \ \ \ \ \ \ \ \ \ \ \ \ \ \ \ \ \ \ \ \ \ \ \ \ \ \ \ \ \ \ \ \ 
\ \ \ \ \ \ \ \ \ \ \ \,
+
z^3\overline{z}^2\overline{\zeta}\,
F_{3,0,2,1,0}^p
+
z^3\overline{z}\overline{\zeta}^2\,
F_{3,0,1,2,0}^p
\\
&
\ \ \ \ \ \ \ \ \ \ \ \ \ \ \ \ \ \ \ \ \ \ \ \ \ \ \ \ \ \ \ \ \ \ \
\ \ \ \ \ \ \ \ \ \ \ \ \ \ \ \ \ \ \ \ \ \ \ \ \ \ \ \ \ \ \ \ \ \ \
\ \ \ \ \
+
z^3\overline{\zeta}^3\,
F_{3,0,0,3,0}^p
\Big\}
+
{\rm O}_{z,\zeta,\overline{z},\overline{\zeta},v}(7),
\endaligned
\]
{\em without changing} the CR-transversal curve 
$0 \in \gamma \subset M$ being the $v$-axis,
hence having {\em flat} $1$-jet at the origin.

\begin{Assertion}
\label{Assertion-F-40110-zero}
Then $F_{4,0,1,1,0}^p = 0$ holds automatically, without
having the needs to perform any further biholomorphism.
\end{Assertion}

\proof
Indeed, from the proof of Lemma~{\ref{Lm-kill-F-30300-F40110}}
we already know that with the choice:
\[
\tau 
\,:=\,
\frac{1}{2}\,
\overline{F_{4,0,1,1,0}^p},
\]
one can continue to normalize
and make $F_{4,0,1,1,0}^{p\,\prime} = 0$ by means of the map:
\[
\aligned
z'
&
\,:=\,
z
+
F_{4,0,1,1,0}^p\,
z^2w
-
\tfrac{1}{2}\,
F_{4,0,1,1,0}^p\,
\zeta w^2
+
\tfrac{1}{2}\,
\overline{F_{4,0,1,1,0}^p}\,
w^2,
\\
\zeta'
&
\,:=\,
\zeta
-
2\,\overline{F_{4,0,1,1,0}^p}\,
zw
+
2\,F_{4,0,1,1,0}^p\,
\zeta w^2,
\\
w'
&
\,:=\,
w
+
F_{4,0,1,1,0}^p\,
zw^2,
\endaligned
\]
which we may call $\Psi \colon (N^p, 0) \longrightarrow
(N_\prime^p, 0)$.
We then reason as in~{\cite[9.5]{Merker-2020}}

If $F_{4,0,1,1,0}^p \neq 0$ would be nonzero, due to the
presence in $z'$ of the quadratic term 
$\frac{1}{2} \overline{F_{4,0,1,1,0}^p}\, w^2$, 
this map $\Psi$ would {\sl not} stabilize the flat order $2$ jet
$j_0^2 = (0, 0, 0, 0, \, 0, 0, 0, 0)$, and so,
this would contradict
Definition~{\ref{Def-2-jet-chain-pullback-flat-M5-C3}} applied
to $(M,p) := (N^p, 0)$, to $\varphi := \Psi$,
and to $(N^p, 0) := (N_\prime^p, 0)$.
\endproof

\Section{\bf Moser-like Normal Form for $\mathfrak{C}_{2,1}$ 
Hypersurfaces $M^5 \subset \C^3$}
\label{Moser-normal-C-2-1-hypersurfaces-M5-C3}
\HEAD{{\ref{Moser-normal-C-2-1-hypersurfaces-M5-C3}}.~{\sf Moser-like 
Normal Form for
$\mathfrak{C}_{2,1}$ Hypersurfaces $M^5 \subset \C^3$}
}{
Wei-Guo {\sc Foo}, Joël {\sc Merker}, The-Anh {\sc Ta}}

Lastly, coming again back to
Proposition~{\ref{Prp-normalization-order-5-v-axis}},
all what precedes showed that,
{\em without changing} any starting order $2$ chain
$0 \in \gamma \subset M$ to be straighgtened to be
the $v$-axis, we have constructed a normalizing map
$(M, 0) \longrightarrow (N, 0)$ 
so that, in the equation of $N$, 
we may (at last!) let appear all the terms of order $6$ in 
$(z, \zeta, \overline{z}, \overline{\zeta})$:
\[
\aligned
u
&
\,=\,
z\overline{z}
+
\tfrac{1}{2}\,\overline{z}^2\zeta
+
\tfrac{1}{2}\,z^2\overline{\zeta}
+
z\overline{z}\zeta\overline{\zeta}
+
\tfrac{1}{2}\,\overline{z}^2\zeta\zeta\overline{\zeta}
+
\tfrac{1}{2}\,z^2\overline{\zeta}\zeta\overline{\zeta}
+
z\overline{z}\zeta\overline{\zeta}\zeta\overline{\zeta}
\\
&
\ \ \ \ \ \ \ \ \ \ \ 
+
z^3\overline{z}^3\,
F_{3,0,3,0}(v)
+
z^3\overline{z}^3\,
{\rm O}_{z,\overline{z}}(1)
\\
&
\ \ \ \ \ \ \ \ \ \ \ \ 
+
2\,\Re\,
\Big\{
0
+
0
+
z^3\overline{\zeta}^2\,
F_{3,0,0,2}(v)
+
\zeta\overline{\zeta}\,
\big(
3\,z^2\overline{z}\overline{\zeta}\,
F_{3,0,0,2}(v)
\big)
\big\}
\\
&
\ \ \ \ \ \ \ \ \ \ \ \ 
+
2\,\Re\,
\Big\{
z^5\overline{\zeta}\,
F_{5,0,0,1}(v)
+
z^4\overline{z}\overline{\zeta}\,
\zero{F_{4,0,1,1}(v)}
+
z^4\overline{\zeta}^2\,
F_{4,0,0,2}(v)
\\
&
\ \ \ \ \ \ \ \ \ \ \ \ \ \ \ \ \ \ \ \ \ \ \ \ \ \ \ \ \ \ \ \ \ \ \ 
\ \ \ \ \ \ \ \ \ \ \ \ \ \ \,
+
z^3\overline{z}^2\overline{\zeta}\,
F_{3,0,2,1}(v)
+
z^3\overline{z}\overline{\zeta}^2\,
F_{3,0,1,2}(v)
\\
&
\ \ \ \ \ \ \ \ \ \ \ \ \ \ \ \ \ \ \ \ \ \ \ \ \ \ \ \ \ \ \ \ \ \ \
\ \ \ \ \ \ \ \ \ \ \ \ \ \ \ \ \ \ \ \ \ \ \ \ \ \ \ \ \ \ \ \ \ \ \
\ \ \ \ \ \ \ \ \ \ \,
+
z^3\overline{\zeta}^3\,
F_{3,0,0,3}(v)
\Big\}
\\
&
\ \ \ \ \ \ \ \ \ \ \ 
+
\overline{z}^3\zeta\,
{\rm O}_{z,\zeta,\overline{z}}(3)
+
z^3\overline{\zeta}\,
{\rm O}_{z,\overline{z},\overline{\zeta}}(3)
+
\zeta\overline{\zeta}\,
{\rm O}_{z,\zeta,\overline{z},\overline{\zeta}}(5).
\endaligned
\]

\begin{Assertion}
The function $F_{4,0,1,1}(v)
\equiv 0$
vanishes for free.
\end{Assertion}

\proof
What was done an instant ago
by Assertion~{\ref{Assertion-F-40110-zero}} 
at the origin $(z, \zeta, w) = (0, 0, 0)$
applies in fact at {\em every} point $(0, 0, iv)$ along the
$v$-axis, thanks to the fact that the above graphed
equation is the same {\em all along} the $v$-axis.
\endproof

\begin{Proposition}
There exists a biholomorphism of the form:
\[
z'
\,:=\,
z\,\varphi(-iw),
\ \ \ \ \ \ \ \ \ \ \ \ \ \ \ \ \ \ \ \
\zeta'
\,:=\,
\zeta
+
\chi(-iw)\,z^2,
\ \ \ \ \ \ \ \ \ \ \ \ \ \ \ \ \ \ \ \
w'
\,:=\,
i\,\psi(-iw),
\]
with $\psi(v) \in \R$ for $v \in \R$, which normalizes 
{\em in addition} $F_{3,0,3,0}'(v') \equiv 0$.

\end{Proposition}

\proof
Left to the reader. Hint: imitate~{\cite[Lm.~12.4]{Merker-2020}}.
\endproof

In summary, we can state 

\begin{Theorem}
\label{Thm-complete-normal-form}
{\bf [Existence of normal form]}
For every $2$-nondegenerate
hypersurface $M^5 \in \mathfrak{C}_{2,1}$ whose
Levi form has constant rank $1$, 
at any point $p \in M$, 
given any order $2$ CR-transversal chain $p \in \gamma \subset M$,
there exist holomorphic coordinates $(z, \zeta, w) \in \C^3$
vanishing at $p$ in which $\gamma$ is the $v$-axis and in which
$M$ has normalized equation:
\reqnomode\usetagform{EngelLie}
\begin{footnotesize}
\begin{align}
u
&
\,=\,
\tfrac{z\overline{z}+\frac{1}{2}\overline{z}^2\zeta
+z^2\overline{\zeta}}{1-\zeta\overline{\zeta}}
\notag
\\
&
\ \ \ \ \
+
z^3\overline{z}^3\,
{\rm O}_{z,\overline{z}}(1)
+
2\,\Resmall\,
\Big\{
z^3\overline{\zeta}^2\,
F_{3,0,0,2}(v)
+
\zeta\overline{\zeta}\,
\big(
3\,z^2\overline{z}\overline{\zeta}\,
F_{3,0,0,2}(v)
\big)
\big\}
\notag
\\
&
\ \ \ \ \
+
2\,\Resmall\,
\Big\{
z^5\overline{\zeta}\,
F_{5,0,0,1}(v)
+
z^4\overline{\zeta}^2\,
F_{4,0,0,2}(v)
+
z^3\overline{z}^2\overline{\zeta}\,
F_{3,0,2,1}(v)
+
z^3\overline{z}\overline{\zeta}^2\,
F_{3,0,1,2}(v)
+
z^3\overline{\zeta}^3\,
F_{3,0,0,3}(v)
\Big\}
\notag
\\
&
\ \ \ \ \
+
\overline{z}^3\zeta\,
{\rm O}_{z,\zeta,\overline{z}}(3)
+
z^3\overline{\zeta}\,
{\rm O}_{z,\overline{z},\overline{\zeta}}(3)
+
\zeta\overline{\zeta}\,
{\rm O}_{z,\zeta,\overline{z},\overline{\zeta}}(5).
\tag{\qed}
\end{align}
\end{footnotesize}
\end{Theorem}

\Section{\bf Consequence of Prenormalization on Dependent Jets}
\label{consequence-prenormalization-dependent-jets}
\HEAD{{\ref{consequence-prenormalization-dependent-jets}}.~{\sf 
Consequence of Prenormalization on Dependent Jets}
}{
Wei-Guo {\sc Foo}, Joël {\sc Merker}, The-Anh {\sc Ta}}

After the prenormalization 
Proposition~{\ref{Prp-prenormalization}},
we know that we have:
\[
u
\,=\,
F
\,=\,
\maux+G
\,=\,
\maux
+
z^3\overline{z}^3\,{\rm O}_{z,\overline{z}}(0)
+
z^3\overline{\zeta}\,{\rm O}_{z,\zeta,\overline{z}}(0)
+
\overline{z}^3\zeta\,{\rm O}_{z,\overline{z},\overline{\zeta}}(0)
+
\zeta\overline{\zeta}\,
{\rm O}_{z,\zeta,\overline{z},\overline{\zeta}}(3).
\]
The next statement shows that the dependent-jets remainder
is {\em in addition} an ${\rm O}_{z,\overline{z}}(3)$.

\begin{Proposition}
\label{Prp-G-O-z-zbar-3}
In prenormalized coordinates, $G = {\rm O}_{z, \overline{z}}(3)$.
\end{Proposition}

This writing means here that $G$ is of order $3$ in
$(z,\overline{z})$, with coefficients being arbitrary functions of
$(z, \zeta, \overline{z}, \overline{\zeta}, v)$, namely that:
\[
G
\,=\,
z^3\,\big(\cdots\big)
+
z^2\overline{z}\,\big(\cdots\big)
+
z\overline{z}^2\,\big(\cdots\big)
+
\overline{z}^3\,\big(\cdots\big).
\]

\proof
Since the coordinates are prenormalized, we have at least:
\[
u
\,=\,
z\overline{z}
+
\tfrac{1}{2}\,
\overline{z}^2\zeta
+
\tfrac{1}{2}\,
z^2\overline{\zeta}
+
{\rm O}_{z,\zeta,\overline{z},\overline{\zeta}}(4)
\,\,=\,\,
\maux
+
G.
\]
Thus if we write:
\[
G
\,=\,
\sum_{\kappa\geqslant 2}\,
\sum_{a+b+c+d=\kappa}\,
G_{a,b,c,d}(v)\,
z^a\zeta^b\overline{z}^c\overline{\zeta}^d
\,\,=:\,\,
\sum_{\kappa\geqslant 2}\,
G^\kappa(v).
\]
we have $0 = G^2 = G^3$, which are certainly both
${\rm O}_{z,\overline{z}}(3)$.

The proof will consist in examining, order by order, the Levi
determinant for $F = \maux + G$:
\[
\left\vert\!
\def\arraystretch{1.5}
\begin{array}{cccc}
0 & F_z & F_\zeta & -\frac{1}{2}+\frac{1}{2i}F_v
\\
F_{\overline{z}} & F_{z\overline{z}} & F_{\zeta\overline{z}} &
\frac{1}{2i}F_{\overline{z}v}
\\
F_{\overline{\zeta}} & F_{z\overline{\zeta}} & 
F_{\zeta\overline{\zeta}} &
\frac{1}{2i}F_{\overline{\zeta}v}
\\
-\frac{1}{2}-\frac{1}{2i}F_v & -\frac{1}{2i}F_{zv} &
-\frac{1}{2i}F_{\zeta v} & \frac{1}{4}F_{vv}
\end{array}
\!\right\vert.
\]

Reasoning by induction, 
assume, for some $\kappa \geqslant 4$, that
$G^2, G^3, \dots, G^{\kappa-1}$ are all 
${\rm O}_{z,\overline{z}}(3)$. For 
all $2 \leqslant \ell \leqslant \kappa-1$, it then follows that:
\[
\aligned
{}
&
{}
&
\ \ \ \ \ \ \ \ \ \ \
G_z^\ell
&
\,=\,
{\rm O}_{z,\overline{z}}(2),
&
\ \ \ \ \ \ \ \ \ \ \
G_\zeta^\ell
&
\,=\,
{\rm O}_{z,\overline{z}}(3),
&
\ \ \ \ \ \ \ \ \ \ \
G_v^\ell
&
\,=\,
{\rm O}_{z,\overline{z}}(3),
\\
G_{\overline{z}}^\ell
&
\,=\,
{\rm O}_{z,\overline{z}}(2),
&
\ \ \ \ \ \ \ \ \ \ \
G_{z\overline{z}}^\ell
&
\,=\,
{\rm O}_{z,\overline{z}}(1),
&
\ \ \ \ \ \ \ \ \ \ \
G_{\zeta\overline{z}}^\ell
&
\,=\,
{\rm O}_{z,\overline{z}}(2),
&
\ \ \ \ \ \ \ \ \ \ \
G_{\overline{z}v}^\ell
&
\,=\,
{\rm O}_{z,\overline{z}}(2),
\\
G_{\overline{\zeta}}^\ell
&
\,=\,
{\rm O}_{z,\overline{z}}(3),
&
\ \ \ \ \ \ \ \ \ \ \
G_{z\overline{\zeta}}^\ell
&
\,=\,
{\rm O}_{z,\overline{z}}(2),
&
\ \ \ \ \ \ \ \ \ \ \
G_{\zeta\overline{\zeta}}^\ell
&
\,=\,
{\rm O}_{z,\overline{z}}(3),
&
\ \ \ \ \ \ \ \ \ \ \
G_{\overline{\zeta}v}^\ell
&
\,=\,
{\rm O}_{z,\overline{z}}(3),
\\
G_{v}^\ell
&
\,=\,
{\rm O}_{z,\overline{z}}(3),
&
\ \ \ \ \ \ \ \ \ \ \
G_{zv}^\ell
&
\,=\,
{\rm O}_{z,\overline{z}}(2),
&
\ \ \ \ \ \ \ \ \ \ \
G_{\zeta v}^\ell
&
\,=\,
{\rm O}_{z,\overline{z}}(3),
&
\ \ \ \ \ \ \ \ \ \ \
G_{vv}^\ell
&
\,=\,
{\rm O}_{z,\overline{z}}(3).
\endaligned
\]

To capture information about $G^\kappa$, we may truncate
modulo ${\rm O}_{z,\zeta,\overline{z},\overline{\zeta}}(\kappa+1)$:
\[
\aligned
\maux
&
\,\equiv\,
\maux^2
+
\maux^3
+\cdots+
\maux^{\kappa-2}
+
\maux^{\kappa-1}
+
\maux^\kappa,
\\
G
&
\,\equiv\,
G^2
+
G^3
+\cdots+
G^{\kappa-2}
+
G^{\kappa-1}
+
G^\kappa,
\endaligned
\]
where, for any formal:
\[
H
\,=\,
\sum_{a,b,c,d\geqslant 0}\,
z^a\zeta^b\overline{z}^c\overline{\zeta}^d\,
H_{a,b,c,d}(v),
\]
and any $\mu \geqslant 0$, we set:
\[
\aligned
H^\mu
&
\,:=\,
\sum_{a+b+c+d=\mu}\,
z^a\zeta^b\overline{z}^c\overline{\zeta}^d\,
H_{a,b,c,d}(v),
\\
\pi^\mu(H)
&
\,:=\,
\sum_{a+b+c+d\leqslant\mu}\,
z^a\zeta^b\overline{z}^c\overline{\zeta}^d\,
H_{a,b,c,d}(v).
\endaligned
\]
We will insert $F = \maux + G$ in the Levi determinant and
apply the projection $\pi^{\kappa-2} (\centersmallbullet)$
in order to capture $G_{\zeta \overline{\zeta}}^\kappa$.

\begin{Assertion}
Under the induction assumption, 
$G_{\zeta\overline{\zeta}}^\kappa = {\rm O}_{z, \overline{z}}(3)$.
\end{Assertion}

\proof
Some further preliminaries are necessary.
At first, for any formal function $L = L(z, \zeta, 
\overline{z}, \overline{\zeta}, v)$ which is an
${\rm O}_{z, \zeta, \overline{z}, \overline{\zeta}} (\lambda)$ 
for some $\lambda \geqslant 0$, it holds, {\em with
a shift}, that:
\leqnomode\usetagform{default}
\begin{align}
\label{truncate-product-L-H}
\pi^{\kappa-2}
\big(
L
\cdot
H
\big)
\,=\,
\pi^{\kappa-2}
\Big(
\pi^{\kappa-2}
\big(
L
\big)
\cdot
\pi^{\kappa-2-\lambda}
\big(
H
\big)
\Big).
\end{align}

Next, with ${}_\smallbullet$ and 
${}_{\smallbullet, \smallbullet}$ denoting partial derivatives
with respect to any of the variables $z$, $\zeta$, $\overline{z}$,
$\overline{\zeta}$, we have:
\[
\footnotesize
\aligned
\pi^{\kappa-2}
(\maux)
&
\,=\,
\maux^2
+\cdots+
\maux^{\kappa-2},
&
\ \ \ \ \ \ \ \ \ \ \ \ \ \ \ \ \ \ \ \
\pi^{\kappa-2}
(G)
&
\,=\,
G^2
+\cdots+
G^{\kappa-2},
\\
\pi^{\kappa-2}
\big(\maux_{\smallbullet}\big)
&
\,=\,
\maux_{\smallbullet}^2
+\cdots+
\maux_{\smallbullet}^{\kappa-2}
+
\maux_{\smallbullet}^{\kappa-1},
&
\ \ \ \ \ \ \ \ \ \ \ \ \ \ \ \ \ \ \ \
\pi^{\kappa-2}
\big(G_{\smallbullet}\big)
&
\,=\,
G_{\smallbullet}^2
+\cdots+
G_{\smallbullet}^{\kappa-2}
+
G_{\smallbullet}^{\kappa-1},
\\
\pi^{\kappa-2}
\big(\maux_{\smallbullet,\smallbullet}\big)
&
\,=\,
\maux_{\smallbullet,\smallbullet}^2
+\cdots+
\maux_{\smallbullet,\smallbullet}^{\kappa-2}
+
\maux_{\smallbullet,\smallbullet}^{\kappa-1}
+
\maux_{\smallbullet,\smallbullet}^{\kappa},
&
\ \ \ \ \ \ \ \ \ \ \ \ \ \ \ \ \ \ \ \
\pi^{\kappa-2}
\big(G_{\smallbullet,\smallbullet}\big)
&
\,=\,
G_{\smallbullet,\smallbullet}^2
+\cdots+
G_{\smallbullet,\smallbullet}^{\kappa-2}
+
G_{\smallbullet,\smallbullet}^{\kappa-1},
+
G_{\smallbullet,\smallbullet}^{\kappa}.
\endaligned
\]
Also, we will be using various values 
$\lambda = 0, 1, 2$ 
of the integer $\lambda
\geqslant 0$ above:
\[
\aligned
\maux_z
&
\,=\,
\tfrac{\overline{z}+z\overline{\zeta}}{1-\zeta\overline{\zeta}}
\,=\,
{\rm O}_{z,\overline{z}}(1),
&
\ \ \ \ \ \ \ \ \ \ \ \ \ \ \ \ \ \ \ \
\maux_{\overline{z}}
&
\,=\,
\tfrac{z+\overline{z}\zeta}{1-\zeta\overline{\zeta}}
\,=\,
{\rm O}_{z,\overline{z}}(1),
\\
\maux_\zeta
&
\,=\,
\tfrac{1}{2}\,
\tfrac{(\overline{z}+z\overline{\zeta})^2}{
(1-\zeta\overline{\zeta})^2}
\,=\,
{\rm O}_{z,\overline{z}}(2),
&
\ \ \ \ \ \ \ \ \ \ \ \ \ \ \ \ \ \ \ \
\maux_{\overline{\zeta}}
&
\,=\,
\tfrac{1}{2}\,
\tfrac{(z+\overline{z}\zeta)^2}{
(1-\zeta\overline{\zeta})^2}
\,=\,
{\rm O}_{z,\overline{z}}(2),
\\
\maux_{z\overline{z}}
&
\,=\,
\tfrac{1}{1-\zeta\overline{\zeta}}
\,=\,
{\rm O}_{z,\overline{z}}(0),
&
\ \ \ \ \ \ \ \ \ \ \ \ \ \ \ \ \ \ \ \
\maux_{\zeta\overline{z}}
&
\,=\,
\tfrac{\overline{z}+z\overline{\zeta}}{(1-\zeta\overline{\zeta})^2}
\,=\,
{\rm O}_{z,\overline{z}}(1),
\\
\maux_{z\overline{\zeta}}
&
\,=\,
\tfrac{z+\overline{z}\zeta}{(1-\zeta\overline{\zeta})^2}
\,=\,
{\rm O}_{z\overline{z}}(1),
&
\ \ \ \ \ \ \ \ \ \ \ \ \ \ \ \ \ \ \ \
\maux_{\zeta\overline{\zeta}}
&
\,=\,
\tfrac{(z+\overline{z}\zeta)^2}{(1-\zeta\overline{\zeta})^3}
\,=\,
{\rm O}_{z,\overline{z}}(2).
\endaligned
\]

Indeed, we start from:
\[
0
\,\equiv\,
\pi^{\kappa-2}
\left(
\left\vert\!
\def\arraystretch{1.5}
\begin{array}{cccc}
0 & 
\maux_z+\smallsum{4\leqslant j\leqslant\kappa-1}G_z^j &
\maux_\zeta+\smallsum{4\leqslant k\leqslant\kappa-1}G_\zeta^k &
-\frac{1}{2}-\frac{i}{2}
\smallsum{4\leqslant l\leqslant\kappa-2}G_v^l 
\\
\maux_{\overline{z}}
+\smallsum{4\leqslant i\leqslant\kappa-1}G_{\overline{z}}^i &
\maux_{z\overline{z}}
+\smallsum{4\leqslant j\leqslant\kappa}G_{z\overline{z}}^j &
\maux_{\zeta\overline{z}}
+\smallsum{4\leqslant k\leqslant\kappa}G_{\zeta\overline{z}}^k &
-\frac{i}{2}
\smallsum{4\leqslant l\leqslant\kappa-1}G_{\overline{z}v}^l
\\
\maux_{\overline{\zeta}}
+\smallsum{4\leqslant i\leqslant\kappa-1}G_{\overline{\zeta}}^i &
\maux_{z\overline{\zeta}}
+\smallsum{4\leqslant j\leqslant\kappa}G_{z\overline{\zeta}}^j &
\maux_{\zeta\overline{\zeta}}
+\smallsum{4\leqslant k\leqslant\kappa}G_{\zeta\overline{\zeta}}^k &
-\frac{i}{2}
\smallsum{4\leqslant l\leqslant\kappa-1}G_{\overline{\zeta}v}^l
\\
-\frac{1}{2}
+\frac{i}{2}\smallsum{4\leqslant i\leqslant\kappa-2}G_v^i &
\frac{i}{2}
\smallsum{4\leqslant j\leqslant\kappa-1}G_{zv}^j &
\frac{i}{2}
\smallsum{4\leqslant k\leqslant\kappa-1}G_{\zeta v}^k &
\frac{1}{4}
\smallsum{4\leqslant l\leqslant\kappa-2}G_{vv}^l
\end{array}
\!\right\vert
\right).
\]
Let us expand this determinant along its first row,
using~({\ref{truncate-product-L-H}}) in order to take account
of various useful {\em negative shifts} for the summations in
the entries of the obtained
$3 \times 3$ determinants:
\[
\aligned
0
&
\,\equiv\,
\pi^{\kappa-2}
\left(
-\,
\Big(
\maux_z
+
\smallsum{4\leqslant j\leqslant\kappa-1}G_z^j
\Big)\,
\left\vert\!
\def\arraystretch{1.5}
\begin{array}{ccc}
\maux_{\overline{z}}
+\smallsum{4\leqslant i\leqslant\kappa-2}G_{\overline{z}}^i &
\maux_{\zeta\overline{z}}
+\smallsum{4\leqslant k\leqslant\kappa-1}G_{\zeta\overline{z}}^k &
-\frac{i}{2}
\smallsum{4\leqslant l\leqslant\kappa-2}G_{\overline{z}v}^l
\\
\maux_{\overline{\zeta}}
+\smallsum{4\leqslant i\leqslant\kappa-2}G_{\overline{\zeta}}^i &
\maux_{\zeta\overline{\zeta}}
+\smallsum{4\leqslant k\leqslant\kappa-1}G_{\zeta\overline{\zeta}}^k &
-\frac{i}{2}
\smallsum{4\leqslant l\leqslant\kappa-2}G_{\overline{\zeta}v}^l
\\
-\frac{1}{2}
+\frac{i}{2}\smallsum{4\leqslant i\leqslant\kappa-3}G_v^i &
\frac{i}{2}
\smallsum{4\leqslant k\leqslant\kappa-2}G_{\zeta v}^k &
\frac{1}{4}
\smallsum{4\leqslant l\leqslant\kappa-3}G_{vv}^l
\end{array}
\!\right\vert
\right.
\\
&
\ \ \ \ \ \ \ \ \ \ \ \ \ \ \ \ \ \ 
+
\Big(
\maux_\zeta
+
\smallsum{4\leqslant k\leqslant\kappa-1}G_\zeta^k
\Big)
\left\vert\!
\def\arraystretch{1.5}
\begin{array}{ccc}
\\
\maux_{\overline{z}}
+\smallsum{4\leqslant i\leqslant\kappa-3}G_{\overline{z}}^i &
\maux_{z\overline{z}}
+\smallsum{4\leqslant j\leqslant\kappa-2}G_{z\overline{z}}^j &
-\frac{i}{2}
\smallsum{4\leqslant l\leqslant\kappa-3}G_{\overline{z}v}^l
\\
\maux_{\overline{\zeta}}
+\smallsum{4\leqslant i\leqslant\kappa-3}G_{\overline{\zeta}}^i &
\maux_{z\overline{\zeta}}
+\smallsum{4\leqslant j\leqslant\kappa-2}G_{z\overline{\zeta}}^j &
-\frac{i}{2}
\smallsum{4\leqslant l\leqslant\kappa-3}G_{\overline{\zeta}v}^l
\\
-\frac{1}{2}
+\frac{i}{2}\smallsum{4\leqslant i\leqslant\kappa-4}G_v^i &
\frac{i}{2}
\smallsum{4\leqslant j\leqslant\kappa-3}G_{zv}^j &
\frac{1}{4}
\smallsum{4\leqslant l\leqslant\kappa-4}G_{vv}^l
\end{array}
\!\right\vert
\\
&
\ \ \ \ \ \ \ \ \ \ \ \ \ \ \ \ \ \ \ \ \ \,
\left.
-\,
\Big(
-\tfrac{1}{2}
+
{\rm O}_{z,\overline{z}}(3)
\Big)
\left\vert\!
\def\arraystretch{1.5}
\begin{array}{ccc}
\maux_{\overline{z}}
+\smallsum{4\leqslant i\leqslant\kappa-1}G_{\overline{z}}^i &
\maux_{z\overline{z}}
+\smallsum{4\leqslant j\leqslant\kappa}G_{z\overline{z}}^j &
\maux_{\zeta\overline{z}}
+\smallsum{4\leqslant k\leqslant\kappa}G_{\zeta\overline{z}}^k 
\\
\maux_{\overline{\zeta}}
+\smallsum{4\leqslant i\leqslant\kappa-1}G_{\overline{\zeta}}^i &
\maux_{z\overline{\zeta}}
+\smallsum{4\leqslant j\leqslant\kappa}G_{z\overline{\zeta}}^j &
\maux_{\zeta\overline{\zeta}}
+\smallsum{4\leqslant k\leqslant\kappa}G_{\zeta\overline{\zeta}}^k 
\\
-\frac{1}{2}
+\frac{i}{2}\smallsum{4\leqslant i\leqslant\kappa-2}G_v^i &
\frac{i}{2}
\smallsum{4\leqslant j\leqslant\kappa-1}G_{zv}^j &
\frac{i}{2}
\smallsum{4\leqslant k\leqslant\kappa-1}G_{\zeta v}^k 
\end{array}
\!\right\vert
\,\right).
\endaligned
\]

Now, apply the induction assumption, and simultaneously also,
expand the last determinant along its first column:
\[
\aligned
0
&
\,\equiv\,
\pi^{\kappa-2}
\left(
-\,{\rm O}_{z,\overline{z}}(1)\,
\left\vert\!
\begin{array}{ccc}
{\rm O}_{z,\overline{z}}(1) & 
{\rm O}_{z,\overline{z}}(1) & 
{\rm O}_{z,\overline{z}}(2) 
\\
{\rm O}_{z,\overline{z}}(2) & 
{\rm O}_{z,\overline{z}}(2) & 
{\rm O}_{z,\overline{z}}(3) 
\\
{\rm O}_{z,\overline{z}}(0) & 
{\rm O}_{z,\overline{z}}(3) & 
{\rm O}_{z,\overline{z}}(3) 
\end{array}
\!\right\vert
+
{\rm O}_{z,\overline{z}}(2)\,
\left\vert\!
\begin{array}{ccc}
{\rm O}_{z,\overline{z}}(1) & 
{\rm O}_{z,\overline{z}}(0) & 
{\rm O}_{z,\overline{z}}(2) 
\\
{\rm O}_{z,\overline{z}}(2) & 
{\rm O}_{z,\overline{z}}(1) & 
{\rm O}_{z,\overline{z}}(3) 
\\
{\rm O}_{z,\overline{z}}(0) & 
{\rm O}_{z,\overline{z}}(2) & 
{\rm O}_{z,\overline{z}}(3) 
\end{array}
\!\right\vert
\right.
\\
&
\ \ \ \ \
+
\Big(
\tfrac{1}{2}
+
{\rm O}_{z,\overline{z}}(3)
\Big)
\left\{
\Big(
\maux_{\overline{z}}
+
\smallsum{4\leqslant i\leqslant\kappa-1}G_{\overline{z}}^i
\Big)
\left\vert\!
\def\arraystretch{1.5}
\begin{array}{cc}
\maux_{z\overline{\zeta}}
+\smallsum{4\leqslant j\leqslant\kappa-1}G_{z\overline{\zeta}}^j &
\maux_{\zeta\overline{\zeta}}
+\smallsum{4\leqslant k\leqslant\kappa-1}\,
G_{\zeta\overline{\zeta}}^k 
\\
\frac{i}{2}
\smallsum{4\leqslant j\leqslant\kappa-2}G_{zv}^j &
\frac{i}{2}
\smallsum{4\leqslant k\leqslant\kappa-2}G_{\zeta v}^k
\end{array}
\!\right\vert
\right.
\\
&
\ \ \ \ \ \ \ \ \ \ \ \ \ \ \ \ \ \ \ \ \ \ \ \ \ \ \ \ \ \ \ \   
-\,
\Big(
\maux_{\overline{\zeta}}
+
\smallsum{4\leqslant i\leqslant\kappa-1}G_{\overline{\zeta}}^i
\Big)
\left\vert\!
\def\arraystretch{1.5}
\begin{array}{cc}
\maux_{z\overline{z}}
+\smallsum{4\leqslant j\leqslant\kappa-2}G_{z\overline{z}}^j &
\maux_{\zeta\overline{z}}
+\smallsum{4\leqslant k\leqslant\kappa-2}G_{\zeta\overline{z}}^k
\\
\frac{i}{2}
\smallsum{4\leqslant j\leqslant\kappa-3}G_{zv}^j &
\frac{i}{2}
\smallsum{4\leqslant k\leqslant\kappa-3}G_{\zeta v}^k
\end{array}
\!\right\vert
\\
&
\ \ \ \ \ \ \ \ \ \ \ \ \ \ \ \ \ \ \ \ \ \ \ \ \ \ \ \ \ \ \ \ 
\left.
\left.
+
\Big(
-\tfrac{1}{2}
+
{\rm O}_{z,\overline{z}}(3)
\Big)
\left\vert\!
\def\arraystretch{1.5}
\begin{array}{cc}
\maux_{z\overline{z}}
+\smallsum{4\leqslant j\leqslant\kappa}G_{z\overline{z}}^j &
\maux_{\zeta\overline{z}}
+\smallsum{4\leqslant k\leqslant\kappa}G_{\zeta\overline{z}}^k
\\
\maux_{z\overline{\zeta}}
+\smallsum{4\leqslant j\leqslant\kappa}G_{z\overline{\zeta}} &
\maux_{\zeta\overline{\zeta}}
+\smallsum{4\leqslant k\leqslant\kappa}
G_{\zeta\overline{\zeta}}^k
\end{array}
\!\right\vert
\right\}
\right).
\endaligned
\]
Taking account of $0 \equiv \big\vert
\begin{smallmatrix} \maux_{z\overline{z}} & \maux_{\zeta\overline{z}}
\\ \maux_{z\overline{\zeta}} & \maux_{\zeta\overline{\zeta}} 
\end{smallmatrix} \big\vert$ in the 
last $2 \times 2$ determinant, we may continue to expand:
\[
\aligned
0
&
\,\equiv\,
{\rm O}_{z,\overline{z}}(3)
+
{\rm O}_{z,\overline{z}}(1)\,
\left\vert\!
\begin{array}{cc}
{\rm O}_{z,\overline{z}}(1) &
{\rm O}_{z,\overline{z}}(2)
\\
{\rm O}_{z,\overline{z}}(2) &
{\rm O}_{z,\overline{z}}(3)
\end{array}
\!\right\vert
-
{\rm O}_{z,\overline{z}}(2)\,
\left\vert\!
\begin{array}{cc}
{\rm O}_{z,\overline{z}}(0) &
{\rm O}_{z,\overline{z}}(1)
\\
{\rm O}_{z,\overline{z}}(2) &
{\rm O}_{z,\overline{z}}(3)
\end{array}
\!\right\vert
\\
&
\ \ \ \ \ \ 
+
\Big(
-\frac{1}{4}
+
{\rm O}_{z,\overline{z}}(3)
\Big)
\left\{
\aligned
{}&
\maux_{z\overline{z}}
\smallsum{4\leqslant k\leqslant\kappa}G_{\zeta\overline{\zeta}}^k
+
\maux_{\zeta\overline{\zeta}}
\smallsum{4\leqslant j\leqslant\kappa-2}G_{z\overline{z}}^j
+
\Big(
\smallsum{4\leqslant j\leqslant\kappa-2}
G_{z\overline{z}}^j
\Big)
\Big(
\smallsum{4\leqslant k\leqslant\kappa-2}
G_{\zeta\overline{\zeta}}^k
\Big)
\\
{}&
-\,\maux_{\zeta\overline{z}}
\smallsum{4\leqslant j\leqslant\kappa-1}G_{z\overline{\zeta}}^j
-
\maux_{z\overline{\zeta}}
\smallsum{4\leqslant k\leqslant\kappa-1}G_{\zeta\overline{z}}^k
-
\Big(
\smallsum{4\leqslant k\leqslant\kappa-2}G_{\zeta\overline{z}}^k
\Big)
\Big(
\smallsum{4\leqslant j\leqslant\kappa-2}G_{z\overline{\zeta}}^j
\Big)
\endaligned
\right\},
\endaligned
\]
that is:
\[
\aligned
{\rm O}_{z,\overline{z}}(3)
\,\equiv\,
\maux_{z\overline{z}}\,
\Big(
\smallsum{4\leqslant k\leqslant\kappa-1}
G_{\zeta\overline{\zeta}}^k
+
G_{\zeta\overline{\zeta}}^\kappa
\Big)
&
+
{\rm O}_{z,\overline{z}}(2)\,
{\rm O}_{z,\overline{z}}(1)
+
{\rm O}_{z,\overline{z}}(1)\,
{\rm O}_{z,\overline{z}}(3)
\\
&
-\,
{\rm O}_{z,\overline{z}}(1)\,
{\rm O}_{z,\overline{z}}(2)
-
{\rm O}_{z,\overline{z}}(1)\,
{\rm O}_{z,\overline{z}}(2)
-
{\rm O}_{z,\overline{z}}(2)\,
{\rm O}_{z,\overline{z}}(2),
\endaligned
\]
and reminding $\maux_{z\overline{z}} = \frac{1}{1 - \zeta
\overline{\zeta}}$, this gives the concluding identity:
\[
{\rm O}_{z,\overline{z}}(3)
\,=\,
\tfrac{1}{1-\zeta\overline{\zeta}}\,
G_{\zeta\overline{\zeta}}^\kappa.
\qedhere
\]
\endproof

By integration, $G^\kappa = \lambda^\kappa (z, \zeta, 
\overline{z}, v) 
+ \overline{\lambda}^\kappa (\overline{z}, \overline{\zeta}, z, v) 
+ {\rm O}_{z,\overline{z}}(3)$. After absorption
in ${\rm O}_{z, \overline{z}}(3)$, we can 
assume that $\lambda^\kappa$ is of degree $\leqslant 2$ in
$(z, \overline{z})$, hence contains only monomials 
$z^a \zeta^b \overline{z}^c v^e$ with $a + c \leqslant 2$
and $a + b + c = \kappa$. So $b \geqslant
\kappa-2$.

Further, $G^\kappa(z, \zeta, 0, 0, v) \equiv 0$ 
imposes $\lambda^\kappa (z, \zeta, 0, v) \equiv 0$.
So $1 \leqslant c \leqslant 2$. 
Consequently, $\lambda^\kappa$ can contain only three monomials:
\[
\lambda^\kappa(z,\zeta,\overline{z},v)
\,=\,
a(v)\,\overline{z}\zeta^{\kappa-1}
+
b(v)\,z\overline{z}\,\zeta^{\kappa-2}
+
c(v)\,\overline{z}^2\zeta^{\kappa-2}.
\]
Since $\kappa \geqslant 4$, we see that the conjugate
$\overline{\lambda}^\kappa(\overline{z}, \overline{\zeta}, z, v)$
is multiple of $\overline{\zeta}^{\kappa-2 \geqslant 2}$, hence:
\[
G^\kappa\big(z,\zeta,\overline{z},0,v\big)
\,=\,
\lambda^\kappa(z,\zeta,\overline{z},v)
+
\zero{
\overline{\lambda}^\kappa(\overline{z},0,z,v)}
+
{\rm O}_{z,\overline{z}}(3).
\]

Finally, because the prenormalized
coordinates of Proposition~{\ref{Prp-prenormalization}} 
require $G^\kappa(z,\zeta,\overline{z},0,v) 
= {\rm O}_{z,\overline{z}}(3)$,
we reach $\lambda^\kappa(z,\zeta,\overline{z},v) = {\rm O}_{z,
\overline{z}}(3)$, which forces $a = b = c = 0 = \lambda^\kappa$, 
so as asserted $G^\kappa = {\rm O}_{z,\overline{z}}(3)$.
\endproof

\Section{\bf Consequence of Prenormalization on Equivalences}
\label{consequence-prenormalization-equivalences}
\HEAD{{\ref{consequence-prenormalization-equivalences}}.~{\sf 
Consequence of Prenormalization on Equivalences}
}{
Wei-Guo {\sc Foo}, Joël {\sc Merker}, The-Anh {\sc Ta}}

Thanks to Proposition~{\ref{Prp-prenormalization}},
if we are given a holomorphic map 
$H \colon (z, \zeta, w) \longmapsto (z', \zeta', w')$ between two
$\mathfrak{C}_{2,1}$ 
hypersurfaces $M^5 \subset \C^3$ 
and ${M'}^5 \subset {\C'}^3$, we can assume that
both hypersurfaces are prenormalized.
In particular, Proposition~{\ref{Prp-G-O-z-zbar-3}}
tells us that the
{\em whole} remainders after the GM-model part
of their graphing functions is of order $3$ in $(z, \overline{z})$: 
\[
u
\,=\,
\maux
+
G
\,=\,
\maux
+
{\rm O}_{z,\overline{z}}(3)
\ \ \ \ \ \ \ \ \ \ \ \ \ \ \ \ \ \ \ \
\text{and}
\ \ \ \ \ \ \ \ \ \ \ \ \ \ \ \ \ \ \ \
u'
\,=\,
\maux'
+
G'
\,=\,
\maux
+
{\rm O}_{z',\overline{z}'}(3).
\]

\begin{Observation}
Complex scalings $(z, \zeta, w) \longmapsto 
\big( \lambda z,\, \frac{\lambda}{\overline{\lambda}} \zeta,\,
\lambda\overline{\lambda} w \big)$ with $\lambda \in \C^\ast$ preserve
prenormalizations as in 
Proposition~{\ref{Prp-prenormalization}}.\qed
\end{Observation}

With $\lambda := \rho \in \R^\ast$, this is
$\big( \rho^1 z,\, \rho^0 \zeta,\, \rho^2 w\big)$. Hence
this observation suggests naturally to assign the following weights
to the three complex variables and their real and imaginary parts:
\[
[z]
\,:=\,
1
\,=:\,
[\overline{z}],
\ \ \ \ \ \ \ \ \ \ \ \ \ \ \ \ \ \ \ \
[\zeta]
\,:=\,
0
\,=:\,
[\overline{\zeta}],
\ \ \ \ \ \ \ \ \ \ \ \ \ \ \ \ \ \ \ \
[w]
\,:=\,
2
\,=:\,
[\overline{w}].
\]

Accordingly, 
let us decompose the components $(f,g,h)$ of $H$
in weighted homogeneous parts:
\[
f
\,=\,
f_0+f_1+f_2+f_3\cdots,
\ \ \ \ \ \ \ \ \ \ \ \ \ \ \ \ \ \ \ \
g
\,=\,
g_0+g_1+g_2+\cdots,
\ \ \ \ \ \ \ \ \ \ \ \ \ \ \ \ \ \ \ \
h
\,=\,
h_0+h_1+h_2+h_3+h_4+\cdots.
\]

\begin{Proposition}
\label{Prp-f-g-h-nu-0-1-2}
If both $M$ and $M'$ are prenormalized, 
possibly after composing
with a complex dilation $(z', \zeta', w') \longmapsto 
\big( \lambda z', \frac{\lambda}{\overline{\lambda}}
\zeta', \lambda\overline{\lambda} w')$, 
one has $f_0 = 0$, $f_1 = z$, $g_0 = \zeta$, 
$h_0 = 0$, $h_1 = 0$, $h_2 = w$, and
the weighted homogeneous components of
$f$, $g$, $h$ are:
\[
f
\,=\,
z
+
f_2
+
f_3
+\cdots,
\ \ \ \ \ \ \ \ \ \ \ \ \ \ \ \ \ \ \ \
g
\,=\,
\zeta
+
g_1
+
g_2
+\cdots,
\ \ \ \ \ \ \ \ \ \ \ \ \ \ \ \ \ \ \ \
h
\,=\,
w
+
h_3
+
h_4
+\cdots.
\]
\end{Proposition}

Mind the fact that this does not mean that the map
is $\Id + {\rm O}_{z,w,\zeta}(2)$, 
since in $f_2$, there can still be the linear term $f_{0,0,2}\, w$, 
and in $g_1 + g_2$, there can still be the
linear terms $g_{1,0,0}\, z + g_{0,0,1}\, w$.

\proof
The fundamental identity expressing that we have a
map $M \longrightarrow M'$ reads:
\leqnomode\usetagform{default}
\begin{footnotesize}
\begin{align}
\label{F-I-0-1-2}
h_0+h_1+\cdots
+
\overline{h}_0+\overline{h}_1+\cdots
&
\,=\,
2\,F'
\Big(
f_0+f_1+\cdots,\,
g_0+g_1+\cdots,\,
\notag
\\
{}
&
\ \ \ \ \ \ \ \ \ \ \ \ \ \ \ \ \ 
\overline{f}_0+\overline{f}_1+\cdots,\,
\overline{g}_0+\overline{g}_1+\cdots,\,
\tfrac{1}{2i}
\big(
h_0+h_1+\cdots
-
\overline{h}_0-\overline{h}_1
-\cdots
\big)
\Big).
\end{align}
\end{footnotesize}
Observe that $f_0 = f_0(\zeta)$, 
$g_0 = g_0(\zeta)$, $h_0 = h_0(\zeta)$ depend only on $\zeta$.
This identity projected to weight $0$ becomes:
\[
h_0(\zeta)
+
\overline{h}_0(\overline{\zeta})
\,\equiv\,
2\,F'
\Big(
f_0(\zeta),\,
g_0(\zeta),\,
\overline{f}_0(\overline{\zeta}),\,
\overline{g}_0(\overline{\zeta}),\,\,
\tfrac{1}{2i}h_0(\zeta)
-
\tfrac{1}{2i}\overline{h}_0(\overline{\zeta})
\Big).
\]
Put $\overline{\zeta} := 0$, use the assumption 
that there are no pluriharmonic
terms (coordinates are prenormalized), namely that
$0 \equiv F' (z', \zeta', 0, 0, v')$, and get $h_0(\zeta) \equiv 0$.

Once again, look at~({\ref{F-I-0-1-2}}), and get from 
$F' = \maux' + G' = \maux' + {\rm O}_{z', \overline{z}'}(3)$:
\[
0
\,\equiv\,
\frac{2\,f_0(\zeta)\overline{f}_0(\overline{\zeta})
+f_0(\zeta)^2\,\overline{g}_0(\overline{\zeta})
+\overline{f}_0(\overline{\zeta})^2\,g_0(\zeta)}{
1-g_0(\zeta)\overline{g}_0(\overline{\zeta})}
+
{\rm O}_{f_0(\zeta),\overline{f}_0(\zeta)}(3).
\]
We claim that $f_0(\zeta) \equiv 0$. Otherwise, $f_0 = e\, \zeta^\tau
+ {\rm O}_\zeta(\tau+1)$ with $e \neq 0$
and $\tau \in \N_{\geqslant 1}$. Hence:
\[
0
\,\equiv\,
2\,e\overline{e}\,
\zeta^\tau
\overline{\zeta}^\tau\,
\big(
1
+
{\rm O}_{\zeta,\overline{\zeta}}(1)
\big)
+
\zeta^{2\tau}\,
(\cdots)
+
\overline{\zeta}^{2\tau}
(\cdots)
+
{\rm O}_{\zeta,\overline{\zeta}}\big(3\,\tau\big),
\]
and this forces $e\overline{e} = 0$. So $f_0(\zeta) \equiv 0$,
and~({\ref{F-I-0-1-2}}) at weight $0$, namely the identity above,
reduces to $0 = 0$.

Next, examine weight $1$. Certainly, $f_1 = z f_1(\zeta)$ and
$h_1 = z h_1(\zeta)$, while $g$ will not participate here.
Since $\maux'$ is weighted $2$-homogeneous, 
as it contains $z\overline{z}$, $\overline{z}^2$, $z^2$
times functions of $(\zeta, \overline{\zeta})$, 
we have 
$F' = {\rm O}_{z',\overline{z}'}(2)$, so the identity:
\[
z\,h_1(\zeta)
+
\overline{z}\,\overline{h}_1(\overline{\zeta})
\,\equiv\,
{\rm O}_{zf_1(\zeta),\overline{z}\overline{f}_1(\overline{\zeta})}(2)
\,\,=\,\,
{\rm O}_{z,\overline{z}}(2),
\]
forces $h_1(\zeta) \equiv 0$.

Next, expand in powers of $z$, $w$:
\[
\!\!\!\!\!\!\!\!\!\!\!\!\!\!\!
\aligned
f
\,=\,
z\,f_1(\zeta)
+
z^2(\cdots)
+
w(\cdots),
\ \ \ \ \ \ \
g
\,=\,
g_0(\zeta)
+
z\,g_1(\zeta)
+
z^2(\cdots)
+
w(\cdots),
\ \ \ \ \ \ \ 
h
&
\,=\,
h_2+h_3+\cdots,
\\
h_2
&
\,=\,
z^2\,\varphi(\zeta)
+
w\,\psi(\zeta).
\endaligned
\]
The holomorphic Jacobian at the origin is assumed to be invertible:
\[
0
\,\neq\,
\left\vert\!
\begin{array}{ccc}
f_z(0) & f_\zeta(0) & f_w(0)
\\
g_z(0) & g_\zeta(0) & g_w(0)
\\
h_z(0) & h_\zeta(0) & h_w(0)
\end{array}
\!\right\vert
\,\,=\,\,
\left\vert\!
\begin{array}{ccc}
f_1(0) & 0 & f_w(0)
\\
g_1(0) & g_0'(0) & g_w(0)
\\
0 & 0 & h_w(0)
\end{array}
\!\right\vert,
\]
whence $h_w(0) \neq 0$ and $g_0'(0) \neq 0$ and also
$f_1(0) \neq 0$.
Then the fundamental identity~({\ref{F-I-0-1-2}}) 
becomes in weight $2$:
\[
h_2(\zeta)
+
\overline{h}_2(\overline{\zeta})
\,\equiv\,
2\,\maux'
\big(
zf_1(\zeta),\,
g_0(\zeta),\,
\overline{z}\overline{f}_1(\overline{\zeta}),\,
\overline{g}_0(\overline{\zeta})
\big),
\]
that is, after replacing $w = \maux + iv$ in $h_2$:
\[
\aligned
z^2\,\varphi(\zeta)
+
\overline{z}^2\,\overline{\varphi}(\overline{\zeta})
&
+
\maux(z,\zeta,\overline{z},\overline{\zeta})\,
\big[
\psi(\zeta)
+
\overline{\psi}(\overline{\zeta})
\big]
+
i\,v\,
\big[
\psi(\zeta)
-
\overline{\psi}(\overline{\zeta})
\big]
\,\equiv\,
\\
&
\,\equiv\,
\frac{2\,z\,f_1(\zeta)\,\overline{z}\,
\overline{f}_1(\overline{\zeta})
+
\overline{z}^2\,\overline{f}_1(\overline{\zeta})\,g_0(\zeta)
+
z^2\,f_1(\zeta)^2\,\overline{g}_0(\overline{\zeta})}{
1-g_0(\zeta)\,\overline{g}_0(\overline{\zeta})},
\endaligned
\]
this holding identically in $\C\{ z, \zeta, \overline{z}, 
\overline{\zeta}, v \}$. This forces 
$\psi(\zeta) \equiv \rho$ to be constant, with $\rho \in \R^\ast$,
and then $\varphi(\zeta) \equiv 0$ necessarily. 

It remains an identity:
\[
\maux(z,\zeta,\overline{z},\overline{\zeta})\,
2\,\rho
\,\equiv\,
2\,\maux'
\big(
zf_1(\zeta),\,
g_0(\zeta),\,
\overline{z}\overline{f}_1(\overline{\zeta}),\,
\overline{g}_0(\overline{\zeta})
\big),
\]
which expresses that the map $(z, \zeta, w) 
\longmapsto \big(zf_1(\zeta),\, g_0(\zeta),\,
\rho\, w\big)$ is an {\em automorphism}\,\,---\,\,in 
fact a {\em rigid} automorphism, 
{\em cf.}~{\cite{Chen-Foo-Merker-Ta-2019}}\,\,---\,\,of
the Gaussier-Merker model. But we know from
Section~{\ref{fractional-representation-isotropy-group}},
see the fractional expression of $w'$ there, that this requires
$\alpha = 0$ and $r = 0$, while
only $\lambda \in \C^\ast$ is free.
Consequently, the map is of the form 
$(f_1, g_0, h_2) = \big(\lambda z,\,
\frac{\lambda}{\overline{\lambda}} \zeta,\,
\lambda\overline{\lambda} w\big)$. Post-composing
by the inverse map yields the conclusion.
\endproof

\Section{\bf Uniqueness of Normal Form}
\label{uniqueness-normal-form}
\HEAD{{\ref{uniqueness-normal-form}}.~{\sf Uniqueness of Normal Form}
}{
Wei-Guo {\sc Foo}, Joël {\sc Merker}, The-Anh {\sc Ta}}

Starting with a $\mathcal{C}^\omega$ hypersurface
$M^5 \subset \C^3$ which is
$2$-nondegenerate and of constant Levi rank $1$, 
at any point $p \in M$, it is elementary to find
holomorphic coordinates $(z, \zeta,
w)$ vanishing at $p$ in which $M$ has
equation: 
\leqnomode\usetagform{default}
\begin{align}
\label{u-5}
u
\,=\,
F
\,=\,
z\overline{z}
+
\tfrac{1}{2}\,\overline{z}^2\zeta
+
\tfrac{1}{2}\,z^2\overline{\zeta}
+
z\overline{z}\zeta\overline{\zeta}
+
{\rm O}_{z,\zeta,\overline{z},\overline{\zeta},v}(5).
\end{align}
Such an equation can hence freely 
be taken as the starting point towards a complete
normalization of $F(z, \zeta, \overline{z}, \overline{\zeta}, v)$.

In the preceding sections, we have in fact established the {\em
existence} of a {\sl normal form} for $M$. We can now present 
a final {\sl uniqueness} statement which will terminate our article.

\begin{Theorem}
Given $M^5 \subset \C^3$ in the class $\mathfrak{C}_{2,1}$
with $0 \in M$ of the form:
\[
u
\,=\,
z\overline{z}
+
\tfrac{1}{2}\,\overline{z}^2\zeta
+
\tfrac{1}{2}\,z^2\overline{\zeta}
+
z\overline{z}\zeta\overline{\zeta}
+
{\rm O}_{z,\zeta,\overline{z},\overline{\zeta},v}(5),
\]
there exists a biholomorphism $(z,\zeta,w) \longmapsto (z', \zeta',
w')$ fixing $0$ which maps $(M, 0)$ into $(M',0)$ of normalized
equation:
\[
\aligned
u'
&
\,=\,
z'\overline{z}'
+
\tfrac{1}{2}\,{\overline{z}'}^2\zeta'
+
\tfrac{1}{2}\,{z'}^2\overline{\zeta}'
+
z'\overline{z}'\zeta'\overline{\zeta}'
+
\tfrac{1}{2}\,{\overline{z}'}^2\zeta'\zeta'\overline{\zeta}'
+
\tfrac{1}{2}\,{z'}^2\overline{\zeta}'\zeta'\overline{\zeta}'
+
z'\overline{z}'\zeta'\overline{\zeta}'\zeta'\overline{\zeta}'
\notag
\\
&
\ \ \ \ \ \ \ \ \ \ \ 
+
0
+
{z'}^3{\overline{z}'}^3\,
{\rm O}_{z',\overline{z}'}(1)
\notag
\\
&
\ \ \ \ \ \ \ \ \ \ \ \ 
+
2\,\Re\,
\Big\{
0
+
0
+
{z'}^3{\overline{\zeta}'}^2\,
F_{3,0,0,2}'(v')
+
\zeta'\overline{\zeta}'\,
\big(
3\,{z'}^2\overline{z}'\overline{\zeta}'\,
F_{3,0,0,2}'(v')
\big)
\big\}
\notag
\\
&
\ \ \ \ \ \ \ \ \ \ \ \ 
+
2\,\Re\,
\Big\{
{z'}^5\overline{\zeta}'\,
F_{5,0,0,1}'(v')
+
\ \ \ \ \ \ \ \ \ \ \ \ \ \
0
\ \ \ \ \ \ \ \ \ \ \ \ \ \ 
+
{z'}^4{\overline{\zeta}'}^2\,
F_{4,0,0,2}'(v')
\notag
\\
&
\ \ \ \ \ \ \ \ \ \ \ \ \ \ \ \ \ \ \ \ \ \ \ \ \ \ \ \ \ \ \ \ \ \ \ 
\ \ \ \ \ \ \ \ \ \ \ \ \ \ \ \ \ \,
+
{z'}^3{\overline{z}'}^2\overline{\zeta}'\,
F_{3,0,2,1}'(v')
+
{z'}^3\overline{z}'{\overline{\zeta}'}^2\,
F_{3,0,1,2}'(v')
\notag
\\
&
\ \ \ \ \ \ \ \ \ \ \ \ \ \ \ \ \ \ \ \ \ \ \ \ \ \ \ \ \ \ \ \ \ \ \
\ \ \ \ \ \ \ \ \ \ \ \ \ \ \ \ \ \ \ \ \ \ \ \ \ \ \ \ \ \ \ \ \ \ \
\ \ \ \ \ \ \ \ \ \ \ \ \ \ \ \ \
+
{z'}^3{\overline{\zeta}'}^3\,
F_{3,0,0,3}'(v')
\Big\}
\notag
\\
&
\ \ \ \ \ \ \ \ \ \ \ 
+
{\overline{z}'}^3\zeta'\,
{\rm O}_{z',\zeta',\overline{z}'}(3)
+
{z'}^3\overline{\zeta}'\,
{\rm O}_{z',\overline{z}',\overline{\zeta}'}(3)
+
\zeta'\overline{\zeta}'\,
{\rm O}_{z',\zeta',\overline{z}'}(3)\,
{\rm O}_{z',\zeta',\overline{z}',\overline{\zeta}'}(2).
\endaligned
\]

Furthermore, the map exists and is {\em unique} if it is assumed
to be of the form:
\[
\aligned
z'
\,:=\,
&\,
z
+
f_{\geqslant 2}(z,\zeta,w)
\ \ \ \ \ \ \ \ \ \ \ \ \ \ \ \ \ \ \ \
\zeta'
\,:=\,
\zeta
+
g_{\geqslant 1}(z,\zeta,w),
&
\ \ \ \ \ \ \ \ \ \ \ \ \ \ \ \ \ \ \ \
w'
\,:=\,
&\,
w
+
h_{\geqslant 3}(z,\zeta,w),
\\
0
\,=\,
&\,
f_w(0),
&
\ \ \ \ \ \ \ \ \ \ \ \ \ \ \ \ \ \ \ \
0
\,=\,
&\,
\Im\,h_{ww}(0).
\endaligned
\]
\end{Theorem}

Here of course, $f_{\geqslant 2}$ 
is of weight $\geqslant 2$, while $g_{\geqslant 1}$ is of weight
$\geqslant 1$, and $h_{\geqslant 3}$ is
of weight $\geqslant 3$ for the currently useful weighting $[z]
:= 1$, $[\zeta] := 0$, $[w] := 2$.

\proof
By choosing a chain at $0 \in M$ whose first jet is flat,
directed along the $v$-axis, one can verify (exercise) that
all the constructions done in the preceding sections do 
indeed give a biholomorphism of this specific form.
So our job is to establish uniqueness.

\smallskip

Suppose therefore that two such normalizations $H_\iota \colon 
(z,\zeta, w)
\longmapsto (z+f_\iota,\, \zeta+g_\iota, w + h_\iota)$, 
$\iota = 1, 2$, 
are given:
\[
\xymatrix{
& & &
M_1'
\ar[dd]^{H_2\circ H_1^{-1}}
\\
\ar[urrr]^{H_1}
M
\ar[drrr]_{H_2}
\\
& & &
M_2',}
\]
with $0 = f_{\iota,w}(0)$
and $0 = {\rm Re}\, 
h_{\iota,ww}(0)$ for $\iota = 1, 2$.
We leave to the reader to verify
that, then, 
$H := H_2 \circ H_1^{-1}$ is also of the form
$(z, \zeta, w) \longmapsto \big( z + f_{\geqslant 2},\,
\zeta + g_{\geqslant 1},\, w + h_{\geqslant 3} \big)$
also with $0 = f_w(0)$ and $0 = \Im\, h_{ww}(0)$.
For this, one has to take account of~({\ref{u-5}}).

The theorem asserts that $H_1 = H_2$. Equivalently, 
$H_2 \circ H_1^{-1} = \Id$. This will be offered by 
the next independent key uniqueness statement.
\endproof

\begin{Theorem}
\label{Thm-uniqueness-complete-normal-form-M5-C3}
For a given $M^5 \subset \C^3$ in the class $\mathfrak{C}_{2,1}$,
if two normal forms $N$ and $N_\prime$ 
at some point $p \in M$ are constructed, with $N$ having 
normalized equation:
\[
\aligned
u
&
\,=\,
z\overline{z}
+
\tfrac{1}{2}\,\overline{z}^2\zeta
+
\tfrac{1}{2}\,z^2\overline{\zeta}
+
z\overline{z}\zeta\overline{\zeta}
+
\tfrac{1}{2}\,\overline{z}^2\zeta\zeta\overline{\zeta}
+
\tfrac{1}{2}\,z^2\overline{\zeta}\zeta\overline{\zeta}
+
z\overline{z}\zeta\overline{\zeta}\zeta\overline{\zeta}
\notag
\\
&
\ \ \ \ \ \ \ \ \ \ \ 
+
0
+
z^3\overline{z}^3\,
{\rm O}_{z,\overline{z}}(1)
\notag
\\
&
\ \ \ \ \ \ \ \ \ \ \ \ 
+
2\,\Re\,
\Big\{
0
+
0
+
z^3\overline{\zeta}^2\,
F_{3,0,0,2}(v)
+
\zeta\overline{\zeta}\,
\big(
3\,z^2\overline{z}\overline{\zeta}\,
F_{3,0,0,2}(v)
\big)
\big\}
\notag
\\
&
\ \ \ \ \ \ \ \ \ \ \ \ 
+
2\,\Re\,
\Big\{
z^5\overline{\zeta}\,
F_{5,0,0,1}(v)
+
\ \ \ \ \ \ \ \ \ \ \ \ 
0
\ \ \ \ \ \ \ \ \ \ \ \ 
+
z^4\overline{\zeta}^2\,
F_{4,0,0,2}(v)
\notag
\\
&
\ \ \ \ \ \ \ \ \ \ \ \ \ \ \ \ \ \ \ \ \ \ \ \ \ \ \ \ \ \ \ \ \ \ \ 
\ \ \ \ \ \ \ \ \ \ \ \ \ \ \,
+
z^3\overline{z}^2\overline{\zeta}\,
F_{3,0,2,1}(v)
+
z^3\overline{z}\overline{\zeta}^2\,
F_{3,0,1,2}(v)
\notag
\\
&
\ \ \ \ \ \ \ \ \ \ \ \ \ \ \ \ \ \ \ \ \ \ \ \ \ \ \ \ \ \ \ \ \ \ \
\ \ \ \ \ \ \ \ \ \ \ \ \ \ \ \ \ \ \ \ \ \ \ \ \ \ \ \ \ \ \ \ \ \ \
\ \ \ \ \ \ \ \ \ \ \,
+
z^3\overline{\zeta}^3\,
F_{3,0,0,3}(v)
\Big\}
\notag
\\
&
\ \ \ \ \ \ \ \ \ \ \ 
+
\overline{z}^3\zeta\,
{\rm O}_{z,\zeta,\overline{z}}(3)
+
z^3\overline{\zeta}\,
{\rm O}_{z,\overline{z},\overline{\zeta}}(3)
+
\zeta\overline{\zeta}\,
{\rm O}_{z,\zeta,\overline{z}}(3)\,
{\rm O}_{z,\zeta,\overline{z},\overline{\zeta}}(2),
\endaligned
\]
and with $N_\prime$ having similarly normalized equation,
and if the map $(z, \zeta, w) \longmapsto (z', \zeta', w')$ 
between them is of the form:
\[
\aligned
z'
\,:=\,
&\,
z
+
f_{\geqslant 2}(z,\zeta,w)
\ \ \ \ \ \ \ \ \ \ \ \ \ \ \ \ \ \ \ \
\zeta'
\,:=\,
\zeta
+
g_{\geqslant 1}(z,\zeta,w),
&
\ \ \ \ \ \ \ \ \ \ \ \ \ \ \ \ \ \ \ \
w'
\,:=\,
&\,
w
+
h_{\geqslant 3}(z,\zeta,w),
\\
0
\,=\,
&\,
f_w(0),
&
\ \ \ \ \ \ \ \ \ \ \ \ \ \ \ \ \ \ \ \
0
\,=\,
&\,
\Im\,h_{ww}(0),
\endaligned
\]
then the map $(z',\zeta',w') = (z, \zeta, w)$
is the identity, and the two normal forms 
$N = N_\prime$ coincide.
\end{Theorem}

\proof
Equivalently, the graphing function $F = \sum_{a,b,c,d}\,
z^a\zeta^b \overline{z}^c \overline{\zeta}^d\, F_{a,b,c,d}(v)$
of $N$ satisfies the general {\sl prenormalization conditions:}
\[
\aligned
0
\,\equiv\,
&
F_{a,b,0,0}(v)
\equiv\,
F_{0,0,c,d}(v),
\\
0
\,\equiv\,
&
F_{a,b,1,0}(v)
\,\equiv\,
F_{1,0,c,d}(v),
\\
0
\,\equiv\,
&
F_{a,b,2,0}(v)
\,\equiv\,
F_{2,0,c,d}(v),
\endaligned
\]
with the obvious two exceptions 
$F_{1,0,1,0}(v) \equiv 1$ and $F_{0,1,2,0}(v) \equiv
\frac{1}{2} \equiv F_{2,0,0,1}(v)$, together with the
{\sl sporadic normalization conditions}, listed by 
increasing order $4$, $5$, $6$:
\[
\aligned
0
&
\,\equiv\,
F_{3,0,0,1}(v)
\,\equiv\,
F_{0,1,3,0}(v),
&
\ \ \ \ \ \ \ \ \ \ \ \ \ \ \ \ \ \ \ \ \ \ \ \ \ \
&
\\
0
&
\,\equiv\,
F_{4,0,0,1}(v)
\,\equiv\,
F_{0,1,4,0}(v),
&
\ \ \ \ \ \ \ \ \ \ \ \ \ \ \ \ \ \ \ \ \ \ \ \ \ \
0
&
\,\equiv\,
F_{3,0,1,1}(v)
\,\equiv\,
F_{1,1,3,0}(v),
\\
0
&
\,\equiv\,
F_{4,0,1,1}(v)
\,\equiv\,
F_{1,1,4,0}(v),
&
\ \ \ \ \ \ \ \ \ \ \ \ \ \ \ \ \ \ \ \ \ \ \ \ \ \
0
&
\,\equiv\,
F_{3,0,3,0}(v),
\endaligned
\]
and the same holds about $F'$.

Accordingly, let us introduce:
\[
\footnotesize
\aligned
S
\,:=\,
&\,
\Big\{
(a,b,0,0),\,
(0,0,c,d),\,
(a,b,1,0),\,
(1,0,c,d),\,
(a,b,2,0),\,
(2,0,c,d)
\Big\}
\\
&\,
\medcup
\Big\{
(3,0,0,1),\,
(0,1,3,0),\,\,
(4,0,0,1),\,
(0,1,4,0),\,
(3,0,1,1),\,
(1,1,3,0),\,\,
(4,0,1,1),\,
(1,1,4,0),\,
(3,0,3,0)
\Big\}.
\endaligned
\]
Notice that $S$ takes {\em no} dependent 
derivatives $\zeta \overline{\zeta}
(\cdots)$, namely one always has
$b + d \leqslant 1$ for any
$(a,b,c,d) \in S$.

For a general real converging power series vanishing at
$(z,\zeta, \overline{z}, \overline{\zeta}, v) = (0,0,0,0,0)$:
\[
H
\,=\,
\sum_{a,b,c,d,e}\,
H_{a,b,c,d,e}\,
z^a\zeta^b\overline{z}^c\overline{\zeta}^dv^e
\eqno
{\scriptstyle{(\overline{H_{c,d,a,b,e}}\,=\,H_{a,b,c,d,e})}},
\]
{\em i.e.} with $H_{0,0,0,0,0} = 0$, introduce the projection:
\[
\Pi_S(H)
\,:=\,
\sum_{(a,b,c,d)\in S}\,
\sum_{e=0}^\infty\,
H_{a,b,c,d,e}\,
z^a\zeta^b\overline{z}^c\overline{\zeta}^dv^e,
\]
so that:
\[
\Pi_S(F)
\,=\,
z\overline{z}
+
\tfrac{1}{2}\,\overline{z}^2\zeta
+
\tfrac{1}{2}\,z^2\overline{\zeta}
\ \ \ \ \ \ \ \ \ \ \ \ \ \ \ \ \ \ \ \
\text{and}
\ \ \ \ \ \ \ \ \ \ \ \ \ \ \ \ \ \ \ \
\Pi_S(F')
\,=\,
z'\overline{z}'
+
\tfrac{1}{2}\,{\overline{z}'}^2\zeta'
+
\tfrac{1}{2}\,{z'}^2\overline{\zeta}'.
\]

By assumption (or because of 
Proposition~{\ref{Prp-f-g-h-nu-0-1-2}}), the map is of the form:
\[
z'
\,=\,
z
+
f_2
+
f_3
+\cdots,
\ \ \ \ \ \ \ \ \ \ \ \ \ \ \ \ \ \ \ \
\zeta'
\,=\,
\zeta
+
g_1
+
g_2
+\cdots,
\ \ \ \ \ \ \ \ \ \ \ \ \ \ \ \ \ \ \ \
w'
\,=\,
w
+
h_3
+
h_4
+\cdots,
\]
that is, more precisely:
\[
\aligned
f
&
\,=\,
\sum_{\nu\geqslant 3}\,
f_{\nu-1}
\,\,=\,\,
\sum_{\nu\geqslant 3}\,
\bigg(
\sum_{a+b+2e=\nu-1}\,
f_{a,b,e}\,z^a\zeta^bw^e
\bigg),
\\
g
&
\,=\,
\sum_{\nu\geqslant 3}\,
g_{\nu-2}
\,\,=\,\,
\sum_{\nu\geqslant 3}\,
\bigg(
\sum_{a+b+2e=\nu-2}\,
g_{a,b,e}\,z^a\zeta^bw^e
\bigg),
\\
h
&
\,=\,
\sum_{\nu\geqslant 3}\,
h_\nu
\,\,=\,\,
\sum_{\nu\geqslant 3}\,
\bigg(
\sum_{a+b+2e=\nu}\,
h_{a,b,e}\,z^a\zeta^bw^e
\bigg).
\endaligned
\]

Let us introduce the projections:
\[
\pi_{\nu-1}(f)
\,:=\,
f_{\nu-1},
\ \ \ \ \ \ \ \ \ \ \ \ \ \ \ \ \ \ \ \
\pi_{\nu-2}(g)
\,:=\,
g_{\nu-2},
\ \ \ \ \ \ \ \ \ \ \ \ \ \ \ \ \ \ \ \
\pi_\nu(h)
\,:=\,
h_\nu,
\]
and also:
\[
\pi_\nu(H)
\,:=\,
\sum_{a+b+c+d+2e=\nu}\,
H_{a,b,c,d,e}\,
z^a\zeta^b\overline{z}^c\overline{\zeta}^dv^e,
\]
so that:
\[
\Pi_S\big(\pi_\nu(F)\big)
\,\,=\,\,
0
\,\,=\,\,
\Pi_S\big(\pi_\nu(F')\big)
\eqno
{\scriptstyle{(\forall\,\nu\,\geqslant\,3)}}.
\]
Also, let us introduce:
\[
\pi^\nu
\,:=\,
\pi_2
+\cdots+
\pi_\nu.
\]

Now, 
remind that $\maux = \frac{z\overline{z} + \frac{1}{2}
\overline{z}^2 \zeta + z^2 \overline{\zeta}}{
1 - \zeta \overline{\zeta}}$ is homogeneous of weight $2$.
Thanks to Proposition~{\ref{Prp-G-O-z-zbar-3}}, we may write:
\[
u
\,=\,
F
\,=\,
\maux
+
\smallsum{\nu\geqslant 3}\,
G_\nu.
\]
Then
for a holomorphic function $e_\mu = e_\mu(z, \zeta, w)$ which 
is weighed $\mu$-homogeneous, it holds (exercise):
\leqnomode\usetagform{default}
\begin{align}
\label{pi-mu-e-mu}
\pi^\mu
\bigg(
e_\mu
\Big(
z,\,\zeta,\,\,
i\,v
+
\maux(z,\zeta,\overline{z},\overline{\zeta})
+
\smallsum{\nu\geqslant 3}\,
G_\nu
\big(
z,\zeta,\overline{z},\overline{\zeta},v
\big)
\Big)
\bigg)
\,\,=\,\,
e_\mu
\big(z,\zeta,\,i\,v+\maux\big).
\end{align}

Now, the fundamental identity expressing that $(z+f, \zeta+g, w+h)$
is a map $N \longrightarrow N_\prime$ writes:
\leqnomode\usetagform{default}
\begin{align}
\label{weighted-F-I}
0
&
\,\equiv\,
-\,\Re\,
\big(
w+h_3+h_4+\cdots
\big)
\notag
\\
&
\ \ \ \ \
+
F'\Big(
z+f_2+f_3+\cdots,\,\,
\zeta+g_1+g_2+\cdots,\,\,
\notag
\\
&
\ \ \ \ \ \ \ \ \ \ \ \ \ \ \ \ \ 
\overline{z}+\overline{f}_2+\overline{f}_3+\cdots,\,\,
\overline{\zeta}+\overline{g}_1+\overline{g}_2+\cdots,\,\,
\Im\,\big(w+h_3+h_4+\cdots\big)
\Big).
\end{align}

In order to prove that $(f, g, h) = (0, 0, 0)$, we may 
proceed progressively, by induction on $\nu \geqslant 3$:

\medskip\noindent$(\bullet_3)$\,
$(f_2, g_1, h_3) = (0, 0, 0)$;

\medskip\noindent$(\bullet_{\nu-1})$\,
$\big(f_{\mu-1}, g_{\mu-2}, h_\mu \big) = (0,0,0)$ for
$\mu = 3, \dots, \nu-1$ and some $\nu \geqslant 4$
implies that $\big( f_{\nu-1}, g_{\nu-2}, h_\nu \big) = (0, 0, 0)$.

\medskip

Therefore, let us examine first the fundamental identity
in weight $\nu = 3$, remembering that this identity
already holds true in weights $0$, $1$, $2$\,\,---\,\,according
to (the proof of) Proposition~{\ref{Prp-f-g-h-nu-0-1-2}},
or according to our hypothesis\,\,---:
\[
\aligned
0
&
\,\equiv\,
\pi^3
\bigg(
-\,\Re\,\big(w+h_3\big)
+
\maux'
\big(
z+f_2,
\zeta+g_1,
\overline{z}+\overline{f}_2,
\overline{\zeta}+\overline{g}_1
\big)
+
F_3'\big(z+f_2,\zeta+g_1,\overline{z}+\overline{f}_2,
\overline{\zeta}+\overline{g}_1\big)
\bigg)
\\
&
\,\equiv\,
\pi^3
\bigg(
-\,\maux
-
F_3
-
\Re\,h_3
+
\maux'
\big(
z+f_2,
\zeta+g_1,
\overline{z}+\overline{f}_2,
\overline{\zeta}+\overline{g}_1
\big)
\bigg)
+
F_3'\big(z,\zeta,\overline{z},\overline{\zeta}\big),
\endaligned
\]
since $\maux'$ is weighted homogeneous of degree $2$,
since we use here~({\ref{pi-mu-e-mu}}).
Equivalently:
\[
F_3\big(z,\zeta,\overline{z},\overline{\zeta}\big)
-
F_3'\big(z,\zeta,\overline{z},\overline{\zeta}\big)
\,\equiv\,
\pi^3
\Big(
\maux'
\big(
z+f_2,
\zeta+g_1,
\overline{z}+\overline{f}_2,
\overline{\zeta}+\overline{g}_1
\big)
-
\maux\big(z,\zeta,\overline{z},\overline{\zeta}\big)
\Big)
-
\Re\,h_3\big(z,\zeta,\maux+iv\big).
\]

Generally,
for any $\nu \geqslant 3$, starting from the
induction assumption
expressed by $(\bullet_{\nu-1})$ above, 
the same reasoning (exercise) conducts to the identity:
\[
\!\!\!\!\!\!\!\!\!\!\!\!\!\!\!
F_\nu\big(z,\zeta,\overline{z},\overline{\zeta}\big)
-
F_\nu'\big(z,\zeta,\overline{z},\overline{\zeta}\big)
\,\equiv\,
\pi^\nu
\Big(
\maux'
\big(
z+f_{\nu-1},
\zeta+g_{\nu-2},
\overline{z}+\overline{f}_{\nu-1},
\overline{\zeta}+\overline{g}_{\nu-2}
\big)
-
\maux\big(z,\zeta,\overline{z},\overline{\zeta}\big)
\Big)
-
\Re\,h_\nu\big(z,\zeta,\maux+iv\big).
\]
Observe that:
\[
\maux_z
\,=\,
\tfrac{\overline{z}+z\overline{\zeta}}{1-\zeta\overline{\zeta}}
\ \ \ \ \ \ \ \ \ \ \ \ \ \ \ \ \ \ \ \
\text{and}
\ \ \ \ \ \ \ \ \ \ \ \ \ \ \ \ \ \ \ \
\maux_\zeta
\,=\,
\tfrac{1}{2}\,
\tfrac{(\overline{z}+z\overline{\zeta})^2}{
(1-\zeta\overline{\zeta})^2}.
\]

\begin{Lemma}
One has:
\[
\aligned
\pi^\nu
&
\Big(
\maux'
\big(
z+f_{\nu-1},
\zeta+g_{\nu-2},
\overline{z}+\overline{f}_{\nu-1},
\overline{\zeta}+\overline{g}_{\nu-2}
\big)
-
\maux\big(z,\zeta,\overline{z},\overline{\zeta}\big)
\Big)
\\
&
\ \ \ \ \
\,\,=\,\,
2\,\Re\,
\Big\{
\tfrac{\overline{z}+z\overline{\zeta}}{1-\zeta\overline{\zeta}}\,
f_{\nu-1}
\big(z,\zeta,\,\maux+iv\big)
+
\tfrac{1}{2}\,
\tfrac{(\overline{z}+z\overline{\zeta})^2}{
(1-\zeta\overline{\zeta})^2}\,
g_{\nu-2}
\big(z,\zeta,\,\maux+iv\big)
\Big\}.
\endaligned
\]
\end{Lemma}

\proof
The reader is referred to~{\cite[Prp.~6.2]{Chen-Foo-Merker-Ta-2019}}
which provides all arguments.
\endproof

Next, let us apply $\Pi_S (\centersmallbullet)$ to the above
identity, multiplied by $2$, namely to:
\[
2\,F_\nu
-
2\,F_\nu'
\,\equiv\,
\pi^\nu
\big(
2\,\maux'
-
2\,\maux
\big)
-
2\,\Re\,h_\nu,
\]
so that all monomials in the left-hand side disappear due to 
our assumption that both $N$ and $N_\prime$ are in normal form:
\[
\aligned
0
&
\,\equiv\,
\Pi_S
\bigg(
\pi^\nu
\big(
2\,\maux'
-
2\,\maux
\big)
-
2\,
\Re\,h_\nu
\bigg)
\\
&
\,\equiv\,
\Pi_S
\bigg(
2\,\Re\,
\Big\{
2\,
\tfrac{\overline{z}+z\overline{\zeta}}{1-\zeta\overline{\zeta}}\,
f_{\nu-1}
\big(z,\zeta,\,\maux+iv\big)
+
\tfrac{(\overline{z}+z\overline{\zeta})^2}{
(1-\zeta\overline{\zeta})^2}\,
g_{\nu-2}
\big(z,\zeta,\,\maux+iv\big)
-
h_\nu
\big(z,\zeta,\,\maux+iv\big)
\Big\}
\bigg).
\endaligned
\]

Then for all monomials $z^a \zeta^b \overline{z}^c \overline{\zeta}^d
v^e$ with $(a,b,c,d) \in S$ and $a + b + c + d + 2e = \nu$,
we obtain a system of linear equations:
\[
\aligned
({\sf E}_\nu)
\colon\ \ \ \ \ \ \ \ \ \ 
0
\,=\,
L_{a,b,c,d,e}
\Big(
\big\{
f_{a',b',e'}
\big\}_{a'+b'+2e'=\nu-1},\,\,
\big\{
g_{a',b',e'}
\big\}_{a'+b'+2e'=\nu-2},\,\,
\big\{
h_{a',b',e'}
\big\}_{a'+b'+2e'=\nu}
\Big).
\endaligned
\]

On the other hand, by considering
the complete $f = f_2 + f_3 + \cdots$, the complete
$g = g_1 + g_2 + \cdots$, and the
complete $h = h_3 + h_4 + \cdots$\,\,---\,\,not 
to be confused with the previous
$(z',\zeta',w') = (z,\zeta,w) + (f,g,h)$\,\,---, we can introduce the
analog `complete' linear system:
\[
0
\,\equiv\,
\Pi_S
\bigg(
2\,\Re\,
\Big\{
2\,
\tfrac{\overline{z}+z\overline{\zeta}}{1-\zeta\overline{\zeta}}\,
f
\big(z,\zeta,\,\maux+iv\big)
+
\tfrac{(\overline{z}+z\overline{\zeta})^2}{
(1-\zeta\overline{\zeta})^2}\,
g
\big(z,\zeta,\,\maux+iv\big)
-
h
\big(z,\zeta,\,\maux+iv\big)
\Big\}
\bigg),
\]
which, similarly, after extracting the coefficients of
all monomials $z^a\zeta^b \overline{z}^c \overline{\zeta}^d v^e$ with 
$(a,b,c,d) \in S$ and any $e \in \N$, 
can be abbreviated as:
\[
({\sf E})
\colon\ \ \ \ \ \ \ \ \ \ 
0
\,=\,
L_{a,b,c,d,e}
\big(
f_{\smallbullet,\smallbullet,\smallbullet},\,
g_{\smallbullet,\smallbullet,\smallbullet},\,
h_{\smallbullet,\smallbullet,\smallbullet}
\big)
\eqno
{\scriptstyle{((a,b,c,d)\,\in\,S,\,\,
e\,\in\,\N)}}.
\]

The key and elementary observation 
is that, because $\maux + iv$ is
(weighted) $2$-homogeneous,
the full system $({\sf E})$
{\em splits into the linear subsystems 
$({\sf E}_\nu)$
having separate unknowns $\big(f_{\nu-1},\, g_{\nu-2},\,
h_\nu \big)$}:
\[
({\sf E})
\,=\,
({\sf E}_3)
\cup
({\sf E}_4)
\cup\cdots\cup
({\sf E}_\nu)
\cup
\cdots.
\]
Therefore:
\[
\Big(
({\sf E})
\,\,\,\Longrightarrow\,\,\,
(f,g,h)=(0,0,0)
\Big)
\ \ \ \ \ \ \ \ \ \
\Longleftrightarrow
\ \ \ \ \ \ \ \ \ \
\Big(
({\sf E}_\nu)
\,\,\,\Longrightarrow\,\,\,
\big(f_{\nu-1},g_{\nu-2},h_\nu\big)=(0,0)
\ \ \ \ \
\text{for all}\,\,
\nu\geqslant 3
\Big).
\]
The interest of this equivalence is that one will be able to gather
all powers $v^e$ for $e \in \N$ in order to deal with
functions of the real variable $v \in \R$, and
hence, extract only coefficients of powers $z^a \zeta^b 
\overline{z}^c \overline{\zeta}^d$, as we will see in a while.

Thus, we are left with establishing the following main technical
statement, which will close the proof of
Theorem~{\ref{Thm-uniqueness-complete-normal-form-M5-C3}}.
\endproof

\begin{Theorem}
\label{Thm-S-L-Moser}
In weighted expansions, 
assume that $f = f_2 + f_3 + \cdots$, that $g = g_1 + g_2 + \cdots$,
and that $h = h_3 + h_4 + \cdots$
vanish at the origin and satisfy in addition:
\[
0
\,=\,
f_w(0)
\ \ \ \ \ \ \ \ \ \ \ \ \ \ \ \ \ \ \ \
\text{and}
\ \ \ \ \ \ \ \ \ \ \ \ \ \ \ \ \ \ \ \
0
\,=\,
\Im\,h_{ww}(0).
\]
If, for all $(a,b,c,d) \in S$ and all $e \in \N$:
\[
0
\,=\,
\big[
z^a\zeta^b\overline{z}^c\overline{\zeta}^dv^e
\big]
\bigg(
2\,\Re\,
\Big\{
2\,\tfrac{\overline{z}+z\overline{\zeta}}{1-\zeta\overline{\zeta}}\,
f\big(z,\zeta,\maux+iv\big)
+
\tfrac{(\overline{z}+z\overline{\zeta})^2}{
(1-\zeta\overline{\zeta})^2}\,
g\big(z,\zeta,\maux+iv\big)
-
h\big(z,\zeta,\maux+i\,v\big)
\Big\}
\bigg),
\]
then $(f,g,h) = (0,0,0)$.
\end{Theorem}

\proof
For some reason of technical 
simplification, to be explained in a little
interlude below,
we now decide to `shift' to the representation
$v = F(z, \zeta, \overline{z}, \overline{\zeta}, u)$ 
instead of $u = F(z, \zeta,
\overline{z}, \overline{\zeta}, v)$,
where $u = \Re\, w$ and $v = \Im\, w$ as always.

The hypotheses become (exercise), instead:
\[
0
\,=\,
f_w(0)
\ \ \ \ \ \ \ \ \ \ \ \ \ \ \ \ \ \ \ \
\text{and}
\ \ \ \ \ \ \ \ \ \ \ \ \ \ \ \ \ \ \ \
0
\,=\,
\Re\,h_{ww}(0),
\]
and also, for all $(a,b,c,d) \in S$ and all $e \in \N$:
\[
0
\,=\,
\big[
z^a\zeta^b\overline{z}^c\overline{\zeta}^du^e
\big]
\bigg(
2\,\Re\,
\Big\{
2\,\tfrac{\overline{z}+z\overline{\zeta}}{1-\zeta\overline{\zeta}}\,
f\big(z,\zeta,u+i\,\maux\big)
+
\tfrac{(\overline{z}+z\overline{\zeta})^2}{
(1-\zeta\overline{\zeta})^2}\,
f\big(z,\zeta,u+i\,\maux\big)
+
i\,h\big(z,\zeta,u+i\,\maux\big)
\Big\}
\bigg).
\]

Because $S$ does not contain any dependent-jet
monomial $\zeta \overline{\zeta} (\cdots)$ by
its very definition given above, 
we may compute everything modulo $\zeta 
\overline{\zeta} (\cdots)$, and this will
simplify our task. Thus, by expanding:
\[
\aligned
\maux
&
\,=\,
\frac{z\overline{z}+\frac{1}{2}\,\overline{z}^2\zeta
+\frac{1}{2}z^2\overline{\zeta}}{1-\zeta\overline{\zeta}}
\\
&
\,=\,
z\overline{z}
+
\tfrac{1}{2}\,\overline{z}^2\zeta
+
\tfrac{1}{2}\,z^2\overline{\zeta}
+
z\overline{z}\zeta\overline{\zeta}
+
\tfrac{1}{2}\,\overline{z}^2\zeta\zeta\overline{\zeta}
+
\tfrac{1}{2}\,z^2\overline{\zeta}\zeta\overline{\zeta}
+
z\overline{z}\zeta\overline{\zeta}\zeta\overline{\zeta}
+
\cdots,
\endaligned
\]
we visibly have:
\[
\aligned
\maux
&
\,\equiv\,
z\overline{z}
+
\tfrac{1}{2}\,\overline{z}^2\zeta
+
\tfrac{1}{2}\,z^2\overline{\zeta},
\\
\maux_z
&
\,\equiv\,
\overline{z}
+
z\overline{\zeta},
\\
\maux_\zeta
&
\,\equiv\,
\tfrac{1}{2}\,
\overline{z}^2
+
z\overline{z}\overline{\zeta}
+
\tfrac{1}{2}\,
z^2\overline{\zeta}^2,
\\
\maux_{\overline{z}}
&
\,\equiv\,
z
+
\overline{z}\zeta,
\\
\maux_{\overline{\zeta}}
&
\,\equiv\,
\tfrac{1}{2}\,z^2
+
z\zeta\overline{z}
+
\tfrac{1}{2}\,
\overline{z}^2\zeta^2.
\endaligned
\]
We will also need (little exercise), 
still modulo $\zeta \overline{\zeta} (\cdots)$:
\[
\aligned
\maux^2
&
\,\equiv\,
z^2\overline{z}^2
+
z\zeta\overline{z}^3
+
z^3\overline{z}\overline{\zeta}
+
\tfrac{1}{4}\,
z^4\overline{\zeta}^2
+
\tfrac{1}{4}\,\zeta^2\overline{z}^4,
\\
\maux^3
&
\,\equiv\,
z^3\overline{z}^3
+
\tfrac{3}{2}\,
z^2\zeta\overline{z}^4
+
\tfrac{3}{2}\,
z^4\overline{z}^2\overline{\zeta}^2
+
\tfrac{3}{4}\,
z^5\overline{z}\overline{\zeta}^2
+
\tfrac{3}{4}\,z\overline{z}^5\overline{\zeta}^2
+
\tfrac{1}{8}\,z^6\overline{\zeta}^3
+
\tfrac{1}{8}\,\zeta^3\overline{z}^6,
\\
\maux\,\maux_\zeta
&
\,\equiv\,
\tfrac{1}{2}\,
z\overline{z}^3
+
\tfrac{1}{4}\,\zeta\overline{z}^4
+
\tfrac{5}{4}\,z^2\overline{z}\overline{\zeta}
+
z^3\overline{z}\overline{\zeta}^2
+
\tfrac{1}{4}\,z^4\overline{\zeta}^3,
\\
\maux\,\maux_z
&
\,\equiv\,
z\overline{z}^2
+
\tfrac{1}{2}\,\zeta\overline{z}^3
+
\tfrac{3}{2}\,z^2\overline{z}\overline{\zeta}
+
\tfrac{1}{2}\,z^3\overline{\zeta}^2,
\\
\maux^2\,\maux_z
&
\,\equiv\,
z^2\overline{z}^3
+
z\zeta\overline{z}^4
+
2\,z^3\overline{z}^2\overline{\zeta}
+
\tfrac{1}{4}\,\zeta^2\overline{z}^5
+
\tfrac{5}{4}\,
z^4\overline{z}\overline{\zeta}^2
+
\tfrac{1}{4}\,z^5\overline{\zeta}^3.
\endaligned
\]

Assuming therefore that the graphing equation
$v = F = \maux + G$ is solved with respect to $v$, not to $u$,
with arguments $(z,\zeta,w) = (z,\zeta,\, u+i\,\maux)$,
the Moser (linear) operator is defined as:
\[
L(f,g,h)
\,:=\,
2\,\maux_z\,f
+
2\,\maux_{\overline{z}}\,\overline{f}
+
2\,\maux_\zeta\,g
+
2\,\maux_{\overline{\zeta}}\,\overline{g}
+
i\,h
-
i\,\overline{h}.
\]
Given a holomorphic function $e = e(w)$, we may Taylor expand
at $u$:
\[
\aligned
e\big(u+i\,\maux\big)
\,=\,
&\,
e(u)
+
e_w(u)\,
[i\,\maux]
+
e_{ww}(u)\,
\big[
-\tfrac{\maux^2}{2}
\big]
+
e_{www}(u)\,
\big[
-i\,\tfrac{\maux^3}{6}
\big]
+\cdots
\\
\,=:\,
&\,
e
+
e'\,
[i\,\maux]
+
e''\,
\big[
-\tfrac{\maux^2}{2}
\big]
+
e'''\,
\big[
-i\,\tfrac{\maux^3}{6}
\big]
+\cdots,
\endaligned
\]
and we can abbreviate derivatives using primes,
even without writing
the argument $u$. Let us now make 
the promised little interlude.

The other choice of graphing $u = F = \maux + G$ leads to
$e(w) = e(iv + \maux)$ which expands as:
\[
\aligned
e\big(iv+\maux\big)
&
\,=\,
e(iv)
+
e_w(iv)\,
[\maux]
+
e_{ww}(iv)\,
\big[
\tfrac{\maux^2}{2}
\big]
+
e_{www}(iv)\,
\big[
\tfrac{\maux^3}{6}
\big]
+\cdots.
\endaligned
\]
It is then convenient to consider the {\em composed} function
of one real variable:
\[
v
\,\,\longmapsto\,\,
e(iv)
\,=:\,
E(v),
\]
which satisfies:
\[
\aligned
\tfrac{d}{dv}\,
E(v)
&\,
\,=\,
i\,e_w(iv)
\ \ \ \ \ \ \ \ \ \ \ \ \ \ \ \ \ \ \ \
&
\Longleftrightarrow
\ \ \ \ \ \ \ \ \ \ \ \ \ \ \ \ \ \ \ \
-\,i\,E'(v)
&
\,=\,
e_w(iv),
\\
\tfrac{d^2}{dv^2}\,
E(v)
&\,
\,=\,
-\,e_{ww}(iv)
\ \ \ \ \ \ \ \ \ \ \ \ \ \ \ \ \ \ \ \
&
\Longleftrightarrow
\ \ \ \ \ \ \ \ \ \ \ \ \ \ \ \ \ \ \ \
-\,E''(v)
&
\,=\,
e_{ww}(iv),
\\
\tfrac{d^3}{dv^3}\,
E(v)
&\,
\,=\,
-\,i\,e_{www}(iv)
\ \ \ \ \ \ \ \ \ \ \ \ \ \ \ \ \ \ \ \
&
\Longleftrightarrow
\ \ \ \ \ \ \ \ \ \ \ \ \ \ \ \ \ \ \ \
i\,E'''(v)
&
\,=\,
e_{www}(iv).
\endaligned
\]
Thus:
\[
\aligned
e\big(iv+\maux\big)
&
\,=\,
E(v)
-
i\,E'(v)\,[\maux]
-
E''(v)\,
\big[\tfrac{\maux^2}{2}\big]
+
i\,E'''(v)\,
\big[\tfrac{\maux^3}{6}\big]
+\cdots,
\\
\overline{e}\big(-iv+\maux\big)
&
\,=\,
\overline{E}(v)
+
i\,\overline{E}'(v)\,[\maux]
-
\overline{E}''(v)\,
\big[\tfrac{\maux^2}{2}\big]
-
i\,\overline{E}'''(v)\,
\big[\tfrac{\maux^3}{6}\big]
+\cdots,
\endaligned
\]
and similarly for the conjugate.
If by convention, we then make the abuse of notation 
to denote $e$ instead of $E$, that is $e(v)$ instead of
$E(v) = e(iv)$, we can abbreviate, without writing
the arguments $iv$ or $-iv$:
\[
\aligned
e\big(iv+\maux\big)
&
\,=\,
e
+
e'\,\big[-i\,\maux\big]
+
e''\,\big[-\tfrac{\maux^2}{2}\big]
+
e'''\,\big[i\,\tfrac{\maux^3}{6}\big]
+\cdots,
\\
\overline{e}\big(-iv+\maux\big)
&
\,=\,
\overline{e}
+
\overline{e}'\,\big[i\,\maux\big]
+
\overline{e}''\,\big[-\tfrac{\maux^2}{2}\big]
+
\overline{e}'''\,\big[-i\,\tfrac{\maux^3}{6}\big]
+\cdots.
\endaligned
\]
This can be applied to functions $e = f_{j,k}$ or
$e = g_{j,k}$ or $e = h_{j,k}$ in the useful expansions:
\[
f
\,=\,
\smallsum{j}\,
\smallsum{k}\,
z^j\zeta^k\,
f_{j,k}(w),
\ \ \ \ \ \ \ \ \ \ \ \ \ \ \ \ \ \ \ \
g
\,=\,
\smallsum{j}\,
\smallsum{k}\,
z^j\zeta^k\,
g_{j,k}(w),
\ \ \ \ \ \ \ \ \ \ \ \ \ \ \ \ \ \ \ \
h
\,=\,
\smallsum{j}\,
\smallsum{k}\,
z^j\zeta^k\,
h_{j,k}(w).
\]
But in these last 
paragraphs of our paper, we 
decided to choose $v = F$
in order to simplify a bit the
presentation, so that $e = e(u) = E(u)$ and there 
will be no abuse
of notation.

We can write the Moser operator as:
\[
L(f,g,h)
\,=\,
{\sf T}_1
+
\overline{\sf T}_1
+
{\sf T}_2
+
\overline{\sf T}_2
+
{\sf T}_3
+
\overline{\sf T}_3.
\]
Computing modulo $\zeta \overline{\zeta} (\cdots)$, 
start with:
\[
\!\!\!\!\!\!\!\!\!\!\!\!\!\!\!\!\!\!\!\!\!\!\!\!\!
\scriptsize
\aligned
{\sf T}_3
&
\,\equiv\,
\smallsum{j}\,\smallsum{k}\,
z^j\zeta^k\,
i\,h_{j,k}\big(u+i\,\maux)
\\
&
\,\equiv\,
\smallsum{j}\,\smallsum{k}\,
z^j\zeta^k\,
\Big\{
i\,h_{j,k}
+
h_{j,k}'\,\big[-\maux\big]
+
h_{j,k}''\,\big[-\tfrac{i}{2}\,\maux^2\big]
+
h_{j,k}'''\,\big[\tfrac{1}{6}\,\maux^3\big]
+
\cdots
\Big\}
\\
&
\,\equiv\,
\smallsum{j}\,\smallsum{k}\,
z^j\zeta^k\,
\Big\{
h_{j,k}\,[i]
+
h_{j,k}'
\big[
-z\overline{z}-\tfrac{1}{2}\overline{z}^2\zeta
-\tfrac{1}{2}z^2\overline{\zeta}
\big]
+
h_{j,k}''
\big[
-\tfrac{i}{2}z^2\overline{z}^2
-
\tfrac{i}{2}z\zeta\overline{z}^3
-
\tfrac{i}{2}z^3\overline{z}\overline{\zeta}
-
\tfrac{i}{8}z^4\overline{\zeta}^2
-
\tfrac{i}{8}\zeta^2\overline{z}^4
\big]
\\
&
\ \ \ \ \ \ \ \ \ \ \ \ \ \ \ \ \ \ \ \ \ \ \ \ \ \ \ \ \ \ \ \ \ \ \
\ \ \ \ \ \ \ \ \ \ \ \ \ \ \ \ \ \ \ \ \ \ \ \ \ \ \ \ \ \ \ \ \ \ \
\ \ \ \ \ \ \ \ \ \ \ \ \ \ \ \ \ \ 
+
h_{j,k}'''
\big[
\tfrac{1}{6}z^3\overline{z}^3
+
\tfrac{1}{4}z^2\zeta\overline{z}^4
+
\tfrac{1}{4}z^4\overline{z}^2\overline{\zeta}^2
+
\tfrac{1}{8}z^5\overline{z}\overline{\zeta}^2
+
\tfrac{1}{8}z\overline{z}^5\overline{\zeta}^2
+
\tfrac{1}{48}z^6\overline{\zeta}^3
+
\tfrac{1}{48}\zeta^3\overline{z}^6
\big]
+\cdots
\Big\}
\\
&
\,\equiv\,
\smallsum{j}\,\smallsum{k}\,
\Big\{
h_{j,k}
\big[
iz^j\zeta^k
\big]
+
h_{j,k}'
\big[
-z^{j+1}\zeta^k\overline{z}
-
\tfrac{1}{2}z^j\zeta^{k+1}\overline{z}^2
-
\tfrac{1}{2}z^{j+2}\zeta^k\overline{\zeta}
\big]
\\
&
\ \ \ \ \ \ \ \ \ \ \ \ \ \ \ \ \ \ \ \ \ \ \ \ \ \ \ \ \ \ \ \ \ \ \
\ \ \ \ \ 
+
h_{j,k}''
\big[
-\tfrac{i}{2}z^{j+2}\zeta^k\overline{z}^2
-
\tfrac{i}{2}z^{j+1}\zeta^{k+1}\overline{z}^3
-
\tfrac{i}{2}z^{j+3}\zeta^k\overline{z}\overline{\zeta}
-
\tfrac{i}{8}z^{j+4}\zeta^k\overline{\zeta}^2
-
\tfrac{i}{8}z^j\zeta^{k+2}\overline{z}^4
\big]
\\
&
\ \ \ \ \ \ \ \ \ \ \ \ \ \ \ \ \ \ \ \ \ \ \ \ \ \ \ \ \ \ \ \ \ \ \
\ \ \ \ \ 
+
h_{j,k}'''
\big[
\tfrac{1}{6}z^{j+3}\zeta^k\overline{z}^3
+
\tfrac{1}{4}z^{j+2}\zeta^{k+1}\overline{z}^4
+
\tfrac{1}{4}z^{j+4}\zeta^k\overline{z}^2\overline{\zeta}^2
+
\tfrac{1}{8}z^{j+5}\zeta^k\overline{z}\overline{\zeta}^2
+
\tfrac{1}{8}z^{j+1}\zeta^k\overline{z}^5\overline{\zeta}^2
+
\tfrac{1}{48}z^{j+6}\zeta^k\overline{\zeta}^3
+
\tfrac{1}{48}z^j\zeta^{k+3}\overline{z}^6
\big]
+
\cdots
\Big\}.
\endaligned
\]
The useful expression of $\overline{\sf T}_3$ is obtained
by plain complex conjugation. 

Next, going
only to derivatives of $g_{j,k}$ up to order $1$, which
will be enough, we obtain, without intermediate explanations:
\[
\footnotesize
\aligned
{\sf T}_2
&
\,\equiv\,
\smallsum{j}\,\smallsum{k}\,
\Big\{
g_{j,k}
\big[
z^j\zeta^k\overline{z}^2
+
2z^{j+1}\zeta^k\overline{z}\overline{\zeta}
+
z^{j+2}\zeta^k\overline{\zeta}^2
\big]
\\
&
\ \ \ \ \ \ \ \ \ \ \ \ \ \ \ \ \ \ \ \ \ \ \ \ \ \
+
g_{j,k}'
\big[
iz^{j+1}\zeta^k\overline{z}^3
+
\tfrac{i}{2}z^j\zeta^{k+1}\overline{z}^4
+
\tfrac{5i}{2}z^{j+2}\zeta^k\overline{z}^2\overline{\zeta}
+
2iz^{j+3}\zeta^k\overline{z}\overline{\zeta}^2
+
\tfrac{i}{2}z^{j+4}\zeta^k\overline{\zeta}^3
\big]
+\cdots
\Big\}.
\endaligned
\]
Lastly, going to derivatives of order $2$ of the $f_{j,k}$:
\[
\!\!\!\!\!\!\!\!\!\!\!\!\!\!\!\!\!\!\!\!
\footnotesize
\aligned
{\sf T}_1
&
\,\equiv\,
\smallsum{j}\,\smallsum{k}\,
\Big\{
f_{j,k}
\big[
2z^j\zeta^k\overline{z}
+
2z^{j+1}\zeta^k\overline{\zeta}
\big]
+
f_{j,k}'
\big[
2iz^{j+1}\zeta^k\overline{z}^2
+
iz^j\zeta^{k+1}\overline{z}^3
+
3iz^{j+2}\zeta^k\overline{z}\overline{\zeta}
+
iz^{j+3}\zeta^k\overline{z}^2
\big]
\\
&
\ \ \ \ \ \ \ \ \ \ \ \ \ \ \ \ \ \ \ \ \ \ \ \ \ \ \ \ \ \ \ \ \ \ \
\ \ \ \ \ 
+
f_{j,k}''
\big[
-\tfrac{1}{2}z^{j+2}\zeta^k\overline{z}^3
-
\tfrac{1}{2}z^{j+1}\zeta^{k+1}\overline{z}^4
-
z^{j+3}\zeta^k\overline{z}^2\overline{\zeta}
-
\tfrac{1}{8}z^j\zeta^{k+2}\overline{z}^5
-
\tfrac{5}{8}z^{j+4}\zeta^k\overline{z}\overline{\zeta}^2
-
\tfrac{1}{8}z^{j+5}\zeta^k\overline{\zeta}^3
\big]
+
\cdots
\Big\}.
\endaligned
\]

Now, patiently, in $0 = {\sf T}_1 + {\sf T}_2 + {\sf T}_3 + 
\overline{\sf T}_1 + \overline{\sf T}_2 + \overline{\sf T}_3$, 
we chase coefficients 
$z^a \zeta^b \overline{z}^c \overline{\zeta}^d$
for all $(a, b, c, d) \in S$,
and each time, we obtain linear combinations of
(differentiated) 
functions of $u$.
Using a computer helps to avoid mistakes.

We hence obtain several groups of linear differential equations
in the functions $f_{j,k}(u)$, 
$g_{j,k}(u)$, $h_{j,k}(u)$. We begin with three major groups
coming from (part of) the prenormalization assumption and
which imply a certain agreeable {\sl `nilpotency phenomenon'},
well known to also hold for Levi nondegenerate
hypersurfaces ({\cite{Chern-Moser-1974,
Jacobowitz-1990, Merker-2020}}).
Figures help to grasp the inequalities we are stating below,
which show certain {\sl regions} 
${\bf R}_h^\ast$, ${\bf R}_f^\ast$, ${\bf R}_g^\ast$.

\begin{center}
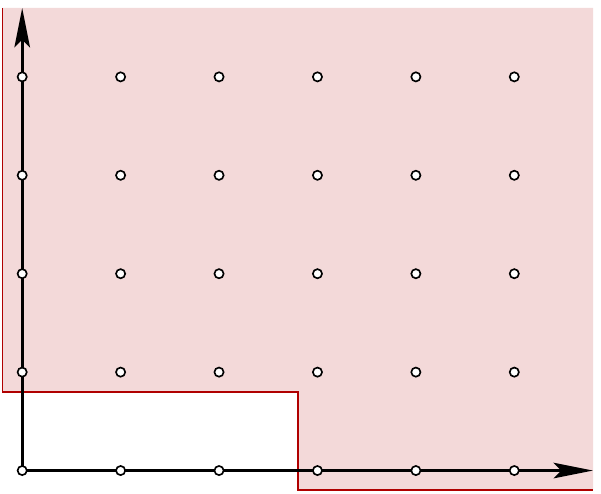
\end{center}

\medskip\noindent$({\bf R}_h^1)$\,
\fbox{$0 \,=\, i\,h_{j,k}(u)$}\,\,
for $(j,k,0,0) \in S$ with $j \geqslant 3$ or with $k \geqslant 1$.
This yields, without writing the argument $u$
of the $h_{j,k}(u)$, that $h$
is a relative polynomial in $(z, \zeta)$:
\[
h
\,=\,
h_{0,0}
+
h_{1,0}\,z
+
h_{2,0}\,z^2.
\]

\begin{center}
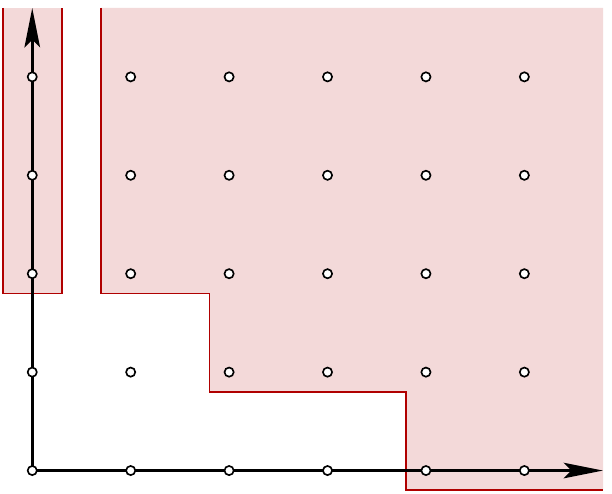
\end{center}

\medskip\noindent$({\bf R}_f^1)$\,
\fbox{$0=2\,f_{0,k}(u)$}\,\,
for $(j,k,1,0) \in S$
with $j = 0$ and $k \geqslant 2$.

\medskip\noindent$({\bf R}_f^2)$\,
\fbox{$0 \,=\, 2\, f_{j,k}(u) - h_{j-1,k}'(u)$}\,\,
for $(j,k,1,0) \in S$ with $j \geqslant 1$ 
and: with $k \geqslant 2$ when $j = 1$; with 
$k \geqslant 1$ when $j = 2, 3$; with $k \geqslant 0$ when $j
\geqslant 4$.
This yields relative polynomialness of:
\[
\aligned
f
&
\,=\,
f_{0,1}\,\zeta
+
f_{1,1}\,z\zeta
\\
&
\ \ \ \ \
+
f_{0,0}
+
f_{1,0}\,z
+
f_{2,0}\,z^2
+
f_{3,0}\,z^3.
\endaligned
\]

\begin{center}
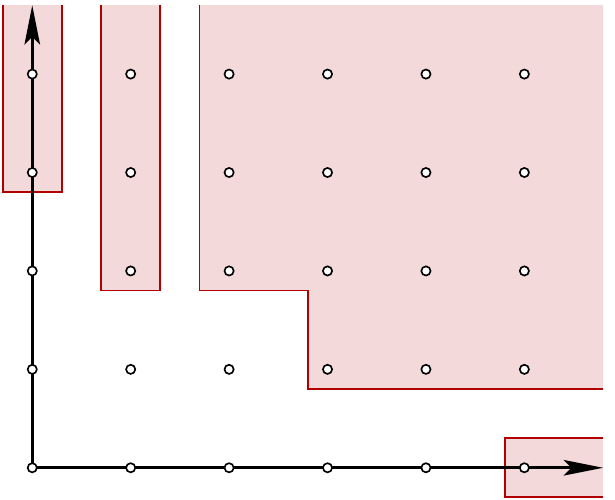
\end{center}

\medskip\noindent$({\bf R}_g^1)$\,
\fbox{$0=g_{0,k} - \frac{1}{2} h_{0,k-1}'$}\,\,
for $(j,k,2,0) \in S$ with $j = 0$ and $k \geqslant 3$.

\medskip\noindent$({\bf R}_g^2)$\,
\fbox{$0 = g_{1,k} - \frac{1}{2} h_{1, k-1}' + 2i\,
f_{0,k}'$}\,\,
for $(j,k,2,0) \in S$ with 
$j = 1$ and $k \geqslant 2$.

\medskip\noindent$({\bf R}_g^3)$\,
\fbox{$0 = g_{j,k} - \frac{1}{2} h_{j,k-1}' + 2i f_{j-1,k}'
- \frac{i}{2} h_{j-2,k}''$}\,\,
for $(j,k,2,0) \in S$ with $j \geqslant 2$ and $k \geqslant 1$
excepting $(j,k) = (2,1)$.

\medskip\noindent$({\bf R}_g^4)$\,
\fbox{$0=g_{j,0} + 2if_{j-1,0}' - \frac{i}{2} h_{j-2,0}''$},
for $(j,k,0,0) \in S$ with $j \geqslant 5$ and $k = 0$.

\medskip

All this also yields relative polynomialness of:
\[
\aligned
g
&
\,=\,
g_{0,2}\,\zeta^2
\\
&
\ \ \ \ \
+
g_{0,1}\,\zeta
+
g_{1,1}\,z\zeta
+
g_{2,1}\,z^2\zeta
\\
&
\ \ \ \ \
+
g_{0,0}
+
g_{1,0}\,z
+
g_{2,0}\,z^2
+
g_{3,0}\,z^3
+
g_{4,0}\,z^4.
\endaligned
\]

To prove that $(f, g, h) = (0, 0, 0)$, it suffices to 
prove that the $3 + 6 + 9$ remaining functions of $u$, namely
$h_{0,0}$, $h_{1,0}$, $h_{2,0}$ 
and $f_{0,1}$, $f_{1,1}$, $f_{0,0}$, $f_{1,0}$, $f_{2,0}$, $f_{3,0}$,
and $g_{0,2}$, $g_{0,1}$, $g_{1,1}$, $g_{2,1}$,
$g_{0,0}$, $g_{1,0}$, $g_{2,0}$, $g_{3,0}$, $g_{4,0}$ 
are identically zero.

For this, we have to examine the remaining groups of
linear ordinary differential equations
with $(a,b,c,d) \in S$.

Firstly (first group), 
the equations for $(j,k,0,0) \in S$ outside the region
${\bf R}_h^1$ are:
\leqnomode\usetagform{default}
\begin{align}
0
&
\,=\,
i\,h_{0,0}
-
i\,\overline{h}_{0,0},
\tag{$0,0,0,0$}
\\
0
&
\,=\,
2\,\overline{f}_{0,0}
+
i\,h_{1,0},
\tag{$1,0,0,0$}
\\
0
&
\,=\,
\overline{g}_{0,0}
+
i\,h_{2,0}.
\tag{$2,0,0,0$}
\end{align}
The conjugate equations are not written, should be understood,
and will in fact be considered later.

Secondly (second group), 
the equations for $(j,k,1,0) \in S$ outside 
${\bf R}_f^1 \cup {\bf R}_f^2$ are:
\leqnomode\usetagform{default}
\begin{align}
0
&
\,=\,
2\,f_{0,0}
-
i\,\overline{h}_{1,0}
\ \ \ \ \ \ \ \ \ \ 
\text{\footnotesize[Already seen]},
\tag{$0,0,1,0$}
\\
0
&
\,=\,
2\,f_{1,0}
-
h_{0,0}'
-
\overline{h}_{0,0}'
+
2\,\overline{f}_{1,0},
\tag{$1,0,1,0$}
\\
0
&
\,=\,
2\,f_{2,0}
+
\overline{g}_{1,0}
-
h_{1,0}'
-
2i\,\overline{f}_{0,0}',
\tag{$2,0,1,0$}
\\
0
&
\,=\,
2\,f_{3,0}
-
h_{2,0}'
-
i\,\overline{g}_{0,0}',
\tag{$3,0,1,0$}
\\
0
&
\,=\,
2\,f_{0,1}
+
2\,\overline{f}_{0,0},
\tag{$0,1,1,0$}
\\
0
&
\,=\,
2\,f_{1,1}
-
\zero{h_{0,1}'}
+
2\,\overline{g}_{0,0}
\tag{$1,1,1,0$}.
\end{align}
Notice that the last equation let appear
$h_{0,1}(u)$, which we already know is identically zero.
Again, the conjugate equations are understood.

Thirdly (third group), the equations for $(j,k,2,0)$ outside 
${\bf R}_g^1 \cup {\bf R}_g^2 \cup {\bf R}_g^3 \cup 
{\bf R}_g^4$ are:
\leqnomode\usetagform{default}
\begin{align}
0
&
\,=\,
g_{0,0}
-
i\,\overline{h}_{2,0}
\ \ \ \ \ \ \ \ \ \ 
\text{\footnotesize[Already seen]},
\tag{$0,0,2,0$}
\\
0
&
\,=\,
g_{1,0}
+
2\,\overline{f}_{2,0}
-
\overline{h}_{1,0}'
+
2i\,f_{0,0}'
\ \ \ \ \ \ \ \ \ \ 
\text{\footnotesize[Already seen]},
\tag{$1,0,2,0$}
\\
0
&
\,=\,
g_{2,0}
+
\overline{g}_{2,0}
-
2i\,\overline{f}_{1,0}'
+
2i\,f_{1,0}'
-
\tfrac{i}{2}\,h_{0,0}''
+
\tfrac{i}{2}\,
\overline{h}_{0,0}'',
\tag{$2,0,2,0$}
\\
0
&
\,=\,
g_{3,0}
-
i\,\overline{g}_{1,0}'
+
2i\,f_{2,0}'
-
\tfrac{i}{2}\,h_{1,0}''
-
\overline{f}_{0,0}'',
\tag{$3,0,2,0$}
\\
0
&
\,=\,
g_{4,0}
+
2i\,f_{3,0}'
-
\tfrac{i}{2}\,h_{2,0}''
-
\tfrac{1}{2}\,\overline{g}_{0,0}'',
\tag{$4,0,2,0$}
\\
0
&
\,=\,
g_{0,1}
+
2\,\overline{f}_{1,0}
-
\tfrac{1}{2}\,
\overline{h}_{0,0}'
-
\tfrac{1}{2}\,
h_{0,0}',
\tag{$0,1,2,0$}
\\
0
&
\,=\,
g_{1,1}
+
2\,\overline{g}_{1,0}
-
\tfrac{1}{2}\,
h_{1,0}'
-
3i\,\overline{f}_{0,0}'
+
2i\,f_{0,1}',
\tag{$1,1,2,0$}
\\
0
&
\,=\,
g_{2,1}
-
\tfrac{1}{2}\,h_{2,0}'
-
\tfrac{5i}{2}\,
\overline{g}_{0,0}'
+
2i\,f_{1,1}'
-
\tfrac{i}{2}\,\zero{h_{0,1}''},
\tag{$2,1,2,0$}
\\
0
&
\,=\,
g_{0,2}
+
\overline{g}_{0,0}
-
\tfrac{1}{2}\,\zero{h_{0,1}'}.
\tag{$0,2,2,0$}
\end{align}
Notice that the last two equations let appear
$h_{0,1}(u)$, which we already know is identically zero.

Fourthly (fourth group) and lastly, we list the sporadic equations:
\leqnomode\usetagform{default}
\begin{align}
0
&
\,\equiv\,
2\,f_{2,0}
-
\tfrac{1}{2}\,
h_{1,0}'
-
i\,\overline{f}_{0,0}',
\tag{$3,0,0,1$}
\\
0
&
\,\equiv\,
\tfrac{1}{6}\,
h_{0,0}'''
+
\tfrac{1}{6}\,
\overline{h}_{0,0}'''
-
f_{1,0}''
-
\overline{f}_{1,0}''
+
i\,g_{2,0}'
-
i\,\overline{g}_{2,0}',
\tag{$3,0,3,0$}
\\
0
&
\,\equiv\,
2\,f_{3,0}
-
\tfrac{1}{2}\,
h_{2,0}'
-
\tfrac{i}{2}\,
\overline{g}_{0,0}',
\tag{$4,0,0,1$}
\\
0
&
\,\equiv\,
2\,g_{2,0}
-
i\,\overline{g}_{0,1}'
-
i\,\overline{f}_{1,0}'
+
3i\,f_{1,0}'
+
\tfrac{i}{2}\,\overline{h}_{0,0}''
-
\tfrac{i}{2}\,h_{0,0}'',
\tag{$3,0,1,1$}
\\
0
&
\,\equiv\,
2\,g_{3,0}
-
\tfrac{i}{2}\,
\overline{g}_{1,0}'
+
3i\,f_{2,0}'
-
\tfrac{i}{2}\,h_{1,0}''
-
\overline{f}_{0,0}''.
\tag{$4,0,1,1$}
\end{align}

Now, the assumptions of 
Theorem~{\ref{Thm-S-L-Moser}}
can be reformulated by comparing the two representations:
\[
f_{\geqslant 2}
\,=\,
\smallsum{j}\,\smallsum{k}\,
z^j\zeta^k\,
f_{j,k}(u),
\ \ \ \ \ \ \ \ \ \ \ \ \ \ \ \ \ \ \ \
g_{\geqslant 1}
\,=\,
\smallsum{j}\,\smallsum{k}\,
z^j\zeta^k\,
g_{j,k}(u),
\ \ \ \ \ \ \ \ \ \ \ \ \ \ \ \ \ \ \ \
h_{\geqslant 3}
\,=\,
\smallsum{j}\,\smallsum{k}\,
z^j\zeta^k\,
h_{j,k}(u),
\]
and one realizes that:
\[
\aligned
0
\,=\,
f(0,0,0)
\,=\,
g(0,0,0)
\,=\,
h(0,0,0)
\ \ \ \ \ \ \ 
&
\Longleftrightarrow
\ \ \ \ \ \ \ 
0
\,=\,
f_{0,0}(0)
\,=\,
g_{0,0}(0)
\,=\,
h_{0,0}(0),
\\
f
\,=\,
f_2+f_3+\cdots
\ \ \ \ \ \ \ 
&
\Longrightarrow
\ \ \ \ \ \ \ 
f_{1,0}(0)
\,=\,
0,
\\
h
\,=\,
h_3+h_4+\cdots
\ \ \ \ \ \ \ 
&
\Longrightarrow
\ \ \ \ \ \ \ 
h_{0,0}'(0)
\,=\,
0,
\\
f_w(0)
\,=\,
0
\ \ \ \ \ \ \ 
&
\Longleftrightarrow
\ \ \ \ \ \ \ 
f_{0,0}'(0)
\,=\,
0,
\\
\Re\,h_{ww}(0)
\,=\,
0
\ \ \ \ \ \ \ 
&
\Longleftrightarrow
\ \ \ \ \ \ \ 
\Re\,h_{0,0}''(0)
\,=\,
0.
\endaligned
\]
The proof of Theorem~{\ref{Thm-S-L-Moser}} will hence be finished
with the next statement.
\endproof

\begin{Proposition}
If $3 + 6 + 9$ analytic functions 
$h_{0,0}$, $h_{1,0}$, $h_{2,0}$ 
and $f_{0,1}$, $f_{1,1}$, $f_{0,0}$, $f_{1,0}$, $f_{2,0}$, $f_{3,0}$,
and $g_{0,2}$, $g_{0,1}$, $g_{1,1}$, $g_{2,1}$,
$g_{0,0}$, $g_{1,0}$, $g_{2,0}$, $g_{3,0}$, $g_{4,0}$ 
of the real variable $u \in \R$ 
with:
\[
\aligned
0
&
\,=\,
f_{0,0}(0)
\,=\,
f_{1,0}(0),
\ \ \ \ \ \ \ \ \ \ \ \ \ \ \ \ \ \ \ \
0
&
\,=\,
g_{0,0}(0),
\ \ \ \ \ \ \ \ \ \ \ \ \ \ \ \ \ \ \ \
0
&
\,=\,
h_{0,0}(0),
\\
0
&
\,=\,
f_{0,0}'(0),
\ \ \ \ \ 
&
\ \ \ \ \ 
0
&
\,=\,
h_{0,0}'(0)
\,=\,
\Re\,h_{0,0}''(0),
\endaligned
\]
satisfy the above system of four groups of
linear ordinary differential equations,
then they all vanish identically.
\end{Proposition}

\proof
From the first two groups of equations {\em and conjugate equations}, 
we may solve:
\[
\aligned
{}&
\aligned
\overline{h}_{0,0}
&
\,:=\,
h_{0,0},
\\
h_{1,0}
&
\,:=\,
2i\,\overline{f}_{0,0},
&
\ \ \ \ \ \ \ \ \ \ \ \ \ \ \ \ \ \ \ \
\ \ \ \ \ \ \ \ \ \ \ \ \ \ \ \ \ \,
\overline{h}_{1,0}
&
\,:=\,
-\,2i\,f_{0,0},
\\
h_{2,0}
&
\,:=\,
i\,\overline{g}_{0,0},
&
\ \ \ \ \ \ \ \ \ \ \ \ \ \ \ \ \ \ \ \
\ \ \ \ \ \ \ \ \ \ \ \ \ \ \ \ \ \,
\overline{h}_{2,0}
&
\,:=\,
-\,i\,g_{0,0},
\endaligned
\\
{}
&
\aligned
\overline{f}_{1,0}
&
\,:=\,
-\,f_{1,0}
+
h_{0,0}',
&
\ \ \ \ \ \ \ \ \ \ \ \ \ \ \ \ \ \ \ \
&
\\
f_{2,0}
&
\,:=\,
-\,\tfrac{1}{2}\,
\overline{g}_{1,0}
+
2i\,\overline{f}_{0,0}',
&
\ \ \ \ \ \ \ \ \ \ \ \ \ \ \ \ \ \ \ \
\overline{f}_{2,0}
&
\,:=\,
-\,\tfrac{1}{2}\,
g_{1,0}
-
2i\,f_{0,0}',
\\
f_{3,0}
&
\,:=\,
i\,\overline{g}_{0,0}',
&
\ \ \ \ \ \ \ \ \ \ \ \ \ \ \ \ \ \ \ \
\overline{f}_{3,0}
&
\,:=\,
-\,i\,g_{0,0}',
\\
f_{0,1}
&
\,:=\,
-\,\overline{f}_{0,0},
&
\ \ \ \ \ \ \ \ \ \ \ \ \ \ \ \ \ \ \ \
\overline{f}_{0,1}
&
\,:=\,
-\,f_{0,0},
\\
f_{1,1}
&
\,:=\,
-\,\overline{g}_{0,0},
&
\ \ \ \ \ \ \ \ \ \ \ \ \ \ \ \ \ \ \ \
\overline{f}_{1,1}
&
\,:=\,
-\,g_{0,0}.
\endaligned
\endaligned
\]
Once this is done, these first two groups of equations become just
$0 = 0$, while the third group 
becomes\footnote{\,\,---\,\,mind 
the fact that because we have sometimes
solved $\overline{e}$ in terms of $e$ for certain
functions $e = e(u)$, the obtained equations
are {\em not} all pairwise conjugates on certain lines,
and this is normal\,\,---}:
\[
\aligned
0
&
\overset{2020}{\,\,=\,\,}
g_{2,0}
+
\overline{g}_{2,0}
-
2i\,h_{0,0}''
+
4i\,f_{1,0}',
&
\ \ \ \ \ \ \ \ \ \ \ \ \ \ \ \ \ \ \ \
&
\\
0
&
\overset{3020}{\,\,=\,\,}
g_{3,0}
-
2i\,\overline{g}_{1,0}'
-
4\,\overline{f}_{0,0}'',
&
\ \ \ \ \ \ \ \ \ \ \ \ \ \ \ \ \ \ \ \
0
&
\overset{2030}{\,\,=\,\,}
\overline{g}_{3,0}
+
2i\,g_{1,0}'
-
4\,f_{0,0}'',
\\
0
&
\overset{4020}{\,\,=\,\,}
g_{4,0}
-
2\,\overline{g}_{0,0}'',
&
\ \ \ \ \ \ \ \ \ \ \ \ \ \ \ \ \ \ \ \
0
&
\overset{2040}{\,\,=\,\,}
\overline{g}_{4,0}
-
2\,g_{0,0}'',
\\
0
&
\overset{0120}{\,\,=\,\,}
g_{0,1}
-
2\,f_{1,0}
+
h_{0,0}',
&
\ \ \ \ \ \ \ \ \ \ \ \ \ \ \ \ \ \ \ \
0
&
\overset{2001}{\,\,=\,\,}
\overline{g}_{0,1}
+
2\,f_{1,0}
-
h_{0,0}',
\\
0
&
\overset{1120}{\,\,=\,\,}
g_{1,1}
+
2\,\overline{g}_{1,0}
-
6i\,\overline{f}_{0,0}',
&
\ \ \ \ \ \ \ \ \ \ \ \ \ \ \ \ \ \ \ \
0
&
\overset{2011}{\,\,=\,\,}
\overline{g}_{1,1}
+
2\,g_{1,0}
+
6i\,f_{0,0}',
\\
0
&
\overset{2120}{\,\,=\,\,}
g_{2,1}
-
5i\,\overline{g}_{0,0}',
&
\ \ \ \ \ \ \ \ \ \ \ \ \ \ \ \ \ \ \ \
0
&
\overset{2021}{\,\,=\,\,}
\overline{g}_{2,1}
+
5i\,g_{0,0}',
\\
0
&
\overset{0220}{\,\,=\,\,}
g_{0,2}
+
\overline{g}_{0,0},
&
\ \ \ \ \ \ \ \ \ \ \ \ \ \ \ \ \ \ \ \
0
&
\overset{2002}{\,\,=\,\,}
\overline{g}_{0,2}
+
g_{0,0},
\endaligned
\]
and the fourth, last, sporadic group becomes:
\[
\aligned
0
&
\overset{3001}{\,\,=\,\,}
2i\,\overline{f}_{0,0}'
-
\overline{g}_{1,0},
&
\ \ \ \ \ \ \ \ \ \ \ \ \ \ \ \ \ \ \ \
0
&
\overset{0130}{\,\,=\,\,}
-\,2i\,f_{0,0}'
-
g_{1,0},
\\
0
&
\overset{3030}{\,\,=\,\,}
-\,\tfrac{2}{3}\,h_{0,0}'''
+
i\,g_{2,0}'
-
i\,\overline{g}_{2,0}',
&
\ \ \ \ \ \ \ \ \ \ \ \ \ \ \ \ \ \ \ \
&
\\
0
&
\overset{4001}{\,\,=\,\,}
i\,\overline{g}_{0,0}',
&
\ \ \ \ \ \ \ \ \ \ \ \ \ \ \ \ \ \ \ \
0
&
\overset{0140}{\,\,=\,\,}
-\,i\,g_{0,0}',
\\
0
&
\overset{3011}{\,\,=\,\,}
2\,g_{2,0}
-
i\,\overline{g}_{0,1}'
-
i\,h_{0,0}''
+
4i\,f_{1,0}',
&
\ \ \ \ \ \ \ \ \ \ \ \ \ \ \ \ \ \ \ \
0
&
\overset{1130}{\,\,=\,\,}
2\,\overline{g}_{2,0}
+
i\,g_{0,1}'
-
3i\,h_{0,0}''
+
4i\,f_{1,0}',
\\
0
&
\overset{4011}{\,\,=\,\,}
2\,g_{3,0}
-
6\,\overline{f}_{0,0}''
-
2i\,\overline{g}_{1,0}',
&
\ \ \ \ \ \ \ \ \ \ \ \ \ \ \ \ \ \ \ \
0
&
\overset{1140}{\,\,=\,\,}
2\,\overline{g}_{3,0}
-
6\,f_{0,0}''
+
2i\,g_{1,0}'.
\endaligned
\]

Hence, from the third group, we can solve:
\[
\aligned
\overline{g}_{2,0}
&
\,:=\,
-\,g_{2,0}
-
4i\,f_{1,0}'
+
2i\,h_{0,0}'',
&
\ \ \ \ \ \ \ \ \ \ \ \ \ \ \ \ \ \ \ \
&
\\
g_{3,0}
&
\,:=\,
2i\,\overline{g}_{1,0}'
+
4\,\overline{f}_{0,0}'',
&
\ \ \ \ \ \ \ \ \ \ \ \ \ \ \ \ \ \ \ \
\overline{g}_{3,0}
&
\,:=\,
-\,2i\,g_{1,0}'
+
4\,f_{0,0}'',
\\
g_{4,0}
&
\,:=\,
2\,\overline{g}_{0,0}'',
&
\ \ \ \ \ \ \ \ \ \ \ \ \ \ \ \ \ \ \ \
\overline{g}_{4,0}
&
\,:=\,
2\,g_{0,0}'',
\\
g_{0,1}
&
\,:=\,
2\,f_{1,0}
-
h_{0,0}',
&
\ \ \ \ \ \ \ \ \ \ \ \ \ \ \ \ \ \ \ \
\overline{g}_{0,1}
&
\,:=\,
-\,2\,f_{1,0}
+
h_{0,0}',
\\
g_{1,1}
&
\,:=\,
6i\,\overline{f}_{0,0}'
-
2\,\overline{g}_{1,0},
&
\ \ \ \ \ \ \ \ \ \ \ \ \ \ \ \ \ \ \ \
\overline{g}_{1,1}
&
\,:=\,
-\,6i\,f_{0,0}'
-
2\,g_{1,0},
\\
g_{2,1}
&
\,:=\,
5i\,\overline{g}_{0,0}',
&
\ \ \ \ \ \ \ \ \ \ \ \ \ \ \ \ \ \ \ \
\overline{g}_{2,1}
&
\,:=\,
-\,5i\,g_{0,0}',
\\
g_{0,2}
&
\,:=\,
-\,\overline{g}_{0,0},
&
\ \ \ \ \ \ \ \ \ \ \ \ \ \ \ \ \ \ \ \
\overline{g}_{0,2}
&
\,:=\,
-\,g_{0,0},
\endaligned
\]
and after that, all equations of the third group reduce to $0 = 0$.
Then the equations of the fourth group become:
\[
\aligned
0
&
\overset{3001}{\,\,=\,\,}
2i\,\overline{f}_{0,0}'
-
\overline{g}_{1,0},
&
\ \ \ \ \ \ \ \ \ \ \ \ \ \ \ \ \ \ \ \
0
&
\overset{0130}{\,\,=\,\,}
-\,2i\,f_{0,0}'
-
g_{1,0},
\\
0
&
\overset{3030}{\,\,=\,\,}
\tfrac{4}{3}\,
h_{0,0}'''
+
2i\,g_{2,0}'
-
4\,f_{1,0}'',
&
\ \ \ \ \ \ \ \ \ \ \ \ \ \ \ \ \ \ \ \
&
\\
0
&
\overset{4001}{\,\,=\,\,}
i\,\overline{g}_{0,0}',
&
\ \ \ \ \ \ \ \ \ \ \ \ \ \ \ \ \ \ \ \
0
&
\overset{0140}{\,\,=\,\,}
-\,i\,g_{0,0}',
\\
0
&
\overset{3011}{\,\,=\,\,}
2\,g_{2,0}
-
2i\,h_{0,0}''
+
6i\,f_{1,0}',
&
\ \ \ \ \ \ \ \ \ \ \ \ \ \ \ \ \ \ \ \
0
&
\overset{1130}{\,\,=\,\,}
-\,2\,g_{2,0}
-
2i\,f_{1,0}',
\\
0
&
\overset{4011}{\,\,=\,\,}
2\,\overline{f}_{0,0}''
+
2i\,\overline{g}_{1,0}',
&
\ \ \ \ \ \ \ \ \ \ \ \ \ \ \ \ \ \ \ \
0
&
\overset{1140}{\,\,=\,\,}
2\,f_{0,0}''
-
2i\,g_{1,0}'.
\endaligned
\]
From this, we can solve, thanks to the assumption $g_{0,0}(0) = 0$:
\[
\aligned
g_{0,0}
&
\,=\,
0,
&
\ \ \ \ \ \ \ \ \ \ \ \ \ \ \ \ \ \ \ \
&
\overline{g}_{0,0}
&
\,=\,
0,
\\
g_{1,0}
&
\,=\,
-\,2i\,f_{0,0}',
&
\ \ \ \ \ \ \ \ \ \ \ \ \ \ \ \ \ \ \ \
&
\overline{g}_{1,0}
&
\,=\,
2i\,\overline{f}_{0,0}',
\\
g_{2,0}
&
\,:=\,
-\,2i\,f_{1,0}'.
\endaligned
\]
The remaining equations become:
\[
\aligned
0
&
\overset{3030}{\,\,=\,\,}
\tfrac{4}{3}\,
h_{0,0}'''
-
2\,f_{1,0}'',
&
\ \ \ \ \ \ \ \ \ \ \ \ \ \ \ \ \ \ \ \
&
\\
0
&
\overset{3011}{\,\,=\,\,}
-\,2i\,h_{0,0}''
+
4i\,f_{1,0}',
&
\ \ \ \ \ \ \ \ \ \ \ \ \ \ \ \ \ \ \ \
&
\\
0
&
\overset{4011}{\,\,=\,\,}
-\,2\,\overline{f}_{0,0}''
&
\ \ \ \ \ \ \ \ \ \ \ \ \ \ \ \ \ \ \ \
0
&
\overset{1140}{\,\,=\,\,}
-2\,f_{0,0}''.
\endaligned
\]
Differentiating once the second equation, using
$0 \neq \big\vert \begin{smallmatrix} \frac{4}{3} & -2 \\ -2i & 4i 
\end{smallmatrix} \big\vert$, we get:
\[
h_{0,0}'''
\,=\,0,
\ \ \ \ \ \ \ \ \ \ \ \ \ \ \ \ \ \ \ \ \ \ \ \ \ \
f_{1,0}''
\,=\,
0.
\]
But we have assumed $0 = h_{0,0}(0) = h_{0,0}'(0) = 
\Re\, h_{0,0}''(0)$, and we know from the beginning
that $\overline{h}_{0,0} = h_{0,0}$ is real.
So $h_{0,0} = 0$.

Back to $\overset{3011}{\,\,=\,\,}$ above, we get $f_{1,0}' = 0$.
Also, we have assumed that $f_{1,0}(0) = 0$. So
$f_{1,0} = 0$. 

Lastly, $f_{0,0}'' = 0$ together with $f_{0,0}'(0) = 0$ 
gives $f_{0,0} = 0$.
This concludes everything.
\endproof


\input references.tex


\vfill\end{document}

%% file: macros.tex
\newtheorem{Theorem}[equation]{Theorem}

\newtheorem{Proposition}[equation]{Proposition}

\newtheorem{Lemma}[equation]{Lemma}

\newtheorem{Corollary}[equation]{Corollary}
\newtheorem{Assertion}[equation]{Assertion}

\newtheorem{Observation}[equation]{Observation}

\theoremstyle{definition}

\newtheorem{Definition}[equation]{Definition}

\newtheorem{Terminology}[equation]{Terminology}

\newtheorem{Question}[equation]{Question}

\newcommand{\C}{\mathbb{C}}

\newcommand{\N}{\mathbb{N}}

\newcommand{\R}{\mathbb{R}}

\newcommand{\NN}{\text{\sc n}}

\newcommand{\TT}{\text{\sc t}}

\newcommand{\kaux}{{\text{\usefont{T1}{qcs}{m}{sl}k}}}

\newcommand{\maux}{{\text{\usefont{T1}{qcs}{m}{sl}m}}}

\newcommand{\Iaux}{{\text{\usefont{T1}{qcs}{m}{sl}I}}}
\newcommand{\Jaux}{{\text{\usefont{T1}{qcs}{m}{sl}J}}}

\newcommand{\Paux}{{\text{\usefont{T1}{qcs}{m}{sl}P}}}
\newcommand{\Qaux}{{\text{\usefont{T1}{qcs}{m}{sl}Q}}}
\newcommand{\Raux}{{\text{\usefont{T1}{qcs}{m}{sl}R}}}

\newcommand{\Vaux}{{\text{\usefont{T1}{qcs}{m}{sl}V}}}
\newcommand{\Waux}{{\text{\usefont{T1}{qcs}{m}{sl}W}}}

\definecolor{blue}{cmyk}{1.,1.,0.,0.63}
\definecolor{red}{cmyk}{0.,1.,1.,0.63}
\definecolor{green}{cmyk}{1.,0.,1.,0.63}
\definecolor{black}{cmyk}{1.,1.,1.,1.}

\newcommand{\blue}{\textcolor{blue}}
\newcommand{\green}{\textcolor{green}}
\newcommand{\red}{\textcolor{red}}

\makeatletter
\renewcommand{\@fnsymbol}[1]
{\ensuremath{\ifcase#1\or $*$\or $**$\or $***$\or $****$\or $*****$
\else\@ctrerr\fi}}
\makeatother

\newcommand{\HEAD}[2]{%
\pagestyle{fancy}
\fancyhead[RO]{\tiny\sf\thepage}
\fancyhead[CO]{{\tiny\sf #1}}
\fancyhead[LE]{\tiny\sf\thepage}
\fancyhead[CE]{{\tiny\sf #2}}
\fancyfoot{}}

\numberwithin{equation}{section}

\newcommand{\Section}[1]{
\renewcommand{\thesection}{\bf\arabic{section}}
\section{#1}
\renewcommand{\thesection}{\arabic{section}}}

\newcommand{\style}[1]{\text{\footnotesize{\sf #1}}}

\newcommand{\stylesmall}[1]{{\sf #1}}

\renewcommand{\dim}{\style{dim}}

\newcommand{\Id}{\style{Id}}

\renewcommand{\Im}{\style{Im}}
\newcommand{\Imsmall}{\stylesmall{Im}}

\newcommand{\Levi}{\style{Levi}}

\renewcommand{\lim}{\style{lim}}

\newcommand{\rank}{\style{rank}}

\renewcommand{\Re}{\style{Re}}
\newcommand{\Resmall}{\stylesmall{Re}}

\newcommand{\Span}{\style{Span}}

\newcommand{\centersmallbullet}{{}_{{}^{{}^{
\scriptscriptstyle{\bullet\!}}}}}

\newcommand{\Hall}{\Hall}

\newcommand{\medcup}{\mathbin{\scalebox{1.5}{\ensuremath{\cup}}}}

\newcommand{\smallbullet}{{\scriptscriptstyle{\bullet}}}

\newcommand{\smallsum}[1]{
\underset{#1}{\raisebox{1pt}{$\sum$\,}}
}

\newcommand{\vf}{\vfill\end{document}}

\newcommand{\zero}[1]{\underline{#1}_{\red{\circ}}}

\newcommand{\zerozero}[1]{\underline{#1}_{\red{\circ\circ}}}


%% file: prolong-jet-1-nontransitive.pdf_t
\begin{picture}(0,0)%
\includegraphics{prolong-jet-1-nontransitive.pdf}%
\end{picture}%
\setlength{\unitlength}{4144sp}%
\begingroup\makeatletter\ifx\SetFigFont\undefined%
\gdef\SetFigFont#1#2#3#4#5{%
  \reset@font\fontsize{#1}{#2pt}%
  \fontfamily{#3}\fontseries{#4}\fontshape{#5}%
  \selectfont}%
\fi\endgroup%
\begin{picture}(2294,2928)(879,-2967)
\put(2669,-735){\makebox(0,0)[lb]{\smash{{\SetFigFont{9}{10.8}{\familydefault}{\mddefault}{\updefault}{\color[rgb]{0,0,0}\green{$\vec{\bf v}^{(1)}$}}%
}}}}
\put(1994,-2903){\makebox(0,0)[lb]{\smash{{\SetFigFont{9}{10.8}{\familydefault}{\mddefault}{\updefault}{\color[rgb]{0,0,0}\blue{$0$}}%
}}}}
\put(1342,-2676){\makebox(0,0)[lb]{\smash{{\SetFigFont{9}{10.8}{\familydefault}{\mddefault}{\updefault}{\color[rgb]{0,0,0}\green{$\vec{\bf v}$}}%
}}}}
\put(1332,-1076){\makebox(0,0)[lb]{\smash{{\SetFigFont{9}{10.8}{\familydefault}{\mddefault}{\updefault}{\color[rgb]{0,0,0}\green{$\vec{\bf v}^{(1)}$}}%
}}}}
\put(2070,-965){\makebox(0,0)[lb]{\smash{{\SetFigFont{9}{10.8}{\familydefault}{\mddefault}{\updefault}{\color[rgb]{0,0,0}\blue{$0$}}%
}}}}
\put(2061,-513){\makebox(0,0)[lb]{\smash{{\SetFigFont{9}{10.8}{\familydefault}{\mddefault}{\updefault}{\color[rgb]{0,0,0}\blue{$\R^4$}}%
}}}}
\put(956,-230){\makebox(0,0)[lb]{\smash{{\SetFigFont{9}{10.8}{\familydefault}{\mddefault}{\updefault}{\color[rgb]{0,0,0}\blue{$J_{1,4}^1$}}%
}}}}
\put(2481,-2675){\makebox(0,0)[lb]{\smash{{\SetFigFont{9}{10.8}{\familydefault}{\mddefault}{\updefault}{\color[rgb]{0,0,0}\green{$\vec{\bf v}$}}%
}}}}
\put(2970,-2709){\makebox(0,0)[lb]{\smash{{\SetFigFont{9}{10.8}{\familydefault}{\mddefault}{\updefault}{\color[rgb]{0,0,0}\blue{$M_{\sf LC}$}}%
}}}}
\end{picture}%

%% file: 2-chains-M5-C3.pdf_t
\begin{picture}(0,0)%
\includegraphics{2-chains-M5-C3.pdf}%
\end{picture}%
\setlength{\unitlength}{4144sp}%
\begingroup\makeatletter\ifx\SetFigFont\undefined%
\gdef\SetFigFont#1#2#3#4#5{%
  \reset@font\fontsize{#1}{#2pt}%
  \fontfamily{#3}\fontseries{#4}\fontshape{#5}%
  \selectfont}%
\fi\endgroup%
\begin{picture}(4881,2844)(654,-2793)
\put(5520,-893){\makebox(0,0)[lb]{\smash{{\SetFigFont{9}{10.8}{\familydefault}{\mddefault}{\updefault}{\color[rgb]{0,0,0}$x_1,\!y_1,\!s_1,\!t_1$}%
}}}}
\put(4832,-475){\makebox(0,0)[lb]{\smash{{\SetFigFont{9}{10.8}{\familydefault}{\mddefault}{\updefault}{\color[rgb]{0,0,0}$x_2,\!y_2,\!s_2,\!t_2$}%
}}}}
\put(1687,-450){\makebox(0,0)[lb]{\smash{{\SetFigFont{7}{8.4}{\familydefault}{\mddefault}{\updefault}{\color[rgb]{0,0,0}$\R^4$}%
}}}}
\put(1466,-105){\makebox(0,0)[lb]{\smash{{\SetFigFont{7}{8.4}{\familydefault}{\mddefault}{\updefault}{\color[rgb]{0,0,0}$\R^4$}%
}}}}
\put(1628,-1985){\makebox(0,0)[lb]{\smash{{\SetFigFont{9}{10.8}{\familydefault}{\mddefault}{\updefault}{\color[rgb]{0,0,0}$0$}%
}}}}
\put(1541,-2699){\makebox(0,0)[lb]{\smash{{\SetFigFont{9}{10.8}{\familydefault}{\mddefault}{\updefault}{\color[rgb]{0,0,0}$0$}%
}}}}
\put(1611,-1775){\makebox(0,0)[lb]{\smash{{\SetFigFont{9}{10.8}{\familydefault}{\mddefault}{\updefault}{\color[rgb]{0,0,0}$\R^4$}%
}}}}
\put(725,-1799){\makebox(0,0)[lb]{\smash{{\SetFigFont{9}{10.8}{\familydefault}{\mddefault}{\updefault}{\color[rgb]{0,0,0}\blue{$J_{1,4}^1$}}%
}}}}
\put(976,-121){\makebox(0,0)[lb]{\smash{{\SetFigFont{9}{10.8}{\familydefault}{\mddefault}{\updefault}{\color[rgb]{0,0,0}\blue{$J_{1,4}^2$}}%
}}}}
\put(691,-2718){\makebox(0,0)[lb]{\smash{{\SetFigFont{9}{10.8}{\familydefault}{\mddefault}{\updefault}{\color[rgb]{0,0,0}\blue{$M_{\sf LC}$}}%
}}}}
\put(5514,-1411){\makebox(0,0)[lb]{\smash{{\SetFigFont{11}{13.2}{\familydefault}{\mddefault}{\updefault}{\color[rgb]{0,0,0}\red{$\Sigma_0^2$}}%
}}}}
\end{picture}%

%% file: M5-C3-translation-normalization.pdf_t
\begin{picture}(0,0)%
\includegraphics{M5-C3-translation-normalization.pdf}%
\end{picture}%
\setlength{\unitlength}{4144sp}%
\begingroup\makeatletter\ifx\SetFigFont\undefined%
\gdef\SetFigFont#1#2#3#4#5{%
  \reset@font\fontsize{#1}{#2pt}%
  \fontfamily{#3}\fontseries{#4}\fontshape{#5}%
  \selectfont}%
\fi\endgroup%
\begin{picture}(6341,1405)(661,-1943)
\put(4994,-688){\makebox(0,0)[lb]{\smash{{\SetFigFont{8}{9.6}{\familydefault}{\mddefault}{\updefault}{\color[rgb]{0,0,0}Normalization}%
}}}}
\put(5311,-846){\makebox(0,0)[lb]{\smash{{\SetFigFont{9}{10.8}{\familydefault}{\mddefault}{\updefault}{\color[rgb]{0,0,0}$\Phi_p$}%
}}}}
\put(2004,-1271){\makebox(0,0)[lb]{\smash{{\SetFigFont{9}{10.8}{\familydefault}{\mddefault}{\updefault}{\color[rgb]{0,0,.69}\blue{$p$}}%
}}}}
\put(826,-1411){\makebox(0,0)[lb]{\smash{{\SetFigFont{9}{10.8}{\familydefault}{\mddefault}{\updefault}{\color[rgb]{0,0,.69}\blue{$M$}}%
}}}}
\put(2335,-922){\makebox(0,0)[lb]{\smash{{\SetFigFont{9}{10.8}{\familydefault}{\mddefault}{\updefault}{\color[rgb]{0,0,0}Translation}%
}}}}
\put(2572,-1045){\makebox(0,0)[lb]{\smash{{\SetFigFont{8}{9.6}{\familydefault}{\mddefault}{\updefault}{\color[rgb]{0,0,0}$\tau_p$}%
}}}}
\put(4258,-1307){\makebox(0,0)[lb]{\smash{{\SetFigFont{9}{10.8}{\familydefault}{\mddefault}{\updefault}{\color[rgb]{0,0,.69}\blue{$0$}}%
}}}}
\put(3069,-1414){\makebox(0,0)[lb]{\smash{{\SetFigFont{9}{10.8}{\familydefault}{\mddefault}{\updefault}{\color[rgb]{0,0,.69}\blue{$M^p$}}%
}}}}
\put(6508,-1307){\makebox(0,0)[lb]{\smash{{\SetFigFont{9}{10.8}{\familydefault}{\mddefault}{\updefault}{\color[rgb]{0,0,.69}\blue{$0$}}%
}}}}
\put(5416,-1879){\makebox(0,0)[lb]{\smash{{\SetFigFont{9}{10.8}{\familydefault}{\mddefault}{\updefault}{\color[rgb]{0,0,0}$N^p$}%
}}}}
\put(3399,-1797){\makebox(0,0)[lb]{\smash{{\SetFigFont{8}{9.6}{\familydefault}{\mddefault}{\updefault}{\color[rgb]{0,0,0}$\Phi_p\circ\tau_p=:\varphi$}%
}}}}
\put(4322,-769){\makebox(0,0)[lb]{\smash{{\SetFigFont{9}{10.8}{\familydefault}{\mddefault}{\updefault}{\color[rgb]{0,0,.69}\blue{$u$}}%
}}}}
\put(6572,-769){\makebox(0,0)[lb]{\smash{{\SetFigFont{9}{10.8}{\familydefault}{\mddefault}{\updefault}{\color[rgb]{0,0,.69}\blue{$u$}}%
}}}}
\put(4617,-1108){\makebox(0,0)[lb]{\smash{{\SetFigFont{9}{10.8}{\familydefault}{\mddefault}{\updefault}{\color[rgb]{0,0,.69}\blue{$z,\zeta,\overline{z},\overline{\zeta},v$}}%
}}}}
\put(6872,-1108){\makebox(0,0)[lb]{\smash{{\SetFigFont{9}{10.8}{\familydefault}{\mddefault}{\updefault}{\color[rgb]{0,0,.69}\blue{$z,\zeta,\overline{z},\overline{\zeta},v$}}%
}}}}
\end{picture}%

%% file: h-Regions.pdf_t
\begin{picture}(0,0)%
\includegraphics{h-Regions.pdf}%
\end{picture}%
\setlength{\unitlength}{4144sp}%
\begingroup\makeatletter\ifx\SetFigFont\undefined%
\gdef\SetFigFont#1#2#3#4#5{%
  \reset@font\fontsize{#1}{#2pt}%
  \fontfamily{#3}\fontseries{#4}\fontshape{#5}%
  \selectfont}%
\fi\endgroup%
\begin{picture}(2734,2239)(799,-2413)
\put(3430,-2251){\makebox(0,0)[lb]{\smash{{\SetFigFont{9}{10.8}{\familydefault}{\mddefault}{\updefault}{\color[rgb]{0,0,0}$j$}%
}}}}
\put(948,-303){\makebox(0,0)[lb]{\smash{{\SetFigFont{9}{10.8}{\familydefault}{\mddefault}{\updefault}{\color[rgb]{0,0,0}$k$}%
}}}}
\put(2390,-785){\makebox(0,0)[lb]{\smash{{\SetFigFont{9}{10.8}{\familydefault}{\mddefault}{\updefault}{\color[rgb]{0,0,0}\red{${\bf R}_h^1$}}%
}}}}
\end{picture}%

%% file: f-Regions.pdf_t
\begin{picture}(0,0)%
\includegraphics{f-Regions.pdf}%
\end{picture}%
\setlength{\unitlength}{4144sp}%
\begingroup\makeatletter\ifx\SetFigFont\undefined%
\gdef\SetFigFont#1#2#3#4#5{%
  \reset@font\fontsize{#1}{#2pt}%
  \fontfamily{#3}\fontseries{#4}\fontshape{#5}%
  \selectfont}%
\fi\endgroup%
\begin{picture}(2779,2239)(754,-2413)
\put(3430,-2251){\makebox(0,0)[lb]{\smash{{\SetFigFont{9}{10.8}{\familydefault}{\mddefault}{\updefault}{\color[rgb]{0,0,0}$j$}%
}}}}
\put(948,-303){\makebox(0,0)[lb]{\smash{{\SetFigFont{9}{10.8}{\familydefault}{\mddefault}{\updefault}{\color[rgb]{0,0,0}$k$}%
}}}}
\put(2377,-768){\makebox(0,0)[lb]{\smash{{\SetFigFont{9}{10.8}{\familydefault}{\mddefault}{\updefault}{\color[rgb]{0,0,0}\red{${\bf R}_f^2$}}%
}}}}
\put(777,-758){\makebox(0,0)[lb]{\smash{{\SetFigFont{9}{10.8}{\familydefault}{\mddefault}{\updefault}{\color[rgb]{0,0,0}\red{${\bf R}_f^1$}}%
}}}}
\end{picture}%

%% file: g-Regions.pdf_t
\begin{picture}(0,0)%
\includegraphics{g-Regions.pdf}%
\end{picture}%
\setlength{\unitlength}{4144sp}%
\begingroup\makeatletter\ifx\SetFigFont\undefined%
\gdef\SetFigFont#1#2#3#4#5{%
  \reset@font\fontsize{#1}{#2pt}%
  \fontfamily{#3}\fontseries{#4}\fontshape{#5}%
  \selectfont}%
\fi\endgroup%
\begin{picture}(2804,2284)(754,-2458)
\put(2833,-783){\makebox(0,0)[lb]{\smash{{\SetFigFont{9}{10.8}{\familydefault}{\mddefault}{\updefault}{\color[rgb]{0,0,0}\red{${\bf R}_g^3$}}%
}}}}
\put(3430,-2251){\makebox(0,0)[lb]{\smash{{\SetFigFont{9}{10.8}{\familydefault}{\mddefault}{\updefault}{\color[rgb]{0,0,0}$j$}%
}}}}
\put(948,-303){\makebox(0,0)[lb]{\smash{{\SetFigFont{9}{10.8}{\familydefault}{\mddefault}{\updefault}{\color[rgb]{0,0,0}$k$}%
}}}}
\put(777,-751){\makebox(0,0)[lb]{\smash{{\SetFigFont{9}{10.8}{\familydefault}{\mddefault}{\updefault}{\color[rgb]{0,0,0}\red{${\bf R}_g^1$}}%
}}}}
\put(1256,-750){\makebox(0,0)[lb]{\smash{{\SetFigFont{9}{10.8}{\familydefault}{\mddefault}{\updefault}{\color[rgb]{0,0,0}\red{${\bf R}_g^2$}}%
}}}}
\put(3543,-2335){\makebox(0,0)[lb]{\smash{{\SetFigFont{9}{10.8}{\familydefault}{\mddefault}{\updefault}{\color[rgb]{0,0,0}\red{${\bf R}_g^4$}}%
}}}}
\end{picture}%